UNIVERSITÉ DU QUÉBEC À MONTRÉAL

MÉMOIRE
PRÉSENTÉ
COMME EXIGENCE PARTIELLE
DE LA MAÎTRISE EN MATHÉMATIQUES

par

SIMON PLOUFFE

# APPROXIMATIONS
# DE SÉRIES GÉNÉRATRICES
# ET QUELQUES CONJECTURES

AOÛT 1992

# AVANT-PROPOS

Un livre très intéressant a été publié en 1973 par N.J.A. Sloane. Il portait le titre *"A Handbook of Integer Sequences* " et comporte plus de 2372 suites d'entiers prises dans tous les domaines des mathématiques et des sciences en général. Depuis sa publication, des milliers de suites nouvelles ont été trouvées, spécialement en combinatoire. L'auteur invitait ses lecteurs à lui communiquer toute correction ou information nouvelle concernant une suite. Il a reçu environ un mètre cube de lettres depuis.

J'entrepris de taper le livre au complet à la main dans un ordinateur au début de 1990; cela m'a demandé 6 mois de travail. Je n'étais pas au bout de mes peines, car une fois cette tâche terminée, j'envoyai une lettre à l'auteur lui indiquant les erreurs dans certaines d'entre elles et que j'avais commencé à constituer une banque de données avec ses suites, etc. Je reçus un coup de téléphone environ deux semaines plus tard. L'auteur était un peu surpris (et moi donc) que quelqu'un se soit donné la peine de taper tout le livre alors que lui avait un fichier sur ordinateur qui contenait toutes les suites. Après une heure de discussion, l'auteur disait qu'il était temps qu'il fasse la 2ème édition de ce livre. Moi je lui disais qu'il était temps que je complète mes études, etc. C'est là que tout a commencé. C'est en essayant de vérifier les suites d'entiers avec un programme que ce projet est né. Je voulais pouvoir vérifier les chiffres des suites pour qu'il n'y ait pas d'erreurs.

C'est également avec l'encouragement et la vision de mon directeur, Gilbert Labelle, que ce mémoire a vu le jour, à la confiance de Pierre Leroux, aux idées génératrices de mon co-directeur, François Bergeron. Tous les autres aussi, qui sont en France, à Bordeaux au LaBRI avec leurs chauds encouragements. Je pense à Xavier Viennot qui m'impressionnait tellement avec ses conférences en 1985, à Maylis Delest, Serge Dulucq, Jean-Guy Penaud, Jean-Marc Fedou, Mireille Bousquet-Mélou, etc. Ceux de Paris à l'INRIA qui m'ont invité à leur en parler et qui ont contribué grandement à faire que le programme *gfun* soit une réalité. Je pense à Paul Zimmermann , “en possession tranquille de la vérité”, Bruno Salvy “le fou de Maple” à qui je dois de vraies belles formules trouvées grâce à ses méthodes (elles sont dans la table en appendice), Philippe Flajolet, “le bon maître”. Je leur dois des discussions fort enrichissantes.

A Neil Sloane évidemment, mon guide et mon maître à penser, qui m’a fait l’honneur de bien vouloir être mon “advisor” comme il se plait lui-même à le dire. Je lui dois de précieux conseils.

A ma mère, qui sera pas mal fière et contente que son garçon fasse une maîtrise en mathématiques.

A ma compagne Danièle, qui m’a  beaucoup aidé au tout début pour la vérification des suites et qui m’a soutenu jusqu’à la fin. Je lui dois et lui dédie ce mémoire.

# RÉSUMÉ

Le présent mémoire tente de répondre à une question simple : Étant donné une suite numérique, comment trouve t-on la fonction génératrice de cette suite? Il s'agit donc de prendre les termes d'une suite et de proposer une façon de les générer à l'aide d'une formule quelconque (simple si possible). Pour ce faire nous avons utilisé des programmes de calcul symbolique couramment disponibles, soit MapleV de l'Université de Waterloo et Pari-GP, un programme développé à l'Université Bordeaux I. Le jeu d'essai des suites est le livre bien connu de Neil J.A. Sloane, *A Handbook of Integer Sequences*[1] . L'exposé se compose de deux parties principales. La première explique les quatre méthodes qui ont permis de répondre à notre question initiale. La deuxième contient une table des formules trouvées à l'aide de ces méthodes. En tout, 1031 fonctions génératrices forment la table sur un total de 4568 suites que composait le jeu d'essai, soit à peu près 23% des suites.

Ces 1031 formules ont toutes été obtenues expérimentalement. C'est donc dire qu'en fait ce sont autant de conjectures. Mais nous verrons que dans presque tous les cas les méthodes sont suffisamment sophistiquées pour pouvoir affirmer que les formules sont les bonnes.

---

[1] Le jeu d'essai est en fait la 2e édition de ce livre qui est en préparation.

# TABLE DES MATIÈRES

# INTRODUCTION

Dans toute cette étude nous procéderons selon une seule ligne directrice: il s'agira de prendre une suite numérique finie et à l'aide d'un programme informatique spécialisé, d'identifier un bon candidat pour la fonction génératrice. Une telle approche pourrait se limiter simplement à consulter un livre de table de suites. Dans [GKP] p.42 on note: "the best source for questions about sequences is an amazing little book called the Handbook of Integer Sequences, by Sloane, which lists thousands of sequences by their numerical values."; et aussi: "the look-up method is limited to problems that other people have decided". De plus, après avoir donné un exemple de suite numérique, les auteurs [GKP] p.327, ajoutent: "no closed form is evident, and this sequence isn't even listed in Sloane's Handbook".

Nous présentons ici une solution à ce problème: c'est-à-dire une méthode alternative aux méthodes standard connues dans ce domaine. Ces dernières partent des propriétés mathématiques d'une suite et de là en font l'analyse, le tout étant basé sur la connaisssance a priori des ces propriétés. Dans la présente étude nous proposons de procéder en sens inverse: c'est uniquement à partir de la suite numérique que les propriétés sont établies. A cette fin, nous décrirons quatre méthodes d'analyse d'une suite numérique.

Ces quatre méthodes s'appuient sur quatre modèles de fonctions génératrices. Le premier modèle suppose que les termes de la suite peuvent être générés avec le développement en série de Taylor d'un quotient de polynômes. Le deuxième modèle suppose que la suite satisfait à une récurrence linéaire à coefficients polynomiaux, appelée aussi une P-récurrence. Le troisième modèle suppose que la suite est donnée par le développement en série d'un produit infini, comme la suite des partages d'entiers ordinaires. Le quatrième modèle enfin suppose qu'une transformation simple de la suite permet de retrouver une fonction génératrice connue. Cette dernière est en fait une version améliorée de la "look-up method" de [GKP].

## NOTES

Pour éviter les répétitions inutiles tout au long de cette étude, nous emploierons la notation Nxxxx pour désigner la suite numéro xxxx du livre de Sloane [Sl] cité en bibliographie. Par exemple, la suite des nombres de Catalan porte le numéro 577, on y fera référence en écrivant N0577. Les autres suites, celles apparues après 1973 dans la littérature, ont été recataloguées dans une deuxième édition que nous préparons avec Neil J.A. Sloane [PlSl]. Elles portent un numéro séquentiel “absolu” noté Axxxx. Donc quand nous parlerons du numéro de suite A3890, nous entendrons le numéro séquentiel de cette table. C’est cette même numérotation qui apparaît dans la table des résultats en appendice. Pour des raisons évidentes de consistence, il était nécessaire de conserver un numéro qui fasse référence toujours à la même suite sans ambiguité. En résumé :

- Nxxxx :Numéro séquentiel de la suite du livre de Sloane [Sl].
- Axxxx :Numéro séquentiel de la suite du livre [PlSl].

La plupart des algorithmes et méthodes décrites dans cette étude ont été regroupés dans un programme appelé “gfun” qui fait partie de la librairie publique de Maple de l’université de Waterloo. On peut avoir une copie de ce programme par transfert électronique via “ftp/anonymous”. Le programme a été écrit en collaboration avec François Bergeron, professeur au département de Mathématiques/Informatique à l’Université du Québec à Montréal et également avec Paul Zimmermann et Bruno Salvy tous deux chercheurs à L’INRIA/Rocquencourt.

# CHAPITRE 1

# LA MÉTHODE DES APPROXIMANTS DE PADÉ

## 1.1 Les fractions rationnelles.

Une façon de donner les termes d'une suite est de les engendrer à l'aide d'une fraction rationnelle. Par exemple la suite de Fibonacci peut être générée à l'aide du développement en série de Taylor à l'origine de

(1.1) $1/(1- z - z^2) = 1 + z + 2 z^2 + 3 z^3 + 5 z^4 +8 z^5 + 13 z^6 + ... + a_n z^n + O(z^{n+1})$.

De la même façon ces termes peuvent être calculés avec la récurrence $a_n = a_{n-1} + a_{n-2}$. Les deux représentations sont équivalentes. Il y a une correspondance assez simple entre la fraction rationnelle et la récurrence. Le dénominateur de la fraction rationnelle "est" la relation de récurrence. Le numérateur tient lieu de conditions initiales de cette même récurrence. Le lien se trouve en fait dans la réécriture de la récurrence en termes de $z^n$ plutôt que n. Une procédure simple en quatre étapes, permettant de passer d'une récurrence linéaire à coefficients constants à la fraction rationnelle correspondante, est décrite dans [GKP]. On peut montrer que la réécriture se fait dans l'autre sens également. Cette mécanique est très connue mais suppose toujours que l'on connaisse au moins l'une des deux représentations.

Notre seul point de départ est la série S(z) tronquée à l'ordre k. Ce qu'on désire faire est de la représenter par une fraction rationnelle. Alors si on pose k=L+M et

(1.2) $$S(z) = \frac{u_0 + u_1 z^1 + ... + u_L z^L}{v_0 + v_1 z^1 + ... + v_m z^m} + O(z^{L+M+1)})$$

il est toujours possible de trouver une solution à cette équation. La façon de faire est de fixer L et M

d'abord. Puis en multipliant le membre de gauche avec le dénominateur on obtient un système de M équations à M inconnues qui déterminent les constantes en v. Pour poser ces équations, on "identifie" les coefficients de $z^i$ avec $L+1 \leq i \leq L+M$. Une fois trouvés les $v_i$, on peut faire de même avec le numérateur pour déterminer les constantes en $u_j$, en identifiant cette fois les coefficients de $z^j$, $0 \leq j \leq L$. Il est toutefois plus aisé de poser en partant que $v_0 = 1$. On donne le nom d'approximant de Padé [L/M] à l'expression trouvée pour un L et M donnés. Le calcul d'un approximant de Padé se fait en principe de façon mécanique. La plupart des programmes de calcul symbolique sur le marché aujourd'hui effectuent ce calcul automatiquement. On parle ici de la résolution du système d'équations linéaires pour un L et un M donnés. En théorie le problème est clos, mais dans la pratique il en est autrement.

Nous illustrerons les difficultés rencontrées en donnant deux exemples extrêmes d'approximants de Padé.

**Exemple 1.1** La suite des parts de gâteaux en 3 dimensions.

Le premier est la suite N0419 qui porte le nom de: "Slicing a cake with n slices"; elle est plutôt simple et connue. C'est le nombre de parts de gâteaux différents avec n coupes en 3 dimensions. Nous nous en tiendrons uniquement aux termes numériques sans tenir compte du contexte. Considérant une quinzaine de termes, on pose les équations et les degrés des deux polynômes, en supposant que les degrés sont de taille égale, i.e. que L=M. La difficulté réside dans le fait que si le système *peut* se réduire, il faut prévoir un algorithme pour le simplifier d'une façon ou d'une autre. Justement, cette suite N0419 est une fraction rationnelle de degré [2/4]. Elle est complètement décrite par cette fraction rationnelle. C'est donc que, ayant pris notre quinzaine de termes et ayant supposé que L=M=7, on aurait été conduit à réduire le système à un nombre d'inconnues et d'équations plus petit. Donc à moins d'être chanceux, i.e. de prévoir exactement à l'avance le degré de la fraction rationnelle, on n'est pas assuré de trouver la *juste* fraction rationnelle.

La deuxième difficulté vient de la taille des calculs. Si la suite considérée EST une fraction rationnelle, comme la suite de Fibonacci (1.1), cela n'a rien de dramatique si on a fait un choix de L et M heureux. Si la suite N'EST PAS une fraction rationnelle, c'est là que les calculs deviennent énormes. Selon l'équation (1.2), il est quand même possible de trouver une fraction rationnelle qui se juxtaposera aux k premiers termes de toute suite, mais elle ne se simplifiera pas.

**Exemple 1.2** La suite des nombres premiers : 2,3,5,7,11,13,....

Nous prendrons ici la suite des nombres premiers N0241. On le sait, il n'existe pas de fonction rationnelle qui permette de les obtenir successivement. Si on prend les 20 premiers termes , de 2 à 71 et que l'on cherche une expression rationnelle qui se juxtapose à cette suite, l'expression que l'on trouvera sera une fraction rationnelle d'une taille appréciable. La taille, disons en nombre de caractères, dépassera largement celle de la suite. Il ne faut pas oublier que l'on cherche une solution rationnelle, donc exacte à l'ordre d'approximation de la série de départ; ce ne sont pas des calculs en "virgule flottante". Ainsi, avec les 48 premiers termes de la suite des nombres premiers on obtient une fraction rationnelle d'une taille de l'ordre de 10,000 caractères, chaque coefficient étant de l'ordre de 120 chiffres. La taille de la suite de départ avec ses 48 termes, pour sa part, ne dépasse pas 200 caractères.

Il existe une procédure en Maple qui permet de convertir une série (tronquée) en une fraction rationnelle. Elle porte le nom de "ratpoly" pour "rational polynomial". Cette procédure est une véritable perle de programmation (elle a plus de 500 lignes). Non seulement elle fait le calcul exactement à l'ordre d'approximation de la série, mais en plus elle le fait bien. On le sait, Maple est en mesure d'effectuer des calculs symboliquement et en principe avec une précision infinie. Le résultat en est que les deux difficultés rencontrées plus tôt sont complètement tranparentes à l'utilisateur.

Donc avec cet outil presque "magique" qu'est "ratpoly", il est possible assez facilement d'effectuer le calcul fastidieux de représentation d'une suite sous forme de série avec une fraction rationnelle. En fait, deux critères simples nous permettrons de détecter une *bonne* fraction rationnelle. Le premier est le degré de l'expression trouvée: si le degré total (L + M) retourné par le programme est plus petit que le nombre de termes, on est alors potentiellement en présence d'une bonne représentation. Le deuxième critère est la taille (en nombre de caractères) de l'expression: si la taille de l'expression est plus grande que la taille de la suite testée, on rejette alors l'expression rationnelle candidate. En combinant ces deux critères, il est possible de détecter avec une assez grande certitude une suite qui EST une fraction rationnelle simple.

En soumettant toute notre table de 4568 suites à cette simple procédure qu'est "ratpoly" , nous avons détecté 614 fractions rationnelles. De ce nombre, 580 nous semblent bonnes: elles sont répertoriées dans la table en appendice. On peut consulter [BP] à ce sujet également.

## 1.2 La dérivée logarithmique et l'inverse fonctionnel.

Malgré le succès remporté (en nombre de fonctions génératrices trouvées) avec notre méthode des approximants de Padé, une partie du problème demeure. Si 580 suites sur 4568 sont des fractions rationnelles, quelle est alors la nature des quelque 4000 qui restent ? La réponse à cette question est inconnue. Ce que l'on sait, c' est que notre méthode permet de détecter la fraction rationnelle d'une suite comme celle de Fibonacci. Elle permet également de détecter des variantes de celle-ci. Il se trouve que la plupart des opérations simples et connues que l'on peut effectuer sur une suite sont en fait des *transformations rationnelles* . Une TR en plus court. Si S(z) est notre suite sous forme de série tronquée, une TR conservera le caractère rationnel de la fonction génératrice. Par exemple, la différence terme à terme de la suite est une TR puisqu'il suffit d'effectuer S(z)(1-z). La somme de deux termes successifs est également une TR : il suffit de faire S(z)(1+z). La suite des sommes partielles s'obtient en prenant S(z)/(1-z), etc. Il en est de même de l'inverse de ces transformations. Ce point est essentiel.

Donc les suites qui ont une fonction génératrice qui est une fraction rationnelle et toutes les variations usuelles de celles-ci sont détectées avec notre méthode.

L'idée fort simple est alors d'utiliser notre méthode et une transformation qui ne soit pas rationnelle dans les deux sens, dans le but de détecter d'autres types de fonctions génératrices. Par exemple, bien que la dérivée soit une transformation qui conserve le caractère rationnel d'une expression, il suffit de prendre une fraction rationnelle quelconque pour se rendre compte que l'intégrale n'est pas une fraction rationnelle en général. En effectuant une dérivation et en appliquant ensuite notre méthode de détection, on pourra obtenir des fonctions génératrices qui sont en fait des intégrales de fractions rationnelles. En poussant le même raisonnement plus loin, on pourrait effectuer d'autres transformations de ce type comme la dérivée du logarithme ou l'inverse fonctionnel. Si nous pouvons toujours retourner sur nos pas à chaque fois, cela nous donne une façon de détecter des expressions qui font partie d'une classe plus vaste que les fractions rationnelles. Avec la dérivée, il est facile de revenir en arrière: une fois le test effectué, si c'est rationnel, il suffit de faire l'intégrale de l'expression. La dérivée du logarithme est aussi "réversible": il suffit de faire l'exponentielle de l'intégrale de l'expression trouvée. L'inverse fonctionnel d'une série est également "réversible" à condition que la suite débute par 0,1,...:

en effet, l’inverse d’une série à coefficients entiers est aussi à coefficients entiers, si la série s’annule en zéro et son premier terme non nul est 1.

C’est l’expérience qui a orienté le choix des transformations judicieuses à effectuer. Le succès d’une transformation plutôt que d’une autre étant guidé simplement par le nombre de fonctions génératrices trouvées une fois la table complète traitée par le programme. Le rejet ou l’acceptation d’une expression est donné par les deux critères énoncés plus haut. Il y a aussi le fait que plus on tranforme une suite avec de telles opérations, plus précises et strictes sont les conditions imposées à la suite de départ. Par exemple, l’inverse fonctionnel de la dérivée du logarithme d’une suite sous forme de série tronquée doit se faire seulement si les coefficients sont restés entiers et débutent par 0,1, ... , une fois que la première transformation a été effectuée.

Notre choix s’est arrêté sur la dérivée, la dérivée logarithmique et l’inverse fonctionnel. Ce sont ces opérations qui ont remporté le plus de succès. Exactement 120 fonctions génératrices qui ne sont pas des fractions rationnelles ont été isolées de cette façon. En tout 700 fonctions génératrices (incluant les fractions rationnelles) ont été trouvées grâce à la procédure “ratpoly”. Les résultats sont présentés en appendice.

# CHAPITRE 2

# LA MÉTHODE DES P-RÉCURRENCES

## 2.1 Les suites P-récurrentes.

L'hypothèse de travail que nous posons ici sur la suite $a_n$ consiste à dire que chaque terme de celle-ci peut être calculé à partir des termes précédents. Dans [Sta80] on introduit ce genre de dépendance sur les autres termes en disant que la suite $a_n$ est une suite P-récurrente, si elle satisfait l'équation suivante

(2.1) $$a_n P_0(n) = a_{n-1}P_1(n) + a_{n-2}P_2(n) + ..., + a_{n-k}P_k(n)$$

où les $P_i(n)$, $0 \le i \le k$, sont des polynômes à coefficients rationnels. Ce type de relation est une classe plus vaste que les relations de récurrences linéaires ordinaires à coefficients constants rencontrées au chapitre précédent. En effet, il y a équivalence entre les fonctions génératrices rationnelles et les relations de récurrence à coefficients constants. Il n'y a cependant pas d'équivalent en termes de fonctions génératrices pour les P-récurrences en général. A l'heure actuelle, il n'existe pas de méthode pour trouver la fonction génératrice correspondant à une P-récurrence quelconque; seuls certains types de P-récurrences peuvent être résolus. Ce qui peut être fait, par contre, est de vérifier si la suite satisfait *numériquement* une P-récurrence. On ne peut donner qu'une P-récurrence *vraisemblable*.

Il faut donc procéder pas-à-pas en augmentant le degré et le nombre de termes. Nous posons d'abord les équations et, en supposant que la suite satisfasse l'équation (2.1) où les $P_i(n)$ sont des polynômes de degré d, il y aura (d+1)(k+1) équations (il faut tenir compte du terme de rang 0). On dira alors qu'elle satisfait une P-récurrence de type (d,k). On remarque que le système admet toujours une solution nulle. S'il y a solution, il y en aura une infinité, ce qui découle du fait que le système d'équations est non-homogène. Ceci est évident, puisque l'on peut multiplier par une constante C arbitraire de chaque côté sans changer l'équation. On prendra donc soin de garder la solution la plus simple. La résolution d'un système d'équations linéaires est une chose que les programmes de calcul symbolique comme Maple font couramment. Un programme a donc été écrit pour permettre de résoudre le système à (d+1)(k+1) inconnues. Le voici, en entrée il accepte une suite et en sortie il donne soit 0 soit une ou plusieurs constantes, quand le nombre de constantes est 1 on pose la solution comme étant la plus simple en substituant la constante à 1.

```
1) read suite : listesuite:=":
2) nbrdetermes:=nops(listesuite):
3)    rec:=proc(w,n,t) local ff,c,d,i,j,k,ii;
4)     option remember;
5)    termes:=(n+1)*(t+1);
6)     if termes>=nbrdetermes then RETURN ( `impossible de resoudre` ) fi;
7)     for ii from 1 to nbrdetermes do a(ii):=op(ii,w) od:
8)            ens:={seq(c[jj],jj=1..termes)}:
9)    s:={seq(sum(sum(k**d*c[j*n+jn+d],d=0..n)*a(kj+1),j=1..t+1),
              k=t+1..termes+t)}:
10)   solution:=[solve(s,ens)];
11)           if sol=[] then RETURN (0) else
              RETURN(assign(solution),[seq(c[kk],kk=1..termes)])
              fi;
      end:
```

Donnons une courte description du programme.

**1)** On lit la suite provenant d'un fichier.

**2)** On pose que la variable nbrdetermes est égal au nombre d'éléments de la liste qui contient la suite.

**3)** Appel de la procédure et on pose les variables locales.

**4)** On prend l'option "remember" , très importante.

**5)** On prend un nombre de termes suffisant pour résoudre le système d'équations linéaires.

**6)** Si le nombre de termes nécessaires est trop grand, un message d'erreur est imprimé.

**7)** On pose les constantes dans notre système d'équations. Ici ce sont les termes de la suite.

**8)** On pose les inconnues de notre système sous forme d'ensemble.

**9)** On pose les équations linéaires.

**10)** On tente de résoudre.

**11)** Si le système admet une solution nulle ( liste vide ici ) on retourne 0. Sinon on assigne les solutions trouvées.

Donc en entrée le programme accepte une suite numérique et teste si celle-ci satisfait une équation P-récurrente de degré d à k termes.

Le programme qui détermine si une suite satisfait une P-récurrence est une des méthodes les plus rapides et de plus, une fois la P-récurrence candidate trouvée, il est très facile d'obtenir des centaines de termes de la suite. En principe si on veut calculer les termes d'une suite une fois obtenue une P-récurrence, il suffit de la mettre telle quelle dans un programme. Il n'est cependant pas approprié d'utiliser une procédure qui soit purement récursive même si c'est d'abord ce qui vient à l'esprit. Il faut linéariser le temps de calcul d'une procédure qui s'appelle elle-même, sinon celui-ci devient vite exponentiel. L'exemple souvent donné dans les cours de programmation de base est la suite des nombres factoriels, 1,1,2,6,24,120,720,..., définie par $a_0 = 1$ et $a_n = n\, a_{n-1}$. Ce problème est facilement résoluble en Maple, puisque les procédures récursives peuvent être *linéarisées* simplement en écrivant "option remember" dans l'appel de la procédure. Maple se charge alors de ré-écrire la procédure en créant une table d'adressage ( interne) automatiquement.

Comme avec les autres méthodes, nous avons utilisé la table de [PlSl] au complet. A chaque suite, le test a été effectué sur les degrés 1 à 4 et sur un nombre de termes variant de 1 à 5, compte tenu que l'expérience indique que la plupart des suites P-récurrentes ont un degré assez bas. Stanley [Sta80] donne un exemple de suite P-récurrente de degré 3 à 2 termes qui donne les nombres d'une suite de Apéry utilisée dans la preuve de l'irrationalité de ζ(3).

**Exemple 2.1 :** $n^3 a_{(n)} + (n-1)^3 a_{(n-2)} = (34 n^3 - 51 n^2 + 27 n - 5) a_{(n-1)},$

En tout, 250 des 1031 suites que contient la table en appendice, seraient P-récurrentes. De ces 250, 220 ont une fonction génératrice associée trouvée par d'autre méthodes. Il en reste donc 30 dont on ne connait que la P-récurrence. Sont comptées ici les suites P-récurrentes de degré 1 ou plus; les fractions rationnelles, au nombre de 580, sont aussi P-récurrentes mais de degré 0. C'est de loin la méthode la plus puissante, puisque au total, près de 81 % des suites qui ont une fonction génératrice connue sont P-récurrentes à des degrés divers, ce qui représente 18 % de tout le catalogue des suites de [PlSl].

## 2.2 Les suites hypergéométriques.

Dans [GKP] on fait une remarque très simple au sujet des P-récurrences d'un certain type. Si une suite $t_k$ satisfait une P-récurrence de type (d,1), c'est donc que le quotient des termes successifs $t_{k+1}/t_k = P(k)/Q(k)$, où P(k) et Q(k) sont deux polynômes. La fonction hypergéométrique est à peu de chose près la même chose. En effet, la définition de celle-ci étant

$$F\left[\begin{matrix} a_1, a_2, \ldots a_m \\ b_1, b_2, \ldots, b_n \end{matrix} \middle| z\right] = \sum_{k=0}^{\infty} \frac{a_1^{\overline{k}} \ldots a_m^{\overline{k}}}{b_1^{\overline{k}} \ldots b_n^{\overline{k}}} \frac{z^k}{k!} \tag{2.2}$$

où le membre de gauche en est l'écriture avec les paramètres en a et en b et où le membre de droite en est le développement en série sous forme de somme de quotients de produits de polynômes factoriels ascendants. En spécifiant que les termes en b ne s'annulent nulle part, on évite la division par zéro; il sufit simplement pour cela qu'ils soient toujours positifs. Considérons le rapport de deux termes successifs et en posant que le premier terme $t_0=1$,

$$\frac{t_{k+1}}{t_k} = \frac{a_1^{\overline{k+1}} \ldots a_m^{\overline{k+1}}}{a_1^{\overline{k}} \ldots a_m^{\overline{k}}} \frac{b_1^{\overline{k}} \ldots b_n^{\overline{k}}}{b_1^{\overline{k+1}} \ldots b_n^{\overline{k+1}}} \frac{k!}{(k+1)!} \frac{z^{k+1}}{z^k}$$

il est alors facile de simplifier cette expression en revenant à la définition d'un polynôme factoriel ascendant de degré k+1 et de degré k. D'où l'expression:

$$\frac{t_{k+1}}{t_k} = \frac{(k+a_1)\ldots(k+a_m)z}{(k+b_1)\ldots(k+b_m)(k+1)} \quad .$$

On obtient alors une fraction rationnelle en k seulement. Donc si on a une suite qui débute avec 1 et dont le rapport des termes successifs est une fraction rationnelle (une P-récurrence de type (d,1)), elle pourra être "lue" directement comme étant une série hypergéométrique. L'avantage énorme de la représentation d'une suite comme "hypergéométrique" est que le programme de calcul symbolique Maple est en mesure de manipuler et de simplifier de telles séries. Dans sa version 5, Maple utilise les tables d'identités hypergéométriques qui se trouvent dans [AS1]. Ce livre étant une véritable bible de formules mathématiques, nous avons à notre disposition un outil excessivement puissant. En fait, dès que l'on sait qu'une suite satisfait une P-récurrence de type (d,1) nous disposons déjà d'une information très précieuse.

Cette représentation en série hypergéométrique ouvre la porte à d'autres formes de fonctions génératrices. Le programme Maple est en effet capable, dans certains cas, de donner directement la fonction génératrice explicite sous forme simplifiée. Il suffit de faire appel à la procédure “simplify” qui réussit à reconnaître les expressions contenant des termes hypergéométriques. C’est alors que les tables d’identités de [AS1] sont appelées et, si la forme le permet, Maple retourne directement une expression algébrique explicite.

Conformément aux autres méthodes nous avons donc, encore une fois, testé toute la table de [PlSl] en recherchant des P-récurrences de type (d,1). Plus de 94 suites satisfont à ce type de récurrence. Dans certains cas, la forme hypergéométrique a été directement simplifiée automatiquement par le programme Maple. Les résultats sont présentés dans la table de fonctions génératrices en appendice.

## 2.3 L’algorithme LLL[2] .

Nous décrirons ici la méthode qui est la plus complexe et puissante de toute cette étude. On s’intéresse aux suites qui sont P-récurrentes de type (d,k) en général. Cette méthode ne s’applique que si on peut avoir autant de termes de la suite que l’on veut. Comme nous l’avons vu à la section précédente, Maple est en mesure, dans les cas où la P-récurrence est de type (d,1), de donner une forme hypergéométrique et une fois obtenu cette forme, de produire directement la fonction génératrice algébrique lorsqu’elle s’y prête. C’est donc que : les P-récurrences de type (d,1) sont quelquefois algébriques. Il en est de même pour les P-récurrences de d’ordre plus élevé. Ce qui nous manque est la façon d’obtenir la forme close. On ne dispose malheureusement pas de moyen de savoir quel type de P-récurrence représente une suite qui a une fonction génératrice algébrique. D’après Stanley [Sta80], une fonction génératrice algébrique est toujours P-récurrente. Ici c’est l’inverse qu’on cherche, malheureusement ce n’est pas toujours vrai : la fonction exp(x) est P-récurrente mais certainement pas algébrique.

Une suite a une fonction génératrice algébrique si elle satisfait à

(2.3) $$\sum_{0 \le j,k \le m} c_{j,k} S(z)^j z^k = 0$$

[2] **Nommé ainsi à cause des travaux de Lenstra, Lenstra et Lovasz.**

où S(z) est la série qui représente la suite $a_n$ et les $c_{j,k}$ sont constantes. On pourra alors obtenir la fonction génératrice close si on peut isoler S(z). Le problème est double ici: il faut d'abord obtenir l'équation (2.3) et de plus on n'est pas assuré de pouvoir isoler S(z). Ce qui vient à l'esprit est d'essayer de trouver “à tâtons” une équation en S(z) et z qui s'annulera. Il est possible effectivement de faire un programme qui fonctionnerait sur le même principe que les P-récurrences. Mais malheureusement la forme qu'on obtiendra ne sera pas, en général, la plus simple. Par exemple, la suite N0577, les nombres de Catalan, satisfait à une telle équation. Elle est de degré 2 : $S(z)^2 z - S(z) + 1 = 0$. Si on résout cette équation par rapport à S(z), on obtient une fonction génératrice close des nombres de Catalan. On s'aperçoit alors qu'il il y a une infinité de telles équations que l'on peut poser. On pourrait peut-être en obtenir une plus simple.

Il existe un algorithme implanté en Maple qui porte le nom de “minpoly”. Il fait appel à l'algorithme LLL. Disons simplement qu'il permet de résoudre numériquement le problème exactement inverse de trouver une racine d'un polynôme. La recherche numérique des racines d'un polynôme est un problème résolu. Mais nous posons la question suivante: étant donné un nombre réel, de quel polynôme minimal est-il racine ? Mentionnons dès maintenant qu'on parle ici d'un nombre réel donné avec une certaine précision numérique. On ne pourra (une fois l'opération réussie) qu'isoler un polynôme qui *semble* avoir ce nombre réel comme racine. Il serait un peu long de donner tous les détails qui font qu'aujourd'hui ce problème est pour ainsi dire *numériquement résolu.* Mentionnons cependant qu'au moins trois programmes de calcul symbolique ont implanté cet algorithme: soit Maple, Mathematica et Pari-GP. Pour la description de cet algorithme, on pourra consulter [BaKa] ou l'article original de [LLL]. La meilleure version de cet algorithme et de loin la plus rapide est celle qui existe sur Pari-GP [Pari]; elle est au moins 800 fois plus rapide que la version équivalente sur Maple. Quant à Mathematica, disons qu'il est, de façon générale, 4 fois plus lent que Maple dans tous les calculs. Nous ne l'avons pas considéré ici.

Cette procédure accepte donc en entrée un nombre décimal et donne (selon la précision numérique en vigueur) le polynôme minimal dont il serait racine. La précision numérique en vigueur est celle que l'utilisateur demande. Elle devrait idéalement être infinie. Plus raisonnablement, la limite est d'environ 100 chiffres décimaux sur les machines à notre disposition avec Maple et d'environ 500 chiffres décimaux avec Pari-GP. Le degré maximal du polynôme que l'on puisse demander dépend largement de la précision. Dans la pratique, la limite est un polynôme de degré 20. Ceci est quand même suffisant pour obtenir des résultats intéressants.

Evidemment, si la fonction génératrice close qui représente S(z) est algébrique et si z=1/m est un nombre rationnel, le résultat, S(1/m) sera alors un nombre algébrique. C'est précisément ici que l'on utilise l'algorithme LLL. Les centaines de termes que nous donnent la P-récurrence serviront pour évaluer S(z) en un point 1/m "très petit", de telle sorte que le résultat soit un nombre algébrique approché à une grande précision numérique. On ira ensuite chercher avec celui-ci le polynôme dont S(1/m) est racine. Une fois le polynôme candidat trouvé, on réévalue la série S(z) en un autre point rationnel ,1/(m+1), et on répète l'appel à l'algorithme. Il se trouve que la version de LLL sur le programme Pari-GP est extrêmement efficace. Non seulement la procédure (qui s'appelle "algdep") retourne en général le bon polynôme, mais de surcroît il est simplifié au maximum. De plus, les solutions trouvées sont *stables*; elles sont stables au point qu'elles permettent de reconstruire la fonction génératrice algébrique. Une fois ces solutions trouvées en fait, la reconstruction de la fonction génératrice se résume à un calcul d'interpolation assez simple. Comme on l'a mentionné plus tôt, l'appel de la procédure demande un nombre décimal et un degré. Pour arrêter notre choix sur le bon polynôme, il nous suffit de rejeter ceux dont la taille est trop grande (en nombre de caractères). Nous utilisons le même critère que notre méthode des approximants de Padé.

La procédure est la suivante, avec en entrée une suite de la table:

**1)** On teste si la suite est P-récurrente. Si oui on passe à l'étape 2), sinon on arrête.
**2)** On calcule plusieurs centaines de termes de la récurrence (dans la pratique 200 termes suffisent).
**3)** On construit une série S(z) avec ces 200 termes.
**4)** On évalue la série S(z) en des points rationnels 1/m, 1/(m+1),1/(m+2),... . En pratique m=100 et le nombre de termes = 12.
**5)** On appelle la procédure "algdep" de Pari-GP avec les 12 valeurs trouvées.
**6)** On teste avec des polynômes de degré 2,3,4,..., (dans la pratique les degrés 2 à 8 sont suffisants).
**7)** On récupère les bons polynômes, on pose la variable comme étant x.
**8)** On identifie les coefficients de même degré et on calcule le polynôme d'interpolation en t avec la méthode de Newton.
**9)** On substitue t=1/z dans l'expression trouvée.
**10)** On résout (si le degré de l'expression le permet).

Cet algorithme, quoique très technique, fonctionne très rapidement. Il nous a permis de trouver 32 fonctions génératrices algébriques de degré et de complexité assez élevés.

Illustrons cet algorithme en donnant un exemple.

**Exemple 2.3** La suite N0768 des cartes planaires.

Cette suite porte le nom de “Rooted Maps” dans [Sl] mais le titre a été modifié dans [PlSl]. Avec l’étape 1) de notre algorithme, on trouve que la suite satisfait la P-récurrence :

$$(n + 1)\, a_n = (12\, n - 18)\, a_{n-1}.$$

C’est une P-récurrence candidate pour notre méthode hypergéométrique plutôt que pour l’algorithme LLL. On procède donc avec celle-ci et il s’avère que c’est une hypergéométrique:

$${}_2F_1\,([1, 1/2], [3], 12\, z).$$

En demandant à Maple de la simplifier avec “simplify”, celui-ci retourne effectivement une expression algébrique.

Mais cette expression n’est pas très élégante:

$$-1/9\, \frac{\left(1 - 12\, z + 24\, z\, I\, (12\, z - 1)^{1/2} - I\, (12\, z - 1)^{1/2}\right) I}{z\, \left(1 + I\, (12\, z - 1)^{1/2}\right) (12\, z - 1)^{1/2}}$$

On voudrait avoir une expression sans valeurs complexes et simplifiée que l’on obtiendrait de façon automatique. On peut toujours la manipuler à la main, mais notre but est d’obtenir une forme close *automatiquement.* On essaye donc avec une autre méthode: la dérivée et les approximants de Padé.

En dérivant S(z) on obtient une expression qui, mise sous forme d’approximant de Padé, nous donne : (une fois factorisée).

$$-2\, \frac{(81\, z^4 - 648\, z^3 + 234\, z^2 - 27\, z + 1)\,(9\, z^2 - 9\, z + 1)}{(9\, z - 1)\,(27\, z^3 - 81\, z^2 + 18\, z - 1)\,(81\, z^3 - 81\, z^2 + 18\, z - 1)}$$

Si on intègre par rapport à z, on devrait retrouver l’expression, mais il y a des polynômes qui sont du 4è degré et la solution n’est pas élégante non plus.

On s’en remet donc à notre méthode LLL.

**(étape 1)** On reprend la P-récurrence et on recalcule la suite mais avec 200 termes.

**(étape 2 et 3)**. On réévalue la série avec ces mêmes 200 termes et nos points d’interpolation 1/(m+i) avec i=0..4 (5 points d’interpolation devraient suffire)

**(étape 4)**. La première valeur, en m=100, nous donne le premier nombre réel à tester, soit :

1.0209580979488151117686851821900121080607759630492109323339875590733954378833687001578416494132577448905329282269472068...

**(étape 5 et 6)** En appelant la procédure “algdep” avec ce nombre réel bon à 118 décimales et un polynôme de degré 2, on obtient, pour les valeurs 1/m, 1/(m+1),1/(m+2),1/(m+3),1/(m+4)

**(étape 7)** On récupère les bons polynômes:

$$27 x^2 + 8200 x - 8400$$

$$27 x^2 + 8383 x - 8585$$

$$27 x^2 + 8568 x - 8772$$

$$27 x^2 + 8755 x - 8961$$

$$27 x^2 + 8944 x - 9152$$

**(étape 8)** Il nous reste à identifier les coefficients de même degré et à calculer les polynômes d’interpolation correspondants. On aura: pour le coefficient de $x^2$, les valeurs 27,27,27, ... ,. pour le coefficient de x, les valeurs 8200, 8383, 8568, 8755 et 8944, aux points d’interpolation 100,101,102,103 et 104. Enfin, pour le coefficient constant, on aura les valeurs -8400, -8585, -8772, -8961 et -9152 aux mêmes points d’interpolation. On applique alors simplement la formule d’interpolation de Newton pour trouver une expression polynômiale pour chaque degré. On peut faire appel à la procédure de la librairie Maple appelée “interp” qui effectue ce calcul automatiquement. Ce qui nous donnera deux variables, x et t.

**(étape 9)** Il restera à substituer t=1/z. On obtient finalement:

$$\frac{-1 + 16 z + x - 18 x z + 27 x^2 z^2}{z^2}$$

**(étape 10)** Il ne reste qu'à résoudre cette équation par rapport à x: on prendra alors la solution positive. Finalement l'expression algébrique de notre suite de départ serait :

$$\frac{1}{54}\,\frac{-1+18z+\left(-(12z-1)^{3}\right)^{1/2}}{z^{2}}$$

C'est l'expression la plus simple qu'on ait obtenu pour cette suite. La magie de cet algorithme LLL est qu'il trouve une expression polynômiale pour un nombre réel qui est en général minimale. Des expressions de plus haut degré encore ont été obtenues de cette façon, la plus grosse étant de degré 8 et elles sont répertoriées dans notre table en appendice.

# CHAPITRE 3

# LA MÉTHODE D'EULER

Ainsi nommée parce qu'elle semble avoir été développée à l'époque d'Euler. Nous n'avons pas trouvé de références historiques sur cette méthode, bien que Andrews [And] la mentionne.

L'idée en est simple: étant donné une suite $a_n$ dont on suppose la série génératrice de la forme,

(3.1) $$S(z) = 1 + \sum_{n=1}^{\infty} a_n z^n = \prod_{n=1}^{\infty} (1 - z^n)^{-c_n},$$

la question est : comment trouver les $c_n$ en fonction des $a_n$. Comme l'explique Andrews à la page 104, il suffit d'utiliser la formule d'inversion de Möbius. En effet, puisque le membre de droite de (3.1) est un produit infini, c'est en prenant le logarithme ou la dérivée logarithmique que nous retrouvons alors une somme ordinaire. En identifiant le coefficient de degré n (pour exprimer chaque coefficient de $a_n$) et en inversant (par Möbius) par rapport à la somme, nous obtenons les coefficients $c_n$ en fonction des $a_n$. La somme s'exprimera en termes des diviseurs de n. Inversement, si nous connaissons les $c_n$ et que l'on cherche les $a_n$ , l'opération est directe; il suffit de développer le produit en série. En prenant soin de garder le même ordre de grandeur des séries correspondantes, nous obtenons le même nombre de termes pour les $c_n$ que pour les $a_n$. Autrement dit, si les k premiers coefficients de $a_n$ sont connus, il y aura alors k coefficients de bons pour les $c_n$.

On peut donc programmer la transformation dans les deux sens en une vingtaine de lignes. La procédure accepte en entrée une suite et donne du même coup une représentation en "partages", c'est-à-dire qu'elle propose un produit infini. Par exemple, la suite N0244 énumère les partages ordinaires de l'entier n. En effectuant le calcul on trouve la suite 1,1,1,1,1,.... C'est la forme de produit infini de ce type la plus simple. Mais pour détecter un bon candidat de produit infini avec ce type de fonction génératrice, dans un cadre plus général, nous avons utilisé la méthode des approximants de Padé qui permet de détecter les “motifs” dans les exposants. En tout, 94 produits infinis ont ainsi été isolés grâce à cette méthode. Les résultats sont présentés dans la table en appendice.

# CHAPITRE 4

# LA MÉTHODE DES RECOUPEMENTS

## 4.1 Les recoupements indirects.

L'hypothèse que l'on pose ici est que la suite dont on cherche la fonction génératrice est en fait une suite connue mais transformée. Par exemple, la suite des partages d'entiers N0244 de [Sl] est très facile à détecter. Il suffit de prendre la méthode d'Euler et le programme nous propose immédiatement un produit infini très simple. Mais si on effectue la translation $a_n+3$, le programme ne détectera pas ce produit infini. Pour une bonne raison car, si la suite ne commence pas naturellement à 1, alors l'opération d'Euler n'est pas valide et même si on l'effectue, les termes seront des nombres rationnels (non entiers). Donc afin de pouvoir isoler le plus possible de suites, on se sert, comme base de comparaison, de la table des suites qui en contient 4568. En prenant chaque suite transformée de façon élémentaire, on compare avec la table afin de voir s'il n'y aurait pas un croisement. En tout, nous avons répertorié 97 transformations élémentaires d'une suite susceptibles de se retrouver dans la table, soit 54 transformations avec la suite sous forme de série ordinaire et 43 avec la suite sous forme de série exponentielle.

Il serait fastidieux de les énumérer toutes, mais en voici quelques unes. Avec S(z) : la suite sous forme de série ordinaire par exemple, nous avons S(z) + cz/(1-z) ou c=±1,±2,±3, 1/S(z), $S(z)^2$, $S(z)^3$, S(z)/(1-z) ce qui équivaut à considérer la suite des sommes partielles de la suite. D'autres transformations sont plus simples encore, comme $\mathbf{N} \setminus \{a_n\}$, la différence ensembliste des entiers et de la suite. On ne tient pas compte ici de la multiplicité des termes. On a considéré aussi de prendre $a_n/\mathrm{pgcd}(a_0, a_1, a_2,... a_k)$ ou de prendre la suite avec les indices de rang pairs et impairs. L'idée est de prendre des transformations les plus simples possibles. La transformation d'Euler dans les 2 sens complète la liste.

On pourrait les classer en ces quelques catégories :

1) Translations : S(z) +/- cz/(1-z), avec c=1,2,3.

2) Inverses : 1/S(z),1/S(z)$^2$ ,1/S(z)$^3$.

3) Puissances : S(z)$^k$ avec k=1,2,3.

4) Sommes et différences.

5) Transformation de type Euler ( voir chapitre 3).

6) Transformations de type ensembliste comme **N** \ {$a_n$}.

7) Transformations avec le p.g.c.d. .

Les autres sont données en considérant des combinaisons de ces dernières.

Par exemple, de la suite N0577 de [Sl] (les nombres de Catalan), on en obtient 97 autres et en comparant ces 97 suites avec la table, 6 autres suites au moins seraient liées à cette dernière. C'est donc que, si on connait déjà la fonction génératrice des nombres de Catalan obtenue avec d'autres méthodes, alors du même coup on obtient la fonction génératrice de ces 6 autres suites. C'est un avantage, parce que justement avec cet exemple, si l'on prend $a_n$-1 et que l'on compare avec la table, on retrouve la suite N1409 de [Sl]. Cette suite n'est pas hypergéométrique en vertu d'un critère assez simple de [GKP], elle ne commence pas par 1. De plus son inverse fonctionnel est impossible à effectuer pour le même genre de raisons; le premier terme est nul mais le deuxième terme n'est pas 1. Elle est cependant algébrique et c'est avec la méthode LLL (beaucoup plus lourde) que la fonction génératrice a été trouvée. En fait, elle est évidemment de la forme S(z) -1/(1-z) où S(z) est la fonction génératrice des nombres de Catalan. Mais ceci constitue un raisonnement a posteriori. Donc cette méthode des recoupements peut mener à des résultats très intéressants en autant que le traitement informatique des 97 transformations appliquées aux 4568 suites et comparées avec ces dernières à chaque fois ne soit pas trop lourd également.

Un détail ne doit cependant pas être oublié. La comparaison de 2 suites entre elles peut mener à des erreurs. On doit faire la comparaison à partir du deuxième terme, parce que souvent la suite est répertoriée mais les premiers termes peuvent être d'indices 0 ou 1. C'est-à-dire que la suite ne débute pas au terme de rang 0. Également on ne doit pas prendre toute la suite: il ne faut pas oublier que certaines suites de la table sont très courtes et ne contiennent que quelques termes. Elles ne sont pas moins importantes, par exemple la suite N0323 de [Sl]. Il y a un juste milieu et l'expérience montre que les indices de rang 2 à 16 sont suffisants, c’est-à-dire les 15 premiers termes de la suite à partir du rang 2.

A cet effet un programme appelé HIS (Handbook of Integer Sequences) a été mis au point. Il n’est cependant pas public comme le programme gfun. Dans HIS se trouve la table numérique des suites et 2 procédures appelées “find” et “findhard”. La première sert simplement à savoir si une suite se trouve dans la table et la deuxième fait une recherche dans la table après avoir effectué les 97 transformations en question. Le programme et la table sont entièrement contenus en Maple. Le programme est donc de cette façon transportable sur toute machine qui peut recevoir Maple.

La procédure “find” qui en principe ne fait que regarder si une suite se trouve dans la table emploie une procédure de recherche mise au point par Bruno Salvy de l’INRIA. Il était essentiel d’avoir à notre disposition un algorithme de recherche qui soit très rapide étant donné le nombre important de comparaisons à chaque opération. Une structure de données adaptée à ces besoins a été construite sous forme d’arbre binaire. En effectuant une boucle de calcul sur toute la table avec la procédure “findhard”, une banque de données des croisements a été obtenue. En tout il y aurait 3800 croisements. Une proportion appréciable des ces croisements, soit environ 25% selon nous, est fortuite ou accidentelle. Ceci est relié à la décision de ne prendre qu’une partie de chaque suite pour comparer. Donc pour pouvoir retrouver la fonction génératrice, il y a un travail de vérification nécessaire.

Ce travail de vérification est très long, mais il en vaut la peine. Evidemment, beaucoup de croisements ne sont pas surprenants: à titre d’exemple, la suite de Fibonacci, qui est très connue et qui a une fonction génératrice assez simple croise avec une bonne centaine de suites. Aucun de ces croisements n’est vraiment nouveau. C’est lorsque la suite est intrinsèquement plus complexe que le jeu en vaut la chandelle. Sans exagérer, nous avons effectué patiemment des centaines d’heures de calcul et de vérification pour trouver ces résultats et il y en a beaucoup à faire encore, puisque cette table des 3000 bons croisements environ n’a pas été passée en revue au complet. Par ce procédé, 38 fonctions génératrices ont été obtenues. Elles sont répertoriées dans la table en appendice.

## 4.2 Les tableaux.

Le programme Maple manipule des données numériquement aussi bien que symboliquement. La procédure “ratpoly” est capable de trouver une fraction rationnelle de séries à une variable aussi bien qu’à 2 variables comme les tableaux à 2 dimensions. Un bon exemple est le triangle de Pascal. Il suffit de le mettre sous forme de tableau “carré” où chaque rangée sera un polynôme. En prenant les 5 premières rangées, on aura la suite

$$1,\ 1 + t,\ 1 + 2\,t + t^2,\ 1 + 3\,t + 3\,t^2 + t^3,\ 1 + 4\,t + 6\,t^2 + 4\,t^3 + t^4\,.$$

Si cette suite (de polynômes) est maintenant convertie en série de puissances en z et passée à la procédure “ratpoly”, elle retourne immédiatement
$1/(1 - t\,z - z)$. Si on développe en série par rapport à t, on obtient la fonction génératrice de chaque colonne et inversement, en développant par rapport à z, on obtient la fonction génératrice de chaque rangée (qui sont ici des polynômes). Notons que 4 termes suffisent pour trouver la fonction génératrice du tableau.

Contrairement aux autres méthodes, il n’existe pas de livre ou de catalogue de tels tableaux. Il y en a un bon nombre dans la littérature, mais ce qui a été fait plutôt est d’en générer de façon *ad hoc* . Il faut prendre un modèle de tableaux assez général, par exemple dans [GKP] ou [Théo], où on introduit les tableaux $A_{[n,k]}$ définis par la relation de récurrence

$$A_{[n+1,k+1]} = (r\,n+s\,k+t)\,A_{[n,k+1]} + (a\,n+b\,k+c)\,A_{[n,k]}$$

où a,b,c,r,s et t sont entiers. Il se trouve qu’une bonne partie des tableaux étudiés en combinatoire sont de ce type: les coefficients binomiaux, les nombres de Stirling de 1ère et de 2ème espèce, les nombres eulériens, les coefficients des polynômes de Tchébycheff, etc. On peut consulter [Théo] à ce sujet où une étude approfondie de ces tableaux a été menée. Il reste donc à en générer un bon nombre en prenant les entiers a,b,c,r,s et t compris entre -4 et 4 et de *tenter* de trouver la fonction génératrice. Sur des milliers tableaux générés de cette façon, 430 fonctions génératrices à deux variables ont été trouvées, couvrant la plupart des cas simples de ces tableaux. On obtient ainsi un échantillonnage assez important de formules, suffisamment important pour y trouver la fonction génératrice de centaines de suites de notre table. En tout, 20 nouvelles fonctions génératrices ont été isolées. Ces résultats sont présentés dans la table en appendice.

# CONCLUSION

On conclut que nos méthodes peuvent dans 23 % des cas donner la fonction génératrice d'une suite d'entiers "quelconque". Le mot quelconque signifie ici: ce qui est catalogué dans la table de suites [PlSl]. Nous croyons qu'il en est de même avec toute suite d'entiers qui se présente au mathématicien dans ses recherches, quel que soit son domaine. Nous souhaitons que ces méthodes deviennent des outils de travail.

Il reste cependant beaucoup à faire. Il faut trouver une explication raisonnable au fait que nous sommes passés à côté de 77% des suites. On pourrait peut-être étendre encore les méthodes en formulant d'autres modèles de fonctions génératrices. En fait, il en existe déjà. Par exemple, la fonction "plancher" ou partie entière permet de construire des suites très simples que nos méthodes n'ont pas détectées; la suite $[(3/2)^n]$ en est un bon représentant. On pourrait également mettre dans la même catégorie les suites définies avec des nombres irrationnels comme $[\sqrt{2}n]$. Un autre modèle pourrait être basé sur les récurrences quadratiques comme la suite 2,4,16,256,... (en mettant au carré à chaque fois). Elle est extrèmement simple mais indétectable par nos méthodes. Un autre serait basé sur les suites "doublement" récurrentes, là où il y a une fonction de l'indice comme $a_{a(n)}$. On pourrait multiplier les exemples de suites très simplement définies mais indétectables. Ce qui caractérise une table de suites comme [Sl] ou [PlSl], c'est la variété et c'est précisément ce qui nous passionne.

# BIBLIOGRAPHIE

# A.0 NOTES À L'UTILISATEUR DE LA TABLE

Chaque fonction génératrice trouvée à l'aide de l'une de nos méthodes est répertoriée dans la table qui suit sous forme de fiche. Chaque fiche contient les informations pertinentes à cette suite :

- Numéro séquentiel Axxxx et Nxxxx (s'il existe)
- Nom de la suite
- Les références bibliographiques avec dans l'ordre : Périodique  Volume Page Année
- La méthode employée pour trouver la fonction génératrice
- Le type de fonction génératrice
- Commentaires additionnels
- Autres formules connues ou trouvées
- La fonction génératrice
- La suite numérique

Elles apparaissent selon le schéma suivant:

<table>
<tr><td colspan="3">Nom de la suite</td></tr>
<tr><td colspan="3">Références</td></tr>
<tr><td>Numéro Axxxx</td><td>Méthode employée</td><td rowspan="2">Commentaires</td></tr>
<tr><td>Numéro Nxxxx</td><td>Type de fonction génératrice</td></tr>
<tr><td colspan="3">Autre formules</td></tr>
<tr><td colspan="3">Fonction génératrice</td></tr>
<tr><td colspan="3">Suite numérique</td></tr>
</table>

• Les références bibliographiques sont notées exactement comme dans le livre [Sl]. La liste des ouvrages se trouve dans une bibliographie séparée à la fin de la table.

• La fonction génératrice qui apparaît au centre est toujours une fonction génératrice ordinaire à moins qu'il en soit indiqué autrement (exponentielle ou double exponentielle).

• Les fiches ont été triées par ordre numérique sur les numéros Axxxx. Cette table est une pile Hypercard. On peut donc l'utiliser sur tout ordinateur Macintosh et la consulter comme une banque de donnée. Nous prévoyons un accès à Maple. De cette façon l'utilisateur pourra vérifier chaque formule.

• W(z) désigne la fonction Oméga, définie implicitement par W(z) exp(W(z)) = z. On la connait aussi sous

sa forme de série exponentielle dont les coefficients sont donnés, en valeur absolue, $|a_0| = 0$, $|a_n| = n^{n-1}$ pour n>0 (série alternante à terme constant nul). Son rayon de convergence est 1/e et elle est souvent utilisée pour le développement en série de certaines fonctions génératrices de structures arborescentes. Elle est très commode dans les calculs.

• La fonction génératrice qui apparaît au centre de chaque fiche est la plus simple ou plus élégante expression que nous connaissons donnant les termes de la suite.

• Le nom de chaque suite (s'il est présent) est tel qu'il apparaît dans [PlSl]. Quand il est omis c'est qu'il est d'une forme que nous jugeons redondante par rapport à la fonction génératrice.

• Les P-récurrences qui apparaissent dans la case "fonction génératrice" ou "autres formules" ont leur conditions initiales données par les premiers termes de la suite.

• La suite qui apparaît dans la case "suite numérique" est telle qu'elle apparaît dans [PlSl], plus determes peuvent être évidemment obtenus avec la fonction génératrice.

• Certaines fiches ont été imprimées en format pleine grandeur pour plus de lisibilité

# 1031 Generating Functions

*par*
*Simon Plouffe*
*August 1992*
*found using GFUN and other tools with a sample of the Encyclopedia of Integer Sequences (as of 1992)*
*with 4568 sequences.*

## Denumerants

**Réf.** R1 152.

**HIS2 A0008** **Euler** **erreur au 19è terme corrigée avec la formule**
**HIS1 N0099** **Fraction rationnelle**

$$\frac{1}{(1 - z)\,(1 - z^{2})\,(1 - z^{5})\,(1 - z^{10})}$$

**1, 1, 2, 2, 3, 4, 5, 6, 7, 8, 11, 12, 15, 16, 19, 22, 25, 28, 34, 40**

## Partitions n into distinct parts

**Réf.** AS1 836.

**HIS2 A0009** **Euler**

**HIS1 N0100** **Produit infini**

$$\prod_{n\geq 0} (1 - z^{2n+1})$$

**1, 1, 1, 2, 2, 3, 4, 5, 6, 8, 10, 12, 15, 18, 22, 27, 32, 38, 46, 54, 64, 76, 89, 104, 122, 142, 165, 192, 222, 256, 296, 340, 390, 448, 512, 585, 668, 760, 864, 982, 1113, 1260, 1426**

## Related to Latin Rectangles

**Réf.** R1 210.

**HIS2 A0023** **Recoupements** **Suite P-récurrente**

**HIS1 N0140** **exponentielle (rationnelle)**

**a(n) = (3 n - 1) a(n - 1) + (- 4 n + 2) a(n - 2)**

$$\frac{1}{\exp(2z)\,(1-z)}$$

**1, 1, 2, 2, 8, 8, 112, 656, 5504, 49024, 491264**

## The natural numbers

**Réf.**

**HIS2 A0027** **Approximants de Padé**

**HIS1 N0173** **Fraction rationnelle**

$$\frac{1}{(1 - z)^2}$$

**1, 2, 3, 4, 5, 6, 7, 8, 9, 10, 11, 12, 13, 14, 15, 16, 17, 18, 19, 20, 21, 22, 23, 24, 25, 26, 27, 28, 29, 30, 31, 32, 33, 34, 35, 36, 37, 38, 39, 40, 41, 42, 43, 44, 45, 46, 47, 48, 49**

## Partitions of n

**Réf.** RS4 90. R1 122. AS1 836.

**HIS2 A0041** **Euler**

**HIS1 N0244** **Produit infini**

$$\prod_{n \geq 1} \frac{1}{(1 - z^n)}$$

**1, 1, 2, 3, 5, 7, 11, 15, 22, 30, 42, 56, 77, 101, 135, 176, 231, 297, 385, 490, 627, 792, 1002, 1255, 1575, 1958, 2436, 3010, 3718, 4565, 5604, 6842, 8349, 10143, 12310, 14883**

## Dying Rabbits

**Réf.** FQ 2 108 64.

**HIS2 A0044** **Approximants de Padé**

**HIS1 N0255** **Fraction rationnelle**

**a(n+13)=a(n+12)+a(n+11)+a(n)**

$$\frac{1 + z^{2} + z^{4} + z^{6} + z^{8} + z^{10}}{1 - z - z^{3} - z^{5} - z^{7} - z^{9} - z^{11}}$$

**1, 1, 2, 3, 5, 8, 13, 21, 34, 55, 89, 144, 232, 375, 606, 979, 1582, 2556, 4130, 6673, 10782, 17421, 28148, 45480, 73484, 118732, 191841, 309967, 500829, 809214, 1307487**

## Fibonacci numbers

**Réf.** HW1 148. HO69.

**HIS2 A0045** **Approximants de Padé**

**HIS1 N0256** **Fraction rationnelle**

$$\frac{1}{1 - z - z^{2}}$$

**1, 1, 2, 3, 5, 8, 13, 21, 34, 55, 89, 144, 233, 377, 610, 987, 1597, 2584, 4181, 6765, 10946, 17711, 28657, 46368, 75025, 121393, 196418, 317811, 514229, 832040, 1346269**

## 2 ^ n + 1

**Réf.** BA9.

**HIS2 A0051** **Approximants de Padé**

**HIS1 N0266** **Fraction rationnelle**

$$\frac{2 - 3z}{(1 - z)(1 - 2z)}$$

**2, 3, 5, 9, 17, 33, 65, 129, 257, 513, 1025, 2049, 4097, 8193, 16385, 32769, 65537, 131073, 262145, 524289, 1048577, 2097153, 4194305, 8388609, 16777217**

---

## Denumerants

**Réf.** R1 152.

**HIS2 A0064** **Euler** **erreur au 19è terme corrigée avec la formule**

**HIS1 N0375** **Fraction rationnelle**

$$\frac{1}{(1 - z)^2 (1 - z^2)(1 - z^5)(1 - z^{10})}$$

**1, 2, 4, 6, 9, 13, 18, 24, 31, 39, 50, 62, 77, 93, 112, 134, 159, 187, 252, 292**

## n-node trees of height 2

**Réf.** IBMJ 4 475 60. KU64.

**HIS2 A0065** **Euler**

**HIS1 N0379** **Produit infini**

$$\frac{z}{(1-z)} + \prod_{n\geq 1} \frac{1}{(1-z^n)}$$

**1, 2, 4, 6, 10, 14, 21, 29, 41, 55, 76, 100, 134, 175, 230, 296, 384, 489, 626, 791, 1001, 1254, 1574, 1957, 2435, 3009, 3717, 4564, 5603, 6841, 8348, 10142, 12309**

## Partitions of n into parts of 2 kinds

**Réf.** RS4 90. RCI 199. FQ 9 332 71.

**HIS2 A0070** **Euler**

**HIS1 N0396** **Produit infini**

$$\prod_{n\geq 1} \frac{1}{(1-z^n)^{c(n)}}$$

$$c(n) = 2,1,1,1,1,...$$

**1, 2, 4, 7, 12, 19, 30, 45, 67, 97, 139, 195, 272, 373, 508, 684, 915, 1212, 1597, 2087, 2714, 3506, 4508, 5763, 7338, 9296, 11732, 14742, 18460, 23025, 28629, 35471**

## Fibonacci numbers - 1

**Réf.** R1 155. AENS 79 203 62. FQ 3 295 65.

**HIS2 A0071** **Approximants de Padé**

**HIS1 N0397** **Fraction rationnelle**

$$\frac{1}{1 - 2z + z^3}$$

**1, 2, 4, 7, 12, 20, 33, 54, 88, 143, 232, 376, 609, 986, 1596, 2583, 4180, 6764, 10945, 17710, 28656, 46367, 75024, 121392, 196417, 317810, 514228, 832039, 1346268**

## Tribonacci numbers

**Réf.** FQ 1(3) 71 63; 5 211 67.

**HIS2 A0073** **Approximants de Padé**

**HIS1 N0406** **Fraction rationnelle**

$$\frac{z}{1 - z - z^2 - z^3}$$

**0, 1, 1, 2, 4, 7, 13, 24, 44, 81, 149, 274, 504, 927, 1705, 3136, 5768, 10609, 19513, 35890, 66012, 121415, 223317, 410744, 755476, 1389537, 2555757, 4700770, 8646064**

## Tetranacci numbers

**Réf.** AMM 33 232 26. FQ 1(3) 74 63.

**HIS2 A0078** **Approximants de Padé**

**HIS1 N0423** **Fraction rationnelle**

$$\frac{1}{1 - z - z^2 - z^3 - z^4}$$

**1, 1, 2, 4, 8, 15, 29, 56, 108, 208, 401, 773, 1490, 2872, 5536, 10671, 20569, 39648, 76424, 147312, 283953, 547337, 1055026, 2033628, 3919944, 7555935, 14564533**

## Powers of 2

**Réf.** BA9. MOC 23 456 69.

**HIS2 A0079** **Approximants de Padé**

**HIS1 N0432** **Fraction rationnelle**

$$\frac{1}{1 - 2z}$$

**1, 2, 4, 8, 16, 32, 64, 128, 256, 512, 1024, 2048, 4096, 8192, 16384, 32768, 65536, 131072, 262144, 524288, 1048576, 2097152, 4194304, 8388608, 16777216**

## Rooted trees with n nodes

**Réf.** R1 138. HA69 232.

**HIS2 A0081** **Recoupements**

**HIS1 N0454** **Produit infini**

$$\prod_{n\geq 1} \frac{1}{(1 - z^{n})^{c(n)}}$$

**c(n) = a(n) : la suite elle-même.**

**1, 1, 2, 4, 9, 20, 48, 115, 286, 719, 1842, 4766, 12486, 32973, 87811, 235381, 634847, 1721159, 4688676, 12826228, 35221832, 97055181, 268282855, 743724984**

---

**Réf.** LU91 1 221. R1 86. MU60 6. DMJ 35 659 68.

**HIS2 A0085** **Dérivée logarithmique** **Suite P-récurrente**

**HIS1 N0469** **exponentielle**

**a(n) = a(n - 1) + (n - 1) a(n - 2)**

$$\exp(z + 1/2\ z^{2})$$

**1, 1, 2, 4, 10, 26, 76, 232, 764, 2620, 9496, 35696, 140152, 568504, 2390480, 10349536, 46206736, 211799312, 997313824, 4809701440, 23758664096**

## Permutations with no cycles of length 3

**Réf.** R1 85.

**HIS2 A0090** **Dérivée logarithmique** **Suite P-récurrente**

**HIS1 N0496** **exponentielle**

**a(n) = (n^3-n^2)a(n-1)+(6 n^3-5 n^2+n)a(n-3)+(24 n^3-26 n^2+9 n-1)a(n-4)**

$$\frac{1}{\exp(1/3\, z^3)\,(1 - z)}$$

**1, 1, 2, 4, 16, 80, 520, 3640, 29120, 259840, 2598400, 28582400, 343235200, 4462057600, 62468806400, 936987251200, 14991796019200, 254860532326400, 4587501779660800**

---

**Réf.** AS1 797.

**HIS2 A0096** **Approximants de Padé**

**HIS1 N0522** **Fraction rationnelle**

$$\frac{z\,(z - 2)}{(z - 1)^3}$$

**0, 2, 5, 9, 14, 20, 27, 35, 44, 54, 65, 77, 90, 104, 119, 135, 152, 170, 189, 209, 230, 252, 275, 299, 324, 350, 377, 405, 434, 464, 495, 527, 560, 594, 629, 665, 702, 740, 779**

## Partitions of n into parts of 2 kinds

**Réf.** RS4 90. RCI 199.

**HIS2 A0097** **Euler**

**HIS1 N0525** **Produit infini**

$$\prod_{n\geq 1} \frac{1}{(1 - z^n)^{c(n)}}$$

$$c(n) = 2,2,1,1,1,1,1,1,...$$

**1, 2, 5, 9, 17, 28, 47, 73, 114, 170, 253, 365, 525, 738, 1033, 1422, 1948, 2634, 3545, 4721, 6259, 8227, 10767, 13990, 18105, 23286, 29837, 38028, 48297, 61053**

---

## Partitions of n into parts of 2 kinds

**Réf.** RS4 90. RCI 199.

**HIS2 A0098** **Euler**

**HIS1 N0533** **Produit infini**

$$\prod_{n\geq 1} \frac{1}{(1 - z^n)^{c(n)}}$$

$$c(n) = 2,2,2,1,1,1,1,1,1,...$$

**1, 2, 5, 10, 19, 33, 57, 92, 147, 227, 345, 512, 752, 1083, 1545, 2174, 3031, 4179, 5719, 7752, 10438, 13946, 18519, 24428, 32051, 41805, 54265, 70079, 90102, 115318**

## Compositions

**Réf.** R1 155.

**HIS2 A0100** **Approximants de Padé**

**HIS1 N0543** **Fraction rationnelle**

$$\frac{1}{(1 - z - z^2)(1 - z - z^2 - z^3)}$$

**1, 2, 5, 11, 23, 47, 94, 185**

## Compositions

**Réf.** R1 155.

**HIS2 A0102** **Approximants de Padé**

**HIS1 N0551** **Fraction rationnelle**

$$\frac{1}{(1 - z - z^2 - z^3)(1 - z - z^2 - z^3 - z^4)}$$

**1, 2, 5, 12, 27, 59, 127**

## Catalan's Numbers

**Réf.** AMM 72 973 65. RCI 101. C1 53. PLC 2 109 71. MAG 61 211 88.

**HIS2 A0108** **Inverse fonctionnel** **Suite P-récurrente**

**HIS1 N0577** **algébrique**

**2F1 ([1, 1/2], [2], 4 z)**

**n a(n) = (4 n - 6) a(n - 1)**

$$\frac{2}{1 + (1 - 4 z)^{1/2}}$$

**1, 1, 2, 5, 14, 42, 132, 429, 1430, 4862, 16796, 58786, 208012, 742900, 2674440, 9694845, 35357670, 129644790, 477638700, 1767263190, 6564120420, 24466267020**

## Bell Numbers

**Réf.** MOC 16 418 62. AMM 71 498 64. PSPM 19 172 71. GO71.

**HIS2 A0110** **Recoupements**

**HIS1 N0585** **exponentielle**

$$\exp(\exp(z) - 1)$$

**1, 1, 2, 5, 15, 52, 203, 877, 4140, 21147, 115975, 678570, 4213597, 27644437, 190899322, 1382958545, 10480142147, 82864869804, 682076806159, 5832742205057**

## Euler numbers

**Réf.** JDM 7 171 1881. JO61 238. NET 110. DKB 262. C1 259.

**HIS2 A0111** **Inverse fonctionnel**
**HIS1 N0587** **exponentielle (complexe)**

$$\tan(1/4\, \pi + 1/2\, z) - 1$$

**1, 1, 1, 2, 5, 16, 61, 272, 1385, 7936, 50521, 353792, 2702765, 22368256, 199360981, 1903757312, 19391512145, 209865342976, 2404879675441, 29088885112832**

---

## Denumerants

**Réf.** R1 152.

**HIS2 A0115** **Euler** **erreur au 19è terme corrigée avec la formule**
**HIS1 N0098** **Fraction rationnelle**

$$\frac{1}{(1 - z)(1 - z^2)(1 - z^5)}$$

**1, 1, 2, 2, 3, 4, 5, 6, 7, 8, 10, 11, 13, 14, 16, 18, 20, 22, 26, 29**

## Representations of n as a sum of distinct Fibonaccis

**Réf.** FQ 4 305 66. BR72 54.

**HIS2 A0119** **Euler**

**HIS1 N0037** **Produit infini**

$$\prod_{n \geq 1} (1 + z^{c(n)})$$

**c(n) = 1,2,3,5,8,... nombres de Fibonacci**

**1, 1, 1, 2, 1, 2, 2, 1, 3, 2, 2, 3, 1, 3, 3, 2, 4, 2, 3, 3, 1, 4, 3, 3, 5, 2, 4, 4, 2, 5, 3, 3, 4, 1, 4, 4, 3, 6, 3, 5, 5, 2, 6, 4, 4, 6, 2, 5, 5, 3, 6, 3, 4, 4, 1, 5, 4, 4, 7, 3, 6, 6, 3, 8, 5, 5, 7, 2, 6, 6, 4**

## Representations of n as a sum of Fibonacci numbers

**Réf.** FQ 4 304 66.

**HIS2 A0121** **Euler**

**HIS1 N0088** **Produit infini**

$$(1 + z) \prod_{n \geq 1} (1 + z^{c(n)})$$

**c(n) = 1,2,3,5,8,... nombres de Fibonacci**

**1, 2, 2, 3, 3, 3, 4, 3, 4, 5, 4, 5, 4, 4, 6, 5, 6, 6, 5, 6, 4, 5, 7, 6, 8, 7, 6, 8, 6, 7, 8, 6, 7, 5, 5, 8, 7, 9, 9, 8, 10, 7, 8, 10, 8, 10, 8, 7, 10, 8, 9, 9, 7, 8, 5, 6, 9, 8, 11, 10, 9, 12, 9, 11, 13**

## Binary partitions (partitions of 2n into powers of 2)

**Réf.** FQ 4 117 66. PCPS 66 376 69. AB71 400. BIT 17 387 77.

**HIS2 A0123** **Euler**

**HIS1 N0378** **Produit infini**

$$\frac{1}{(1-z)^2(1-z^2)(1-z^4)(1-z^8)(1-z^{16})(1-z^{32})\ldots}$$

**1, 2, 4, 6, 10, 14, 20, 26, 36, 46, 60, 74, 94, 114, 140, 166, 202, 238, 284, 330, 390, 450, 524, 598, 692, 786, 900, 1014, 1154, 1294, 1460, 1626, 1828, 2030, 2268, 2506**

## Central polygonal numbers

**Réf.** MAG 30 150 46. HO50 22. FQ 3 296 65.

**HIS2 A0124** **Approximants de Padé**

**HIS1 N0391** **Fraction rationnelle**

$$\frac{1-z+z^2}{(1-z)^3}$$

**1, 2, 4, 7, 11, 16, 22, 29, 37, 46, 56, 67, 79, 92, 106, 121, 137, 154, 172, 191, 211, 232, 254, 277, 301, 326, 352, 379, 407, 436, 466, 497, 529, 562, 596, 631, 667, 704, 742**

## Slicing a cake with n slices

**Réf.** MAG 30 150 46. FQ 3 296 65.

**HIS2 A0125** **Approximants de Padé**

**HIS1 N0419** **Fraction rationnelle**

**1+C(n,1)+C(n,2)+C(n,3)**

$$\frac{1 - 2z + 2z^2}{(1 - z)^4}$$

**1, 2, 4, 8, 15, 26, 42, 64, 93, 130, 176, 232, 299, 378, 470, 576, 697, 834, 988, 1160, 1351, 1562, 1794, 2048, 2325, 2626, 2952, 3304, 3683, 4090, 4526, 4992, 5489**

## A nonlinear binomial sum

**Réf.** FQ 3 295 65.

**HIS2 A0126** **Approximants de Padé**

**HIS1 N0421** **Fraction rationnelle**

$$\frac{1 - z + z^3}{(1 - z - z^2)\,(z - 1)^2}$$

**1, 2, 4, 8, 15, 27, 47, 80, 134, 222, 365, 597, 973, 1582, 2568, 4164, 6747, 10927, 17691, 28636, 46346, 75002, 121369, 196393, 317785, 514202, 832012, 1346240**

## C(n,4)+C(n,3)+ ... +C(n,0)

**Réf.** MAG 30 150 46. FQ 3 296 65.

**HIS2 A0127** **Approximants de Padé**

**HIS1 N0427** **Fraction rationnelle**

$$\frac{1 - 3z + 4z^2 - 2z^3 + z^4}{(1 - z)^5}$$

**1, 2, 4, 8, 16, 31, 57, 99, 163, 256, 386, 562, 794, 1093, 1471, 1941, 2517, 3214, 4048, 5036, 6196, 7547, 9109, 10903, 12951, 15276, 17902, 20854, 24158, 27841, 31931**

## A nonlinear binomial sum

**Réf.** FQ 3 295 65.

**HIS2 A0128** **Approximants de Padé**

**HIS1 N0428** **Fraction rationnelle**

$$\frac{1 - 2z + z^2 + z^3}{(1 - z - z^2)(1 - z)^3}$$

**1, 2, 4, 8, 16, 31, 58, 105, 185, 319, 541, 906, 1503, 2476, 4058, 6626, 10790, 17537, 28464, 46155, 74791, 121137, 196139, 317508, 513901, 831686, 1345888**

## Pell numbers

**Réf.** FQ 4 373 66. RI89 43.

**HIS2 A0129** **Approximants de Padé**

**HIS1 N0552** **Fraction rationnelle**

**a(n)=2 a(n-1)+a(n-2)**

$$\frac{1}{1 - 2z - z^2}$$

**1, 2, 5, 12, 29, 70, 169, 408, 985, 2378, 5741, 13860, 33461, 80782, 195025, 470832, 1136689, 2744210, 6625109, 15994428, 38613965, 93222358, 225058681**

---

**Réf.** R1 85.

**HIS2 A0138** **Dérivée logarithmique** **Suite P-réccurente**

**HIS1 N0638** **exponentielle**

**a(n) =(n - 1) a(n - 1) - (n^3 - 9 n^2 + 26 n - 24) a(n - 4) +**
**(n^4 - 14 n^3 + 71 n^2 - 154 n + 120) a(n - 5)**

$$\frac{1}{\exp(1/4\, z^4)\,(1 - z)}$$

**1, 1, 2, 6, 18, 90, 540, 3780, 31500, 283500, 2835000, 31185000, 372972600, 4848643800, 67881013200, 1018215198000, 16294848570000, 277012425690000, 4986223662420000**

**Réf.** CJM 15 257 63. AB71 363.

**HIS2 A0139** **Hypergéométrique** **Suite P-récurrente**

**HIS1 N0651** **algébrique** **équation du 3è degré**

**1/2 (n + 1) (2 n + 1) a(n) = 3/4 (3 n - 1) (3 n - 2) a(n - 1)**

$$_3F_2([1, 4/3, 5/3],[3, 5/2],27 z / 4)$$

**1, 2, 6, 22, 91, 408, 1938, 9614, 49335, 260130, 1402440, 7702632, 42975796, 243035536, 1390594458, 8038677054, 46892282815, 275750636070, 1633292229030, 9737153323590**

---

## Factorial numbers

**Réf.** AS1 833. MOC 24 231 70.

**HIS2 A0142** **Dérivée logarithmique** **Suite P-récurrente**

**HIS1 N0659** **Fraction rationnelle**

**a(n) = n a(n-1)**

$$\frac{1}{1 - z}$$

**1, 1, 2, 6, 24, 120, 720, 5040, 40320, 362880, 3628800, 39916800, 479001600, 6227020800, 87178291200, 1307674368000, 20922789888000, 355687428096000**

## Oriented rooted trees with n nodes

**Réf.** R1 138.

**HIS2 A0151** **Euler**

**HIS1 N0701** **Produit infini**

$$\prod_{n \geq 1} \frac{1}{(1 - z^n)^{c(n)}}$$

$$c(n) = 2\ a(n)$$

**1, 2, 7, 26, 107, 458, 2058, 9498, 44947, 216598, 1059952, 5251806, 26297238, 132856766, 676398395, 3466799104, 17873808798, 92630098886, 482292684506**

---

**Réf.** R1 188.

**HIS2 A0153** **Dérivée logarithmique** **Suite P-récurrente**

**HIS1 N0706** **exponentielle**

**a(n) = n a(n-1) + (n - 2) a(n-2)**

$$\frac{1}{(1 - z)^3 \exp(z)}$$

**0, 1, 2, 7, 32, 181, 1214, 9403, 82508, 808393, 8743994, 103459471, 1328953592, 18414450877, 273749755382, 4345634192131, 73362643649444**

## Coefficients of iterated exponentials

**Réf.** SMA 11 353 45.

**HIS2 A0154** **Recoupements**

**HIS1 N0710** **exponentielle (log)**

**L'inverse fonctionnel est exp(exp(z)-1) : Les nombres de Bell.**

$$- \ln(1 + \ln(1 - z)) + 1$$

**1, 1, 2, 7, 35, 228, 1834, 17382, 195866, 2487832, 35499576, 562356672, 9794156448, 186025364016, 3826961710272, 84775065603888, 2011929826983504**

## Double factorials

**Réf.** AMM 55 425 48. MOC 24 231 70.

**HIS2 A0165** **Dérivée logarithmique** **Suite P-récurrente**

**HIS1 N0742** **Fraction rationnelle**

**2^(m-1) Γ(m)**

$$\frac{1}{1 - 2 z}$$

**1, 2, 8, 48, 384, 3840, 46080, 645120, 10321920, 185794560, 3715891200, 81749606400, 1961990553600, 51011754393600, 1428329123020800**

## Subfactorial or rencontres numbers

**Réf.** R1 65. DB1 168. RY63 23. MOC 21 502 67. C1 182.

**HIS2 A0166** **Dérivée logarithmique** **Suite P-récurrente**

**HIS1 N0766** **exponentielle**

**a(n) = (n - 2) a(n-1) + (n - 2) a(n -2)**

$$\frac{1}{(1 - z) \exp(z)}$$

**1, 0, 1, 2, 9, 44, 265, 1854, 14833, 133496, 1334961, 14684570, 176214841, 2290792932, 32071101049, 481066515734, 7697064251745, 130850092279664**

---

**Réf.** CJM 15 254 63; 33 1039 81. JCT 3 121 67.

**HIS2 A0168** **hypergéométrique-LLL** **Suite P-récurrente**

**HIS1 N0768** **algébrique**

**$2F1([1, 1/2], [3], 12 z)$**

**(n + 1) a(n) = (12 n - 18) a(n - 1)**

$$\frac{-1 + 18 z + (-(12 z - 1)^3)^{1/2}}{54 z^2}$$

**1, 2, 9, 54, 378, 2916, 24057, 208494, 1876446, 17399772, 165297834, 1602117468, 15792300756, 157923007560, 1598970451545, 16365932856990**

**Réf.** BA9. R1 128.

**HIS2 A0169** — **Inverse fonctionnel exponentielle** — **L'inverse fonctionnel est z exp(- z)**

**HIS1 N0771**

**n^ (n-1)**

$$- W(- z)$$

**1, 2, 9, 64, 625, 7776, 117649, 2097152, 43046721, 1000000000, 25937424601, 743008370688, 23298085122481, 793714773254144, 29192926025390625**

---

## Card matching

**Réf.** R1 193.

**HIS2 A0172** — **P-récurrences** — **Suite P-récurrente**

**HIS1 N0781** — *** titre modifié**

$$\sum_{k=0}^{n} (n,k)^3 = a(n)$$

$$a(n)\,(n-1)^2 = (7 n^2 - 21 n + 16)\, a(n-1) + (8 n^2 - 32 n + 32)\, a(n-2)$$

**1, 2, 10, 56, 346, 2252, 15184, 104960, 739162, 5280932, 38165260, 278415920, 2046924400, 15148345760, 112738423360, 843126957056, 6332299624282**

## Ménage numbers

**Réf.** CJM 10 478 58. R1 197.

**HIS2 A0179** **P-récurrences** **Suite P-récurrente**

**HIS1 N0815**

$$(n - 39/7)\, a(n) = (n^2 - 47/7\, n + 43/7)\, a(n-1) + (1/7\, n^2 + n - 65/7)\, a(n-2) + (-6/7\, n^2 + 67/7\, n - 26)\, a(n-3) + (-6/7\, n + 36/7)\, a(n-4)$$

**1, 1, 0, 1, 2, 13, 80, 579, 4738, 43387, 439792, 4890741, 59216642, 775596313, 10927434464, 164806435783, 2649391469058, 45226435601207, 817056406224416**

## Permutations with no cycles of length 3

**Réf.** R1 83.

**HIS2 A0180** **Dérivée logarithmique exponentielle** **Suite P-récurrente**

**HIS1 N0816**

**a(n) = (3 n - 4) a(n - 1) + (3 n - 6) a(n - 2)**

$$\frac{1}{(1 - 3z)\exp(z)}$$

**1, 2, 13, 116, 1393, 20894, 376093, 7897952, 189550849, 5117872922, 153536187661, 5066694192812, 182400990941233, 7113638646708086**

## Lucas numbers

**Réf.** HW1 148. HO69. C1 46.

**HIS2 A0204** **Approximants de Padé**

**HIS1 N0924** **Fraction rationnelle**

**a(n) = a(n-1) + a(n-2)**

$$\frac{1 + 2z}{1 - z - z^2}$$

**1, 3, 4, 7, 11, 18, 29, 47, 76, 123, 199, 322, 521, 843, 1364, 2207, 3571, 5778, 9349, 15127, 24476, 39603, 64079, 103682, 167761, 271443, 439204, 710647, 1149851**

---

**Réf.** SMA 20 23 54. R1 233. JCT 7 292 69.

**HIS2 A0211** **Approximants de Padé**

**HIS1 N0953** **Fraction rationnelle**

$$\frac{(1 + z)(4z - 3)}{(1 - z)(1 - z - z^2)}$$

**3, 5, 6, 9, 13, 20, 31, 49, 78, 125, 201, 324, 523, 845, 1366, 2209, 3573, 5780, 9351, 15129, 24478, 39605, 64081, 103684, 167763, 271445, 439206, 710649, 1149853**

**Réf.**
**HIS2 A0212** **Approximants de Padé**
**HIS1 N0966** **Fraction rationnelle**
**Partie entière de (n^2)/3.**

$$\frac{1 - z + 2 z^2 - z^3 + 2 z^4 - z^5}{(z^2 + z + 1)(1 - z)^3}$$

**1, 1, 3, 5, 8, 12, 16, 21, 27, 33, 40, 48, 56, 65, 75, 85, 96, 108, 120, 133, 147, 161, 176, 192, 208, 225, 243, 261, 280, 300, 320, 341, 363, 385, 408, 432, 456, 481, 507, 533**

---

**Réf.** FQ 1(3) 72 63; 2 260 64.
**HIS2 A0213** **Approximants de Padé**
**HIS1 N0975** **Fraction rationnelle**

$$\frac{(z - 1)(1 + z)}{1 - z - z^2 - z^3}$$

**1, 1, 1, 3, 5, 9, 17, 31, 57, 105, 193, 355, 653, 1201, 2209, 4063, 7473, 13745, 25281, 46499, 85525, 157305, 289329, 532159, 978793, 1800281, 3311233, 6090307, 11201821**

## Triangular numbers

**Réf.** D1 2 1. RS3. B1 189. AS1 828.

**HIS2 A0217** **Approximants de Padé**

**HIS1 N1002** **Fraction rationnelle**

$$\frac{1}{(1 - z)^3}$$

**1, 3, 6, 10, 15, 21, 28, 36, 45, 55, 66, 78, 91, 105, 120, 136, 153, 171, 190, 210, 231, 253, 276, 300, 325, 351, 378, 406, 435, 465, 496, 528, 561, 595, 630, 666, 703, 741**

## Planar partitions of n

**Réf.** MA15 2 332. PCPS 63 1099 67. AN76 241.

**HIS2 A0219** **Euler**

**HIS1 N1016** **Produit infini**

$$\prod_{n \geq 1} \frac{1}{(1 - z^n)^{c(n)}}$$

$$c(n) = 1,2,3,4,5,6,7,...$$

**1, 3, 6, 13, 24, 48, 86, 160, 282, 500, 859, 1479, 2485, 4167, 6879, 11297, 18334, 29601, 47330, 75278, 118794, 186475, 290783, 451194, 696033, 1068745, 1632658**

2 ^ ( n - 1)

**Réf.** BA9.

**HIS2 A0225** **Approximants de Padé**

**HIS1 N1059** **fraction rationnelle**

$$\frac{1}{(1 - 2 z) (1 - z)}$$

**1, 3, 7, 15, 31, 63, 127, 255, 511, 1023, 2047, 4095, 8191, 16383, 32767, 65535, 131071, 262143, 524287, 1048575, 2097151, 4194303, 8388607, 16777215, 33554431**

---

**Réf.** R1 65.

**HIS2 A0240** **Dérivée logarithmique** **Suite P-récurrente**

**HIS1 N1111** **exponentielle**

**a(n) = (n - 2) a(n - 1) + (2 n - 3) a(n - 2) + (n - 2) a(n - 3)**

$$\frac{\exp(- z) (z^2 - z + 1)}{(z - 1)^2}$$

**1, 0, 3, 8, 45, 264, 1855, 14832, 133497, 1334960, 14684571, 176214840, 2290792933, 32071101048, 481066515735, 7697064251744, 130850092279665**

## Crossing number of complete graph with n nodes

**Réf.** GU60. AMM 80 53 73.

**HIS2 A0241** **Approximants de Padé** **conjecture connue**

**HIS1 N1115** **Fraction rationnelle**

$$\frac{1 + z + z^2}{(z - 1)^5 (z + 1)^3}$$

**0, 0, 0, 0, 1, 3, 9, 18, 36, 60, 100, 150, 225, 315, 441, 588**

## Powers of 3

**Réf.** BA9.

**HIS2 A0244** **Approximants de Padé**

**HIS1 N1129** **fraction rationnelle**

$$\frac{1}{1 - 3z}$$

**1, 3, 9, 27, 81, 243, 729, 2187, 6561, 19683, 59049, 177147, 531441, 1594323, 4782969, 14348907, 43046721, 129140163, 387420489, 1162261467**

**Réf.** QAM 14 407 56. MOC 29 216 75. FQ 14 397 76.
**HIS2 A0245** **Hypergéométrique** **Suite P-récurrente**
**HIS1 N1130** **algébrique**
**(n + 2) a(n) = (5 n + 2) a(n - 1) + (- 4 n + 6) a(n - 2)**
**2F1([3/2, 2], [4], 4 z)**

$$\frac{8\,z}{\left(1 + (1 - 4\,z)^{1/2}\right)^{3}}$$

**1, 3, 9, 28, 90, 297, 1001, 3432, 11934, 41990, 149226, 534888, 1931540, 7020405, 25662825, 94287120, 347993910, 1289624490, 4796857230, 17902146600**

---

## Permutations of length n with odd cycles

**Réf.** R1 87.
**HIS2 A0246** **Hypergéométrique** **Suite P-récurrente**
**HIS1 N1137** **algébrique**
**a(n) = a(n - 1) + (n^2 - 3 n + 2) a(n - 2)**

$$\frac{1}{(1 - z)^{3/2}\,(1 + z)^{1/2}}$$

**0, 1, 1, 3, 9, 45, 225, 1575, 11025, 99225, 893025, 9823275, 108056025, 1404728325, 18261468225, 273922023375, 4108830350625, 69850115960625**

## Associated Stirling numbers

**Réf.** R1 76. DB1 296. C1 222.

**HIS2 A0247** **Approximants de Padé**

**HIS1 N1141** **fraction rationnelle**

$$\frac{3 - 2 z}{1 - 4 z + 5 z^{2} - 2 z^{3}}$$

**3, 10, 25, 56, 119, 246, 501, 1012, 2035, 4082, 8177, 16368, 32751, 65518, 131053, 262124, 524267, 1048554, 2097129, 4194280, 8388583, 16777190, 33554405**

## Forests with n nodes and height at most 1

**Réf.** JCT 3 134 67; 5 102 68. C1 91.

**HIS2 A0248** **Dérivée logarithmique**

**HIS1 N1148** **exponentielle**

$$\exp(\exp(z)\ z)$$

**1, 1, 3, 10, 41, 196, 1057, 6322, 41393, 293608, 2237921, 18210094, 157329097, 1436630092, 13810863809, 139305550066, 1469959371233**

## Stirling numbers of first kind

**Réf.** AS1 833. DKB 226.

**HIS2 A0254** **équations différentielles** **Suite P-récurrente**

**HIS1 N1165** **exponentielle (log)**

**a(n) = (2 n - 1) a(n - 1) + (- n^2 + 2 n - 1) a(n - 2)**

$$\frac{1 - \ln(1 - z)}{(1 - z)^2}$$

**1, 3, 11, 50, 274, 1764, 13068, 109584, 1026576, 10628640, 120543840, 1486442880, 19802759040, 283465647360, 4339163001600, 70734282393600**

---

**Réf.** R1 188. DKB 263. MAG 52 381 68.

**HIS2 A0255** **Dérivée logarithmique** **Suite P-récurrente**

**HIS1 N1166** **exponentielle**

**a(n) = na(n-1) + (n-1)a(n-2)**

$$\frac{\exp(- z)}{(1 - z)^2}$$

**1, 1, 3, 11, 53, 309, 2119, 16687, 148329, 1468457, 16019531, 190899411, 2467007773, 34361893981, 513137616783, 8178130767479**

Réf. CJM 15 268 63.

HIS2 A0256 LLL Suite P-récurrente

HIS1 N1173 algébrique 3è degré

$$\frac{1}{2}(n-1)(n-3)(2n-1)a(n) = \frac{1}{16}(n-3)(104n^2-430n+414)a(n-1) + \frac{1}{16}(n-3)(27n^2-81n+60)a(n-2)$$

1, 1, 0, 1, 3, 12, 52, 241, 1173, 5929, 30880, 164796, 897380, 4970296, 27930828, 158935761, 914325657, 5310702819, 31110146416, 183634501753, 1091371140915

---

## Rooted bicubic maps

Réf. CJM 15 269 63.

HIS2 A0257 Hypergéométrique Suite P-récurrente

HIS1 N1175 algébrique

$_2F_1([1, 3/2], [4], 8z)$

$(n+2)a(n) = (8n-4)a(n-1)$

$$\frac{3(1-8z)^{1/2} + 8z - 3(1-8z)^{3/2}}{4(1+(1-8z)^{1/2})^3 z}$$

1, 3, 12, 56, 288, 1584, 9152, 54912, 339456

## Coefficients of iterated exponentials

**Réf.** SMA 11 353 45. PRV A32 2342 85.

**HIS2 A0258** **Recoupements**

**HIS1 N1178** **exponentielle**

$$\exp(\exp(\exp(z) - 1) - 1)$$

**1, 1, 3, 12, 60, 358, 2471, 19302, 167894, 1606137, 16733779, 188378402, 2276423485, 29367807524, 402577243425, 5840190914957, 89345001017415**

---

**Réf.** CJM 14 32 62.

**HIS2 A0260** **Hypergéométrique** **Suite P-récurrente**

**HIS1 N1187** **algébrique** **algébrique du 4è degré**

**$_4F_3$ ([1, 1/2, 3/4, 5/4],[2, 5/3, 4/3],(256/27) z)**

$$\frac{1}{9} (3 n - 1) (3 n - 2) n\, a(n) = \frac{8}{27} (4 n - 5) (4 n - 3) (2 n - 3)\, a(n - 1)$$

**1, 1, 3, 13, 68, 399, 2530, 16965, 118668, 857956, 6369883, 48336171, 373537388, 2931682810, 23317105140, 187606350645, 1524813969276, 12504654858828**

**Réf.** R1 188.
**HIS2 A0261** **Dérivée logarithmique** **Suite P-récurrente**
**HIS1 N1189** **exponentielle**
**a(n) = (n + 1) a(n - 1) + (n - 2) a(n - 2)**

$$\frac{\exp(-z)}{(1-z)^4}$$

**0, 1, 3, 13, 71, 465, 3539, 30637, 296967, 3184129, 37401155, 477471021, 6581134823, 97388068753, 1539794649171, 25902759280525, 461904032857319**

---

**Réf.** RCI 194. PSPM 19 172 71.
**HIS2 A0262** **Dérivée logarithmique** **Suite P-récurrente**
**HIS1 N1190** **exponentielle**
**a(n) = (2n-1) a(n-1) - (n-1) (n-2) a(n-2)**

$$\exp(z/(1-z))$$

**1, 1, 3, 13, 73, 501, 4051, 37633, 394353, 4596553, 58941091, 824073141, 12470162233, 202976401213, 3535017524403, 65573803186921, 1290434218669921**

**Réf.** R1 85.

**HIS2 A0266** **Dérivée logarithmique** **Suite P-récurrente**

**HIS1 N1211** **exponentielle**

**a(n) = (n - 1) a(n - 1) + (- n + 2) a(n - 2) + (n^2 - 5 n + 6) a(n - 3)**

$$\frac{1}{\exp(1/2\, z^2)\,(1 - z)}$$

**1, 1, 1, 3, 15, 75, 435, 3045, 24465, 220185, 2200905, 24209955, 290529855, 3776888115, 52876298475, 793144477125, 12690313661025, 215735332237425, 3883235945814225**

---

## Coefficients of iterated exponentials

**Réf.** SMA 11 353 45.

**HIS2 A0268** **Recoupements**

**HIS1 N1218** **exponentielle**

**L'inverse fonctionnel est** exp(exp(exp(z)-1)-1)

$$-\ln(1 + \ln(1 + \ln(1 - z))) + 1$$

**1, 1, 3, 15, 105, 947, 10472, 137337, 2085605, 36017472, 697407850, 14969626900, 352877606716, 9064191508018, 252024567201300, 7542036496650006**

## Sums of ménage numbers

**Réf.** AH21 2 79. CJM 10 478 58. R1 198.

**HIS2 A0271** **P-récurrences** **Suite P-récurrente**

**HIS1 N1222**

$$a(n) = (n + 1)\, a(n - 1) + (n + 1)\, a(n - 2) + a(n - 3)$$

**0, 0, 1, 3, 16, 96, 675, 5413, 48800, 488592, 5379333, 64595975, 840192288, 11767626752, 176574062535, 2825965531593, 48052401132800, 865108807357216**

---

**Réf.** BA9. R1 128.

**HIS2 A0272** **Inverse fonctionnel**

**HIS1 N1227** **exponentielle** **f.g. exponentielle**

**n^(n-2)**

**L'inverse est ln(1+z)/(1+z)**

$$\frac{z + W(- z)}{z}$$

**1, 3, 16, 125, 1296, 16807, 262144, 4782969, 100000000, 2357947691, 61917364224, 1792160394037, 56693912375296, 1946195068359375**

## Permutations of length n by rises

**Réf.** DKB 263. R1 210 (divided by 2).

**HIS2 A0274** **Dérivée logarithmique** **Suite P-récurrente**

**HIS1 N1236** **exponentielle**

**a(n) = (n + 1) a(n - 1) + (n + 3) a(n - 2) + (- n + 3) a(n - 3) + (- n + 2) a(n - 4)**

$$\frac{2 - 5 z^2 + 2 z^3 - z^4}{2 (1 - z)^4 \exp(z)}$$

**1, 3, 18, 110, 795, 6489, 59332, 600732, 6674805, 80765135, 1057289046, 14890154058, 224497707343, 3607998868005**

## Associated Stirling numbers

**Réf.** R1 75. C1 256.

**HIS2 A0276** **équations différentielles** **Suite P-récurrente**

**HIS1 N1248** **exponentielle (log)** **Formule de B. Salvy**

**a(n) = (2 n + 2) a(n - 1) - (n^2 + 1) a(n - 2) - (n^2 + n) a(n - 3)**

$$\frac{2 z - 6 \ln(- z + 1) + 3}{(1 - z)^4}$$

**3, 20, 130, 924, 7308, 64224, 623376, 6636960, 76998240, 967524480, 13096736640, 190060335360, 2944310342400, 48503818137600, 846795372595200**

**Réf.** FQ 3 129 65. BR72 53.

**HIS2 A0285** **Approximants de Padé**

**HIS1 N1309** **Fraction rationnelle**

$$\frac{1 + 3z}{1 - z - z^2}$$

**1, 4, 5, 9, 14, 23, 37, 60, 97, 157, 254, 411, 665, 1076, 1741, 2817, 4558, 7375, 11933, 19308, 31241, 50549, 81790, 132339, 214129, 346468, 560597, 907065, 1467662**

---

## Rooted polyhedral graphs with n edges

**Réf.** CJM 15 265 63.

**HIS2 A0287** **LLL** **suite corrigée avec la formule de**

**HIS1 N1326** **algébrique** **récurrence.**

$$(n + 4)\, a(n) = (3/2\, n - 3)\, a(n - 1) + (8\, n + 4)\, a(n - 2) + (15/2\, n + 6)\, a(n - 3) + (2\, n + 3)\, a(n - 4)$$

$$\frac{(1 + z)\left((-4z + 1)^{3/2} - 1 + 6z - 6z^2 - 4z^3 - 6z^4\right) + 4z^5}{2\,(2 z^5 (z + 2)^3 (1 + z))}$$

**1, 0, 4, 6, 24, 66, 214, 676, 2209, 7296, 24460, 82926, 284068, 981882, 3421318, 12007554, 42416488, 150718770, 538421590, 1932856590, 6969847484**

## Tetranacci numbers

**Réf.** FQ 2 260 64.

**HIS2 A0288** **Approximants de Padé**

**HIS1 N1332** **Fraction rationnelle**

$$\frac{1 - z^{2} - 2 z^{3}}{1 - z - z^{2} - z^{3} - z^{4}}$$

**1, 1, 1, 1, 4, 7, 13, 25, 49, 94, 181, 349, 673, 1297, 2500, 4819, 9289, 17905, 34513, 66526, 128233, 247177, 476449, 918385, 1770244, 3412255, 6577333, 12678217**

## The squares

**Réf.** BA9.

**HIS2 A0290** **Approximants de Padé**

**HIS1 N1350** **Fraction rationnelle**

$$\frac{1 + z}{(1 - z)^{3}}$$

**1, 4, 9, 16, 25, 36, 49, 64, 81, 100, 121, 144, 169, 196, 225, 256, 289, 324, 361, 400, 441, 484, 529, 576, 625, 676, 729, 784, 841, 900, 961, 1024, 1089, 1156, 1225, 1296**

## Tetrahedral numbers

**Réf.** D1 2 4. RS3. B1 194. AS1 828.

**HIS2 A0292** **Approximants de Padé**

**HIS1 N1363** **Fraction rationnelle**

**C(n,3)**

$$\frac{1}{(1-z)^4}$$

**1, 4, 10, 20, 35, 56, 84, 120, 165, 220, 286, 364, 455, 560, 680, 816, 969, 1140, 1330, 1540, 1771, 2024, 2300, 2600, 2925, 3276, 3654, 4060, 4495, 4960, 5456, 5984**

## Related to solid partitions

**Réf.** PNISI 26 135 60. PCPS 63 1100 67.

**HIS2 A0294** **Euler**

**HIS1 N1372** **Produit infini**

$$\prod_{n\geq 1} \frac{1}{(1 - z^{c(n)})}$$

**c(n) = 1,3,6,10,...,nombres triangulaires**

**1, 1, 4, 10, 26, 59, 141, 310, 692, 1483, 3162, 6583, 13602, 27613, 55579, 110445, 217554, 424148, 820294, 1572647, 2992892, 5652954, 10605608, 19765082**

## Eulerian numbers 2^n -n - 1

**Réf.** R1 215. DB1 151.

**HIS2 A0295** **Approximants de Padé**

**HIS1 N1382** **Fraction rationnelle**

$$\frac{1}{(1 - 2z)(1 - z)^2}$$

**0, 1, 4, 11, 26, 57, 120, 247, 502, 1013, 2036, 4083, 8178, 16369, 32752, 65519, 131054, 262125, 524268, 1048555, 2097130, 4194281, 8388584, 16777191, 33554406**

---

**Réf.** FQ 14 69 76. ANY 319 464 79.

**HIS2 A0296** **Dérivée logarithmique exponentielle** **Différences finies des nombres de Bell**

**HIS1 N1387**

$$\exp(\exp(z) - 1 - z)$$

**1, 0, 1, 1, 4, 11, 41, 162, 715, 3425, 17722, 98253, 580317, 3633280, 24011157, 166888165, 1216070380, 9264071767, 73600798037, 608476008122, 5224266196935**

**Réf.** R1 150. FQ 15 194 77.

**HIS2 A0297** **Approximants de Padé**

**HIS1 N1393** **Fraction rationnelle**

$$\frac{(z - 2)^2}{(1 - z)^4}$$

**4, 12, 25, 44, 70, 104, 147, 200, 264, 340, 429, 532, 650, 784, 935, 1104, 1292, 1500, 1729, 1980, 2254, 2552, 2875, 3224, 3600, 4004, 4437, 4900, 5394, 5920, 6479**

---

## Powers of 4

**Réf.** BA9.

**HIS2 A0302** **Approximants de Padé**

**HIS1 N1428** **Fraction rationnelle**

$$\frac{1}{1 - 4 z}$$

**1, 4, 16, 64, 256, 1024, 4096, 16384, 65536, 262144, 1048576, 4194304, 16777216, 67108864, 268435456, 1073741824, 4294967296, 17179869184**

## Coefficients of iterated exponentials

**Réf.** SMA 11 353 45. PRV A32 2342 85.

**HIS2 A0307** **Recoupements**

**HIS1 N1455** **exponentielle**

$$\exp(\exp(\exp(\exp(z) - 1) - 1) - 1)$$

**1, 1, 4, 22, 154, 1304, 12915, 146115, 1855570, 26097835, 402215465, 6734414075, 121629173423, 2355470737637, 48664218965021, 1067895971109199**

## Rooted maps with 2n nodes

**Réf.** CJM 14 416 62.

**HIS2 A0309** **Hypergéométrique** **Suite P-récurrente**

**HIS1 N1460** **algébrique** **Algébrique du 3è degré**

$$\frac{1}{2} (n + 1) (2 n + 1) a(n) = \frac{3}{2} (3 n - 1) (3 n - 2) a(n - 1)$$

$$- \frac{1}{12} \Big( \big(1458 z^2 + 270 z - 1 + 12 (- 2 + 27 z)^{1/2} 3^{1/2} z^{1/2} - 162 (- 2 + 27 z)^{1/2} 3^{1/2} z^{3/2}\big)^{1/3} + \big(1458 z^2 + 270 z - 1 - 12 (- 2 + 27 z)^{1/2} 3^{1/2} z^{1/2} + 162 (- 2 + 27 z)^{1/2} 3^{1/2} z^{3/2}\big)^{1/3} + 12 z + 2 \Big)$$

**1, 4, 24, 176, 1456, 13056, 124032, 1230592, 12629760, 133186560, 1436098560**

## Coefficients of iterated exponentials

**Réf.** SMA 11 353 45.

**HIS2 A0310** **Recoupements**

**HIS1 N1464** **exponentielle (log)**

$$- \ln(1 + \ln(1 + \ln(1 + \ln(1 - z)))) + 1$$

**1, 1, 4, 26, 234, 2696, 37919, 630521, 12111114, 264051201, 6445170229, 174183891471, 5164718385337, 166737090160871, 5822980248613990**

---

## Schroeder's fourth problem

**Réf.** RCI 197. C1 224.

**HIS2 A0311** **Inverse fonctionnel**

**HIS1 N1465** **exponentielle**

**L'inverse fonctionnel de 1 + 2 z - exp(z)**

$$- W(- 1/2* \exp(- 1/2 + 1/2*z)) - 1/2 + 1/2*z$$

**1, 1, 1, 4, 26, 236, 2752, 39208, 660032, 12818912, 282137824, 6939897856, 188666182784, 5617349020544, 181790703209728, 6353726042486272, 238513970965257728**

**Réf.** BA9.
**HIS2 A0312** **Inverse fonctionnel**
**HIS1 N1469** **exponentielle**
**a(n) = n^n**
**L'inverse fonctionnel de z exp(1/(z+1))/(z+1)**

$$\frac{W(-z)}{-1 - W(-z)}$$

**1, 4, 27, 256, 3125, 46656, 823543, 16777216, 387420489, 10000000000, 285311670611, 8916100448256, 302875106592253, 11112006825558016**

---

**Permutations of length n by rises**

**Réf.** DKB 263.
**HIS2 A0313** **Approximants de Padé** **Suite P-récurrente**
**HIS1 N1477** **exponentielle** **Conjecture**

$$\frac{-z^6 + 6z^5 - 18z^4 + 22z^3 - 27z^2 - 6}{(z-1)^5 \exp(z)}$$

**1, 4, 30, 220, 1855, 17304, 177996, 2002440, 24474285, 323060540, 4581585866, 69487385604, 1122488536715**

## Pentanacci numbers

**Réf.** FQ 2 260 64.

**HIS2 A0322** **Approximants de Padé**

**HIS1 N1542** **Fraction rationnelle**

$$\frac{3 z^4 + 2 z^3 + z^2 - 1}{z^5 + z^4 + z^3 + z^2 + z - 1}$$

**1, 1, 1, 1, 1, 5, 9, 17, 33, 65, 129, 253, 497, 977, 1921, 3777, 7425, 14597, 28697, 56417, 110913, 218049, 428673, 842749, 1656801, 3257185, 6403457, 12588865, 24749057**

## Pentagonal numbers

**Réf.** D1 2 1. B1 189. HW1 284. FQ 8 84 70.

**HIS2 A0326** **Approximants de Padé**

**HIS1 N1562** **Fraction rationnelle**

$$\frac{(1 + 2 z)}{(1 - z)^3}$$

**1, 5, 12, 22, 35, 51, 70, 92, 117, 145, 176, 210, 247, 287, 330, 376, 425, 477, 532, 590, 651, 715, 782, 852, 925, 1001, 1080, 1162, 1247, 1335, 1426, 1520, 1617, 1717**

## Square pyramidal numbers

**Réf.** D1 2 2. B1 194. AS1 813.

**HIS2 A0330** **Approximants de Padé**

**HIS1 N1574** **Fraction rationnelle**

$$\frac{1 + z}{(1 - z)^4}$$

**1, 5, 14, 30, 55, 91, 140, 204, 285, 385, 506, 650, 819, 1015, 1240, 1496, 1785, 2109, 2470, 2870, 3311, 3795, 4324, 4900, 5525, 6201, 6930, 7714, 8555, 9455, 10416**

## Figurate numbers C(n,4)

**Réf.** D1 2 7. RS3. B1 196. AS1 828.

**HIS2 A0332** **Approximants de Padé**

**HIS1 N1578** **Fraction rationnelle**

$$\frac{1}{(1 - z)^5}$$

**1, 5, 15, 35, 70, 126, 210, 330, 495, 715, 1001, 1365, 1820, 2380, 3060, 3876, 4845, 5985, 7315, 8855, 10626, 12650, 14950, 17550, 20475, 23751, 27405, 31465**

**Réf.** HB67 16.

**HIS2 A0337** **Approximants de Padé**

**HIS1 N1587** **Fraction rationnelle**

$$\frac{1}{(z-1)(2z-1)^2}$$

**1, 5, 17, 49, 129, 321, 769, 1793, 4097, 9217, 20481, 45057, 98305, 212993, 458753, 983041, 2097153, 4456449, 9437185, 19922945, 41943041, 88080385**

---

**Réf.** SMA 20 23 54.

**HIS2 A0338** **Approximants de Padé**

**HIS1 N1589** **Fraction rationnelle**

$$\frac{(2z-5)(z^2+z+1)}{(z-1)^3}$$

**5, 18, 42, 75, 117, 168, 228, 297, 375, 462, 558, 663, 777, 900, 1032, 1173, 1323, 1482, 1650, 1827, 2013, 2208, 2412, 2625, 2847, 3078, 3318, 3567**

**Réf.** DKB 260.

**HIS2 A0340** **Approximants de Padé**

**HIS1 N1592** **Fraction rationnelle**

$$\frac{1}{(1 - 3 z) (1 - z)^{2}}$$

**1, 5, 18, 58, 179, 543, 1636, 4916, 14757, 44281, 132854, 398574, 1195735, 3587219, 10761672, 32285032, 96855113, 290565357, 871696090, 2615088290, 7845264891**

---

**Réf.** QAM 14 407 56. MOC 29 216 75. FQ 14 397 76.

**HIS2 A0344** **Hypergéométrique** **Suite P-récurrente**

**HIS1 N1602** **algébrique**

**$_3F_2([5/2, 3], [6], 4 z)$**

**(n + 4) (n - 1) a(n) = 2 (n + 1) (2 n + 1) a(n - 1)**

$$\frac{32 z}{(1 + (1 - 4 z)^{1/2})^{5}}$$

**1, 5, 20, 75, 275, 1001, 3640, 13260, 48450, 177650, 653752, 2414425, 8947575, 33266625, 124062000, 463991880, 1739969550, 6541168950, 24647883000**

**Réf.** BAMS 74 74 68. JCT 13 215 72.

**HIS2 A0346** **LLL** **Suite P-récurrente**

**HIS1 N1611** **algébrique**

**n a(n) = (8 n - 6) a(n - 1) + (- 16 n + 24) a(n - 2)**

$$\frac{1 - 4z - (-(-1 + 4z)^3)^{1/2}}{2(z - 8z^2 + 16z^3)}$$

**1, 5, 22, 93, 386, 1586, 6476, 26333, 106762, 431910, 1744436, 7036530, 28354132, 114159428, 459312152, 1846943453, 7423131482, 29822170718, 119766321572, 480832549478**

---

## Powers of 5

**Réf.** BA9.

**HIS2 A0351** **Approximants de Padé**

**HIS1 N1620** **Fraction rationnelle**

$$\frac{1}{1 - 5z}$$

**1, 5, 25, 125, 625, 3125, 15625, 78125, 390625, 1953125, 9765625, 48828125, 244140625, 1220703125, 6103515625, 30517578125, 152587890625**

## Permutations of length n by number of runs

**Réf.** DKB 260.

**HIS2 A0352** **Approximants de Padé**

**HIS1 N1629** **Fraction rationnelle**

$$\frac{5 - 6 z}{(3 z - 1) (2 z - 1) (z - 1)^2}$$

**5, 29, 118, 418, 1383, 4407, 13736, 42236, 128761, 390385, 1179354, 3554454**

---

**Réf.** LU91 1 223. R1 83.

**HIS2 A0354** **Dérivée logarithmique** **Suite P-récurrente**

**HIS1 N1631** **exponentielle**

**1/2 a(n) = (n - 3/2) a(n - 1) + (n - 2) a(n - 2)**

$$\frac{1}{(1 - 2 z) \exp(z)}$$

**1, 1, 5, 29, 233, 2329, 27949, 391285, 6260561, 112690097, 2253801941, 49583642701, 1190007424825, 30940193045449, 866325405272573**

## Hamiltonian rooted maps with 2n nodes

**Réf.** CJM 14 416 62.

**HIS2 A0356** **hypergéométrique** **Suite P-récurrente**

**HIS1 N1647** **Intégrales elliptiques**

**2F1([1/2, -1/2],[2],16 z)**

**1, 5, 35, 294, 2772, 28314, 306735, 3476330, 40831076, 493684828, 6114096716**

---

## Coefficients of iterated exponentials

**Réf.** SMA 11 353 45. PRV A32 2342 85.

**HIS2 A0357** **Recoupements**

**HIS1 N1648** **exponentielle**

**exp(exp(exp(exp(exp(z) - 1) - 1) - 1) - 1)**

**1, 1, 5, 35, 315, 3455, 44590, 660665, 11035095, 204904830, 4183174520, 93055783320, 2238954627848, 57903797748386, 1601122732128779**

## Coefficients of iterated exponentials

**Réf.** SMA 11 353 45.

**HIS2 A0359** **Recoupements**

**HIS1 N1654** **exponentielle (log)**

$$- \ln(1 + \ln(1 + \ln(1 + \ln(1 + \ln(1 - z))))) + 1$$

**1, 1, 5, 40, 440, 6170, 105315, 2120610, 49242470, 1296133195, 38152216495, 1242274374380, 44345089721923, 1722416374173854, 72330102999829054**

---

**Réf.** CMB 4 32 61 (divided by 3).

**HIS2 A0381** **Approximants de Padé**

**HIS1 N1692** **Fraction rationnelle**

$$\frac{2 - z - 2 z^2}{1 - 2 z + z^3}$$

**2, 3, 4, 6, 9, 14, 22, 35, 56, 90, 145, 234, 378, 611, 988, 1598, 2585, 4182, 6766, 10947, 17712, 28658, 46369, 75026, 121394, 196419, 317812, 514230**

## Restricted permutations

**Réf.** CMB 4 32 61 (divided by 4).

**HIS2 A0382** **Approximants de Padé**

**HIS1 N1696** **Fraction rationnelle**

$$\frac{6 - z - 2z^2 - 4z^3 - z^4}{1 - 2z + z^4}$$

**6, 11, 20, 36, 65, 119, 218, 400, 735, 1351, 2484, 4568, 8401, 15451, 28418, 52268, 96135, 176819, 325220, 598172, 1100209, 2023599, 3721978, 6845784**

## Hexanacci numbers

**Réf.** FQ 2 302 64.

**HIS2 A0383** **Approximants de Padé**

**HIS1 N1697** **Fraction rationnelle**

$$\frac{4z^5 + 3z^4 + 2z^3 + z^2 - 1}{z^6 + z^5 + z^4 + z^3 + z^2 + z - 1}$$

**1, 1, 1, 1, 1, 1, 6, 11, 21, 41, 81, 161, 321, 636, 1261, 2501, 4961, 9841, 19521, 38721, 76806, 152351, 302201, 599441, 1189041, 2358561, 4678401, 9279996, 18407641**

## Hexagonal numbers

**Réf.** D1 2 2. B1 189.

**HIS2 A0384** **Approximants de Padé**

**HIS1 N1705** **Fraction rationnelle**

$$\frac{1 + 3 z}{(1 - z)^3}$$

**1, 6, 15, 28, 45, 66, 91, 120, 153, 190, 231, 276, 325, 378, 435, 496, 561, 630, 703, 780, 861, 946, 1035, 1128, 1225, 1326, 1431, 1540, 1653, 1770, 1891, 2016, 2145, 2278**

## Rencontres numbers

**Réf.** R1 65.

**HIS2 A0387** **Dérivée logarithmique** **Suite P-récurrente**

**HIS1 N1716** **exponentielle**

**a(n) = (3 n - 4) a(n - 3) + (n - 2) a(n - 4) + (n - 2) a(n - 1) + (3 n - 3) a(n - 2)**

$$\frac{z^4 - 4 z^3 + 7 z^2 - 4 z + 2}{(z - 1)^3 \exp(z)}$$

**1, 0, 6, 20, 135, 924, 7420, 66744, 667485, 7342280, 88107426, 1145396460, 16035550531, 240533257860, 3848532125880, 65425046139824**

## Binomial coefficients C(n,5)

**Réf.** D1 2 7. RS3. B1 196. AS1 828.

**HIS2 A0389** **Approximants de Padé**

**HIS1 N1719** **Fraction rationnelle**

$$\frac{1}{(1 - z)^{6}}$$

**1, 6, 21, 56, 126, 252, 462, 792, 1287, 2002, 3003, 4368, 6188, 8568, 11628, 15504, 20349, 26334, 33649, 42504, 53130, 65780, 80730, 98280, 118755, 142506**

## Stirling numbers of second kind

**Réf.** AS1 835. DKB 223.

**HIS2 A0392** **Approximants de Padé**

**HIS1 N1734** **Fraction rationnelle**

$$\frac{1}{(1 - z)(1 - 2z)(1 - 3z)}$$

**1, 6, 25, 90, 301, 966, 3025, 9330, 28501, 86526, 261625, 788970, 2375101, 7141686, 21457825, 64439010, 193448101, 580606446, 1742343625, 5228079450**

## Stirling numbers of first kind

**Réf.** AS1 833. DKB 226.

**HIS2 A0399** **Tableaux généralisés** **Suite P-récurrente**

**HIS1 N1762** **exponentielle (log)**

**a(n) = -3 n^2 a(n - 1) + (n^3 - 3 n^2 + 3 n - 1) a(n - 3)**
**+ (n^3 - 3 n^2 - 3 n) a(n - 2)**

$$\frac{\ln(1 - z)^2}{2\,(1 - z)}$$

**1, 6, 35, 225, 1624, 13132, 118124, 1172700, 12753576, 150917976, 1931559552, 26596717056, 392156797824, 6165817614720, 102992244837120**

---

## Powers of 6

**Réf.** BA9.

**HIS2 A0400** **Approximants de Padé**

**HIS1 N1765** **Fraction rationnelle**

$$\frac{1}{1 - 6\,z}$$

**1, 6, 36, 216, 1296, 7776, 46656, 279936, 1679616, 10077696, 60466176, 362797056, 2176782336, 13060694016, 78364164096, 470184984576, 2821109907456**

## Coefficients of iterated exponentials

**Réf.** SMA 11 353 45. PRV A32 2342 85.

**HIS2 A0405** **Recoupements**
**HIS1 N1781** **exponentielle**

exp(exp(exp(exp(exp(exp(z) - 1) - 1) - 1) - 1) - 1)

**1, 1, 6, 51, 561, 7556, 120196, 2201856, 45592666, 1051951026, 26740775306, 742069051906, 22310563733864, 722108667742546, 25024187820786357**

## Coefficients of iterated exponentials

**Réf.** SMA 11 353 45.

**HIS2 A0406** **Recoupements**
**HIS1 N1782** **exponentielle (log)**

- ln(1 + ln(1 + ln(1 + ln(1 + ln(1 + ln(1 - z)))))) + 1

**1, 1, 6, 57, 741, 12244, 245755, 5809875, 158198200, 4877852505, 168055077875, 6400217406500, 267058149580823, 12118701719205803, 594291742526530761**

**Réf.** MOC 3 168 48; 9 174 55. CMA 2 25 70. MAN 191 98 71.
**HIS2 A0407** **Hypergéométrique** **Suite P-récurrente**
**HIS1 N1784** **algébrique**
**(2n)!/(2.n!)**

$$\frac{1}{(1 - 4 z)^{3/2}}$$

**1, 6, 60, 840, 15120, 332640, 8648640, 259459200, 8821612800, 335221286400, 14079294028800, 647647525324800, 32382376266240000**

---

**Powers of 7**

**Réf.** BA9.
**HIS2 A0420** **Approximants de Padé**
**HIS1 N1874** **Fraction rationnelle**

$$\frac{1}{1 - 7 z}$$

**1, 7, 49, 343, 2401, 16807, 117649, 823543, 5764801, 40353607, 282475249, 1977326743, 13841287201, 96889010407, 678223072849, 4747561509943**

## Permutations of length n by number of peaks

**Réf.** DKB 261.

**HIS2 A0431** **Approximants de Padé**

**HIS1 N0824** **Fraction rationnelle**

$$\frac{2}{1 - 8 z + 20 z^{2} - 16 z^{3}}$$

**2, 16, 88, 416, 1824, 7680, 31616, 128512, 518656, 2084864, 8361984, 33497088, 134094848, 536608768, 2146926592, 8588754944, 34357248000, 137433710592**

---

## Powers of rooted tree enumerator

**Réf.** R1 150.

**HIS2 A0439** **Approximants de Padé**

**HIS1 N1965** **Fraction rationnelle**

$$\frac{(3 - 2 z) (z^{2} - 3 z + 3)}{(1 - z)^{5}}$$

**9, 30, 69, 133, 230, 369, 560, 814, 1143, 1560, 2079, 2715, 3484, 4403, 5490, 6764, 8245, 9954, 11913, 14145, 16674, 19525, 22724, 26298, 30275, 34684, 39555**

**Réf.** CC55 742. RCI 217. JO61 7.
**HIS2 A0447** **Approximants de Padé**
**HIS1 N2006** **Fraction rationnelle**

$$\frac{z\,(1 + 6\,z + z^2)}{(z - 1)^4}$$

**0, 1, 10, 35, 84, 165, 286, 455, 680, 969, 1330, 1771, 2300, 2925, 3654, 4495, 5456, 6545, 7770, 9139, 10660, 12341, 14190, 16215, 18424, 20825, 23426, 26235, 29260**

---

## Rencontres numbers

**Réf.** R1 65.
**HIS2 A0449** **Dérivée logarithmique** **Suite P-récurrente**
**HIS1 N2009** **exponentielle**
**(n - 1) a(n) = (n + 2) (n - 2) a(n - 1) + (n + 2) (n + 1) a(n - 2)**

$$\frac{6 - 18\,z + 45\,z^2 - 49\,z^3 + 30\,z^4 - 9\,z^5 + z^6}{6\,(1 - z)^4\,\exp(z)}$$

**1, 0, 10, 40, 315, 2464, 22260, 222480, 2447445, 29369120, 381798846, 5345183480, 80177752655, 1282844041920, 21808348713320, 392550276838944**

## Stirling numbers of second kind

**Réf.** AS1 835. DKB 223.

**HIS2 A0453** **Approximants de Padé**

**HIS1 N2018** **Fraction rationnelle**

$$\frac{1}{(1 - z)(1 - 2z)(1 - 3z)(1 - 4z)}$$

**1, 10, 65, 350, 1701, 7770, 34105, 145750, 611501, 2532530, 10391745, 42355950, 171798901, 694337290, 2798806985, 11259666950, 45232115901**

## Stirling numbers of first kind

**Réf.** AS1 833. DKB 226.

**HIS2 A0454** **Tableaux généralisés**

**HIS1 N2022** **exponentielle (log)**

$$\frac{-\ln(1 - z)^3}{6(1 - z)}$$

**1, 10, 85, 735, 6769, 67284, 723680, 8409500, 105258076, 1414014888, 20313753096, 310989260400, 5056995703824, 87077748875904, 1583313975727488**

**Réf.** TOH 37 259 33. JO39 152. DB1 296. C1 256.

**HIS2 A0457** **Hypergéométrique** **Suite P-récurrente**

**HIS1 N2028** **algébrique** **f.g. exponentielle**

**(n - 1) a(n) = (2 n + 1) n a(n - 1)**

$$\frac{z}{(1 - 2 z)^{5/2}}$$

**1, 10, 105, 1260, 17325, 270270, 4729725, 91891800, 1964187225, 45831035250, 1159525191825, 31623414322500, 924984868933125, 28887988983603750**

---

## Eulerian numbers

**Réf.** R1 215. DB1 151. JCT 1 351 66. DKB 260. C1 243.

**HIS2 A0460** **Approximants de Padé**

**HIS1 N2047** **Fraction rationnelle**

$$\frac{z (1 + z - 4 z^2)}{(1 - z)^3 (1 - 2 z)^2 (1 - 3 z)}$$

**0, 1, 11, 66, 302, 1191, 4293, 14608, 47840, 152637, 478271, 1479726, 4537314, 13824739, 41932745**

## Rencontres numbers

**Réf.** R1 65.

**HIS2 A0475** **Approximants de Padé** **Suite P-récurrente**

**HIS1 N2132** **exponentielle**

$$a(n) = (2n-1)\,a(n-1) - 5\,a(n-2) - 10\,a(n-3) + (5n-10)\,a(n-4)$$
$$(6n-5)\,a(n-5) + (2n-1)\,a(n-6)$$

$$\frac{z^8 - 16z^7 + 94z^6 - 280z^5 + 481z^4 - 496z^3 + 312z^2 - 96z + 24}{24\,(1-z)^5\,\exp(z)}$$

**1, 0, 15, 70, 630, 5544, 55650, 611820, 7342335, 95449640, 1336295961, 20044438050, 320711010620, 5452087178160, 98137569209940, 1864613814984984**

## Associated Stirling numbers

**Réf.** R1 76. DB1 296. C1 222.

**HIS2 A0478** **Approximants de Padé**

**HIS1 N2138** **Fraction rationnelle**

$$\frac{-12z^3 + 40z^2 - 45z + 15}{(3z-1)\,(2z-1)^2\,(z-1)^3}$$

**15, 105, 490, 1918, 6825, 22935, 74316, 235092, 731731, 2252341, 6879678, 20900922, 63259533**

## Stirling numbers of second kind

**Réf.** AS1 835. DKB 223.

**HIS2 A0481** **Approximants de Padé**

**HIS1 N2141** **Fraction rationnelle**

$$\frac{1}{(1 - z)\,(1 - 2 z)\,(1 - 3 z)\,(1 - 4 z)\,(1 - 5 z)}$$

**1, 15, 140, 1050, 6951, 42525, 246730, 1379400, 7508501, 40075035, 210766920, 1096190550, 5652751651, 28958095545, 147589284710, 749206090500**

## Stirling numbers of first kind

**Réf.** AS1 833. DKB 226.

**HIS2 A0482** **Tableaux généralisés**

**HIS1 N2142** **exponentielle (log)**

$$\frac{\ln(1 - z)^{4}}{24\,(1 - z)}$$

**1, 15, 175, 1960, 22449, 269325, 3416930, 45995730, 657206836, 9957703756, 159721605680, 2706813345600, 48366009233424, 909299905844112**

## Restricted permutations

**Réf.** CMB 4 32 61.

**HIS2 A0496** **Approximants de Padé**

**HIS1 N2231** **Fraction rationnelle**

$$\frac{4\,(6 - z - 2\,z^2 - 4\,z^3 - z^4)}{(1 - z)\,(1 - z - z^2 - z^3)}$$

**24, 44, 80, 144, 260, 476, 872, 1600, 2940, 5404, 9936, 18272, 33604, 61804, 113672, 209072, 384540, 707276, 1300880, 2392688, 4400836, 8094396, 14887912**

## Related to remainder in gaussian quadrature

**Réf.** MOC 1 53 43.

**HIS2 A0515** **hypergéométrique** **Suite P-récurrente**

**HIS1 N2087** **Intégrales elliptiques**

**(n - 1)^2 a(n) = 4 (2 n - 1) (2 n - 3) a(n - 1)**

$$_2F_1([1/2,\ 3/2],\ [1],\ 16\,z)$$

**1, 12, 180, 2800, 44100, 698544, 11099088, 176679360, 2815827300, 44914183600, 716830370256, 11445589052352, 182811491808400, 2920656969720000**

**Réf.** R1 16. MAS 31 79 63.

**HIS2 A0522** **Dérivée exponentielle** **Suite P-récurrente**

**HIS1 N0589**

**a(n) = a(n-1) n + (2 - n) a(n-2)**

$$\frac{\exp(z)}{1 - z}$$

**1, 2, 5, 16, 65, 326, 1957, 13700, 109601, 986410, 9864101, 108505112, 1302061345, 16926797486, 236975164805, 3554627472076, 56874039553217**

---

**Powers of rooted tree enumerator**

**Réf.** R1 150.

**HIS2 A0529** **Approximants de Padé**

**HIS1 N2202** **Fraction rationnelle**

$$\frac{(z - 2)\,(3 z^3 - 12 z^2 + 18 z - 10)}{(1 - z)^6}$$

**20, 74, 186, 388, 721, 1236, 1995, 3072, 4554, 6542, 9152, 12516, 16783, 22120, 28713, 36768, 46512, 58194, 72086, 88484, 107709, 130108, 156055, 185952**

## Sums of cubes

**Réf.** AS1 813.

**HIS2 A0537** **Approximants de Padé**

**HIS1 N1972** **Fraction rationnelle**

$$\frac{1 + 4 z + z}{(1 - z)^{5}}$$

**1, 9, 36, 100, 225, 441, 784, 1296, 2025, 3025, 4356, 6084, 8281, 11025, 14400, 18496, 23409, 29241, 36100, 44100, 53361, 64009, 76176, 90000, 105625, 123201**

---

## Sums of fourth powers

**Réf.** AS1 813.

**HIS2 A0538** **Approximants de Padé**

**HIS1 N2179** **Fraction rationnelle**

$$\frac{(1 + z) (z^{2} + 10 z + 1)}{(z - 1)^{6}}$$

**1, 17, 98, 354, 979, 2275, 4676, 8772, 15333, 25333, 39974, 60710, 89271, 127687, 178312, 243848, 327369, 432345, 562666, 722666, 917147, 1151403, 1431244**

## Sums of 5th powers

**Réf.** AS1 813.

**HIS2 A0539** **Approximants de Padé**

**HIS1 N2280** **Fraction rationnelle**

$$\frac{1 + 26 z + 66 z^{2} + 26 z^{3} + z^{4}}{(1 - z)^{7}}$$

**1, 33, 276, 1300, 4425, 12201, 29008, 61776, 120825, 220825, 381876, 630708, 1002001, 1539825, 2299200, 3347776, 4767633, 6657201, 9133300, 12333300**

## Sums of 6th powers

**Réf.** AS1 813.

**HIS2 A0540** **Approximants de Padé**

**HIS1 N2322** **Fraction rationnelle**

$$\frac{(1 + z) (z^{4} + 56 z^{3} + 246 z^{2} + 56 z + 1)}{(z - 1)^{8}}$$

**1, 65, 794, 4890, 20515, 67171, 184820, 446964, 978405, 1978405, 3749966, 6735950, 11562759, 19092295, 30482920, 47260136, 71397705, 105409929, 152455810**

## Sums of 7th powers

**Réf.** AS1 815.

**HIS2 A0541** **Dérivée logarithmique**

**HIS1 N2343** **Fraction rationnelle**

$$\frac{z^6 + 120 z^5 + 1191 z^4 + 2416 z^3 + 1191 z^2 + 120 z + 1}{(z - 1)^9}$$

**1, 129, 2316, 18700, 96825, 376761, 1200304, 3297456, 8080425, 18080425, 37567596, 73399404, 136147921, 241561425, 412420800, 680856256, 1091194929**

## Sums of eighth powers

**Réf.** AS1 815.

**HIS2 A0542** **Recoupements**

**HIS1 N2358** **Fraction rationnelle**

$$\frac{1 + 247 z + 4293 z^2 + 15619 z^3 + 15619 z^4 + 4293 z^5 + 247 z^6 + z^7}{(1 - z)^{10}}$$

**1, 257, 6818, 72354, 462979, 2142595, 7907396, 24684612, 67731333, 167731333, 382090214, 812071910, 1627802631, 3103591687, 5666482312**

## Discordant permutations

**Réf.** SMA 20 23 54.

**HIS2 A0561** **Approximants de Padé**

**HIS1 N1773** **Fraction rationnelle**

$$\frac{4 z^{3} - 5 z^{2} - 20 z - 6}{(1 - z)^{4}}$$

**6, 44, 145, 336, 644, 1096, 1719, 2540, 3586, 4884, 6461, 8344, 10560, 13136, 16099, 19476, 23294, 27580, 32361, 37664, 43516, 49944, 56975, 64636, 72954, 81956**

## Discordant permutations

**Réf.** SMA 20 23 54.

**HIS2 A0562** **Approximants de Padé**

**HIS1 N1994** **Fraction rationnelle**

$$\frac{9 + 50 z + 35 z^{2} - 15 z^{3} + 4 z^{4} - 2 z^{5}}{(1 - z)^{5}}$$

**9, 95, 420, 1225, 2834, 5652, 10165, 16940, 26625, 39949, 57722, 80835, 110260, 147050, 192339, 247342, 313355, 391755, 484000, 591629, 716262, 859600**

## Discordant permutations

**Réf.** SMA 20 23 54.

**HIS2 A0563** **Approximants de Padé**

**HIS1 N2109** **Fraction rationnelle**

$$\frac{8 z^5 + 6 z^4 - 10 z^3 + 128 z^2 + 114 z + 13}{(1 - z)^6}$$

**13, 192, 1085, 3880, 10656, 24626, 50380, 94128, 163943, 270004, 424839, 643568, 944146, 1347606, 1878302, 2564152, 3436881, 4532264, 5890369, 7555800**

## Discordant permutations

**Réf.** SMA 20 23 54.

**HIS2 A0564** **Approximants de Padé**

**HIS1 N2208** **Fraction rationnelle**

$$\frac{2 z^7 + 4 z^6 - 36 z^5 + 29 z^4 + 72 z^3 + 411 z^2 + 231 z + 20}{(1 - z)^7}$$

**20, 371, 2588, 11097, 35645, 94457, 218124, 454220, 872648, 1571715, 2684936, 4388567, 6909867, 10536089, 15624200, 22611330, 32025950, 44499779**

## Discordant permutations

**Réf.** SMA 20 23 54.

**HIS2 A0565** **Approximants de Padé**

**HIS1 N2275** **Fraction rationnelle**

$$\frac{12 z^7 - 6 z^6 + 88 z^5 - 131 z^4 - 548 z^3 - 1123 z^2 - 448 z - 31}{(1 - z)^8}$$

**31, 696, 5823, 29380, 108933, 327840, 848380, 1958004, 4130895, 8107024, 14990889, 26372124, 44470165, 72305160, 113897310, 174496828, 260846703**

## Heptagonal numbers

**Réf.** D1 2 2. B1 189.

**HIS2 A0566** **Approximants de Padé**

**HIS1 N1826** **Fraction rationnelle**

$$\frac{1 + 4 z}{(1 - z)^3}$$

**1, 7, 18, 34, 55, 81, 112, 148, 189, 235, 286, 342, 403, 469, 540, 616, 697, 783, 874, 970, 1071, 1177, 1288, 1404, 1525, 1651, 1782, 1918, 2059, 2205, 2356, 2512, 2673**

## Octagonal numbers

**Réf.** D1 2 1. B1 189.

**HIS2 A0567** **Approximants de Padé**

**HIS1 N1901** **Fraction rationnelle**

$$\frac{1 + 5z}{(1 - z)^3}$$

**1, 8, 21, 40, 65, 96, 133, 176, 225, 280, 341, 408, 481, 560, 645, 736, 833, 936, 1045, 1160, 1281, 1408, 1541, 1680, 1825, 1976, 2133, 2296, 2465, 2640, 2821, 3008**

---

## From expansion (1+x+x ^ 2 )^ n

**Réf.** JCT 1 372 66. C1 78.

**HIS2 A0574** **Approximants de Padé**

**HIS1 N1219** **Fraction rationnelle**

$$\frac{3 - 2z}{(1 - z)^6}$$

**3, 16, 51, 126, 266, 504, 882, 1452, 2277, 3432, 5005, 7098, 9828, 13328, 17748, 23256, 30039, 38304, 48279, 60214, 74382, 91080, 110630, 133380, 159705, 190008**

## Cubes

**Réf.** BA9.

**HIS2 A0578** **Approximants de Padé**

**HIS1 N1905** **Fraction rationnelle**

$$\frac{1 + 4z + z^2}{(z - 1)^4}$$

**1, 8, 27, 64, 125, 216, 343, 512, 729, 1000, 1331, 1728, 2197, 2744, 3375, 4096, 4913, 5832, 6859, 8000, 9261, 10648, 12167, 13824, 15625, 17576, 19683**

## Binomial coefficients C(n,6)

**Réf.** D1 2 7. RS3. B1 196. AS1 828.

**HIS2 A0579** **Approximants de Padé**

**HIS1 N1847** **Fraction rationnelle**

$$\frac{1}{(1 - z)^7}$$

**1, 7, 28, 84, 210, 462, 924, 1716, 3003, 5005, 8008, 12376, 18564, 27132, 38760, 54264, 74613, 100947, 134596, 177100, 230230, 296010, 376740, 475020, 593775**

## Binomial coefficients C(n,7)

**Réf.** D1 2 7. RS3. B1 196. AS1 828.

**HIS2 A0580** **Approximants de Padé**

**HIS1 N1911** **Fraction rationnelle**

$$\frac{1}{(1 - z)^{8}}$$

**1, 8, 36, 120, 330, 792, 1716, 3432, 6435, 11440, 19448, 31824, 50388, 77520, 116280, 170544, 245157, 346104, 480700, 657800, 888030, 1184040, 1560780, 2035800**

## Binomial coefficients C(n,8)

**Réf.** D1 2 7. RS3. B1 196. AS1 828.

**HIS2 A0581** **Approximants de Padé**

**HIS1 N1976** **Fraction rationnelle**

$$\frac{1}{(1 - z)^{9}}$$

**1, 9, 45, 165, 495, 1287, 3003, 6435, 12870, 24310, 43758, 75582, 125970, 203490, 319770, 490314, 735471, 1081575, 1562275, 2220075, 3108105, 4292145**

## Binomial coefficients C(n,9)

**Réf.** D1 2 7. RS3. B1 196. AS1 828.

**HIS2 A0582** **Approximants de Padé**

**HIS1 N2013** **Fraction rationnelle**

$$\frac{1}{(1 - z)^{10}}$$

**1, 10, 55, 220, 715, 2002, 5005, 11440, 24310, 48620, 92378, 167960, 293930, 497420, 817190, 1307504, 2042975, 3124550, 4686825, 6906900, 10015005, 14307150**

## Fourth powers

**Réf.** BA9.

**HIS2 A0583** **Approximants de Padé**

**HIS1 N2154** **Fraction rationnelle**

$$\frac{(1 + z)\,(z^2 + 10\,z + 1)}{(1 - z)^5}$$

**1, 16, 81, 256, 625, 1296, 2401, 4096, 6561, 10000, 14641, 20736, 28561, 38416, 50625, 65536, 83521, 104976, 130321, 160000, 194481, 234256, 279841, 331776**

## 5th powers

**Réf.** BA9.

**HIS2 A0584** **Approximants de Padé**

**HIS1 N2277** **Fraction rationnelle**

$$\frac{1 + 26 z + 66 z^2 + 26 z^3 + z^4}{(1 - z)^6}$$

**1, 32, 243, 1024, 3125, 7776, 16807, 32768, 59049, 100000, 161051, 248832, 371293, 537824, 759375, 1048576, 1419857, 1889568, 2476099, 3200000, 4084101**

## Partitions of n into distinct primes

**Réf.** PNISI 21 186 55. PURB 107 285 57.

**HIS2 A0586** **Euler**

**HIS1 N0004** **Produit infini**

$$\prod_{n \geq 1} (1 + z^{c(n)})$$

**c(n) = 2,3,5,7,11,... Les nombres premiers**

**1, 0, 1, 1, 0, 2, 0, 2, 1, 1, 2, 1, 2, 2, 2, 2, 3, 2, 4, 3, 4, 4, 4, 5, 5, 5, 6, 5, 6, 7, 6, 9, 7, 9, 9, 9, 11, 11, 11, 13, 12, 14, 15, 15, 17, 16, 18, 19, 20, 21, 23, 22, 25, 26, 27, 30, 29, 32, 32, 35, 37, 39, 40, 42**

**Réf.** JIA 76 153 50. FQ 7 448 69.

**HIS2 A0587** **Recoupements** **1/A0296**

**HIS1 N0755** **exponentielle**

$$\frac{1}{\exp(\exp(z) - 1 - z)}$$

**1, 0, 1, 1, 2, 9, 9, 50, 267, 413, 2180, 17731, 50533, 110176, 1966797, 9938669, 8638718, 278475061, 2540956509, 9816860358, 27172288399, 725503033401**

---

**Réf.** QAM 14 407 56. MOC 29 216 75. FQ 14 397 76.

**HIS2 A0588** **Hypergéométrique** **Suite P-récurrente**

**HIS1 N1866** **algébrique**

**$_2F_1$ ([4, 7/2], [8], 4 z)**

$$\frac{128\ z}{\left(1 + (1 - 4\ z)^{1/2}\right)^{7}}$$

**1, 7, 35, 154, 637, 2548, 9996, 38760, 149226, 572033, 2187185, 8351070, 31865925, 121580760, 463991880, 1771605360, 6768687870, 25880277150**

**Réf.** QAM 14 407 56. MOC 29 216 75.
**HIS2 A0589** **Hypergéométrique** **Suite P-récurrente**
**HIS1 N2048** **algébrique**
**2F1 ([6, 11/2], [12], 4 z)**

$$\frac{1}{\left(1/2 + 1/2\,(1 - 4\,z)^{1/2}\right)^{11}}$$

**1, 11, 77, 440, 2244, 10659, 48279, 211508, 904475, 3798795, 15737865, 64512240, 262256280, 1059111900, 4254603804, 17018415216, 67837293986**

---

**Réf.** QAM 14 407 56. MOC 29 216 75.
**HIS2 A0590** **Hypergéométrique** **Suite P-récurrente**
**HIS1 N2104** **algébrique**
**2F1 ([13/2, 7], [14], 4 z)**

$$\frac{1}{\left(1/2 + 1/2\,(1 - 4\,z)^{1/2}\right)^{13}}$$

**1, 13, 104, 663, 3705, 19019, 92092, 427570, 1924065, 8454225, 36463440, 154969620, 650872404, 2707475148, 11173706960, 45812198536, 186803188858**

## Ramanujan τ function

**Réf.** PLMS 51 4 50. MOC 24 495 70.

**HIS2 A0594** **Euler**

**HIS1 N2237** **Produit infini**

$$\prod_{n\geq 1} \frac{1}{(1 - z^n)^{c(n)}}$$

$$c(n) = -24, -24, -24, -24, \ldots$$

**1, 24, 252, 1472, 4830, 6048, 16744, 84480, 113643, 115920, 534612, 370944, 577738, 401856, 1217160, 987136, 6905934, 2727432, 10661420, 7109760, 4219488**

## Central factorial numbers

**Réf.** RCI 217.

**HIS2 A0596** **Approximants de Padé**

**HIS1 N1505** **Fraction rationnelle**

$$\frac{4 + 21 z + 14 z^2 + z^3}{(1 - z)^7}$$

**4, 49, 273, 1023, 3003, 7462, 16422, 32946, 61446, 108031, 180895, 290745, 451269, 679644, 997084, 1429428, 2007768, 2769117, 3757117, 5022787, 6625311**

## Central factorial numbers

**Réf.** RCI 217.

**HIS2 A0597** **Dérivée logarithmique**

**HIS1 N2287** **Fraction rationnelle**

$$\frac{z^5 + 75 z^4 + 603 z^3 + 1065 z^2 + 460 z + 36}{(z - 1)^{10}}$$

**36, 820, 7645, 44473, 191620, 669188, 1999370, 5293970, 12728936, 28285400, 58856655, 115842675, 217378200, 391367064, 679524340, 1142659012**

## A partition function

**Réf.** CAY 2 278. JACS 53 3084 31. AMS 26 304 55.

**HIS2 A0601** **Approximants de Padé** *** titre modifié**

**HIS1 N0392** **Fraction rationnelle**

$$\frac{1}{(1 + z)(z^2 + z + 1)(z - 1)^4}$$

**1, 2, 4, 7, 11, 16, 23, 31, 41, 53, 67, 83, 102, 123, 147, 174, 204, 237, 274, 314, 358, 406, 458, 514, 575, 640, 710, 785, 865, 950, 1041, 1137, 1239, 1347, 1461, 1581**

## Partitions of n into prime parts

**Réf.** PNISI 21 183 55. AMM 95 711 88.

**HIS2 A0607** **Euler**

**HIS1 N0093** **Produit infini**

$$\prod_{n\geq 1} \frac{1}{(1 - z^{c(n)})}$$

**c(n) = 2,3,5,7,...,les nombres premiers**

**1, 0, 1, 1, 1, 2, 2, 3, 3, 4, 5, 6, 7, 9, 10, 12, 14, 17, 19, 23, 26, 30, 35, 40, 46, 52, 60, 67, 77, 87, 98, 111, 124, 140, 157, 175, 197, 219, 244, 272, 302, 336, 372, 413, 456, 504, 557**

## Preferential arrangements of n things

**Réf.** CAY 4 113. PLMS 22 341 1891. AMM 69 7 62. PSPM 19 172 71. DM 48 102 84.

**HIS2 A0670** **Inverse fonctionnel**

**HIS1 N1191** **exponentielle**

$$\frac{1 - \exp(z)}{\exp(z) - 2}$$

**1, 1, 3, 13, 75, 541, 4683, 47293, 545835, 7087261, 102247563, 1622632573, 28091567595, 526858348381, 10641342970443, 230283190977853**

**Réf.** QJM 47 110 16. FMR 1 112. DA63 2 283. PSAM 15 101 63.

**HIS2 A0680** **Hypergéométrique** **Suite P-récurrente**

**HIS1 N1793** **algébrique** **f.g. exponentielle double**

**Γ(2n+1)/2^n**

**a(n) = n (2n-1) a(n-1)**

$$\frac{1}{(1 - 2 z)^{1/2}}$$

**1, 6, 90, 2520, 113400, 7484400, 681080400, 81729648000, 12504636144000, 2375880867360000, 548828480360160000, 151476660579404160000**

---

## Stochastic matrices of integers

**Réf.** PSAM 15 101 63. SS70.

**HIS2 A0681** **équations différentielles** **Suite P-récurrente**

**HIS1 N1250** **exponentielle (algébrique)** **Formule de B. Salvy**

**a(n) = - 1/2 (n - 1) (- 2 n + 2) a(n - 1) - 1/2 (n - 1) (n^2 - 4 n + 4) a(n - 2)**

$$\frac{\exp(z/2)}{(1 - z)^{1/2}}$$

**1, 1, 3, 21, 282, 6210, 202410, 9135630, 545007960, 41514583320, 3930730108200, 452785322266200, 62347376347779600, 10112899541133589200**

## Partitions of n into distinct odd parts

**Réf.** PLMS 42 553 36. CJM 4 383 52.

**HIS2 A0700** **Euler**

**HIS1 N0078** **Produit infini**

$$\prod_{n\geq 1} (1 + z^{c(n)})$$

$$c(n) = 1,3,5,7,9,11,13,...$$

**1, 1, 0, 1, 1, 1, 1, 1, 2, 2, 2, 2, 3, 3, 3, 4, 5, 5, 5, 6, 7, 8, 8, 9, 11, 12, 12, 14, 16, 17, 18, 20, 23, 25, 26, 29, 33, 35, 37, 41, 46, 49, 52, 57, 63, 68, 72, 78, 87, 93, 98, 107, 117, 125, 133, 144**

---

## Degree n even permutations of order dividing 2

**Réf.** CJM 7 168 55.

**HIS2 A0704** **équations différentielles** **Formule de B. Salvy**

**HIS1 N1427** **exponentielle**

$$\exp(z)\ \cosh(z^2 / 2)$$

**1, 1, 1, 1, 4, 16, 46, 106, 316, 1324, 5356, 18316, 63856, 272416, 1264264, 5409496, 22302736, 101343376, 507711376, 2495918224, 11798364736, 58074029056**

## Partitions of n into parts of 2 kinds

**Réf.** RS4 90. RCI 199.

**HIS2 A0710** **Euler**

**HIS1 N0535** **Produit infini**

$$\prod_{n\geq 1} \frac{1}{(1 - z^{n})^{c(n)}}$$

$$c(n) = 2,2,2,2,1,1,1,1,1,1,1,1,...$$

**1, 2, 5, 10, 20, 35, 62, 102, 167, 262, 407, 614, 919, 1345, 1952, 2788, 3950, 5524, 7671, 10540, 14388, 19470, 26190, 34968, 46439, 61275, 80455, 105047, 136541**

## Partitions of n into parts of 3 kinds

**Réf.** RS4 122.

**HIS2 A0711** **Euler**

**HIS1 N1122** **Produit infini**

$$\prod_{n\geq 1} \frac{1}{(1 - z^{n})^{c(n)}}$$

$$c(n) = 3,3,3,3,2,2,2,2,2,2,2,2,....$$

**1, 3, 9, 22, 51, 107, 217, 416, 775, 1393, 2446, 4185, 7028, 11569, 18749, 29908, 47083, 73157, 112396, 170783, 256972, 383003, 565961, 829410, 1206282, 1741592**

## Partitions of n into parts of 2 kinds

**Réf.** RS4 90. RCI 199.

**HIS2 A0712** **Euler**

**HIS1 N0536** **Produit infini**

$$\prod_{n\geq 1} \frac{1}{(1 - z^n)^{c(n)}}$$

$$c(n) = 2,2,2,2,2,2,2,2,....$$

**1, 2, 5, 10, 20, 36, 65, 110, 185, 300, 481, 752, 1165, 1770, 2665, 3956, 5822, 8470, 12230, 17490, 24842, 35002, 49010, 68150, 94235, 129512, 177087, 240840**

## Partitions of n into parts of 3 kinds

**Réf.** RS4 122.

**HIS2 A0713** **Euler** **différences de A0712**

**HIS1 N1096** **Produit infini**

$$\prod_{n\geq 1} \frac{1}{(1 - z^n)^{c(n)}}$$

$$c(n) = 3,2,2,2,2,2,2,2,....$$

**1, 3, 8, 18, 38, 74, 139, 249, 434, 734, 1215, 1967, 3132, 4902, 7567, 11523, 17345, 25815, 38045, 55535, 80377, 115379, 164389, 232539, 326774, 456286, 633373**

## Partitions of n into parts of 3 kinds

**Réf.** RS4 122.

**HIS2 A0714** **Euler**

**HIS1 N1117** **Produit infini**

$$\prod_{n\geq 1} \frac{1}{(1 - z^n)^{c(n)}}$$

$$c(n) = 3, 3, 2, 2, 2, 2, 2, 2, 2, 2, \ldots.$$

**1, 3, 9, 21, 47, 95, 186, 344, 620, 1078, 1835, 3045, 4967, 7947, 12534, 19470, 29879, 45285, 67924, 100820, 148301, 216199, 312690, 448738, 639464, 905024**

## Partitions of n into parts of 3 kinds

**Réf.** RS4 122.

**HIS2 A0715** **Euler**

**HIS1 N1121** **Produit infini**

$$\prod_{n\geq 1} \frac{1}{(1 - z^n)^{c(n)}}$$

$$c(n) = 3, 3, 3, 2, 2, 2, 2, 2, 2, 2, 2, \ldots.$$

**1, 3, 9, 22, 50, 104, 208, 394, 724, 1286, 2229, 3769, 6253, 10176, 16303, 25723, 40055, 61588, 93647, 140875, 209889, 309846, 453565, 658627, 949310, 1358589**

## Partitions of n into parts of 3 kinds

**Réf.** RS4 122.

**HIS2 A0716** **Euler**

**HIS1 N1123** **Produit infini**

$$\prod_{n\geq 1} \frac{1}{(1-z^n)^{c(n)}}$$

$$c(n) = 3,3,3,3,\ldots.$$

**1, 3, 9, 22, 51, 108, 221, 429, 810, 1479, 2640, 4599, 7868, 13209, 21843, 35581, 57222, 90882, 142769, 221910, 341649, 521196, 788460, 1183221, 1762462, 2606604**

## Partitions of n into parts prime to 3

**Réf.** PSPM 8 145 65.

**HIS2 A0726** **Euler**

**HIS1 N0116** **Produit infini**

$$\prod_{n\geq 1} \frac{1}{(1-z^n)^{c(n)}}$$

$$c(n) = 1 \text{ si } n = 1 \text{ ou } 2 \bmod 3.$$

**1, 1, 2, 2, 4, 5, 7, 9, 13, 16, 22, 27, 36, 44, 57, 70, 89, 108, 135, 163, 202, 243, 297, 355, 431, 513, 617, 731, 874, 1031, 1225, 1439, 1701, 1991, 2341, 2731, 3197, 3717**

**Réf.** KNAW 59 207 56.

**HIS2 A0727** **Recoupements**

**HIS1 N1296** **Produit infini**

**La suite est alternée**

$$\prod_{n\geq 1} \frac{1}{(1-z^n)^{c(n)}}$$

**c(n) = -4,-4,-4,-4,...**

**1, 4, 2, 8, 5, 4, 10, 8, 9, 0, 14, 16, 10, 4, 0, 8, 14, 20, 2, 0, 11, 20, 32, 16, 0, 4, 14, 8, 9, 20, 26, 0, 2, 28, 0, 16, 16, 28, 22, 0, 14, 16, 0, 40, 0, 28, 26, 32, 17, 0, 32, 16, 22, 0, 10**

---

**Réf.** KNAW 59 207 56.

**HIS2 A0729** **Recoupements** **La suite est alternée**

**HIS1 N1691** **Produit infini**

$$\prod_{n\geq 1} \frac{1}{(1-z^n)^{c(n)}}$$

**c(n) = -6,-6,-6,-6,-6,...**

**1, 6, 9, 10, 30, 0, 11, 42, 0, 70, 18, 54, 49, 90, 0, 22, 60, 0, 110, 0, 81, 180, 78, 0, 130, 198, 0, 182, 30, 90, 121, 84, 0, 0, 210, 0, 252, 102, 270, 170, 0, 0, 69, 330, 0, 38**

**Réf.** QJM 38 56 07. KNAW 59 207 56. GMJ 8 29 67.

**HIS2 A0735** **Euler** **Inverse de A5758 alternée en signe**

**HIS1 N2069** **Produit infini**

**La suite est alternée**

$$\prod_{n\geq 1} \frac{1}{(1 - z^n)^{c(n)}}$$

$$c(n) = -12, -12, -12, -12, \ldots$$

**1, 12, 54, 88, 99, 540, 418, 648, 594, 836, 1056, 4104, 209, 4104, 594, 4256, 6480, 4752, 298, 5016, 17226, 12100, 5346, 1296, 9063, 7128, 19494, 29160, 10032, 7668**

---

**Réf.** PLMS 31 341 30. SPS 37-40-4 209 66.

**HIS2 A0757** **équations différentielles** **Formule de B. Salvy**

**HIS1 N1915** **exponentielle** **f.g. exponentielle**

**a(n) = 2 n a(n - 2) + n a(n - 3) + (n - 1) a(n - 1)**

$$(- \ln(- z + 1) + 1) \exp(- z)$$

**0, 0, 1, 1, 8, 36, 229, 1625, 13208, 120288, 1214673, 13469897, 162744944, 2128047988, 29943053061, 451123462673, 7245940789072, 123604151490592, 2231697509543361**

## Stirling numbers of second kind

**Réf.** AS1 835. DKB 223.

**HIS2 A0770** **Approximants de Padé**

**HIS1 N2215** **Fraction rationnelle**

$$\frac{1}{(1 - z)(1 - 2z)(1 - 3z)(1 - 4z)(1 - 5z)(1 - 6z)}$$

**1, 21, 266, 2646, 22827, 179487, 1323652, 9321312, 63436373, 420693273, 2734926558, 17505749898, 110687251039, 693081601779, 4306078895384**

## Stirling numbers of second kind

**Réf.** AS1 835. DKB 223.

**HIS2 A0771** **Approximants de Padé**

**HIS1 N2263** **Fraction rationnelle**

$$\frac{1}{(1 - z)(1 - 2z)(1 - 3z)(1 - 4z)(1 - 5z)(1 - 6z)(1 - 7z)}$$

**1, 28, 462, 5880, 63987, 627396, 5715424, 49329280, 408741333, 3281882604, 25708104786, 197462483400, 1492924634839, 11143554045652**

**Réf.** CMB 8 627 65. JRM 4 168 71. FQ 27 16 89.
**HIS2 A0792** **Approximants de Padé**
**HIS1 N0205** **Fraction rationnelle**

$$\frac{1 + 2z + 3z^2 + z^3}{1 - 3z^3}$$

**1, 2, 3, 4, 6, 9, 12, 18, 27, 36, 54, 81, 108, 162, 243, 324, 486, 729, 972, 1458, 2187, 2916, 4374, 6561, 8748, 13122, 19683, 26244, 39366, 59049, 78732, 118098**

---

**Réf.** CMB 7 262 64. JCT 7 315 69.
**HIS2 A0803** **Approximants de Padé**
**HIS1 N2232** **Fraction rationnelle**

$$\frac{24 - 4z - 8z^2 - 16z^3}{1 - 2z + z^4}$$

**24, 44, 80, 144, 264, 484, 888, 1632, 3000, 5516, 10144, 18656, 34312, 63108, 116072, 213488, 392664, 722220, 1328368, 2443248, 4493832, 8265444**

**Réf.** CJM 8 308 56.
**HIS2 A0806** **équations différentielles** **Suite P-récurrente**
**HIS1 N1651** **exponentielle (algébrique)** **Formule de B. Salvy**
**a(n) = (2 n + 1) a(n - 1) + a(n - 2)**

$$-\frac{-4 + 3\,(1 - 2\,z)^{1/2} + 2\,z}{\exp(1 - (1 - 2\,z)^{1/2})\,(1 - 2\,z)^{5/2}}$$

**0, 1, 5, 36, 329, 3655, 47844, 721315, 12310199, 234615096, 4939227215, 113836841041, 2850860253240, 77087063678521, 2238375706930349**

---

**Réf.** LU91 1 221.
**HIS2 A0898** **Dérivée logarithmique** **Suite P-récurrente**
**HIS1 N0645** **exponentielle**
**a(n) = 2 a(n - 1) + (2 n - 4) a(n - 2)**

$$\exp(2\,z + z^2)$$

**1, 2, 6, 20, 76, 312, 1384, 6512, 32400, 168992, 921184, 5222208, 30710464, 186753920, 1171979904, 7573069568, 50305536256, 342949298688, 2396286830080**

## Symmetric permutations

**Réf.** LU91 1 222. LNM 560 201 76.

**HIS2 A0902** **Recoupements exponentielle** **Suite P-récurrente**

**HIS1 N1147**

**a(n) = 2 a(n-1)-(4-2 n) a(n-2)**

$$1/2 \exp(z (2 + z)) + 1/2$$

**1, 3, 10, 38, 156, 692, 3256, 16200**

## Ménage numbers

**Réf.** LU91 1 495.

**HIS2 A0904** **P-récurrences** **Suite P-récurrente**

**HIS1 N1193**

$$a(n) = a(n - 3) + (n + 3) a(n - 2) + (n + 2) a(n - 1)$$

**0, 3, 13, 83, 592, 4821, 43979, 444613, 4934720, 59661255, 780531033, 10987095719, 165586966816, 2660378564777, 45392022568023, 819716784789193**

**Réf.** TOH 37 259 33. JO39 152. DB1 296. C1 256.

**HIS2 A0906** **Hypergéométrique** **Suite P-récurrente**

**HIS1 N0841** **algébrique** **f.g. exponentielle**

**(n - 1) a(n) = (2 n + 1) n a(n - 1)**

$$\frac{z}{(1 - 2z)^{5/2}}$$

**2, 20, 210, 2520, 34650, 540540, 9459450, 183783600, 3928374450, 91662070500, 2319050383650, 63246828645000, 1849969737866250**

## Associated Stirling numbers

**Réf.** TOH 37 259 33. JO39 152. C1 256.

**HIS2 A0907** **Hypergéométrique** **Suite P-récurrente**

**HIS1 N1797** **algébrique** **f.g. exponentielle**

**1/4 a(n) (4 n + 1) (n - 1) = 1/4 a(n - 1) (4 n + 5) (2 n + 1) (n + 1)**

$$\frac{z\,(2z^2 + 33z + 18)}{3\,(1 - 2z)^{9/2}}$$

**6, 130, 2380, 44100, 866250, 18288270, 416215800, 10199989800, 268438920750, 7562120816250, 227266937597700, 7262844156067500**

## Stirling numbers of first kind

**Réf.** AS1 833. DKB 226.

**HIS2 A0914** **Approximants de Padé**

**HIS1 N0789** **Fraction rationnelle**

$$\frac{2 - z}{(1 - z)^5}$$

**2, 11, 35, 85, 175, 322, 546, 870, 1320, 1925, 2717, 3731, 5005, 6580, 8500, 10812, 13566, 16815, 20615, 25025, 30107, 35926, 42550, 50050, 58500, 67977, 78561**

## Stirling numbers of first kind

**Réf.** AS1 833. DKB 226.

**HIS2 A0915** **Dérivée logarithmique**

**HIS1 N2239** **Fraction rationnelle**

$$\frac{z^3 + 22 z^2 + 58 z + 24}{(z - 1)^9}$$

**24, 274, 1624, 6769, 22449, 63273, 157773, 357423, 749463, 1474473, 2749747, 4899622, 8394022, 13896582, 22323822, 34916946, 53327946, 79721796**

## 2 ^ (n-2)

**Réf.** VO11 31. DA63 2 212. R1 33.

**HIS2 A0918** **Approximants de Padé**

**HIS1 N0625** **Fraction rationnelle**

$$\frac{z}{(1 - 2z)(1 - z)}$$

**0, 2, 6, 14, 30, 62, 126, 254, 510, 1022, 2046, 4094, 8190, 16382, 32766, 65534, 131070, 262142, 524286, 1048574, 2097150, 4194302, 8388606, 16777214, 33554430**

## Differences of 0

**Réf.** VO11 31. DA63 2 212. R1 33.

**HIS2 A0919** **Approximants de Padé**

**HIS1 N2235** **Fraction rationnelle**

$$\frac{24}{(1 - z)(1 - 2z)(1 - 3z)(1 - 4z)}$$

**24, 240, 1560, 8400, 40824, 186480, 818520, 3498000, 14676024, 60780720, 249401880, 1016542800, 4123173624, 16664094960, 67171367640**

## Differences of 0

**Réf.** VO11 31. DA63 2 212. R1 33.

**HIS2 A0920** **Recoupements**

**HIS1 N2370** **Fraction rationnelle**

$$\frac{720}{(1 - z)(1 - 2 z)(1 - 3 z)(1 - 4 z)(1 - 5 z)(1 - 6 z)}$$

**720, 15120, 191520, 1905120, 16435440, 129230640, 953029440, 6711344640, 45674188560, 302899156560, 1969147121760, 12604139926560**

---

**Réf.** LA62 13. FQ 2 225 64. JA66 91. MMAG 41 15 68.

**HIS2 A0930** **Approximants de Padé**

**HIS1 N0207** **Fraction rationnelle**

$$\frac{1}{1 - z - z^3}$$

**1, 1, 1, 2, 3, 4, 6, 9, 13, 19, 28, 41, 60, 88, 129, 189, 277, 406, 595, 872, 1278, 1873, 2745, 4023, 5896, 8641, 12664, 18560, 27201, 39865, 58425, 85626, 125491, 183916**

**Réf.** JA66 90. MMAG 41 17 68.
**HIS2 A0931** **Approximants de Padé**
**HIS1 N0102** **Fraction rationnelle**

$$\frac{1 + z}{1 - z^2 - z^3}$$

**1, 1, 1, 2, 2, 3, 4, 5, 7, 9, 12, 16, 21, 28, 37, 49, 65, 86, 114, 151, 200, 265, 351, 465, 616, 816, 1081, 1432, 1897, 2513, 3329, 4410, 5842, 7739, 10252, 13581, 17991, 23833**

---

**Genus of complete graph on n nodes**

**Réf.** PNAS 60 438 68.
**HIS2 A0933** **Approximants de Padé** **conjecture**
**HIS1 N0182** **Fraction rationnelle**

$$\frac{z^4 (1 - z + z^2 - z^3 + z^4)}{(z^2 + z + 1)(1 + z^2)(1 - z)^3}$$

**0, 0, 0, 0, 1, 1, 1, 2, 3, 4, 5, 6, 8, 10, 11, 13, 16, 18, 20, 23, 26, 29, 32, 35, 39, 43, 46, 50, 55, 59, 63, 68, 73, 78, 83, 88, 94, 100, 105, 111, 118, 124, 130, 137, 144, 151, 158, 165, 173, 181**

## Fine's sequence: relations of valence ≥1 on an n-set

**Réf.** IC 16 352 70. JCT A23 90 77. DM 19 101 77.

**HIS2 A0957** **LLL** **Suite P-récurrente**
**HIS1 N0635** **algébrique**

**(n + 2) a(n) = (7/2 n + 1) a(n - 1) + (2 n + 1) a(n - 2)**

$$\frac{1}{2}\,\frac{1 - 2z - 2z^2 - (1 - 4z)^{1/2}}{2z^3 + z^4}$$

**1, 2, 6, 18, 57, 186, 622, 2120, 7338, 25724, 91144, 325878, 1174281, 4260282, 15548694, 57048048, 210295326, 778483932, 2892818244, 10786724388**

---

## A simple recurrence

**Réf.** IC 16 351 70.

**HIS2 A0958** **LLL** **Suite P-récurrente**
**HIS1 N1104** **algébrique**

$$\frac{1 - z - 4z^2 - 2z^3 - (-(4z - 1)(z + 1)^2)^{1/2}}{2(2z^3 + z^4)}$$

**1, 3, 8, 24, 75, 243, 808, 2742, 9458, 33062, 116868, 417022, 1500159, 5434563, 19808976, 72596742, 267343374, 988779258, 3671302176, 13679542632**

## A ternary continued fraction

**Réf.** TOH 37 441 33.

**HIS2 A0962** **Approximants de Padé**

**HIS1 N0582** **Fraction rationnelle**

$$\frac{(1+z)\,(2z^4-7z^3+6z^2+z-1)}{z^6-3z^4+7z^2-1}$$

**1, 0, 0, 1, 2, 5, 15, 32, 99, 210, 650, 1379, 4268, 9055, 28025, 59458, 184021, 390420, 1208340, 2563621, 7934342, 16833545, 52099395, 110534372, 342101079, 725803590**

## A ternary continued fraction

**Réf.** TOH 37 441 33.

**HIS2 A0963** **Approximants de Padé**

**HIS1 N1062** **Fraction rationnelle**

$$\frac{1-4z^2+7z^3-2z^4}{1-7z^2+3z^4-z^6}$$

**0, 1, 0, 3, 7, 16, 49, 104, 322, 683, 2114, 4485, 13881, 29450, 91147, 193378, 598500, 1269781, 3929940, 8337783, 25805227, 54748516, 169445269, 359496044, 1112631142**

## n! never ends in this many 0's

**Réf.** MMAG 27 55 53.

**HIS2 A0966** **Approximants de Padé**

**HIS1 N1557** **Fraction rationnelle**

$$\frac{5 + 6 z + 6 z^2 + 6 z^3 + 6 z^4 + z^5 + z^6}{1 - z - z^6 + z^7}$$

**5, 11, 17, 23, 29, 30, 36, 42, 48, 54, 60, 61, 67, 73, 79, 85, 91, 92, 98, 104, 110, 116, 122, 123, 129, 135, 141, 147, 153, 154, 155**

## Fermat coefficients

**Réf.** MMAG 27 141 54.

**HIS2 A0969** **Approximants de Padé**

**HIS1 N1042** **Fraction rationnelle**

$$\frac{1 + z + 2 z^2}{(z^2 + z + 1) (1 - z)^3}$$

**1, 3, 7, 12, 18, 26, 35, 45, 57, 70, 84, 100, 117, 135, 155, 176, 198, 222, 247, 273, 301, 330, 360, 392, 425, 459, 495, 532, 570, 610, 651, 693, 737, 782, 828, 876, 925, 975**

## Fermat coefficients

**Réf.** MMAG 27 141 54.

**HIS2 A0970** **Approximants de Padé**

**HIS1 N1846** **Fraction rationnelle**

$$\frac{3 z^5 + 2 z^4 + 4 z^3 + 3 z^2 + 3 z + 1}{(z^4 + z^3 + z^2 + z + 1) (1 - z)^5}$$

**1, 7, 25, 66, 143, 273, 476, 775, 1197, 1771, 2530, 3510, 4750, 6293, 8184, 10472, 13209, 16450, 20254, 24682, 29799, 35673, 42375, 49980, 58565, 68211, 79002**

## Fermat coefficients

**Réf.** MMAG 27 141 54.

**HIS2 A0973** **Approximants de Padé**

**HIS1 N2137** **Fraction rationnelle**

$$\frac{(z + 1) (z^2 + 6 z + 1)}{(z - 1)^8}$$

**1, 15, 99, 429, 1430, 3978, 9690, 21318, 43263, 82225, 148005, 254475, 420732, 672452, 1043460, 1577532, 2330445, 3372291, 4790071, 6690585, 9203634**

## Central binomial coefficients

**Réf.** RS3. AS1 828.

**HIS2 A0984** **Hypergéométrique** **Suite P-récurrente**

**HIS1 N0643** **algébrique**

**$_2F_1$([1/2], [ ], 4 z)**

$$\frac{1}{(1 - 4z)^{1/2}}$$

**1, 2, 6, 20, 70, 252, 924, 3432, 12870, 48620, 184756, 705432, 2704156, 10400600, 40116600, 155117520, 601080390, 2333606220, 9075135300, 35345263800**

---

## Stochastic matrices of integers

**Réf.** DMJ 35 659 68.

**HIS2 A0985** **Dérivée logarithmique** **Suite P-récurrente**

**HIS1 N1168** **exponentielle**

**a(n) = (1/2 n^3 - 9/2 n^2 + 13 n - 12) a(n - 4) + (2 n - 3) a(n - 1)**
**+ (- n^2 + 5 n - 6) a(n - 2) + (- n^2 + 5 n - 6) a(n - 3)**

$$\frac{\exp(z\,(z^3 + z^2 - 2)/(4(1-z)))}{(1 - z)^{1/2}}$$

**1, 1, 3, 11, 56, 348, 2578, 22054, 213798, 2313638, 27627434, 360646314, 5107177312, 77954299144, 1275489929604, 22265845018412, 412989204564572**

## Stochastic matrices of integers

**Réf.** DMJ 35 659 68.

**HIS2 A0986** **Dérivée logarithmique** **Suite P-récurrente**

**HIS1 N1437** **exponentielle (algébrique)**

**a(n) = 2 (2 n - 1) n^2 a(n - 1) - 1/2 (2 n - 1) (12 n^2 - 7 n + 1) a(n - 4) - 1/2 (2 n - 1) (- 8 n^2 + 2 n) a(n - 2)**

$$\frac{\exp\left(\dfrac{z^3 + 3 z^2 - 4 z + 2}{4 (1 - z)}\right)}{(z - 1)^{1/2}}$$

**1, 0, 1, 4, 18, 112, 820, 6912, 66178, 708256, 8372754, 108306280, 1521077404, 23041655136, 374385141832, 6493515450688, 119724090206940**

## Stochastic matrices of integers

**Réf.** DMJ 35 659 68.

**HIS2 A0987** **Dérivée logarithmique** **Suite P-récurrente**

**HIS1 N0707** **exponentielle (algébrique)**

$$\frac{\exp(z (z^3 + z^2 - 2)/(4(1 - z)))}{(1 - z)^{3/2}}$$

**0, 1, 1, 2, 7, 32, 184, 1268, 10186, 93356, 960646, 10959452, 137221954, 1870087808, 27548231008, 436081302248, 7380628161076, 132975267434552**

## 2-line partitions of n

**Réf.** DMJ 31 272 64.

**HIS2 A0990** **Euler**

**HIS1 N0978** **Produit infini**

$$\prod_{n\geq 1} \frac{1}{(1-z^n)^{c(n)}}$$

$$c(n) = 1,2,2,2,2,...$$

**1, 3, 5, 10, 16, 29, 45, 75, 115, 181, 271, 413, 605, 895, 1291, 1866, 2648, 3760, 5260, 7352, 10160, 14008, 19140, 26085, 35277, 47575, 63753, 85175, 113175, 149938**

## 3-line partitions of n

**Réf.** DMJ 31 272 64.

**HIS2 A0991** **Euler**

**HIS1 N1011** **Produit infini**

$$\prod_{n\geq 1} \frac{1}{(1-z^n)^{c(n)}}$$

$$c(n) = 1,2,3,3,3,3,3,3,....$$

**1, 3, 6, 12, 21, 40, 67, 117, 193, 319, 510, 818, 1274, 1983, 3032, 4610, 6915, 10324, 15235, 22371, 32554, 47119, 67689, 96763, 137404, 194211, 272939, 381872**

## Dissections of a polygon

**Réf.** EMN 32 6 40. BAMS 54 359 48.

**HIS2 A1002** **Inverse fonctionnel algébrique** **Suite P-récurrente.**

**HIS1 N1146**

$(n - 1)\, n\, a(n) = (22/5\, n^2 - 11\, n + 33/5)\, a(n - 1) + (27/5\, n^2 - 108/5\, n + 21)\, a(n - 2)$

$$\text{Inverse de } z(1 - z - z^2)$$

**1, 1, 3, 10, 38, 154, 654, 2871, 12925, 59345, 276835, 1308320, 6250832, 30142360, 146510216, 717061938, 3530808798, 17478955570, 86941210950, 434299921440**

## Super Catalan numbers

**Réf.** EMN 32 6 40. BAMS 54 359 48. RCI 168. C1 57. VA91 198.

**HIS2 A1003** **Inverse fonctionnel algébrique** **Suite P-récurrente**

**HIS1 N1163**

$n\, a(n) = (6\, n - 9)\, a(n - 1) + (- n + 3)\, a(n - 2)$

$$\frac{1 + z - (1 - 6\,z + z^2)^{1/2}}{4\,z}$$

**1, 1, 3, 11, 45, 197, 903, 4279, 20793, 103049, 518859, 2646723, 13648869, 71039373, 372693519, 1968801519, 10463578353, 55909013009, 300159426963**

## Partitions of points on a circle

**Réf.** BAMS 54 359 48.

**HIS2 A1005** **Inverse fonctionnel algébrique** **Suite P-récurrente algébrique du 3è degré**

**HIS1 N0520**

$$\frac{1}{2}(2n+1)\,n\,a(n) = \left(\frac{193}{4}n^2 - \frac{1015}{4}n + 327\right)a(n-3)$$

$$+\left(-\frac{37}{4}n^2 + \frac{91}{4}n - 9\right)a(n-1) + \left(\frac{9}{4}n^2 - \frac{9}{4}n - 3\right)a(n-2)$$

$$+\left(\frac{279}{4}n^2 - \frac{1953}{4}n + 837\right)a(n-4)$$

**1, 0, 1, 1, 2, 5, 8, 21, 42, 96, 222, 495, 1177, 2717, 6435, 15288, 36374, 87516, 210494, 509694, 1237736, 3014882, 7370860, 18059899, 44379535, 109298070, 269766655**

## Motzkin numbers

**Réf.** BAMS 54 359 48. JSIAM 18 254 69. JCT A23 292 77.

**HIS2 A1006** **LLL** **Suite P-récurrente**

**HIS1 N0456** **algébrique**

$(n+1)\,a(n) = (2n-1)\,a(n-1) + (3n-6)\,a(n-2)$

$$\frac{1 - z - (1 - 2z - 3z^2)^{1/2}}{2z^2}$$

**1, 1, 2, 4, 9, 21, 51, 127, 323, 835, 2188, 5798, 15511, 41835, 113634, 310572, 853467, 2356779, 6536382, 18199284, 50852019, 142547559, 400763223, 1129760415**

## 6th powers

**Réf.** BA9.

**HIS2 A1014** **Approximants de Padé**

**HIS1 N2318** **Fraction rationnelle**

$$\frac{(1 + z)\,(z^4 + 56\,z^3 + 246\,z^2 + 56\,z + 1)}{(1 - z)^7}$$

1, 64, 729, 4096, 15625, 46656, 117649, 262144, 531441, 1000000, 1771561, 2985984, 4826809, 7529536, 11390625, 16777216, 24137569, 34012224, 47045881

## Seventh powers

**Réf.** BA9.

**HIS2 A1015** **Approximants de Padé**

**HIS1 N2341** **Fraction rationnelle**

$$\frac{z^6 + 120\,z^5 + 1191\,z^4 + 2416\,z^3 + 1191\,z^2 + 120\,z + 1}{(z - 1)^8}$$

1, 128, 2187, 16384, 78125, 279936, 823543, 2097152, 4782969, 10000000, 19487171, 35831808, 62748517, 105413504, 170859375, 268435456, 410338673

## Eighth powers

**Réf.** BA9.

**HIS2 A1016** **Recoupements**

**HIS1 N2357** **Fraction rationnelle**

$$\frac{(z+1)\,(z^6 + 246\,z^5 + 4047\,z^4 + 11572\,z^3 + 4047\,z^2 + 246\,z + 1)}{(z-1)^9}$$

**1, 256, 6561, 65536, 390625, 1679616, 5764801, 16777216, 43046721, 100000000, 214358881, 429981696, 815730721, 1475789056, 2562890625, 4294967296**

## Powers of 8

**Réf.** BA9.

**HIS2 A1018** **Approximants de Padé**

**HIS1 N1937** **Fraction rationnelle**

$$\frac{1}{1 - 8\,z}$$

**1, 8, 64, 512, 4096, 32768, 262144, 2097152, 16777216, 134217728, 1073741824, 8589934592, 68719476736, 549755813888, 4398046511104, 35184372088832**

## Powers of 9

**Réf.** BA9.

**HIS2 A1019** **Approximants de Padé**

**HIS1 N1992** **Fraction rationnelle**

$$\frac{1}{1 - 9z}$$

**1, 9, 81, 729, 6561, 59049, 531441, 4782969, 43046721, 387420489, 3486784401, 31381059609, 282429536481, 2541865828329, 22876792454961**

## Powers of 11

**Réf.** BA9.

**HIS2 A1020** **Approximants de Padé**

**HIS1 N2054** **Fraction rationnelle**

$$\frac{1}{1 - 11z}$$

**1, 11, 121, 1331, 14641, 161051, 1771561, 19487171, 214358881, 2357947691, 25937424601, 285311670611, 3138428376721, 34522712143931**

## Powers of 12

**Réf.** BA9.

**HIS2 A1021** **Approximants de Padé**

**HIS1 N2084** **Fraction rationnelle**

$$\frac{1}{1 - 12 z}$$

**1, 12, 144, 1728, 20736, 248832, 2985984, 35831808, 429981696, 5159780352, 61917364224, 743008370688, 8916100448256, 106993205379072**

## Powers of 13

**Réf.** BA9.

**HIS2 A1022** **Approximants de Padé**

**HIS1 N2107** **Fraction rationnelle**

$$\frac{1}{1 - 13 z}$$

**1, 13, 169, 2197, 28561, 371293, 4826809, 62748517, 815730721, 10604499373, 137858491849, 1792160394037, 23298085122481, 302875106592253**

## Powers of 14

**Réf.** BA9.

**HIS2 A1023** **Approximants de Padé**

**HIS1 N2120** **Fraction rationnelle**

$$\frac{1}{1 - 14 z}$$

**1, 14, 196, 2744, 38416, 537824, 7529536, 105413504, 1475789056, 20661046784, 289254654976, 4049565169664, 56693912375296, 793714773254144**

## Powers of 15

**Réf.** BA9.

**HIS2 A1024** **Approximants de Padé**

**HIS1 N2147** **Fraction rationnelle**

$$\frac{1}{1 - 15 z}$$

**1, 15, 225, 3375, 50625, 759375, 11390625, 170859375, 2562890625, 38443359375, 576650390625, 8649755859375, 129746337890625, 1946195068359375**

## Powers of 16

**Réf.** BA9.

**HIS2 A1025** **Approximants de Padé**

**HIS1 N2164** **Fraction rationnelle**

$$\frac{1}{1 - 16 z}$$

**1, 16, 256, 4096, 65536, 1048576, 16777216, 268435456, 4294967296, 68719476736, 1099511627776, 17592186044416, 281474976710656**

## Powers of 17

**Réf.** BA9.

**HIS2 A1026** **Approximants de Padé**

**HIS1 N2182** **Fraction rationnelle**

$$\frac{1}{1 - 17 z}$$

**1, 17, 289, 4913, 83521, 1419857, 24137569, 410338673, 6975757441, 118587876497, 2015993900449, 34271896307633, 582622237229761**

## Powers of 18

**Réf.** BA9.

**HIS2 A1027** **Approximants de Padé**

**HIS1 N2192** **Fraction rationnelle**

$$\frac{1}{1 - 18 z}$$

**1, 18, 324, 5832, 104976, 1889568, 34012224, 612220032, 11019960576, 198359290368, 3570467226624, 64268410079232, 1156831381426176**

## Powers of 19

**Réf.** BA9.

**HIS2 A1029** **Approximants de Padé**

**HIS1 N2198** **Fraction rationnelle**

$$\frac{1}{1 - 19 z}$$

**1, 19, 361, 6859, 130321, 2476099, 47045881, 893871739, 16983563041, 322687697779, 6131066257801, 116490258898219, 2213314919066161**

**Réf.** RCI 217.

**HIS2 A1044** **Hypergéométrique** **Suite P-récurrente**

**HIS1 N1492** **Fraction rationnelle** **double exponentielle**

**a(n) = Γ(n+1)^2**

$$\frac{z}{1 - z}$$

**1, 4, 36, 576, 14400, 518400, 25401600, 1625702400, 131681894400, 13168189440000, 1593350922240000, 229442532802560000, 38775788043632640000**

---

**Réf.** FQ 10 499 72. JCT A26 149 79.

**HIS2 A1045** **Approximants de Padé**

**HIS1 N0983** **Fraction rationnelle**

$$\frac{1}{(1 + z)(1 - 2z)}$$

**1, 1, 3, 5, 11, 21, 43, 85, 171, 341, 683, 1365, 2731, 5461, 10923, 21845, 43691, 87381, 174763, 349525, 699051, 1398101, 2796203, 5592405, 11184811, 22369621**

**Réf.** EUR 24 20 61. CR 268 579 69.
**HIS2 A1047** **Approximants de Padé**
**HIS1 N1596** **Fraction rationnelle**

$$\frac{1}{(1 - 3 z) (1 - 2 z)}$$

**1, 5, 19, 65, 211, 665, 2059, 6305, 19171, 58025, 175099, 527345, 1586131, 4766585, 14316139, 42981185, 129009091, 387158345, 1161737179, 3485735825**

---

**Réf.** CJM 22 26 70.
**HIS2 A1048** **Dérivée logarithmique** **Suite P-récurrente**
**HIS1 N0337** **Fraction rationnelle** **f.g. exponentielle**
**$3F2$([1, 1, 3], [2, 2], z)**

$$\frac{2 - z}{(1 - z)^2}$$

**2, 3, 8, 30, 144, 840, 5760, 45360, 403200, 3991680, 43545600, 518918400, 6706022400, 93405312000, 1394852659200, 22230464256000, 376610217984000**

**Réf.** FQ 3 129 65. BR72 52.

**HIS2 A1060 Approximants de Padé**

**HIS1 N0512 Fraction rationnelle**

$$\frac{2 + 3z}{1 - z - z^2}$$

**2, 5, 7, 12, 19, 31, 50, 81, 131, 212, 343, 555, 898, 1453, 2351, 3804, 6155, 9959, 16114, 26073, 42187, 68260, 110447, 178707, 289154, 467861, 757015, 1224876**

---

**Réf.** NCM 4 167 1878. MMAG 40 78 67. FQ 7 239 69.

**HIS2 A1075 Approximants de Padé**

**HIS1 N0700 Fraction rationnelle**

$$\frac{1 - 2z}{1 - 4z + z^2}$$

**1, 2, 7, 26, 97, 362, 1351, 5042, 18817, 70226, 262087, 978122, 3650401, 13623482, 50843527, 189750626, 708158977, 2642885282, 9863382151, 36810643322**

**Réf.** TH52 282.
**HIS2 A1076** **Approximants de Padé**
**HIS1 N1434** **Fraction rationnelle**

$$\frac{1}{1 - 4z - z^2}$$

**1, 4, 17, 72, 305, 1292, 5473, 23184, 98209, 416020, 1762289, 7465176, 31622993, 133957148, 567451585, 2403763488, 10182505537, 43133785636**

---

**Réf.** TH52 282.
**HIS2 A1077** **Approximants de Padé**
**HIS1 N0764** **Fraction rationnelle**

$$\frac{1 - 2z}{1 - 4z - z^2}$$

**1, 2, 9, 38, 161, 682, 2889, 12238, 51841, 219602, 930249, 3940598, 16692641, 70711162, 299537289, 1268860318, 5374978561, 22768774562, 96450076809**

**Réf.** TH52 281.

**HIS2 A1078** **Approximants de Padé**

**HIS1 N0839** **Fraction rationnelle**

$$\frac{2z}{1 - 10z + z^2}$$

**0, 2, 20, 198, 1960, 19402, 192060, 1901198, 18819920, 186298002, 1844160100, 18255302998, 180708869880, 1788833395802, 17707625088140**

**Réf.** EUL (1) 1 374 11. TH52 281.

**HIS2 A1079** **Approximants de Padé**

**HIS1 N1659** **Fraction rationnelle**

$$\frac{1 - 5z}{1 - 10z + z^2}$$

**1, 5, 49, 485, 4801, 47525, 470449, 4656965, 46099201, 456335045, 4517251249, 44716177445, 442644523201, 4381729054565, 43374646022449**

**Réf.** NCM 4 167 1878. TH52 281.
**HIS2 A1080** **Approximants de Padé**
**HIS1 N1278** **Fraction rationnelle**

$$\frac{3z}{1 - 16z + z^2}$$

**0, 3, 48, 765, 12192, 194307, 3096720, 49353213, 786554688, 12535521795, 199781794032, 3183973182717, 50743789129440, 808716652888323**

---

**Réf.** NCM 4 167 1878. TH52 281.
**HIS2 A1081** **Approximants de Padé**
**HIS1 N1949** **Fraction rationnelle**

$$\frac{1 - 8z}{1 - 16z + z^2}$$

**1, 8, 127, 2024, 32257, 514088, 8193151, 130576328, 2081028097, 33165873224, 528572943487, 8424001222568, 134255446617601, 2139663144659048**

**Réf.** NCM 4 167 1878. MTS 65(4, Supplement) 8 56.
**HIS2 A1084** **Approximants de Padé**
**HIS1 N1284** **Fraction rationnelle**

$$\frac{3z}{1 - 20z + z^2}$$

**0, 3, 60, 1197, 23880, 476403, 9504180, 189607197, 3782639760, 75463188003, 1505481120300, 30034159217997, 599177703239640, 11953519905574803**

---

**Réf.** NCM 4 167 1878. MTS 65(4, Supplement) 8 56.
**HIS2 A1085** **Approximants de Padé**
**HIS1 N2030** **Fraction rationnelle**

$$\frac{1 - 10z}{1 - 20z + z^2}$$

**1, 10, 199, 3970, 79201, 1580050, 31521799, 628855930, 12545596801, 250283080090, 4993116004999, 99612037019890, 1987247624392801**

**Réf.** NCM 4 167 1878.
**HIS2 A1090** **Approximants de Padé**
**HIS1 N1936** **Fraction rationnelle**

$$\frac{1}{1 - 8z + z^2}$$

**1, 8, 63, 496, 3905, 30744, 242047, 1905632, 15003009, 118118440, 929944511, 7321437648, 57641556673, 453811015736, 3572846569215, 28128961537984**

---

**Réf.** NCM 4 167 1878.
**HIS2 A1091** **Approximants de Padé**
**HIS1 N1479** **Fraction rationnelle**

$$\frac{1 - 4z}{1 - 8z + z^2}$$

**1, 4, 31, 244, 1921, 15124, 119071, 937444, 7380481, 58106404, 457470751, 3601659604, 28355806081, 223244789044, 1757602506271, 13837575261124**

## Enneagonal numbers

**Réf.** B1 189.

**HIS2 A1106** **Approximants de Padé**

**HIS1** **Fraction rationnelle**

$$\frac{1 + 6z}{(1 - z)^3}$$

**1, 9, 24, 46, 75, 111, 154, 204, 261, 325, 396, 474, 559, 651, 750, 856, 969, 1089, 1216, 1350, 1491, 1639, 1794, 1956, 2125, 2301, 2484, 2674, 2871, 3075, 3286, 3504, 3729, 3961, 4200**

## Decagonal numbers

**Réf.** B1 189.

**HIS2 A1107** **Approximants de Padé**

**HIS1** **Fraction rationnelle**

$$\frac{1 + 7z}{(1 - z)^3}$$

**1, 10, 27, 52, 85, 126, 175, 232, 297, 370, 451, 540, 637, 742, 855, 976, 1105, 1242, 1387, 1540, 1701, 1870, 2047, 2232, 2425, 2626, 2835, 3052, 3277, 3510, 3751, 4000, 4257, 4522**

**n(n+1)/2 is square**

**Réf.** D1 2 10. MAG 47 237 63. B1 193. FQ 9 95 71.

**HIS2 A1108** **Approximants de Padé**

**HIS1 N1924** **Fraction rationnelle**

$$\frac{1 + z}{(z - 1)(z^2 - 6z + 1)}$$

**1, 8, 49, 288, 1681, 9800, 57121, 332928, 1940449, 11309768, 65918161, 384199200, 2239277041, 13051463048, 76069501249, 443365544448, 2584123765441**

---

**Réf.** D1 2 10. MAG 47 237 63. B1 193. FQ 9 95 71.

**HIS2 A1109** **Approximants de Padé**

**HIS1 N1760** **Fraction rationnelle**

$$\frac{1}{1 - 6z + z^2}$$

**1, 6, 35, 204, 1189, 6930, 40391, 235416, 1372105, 7997214, 46611179, 271669860, 1583407981, 9228778026, 53789260175, 313506783024, 1827251437969**

## Both triangular and square

**Réf.** D1 2 10. MAG 47 237 63. B1 193. FQ 9 95 71.

**HIS2 A1110** **Approximants de Padé**

**HIS1 N2291** **Fraction rationnelle**

$$\frac{1 + z}{(1 - z)\,(z^2 - 34\,z + 1)}$$

**1, 36, 1225, 41616, 1413721, 48024900, 1631432881, 55420693056, 1882672131025, 63955431761796, 2172602007770041, 73804512832419600**

## Differences of 0

**Réf.** VO11 31. DA63 2 212. R1 33.

**HIS2 A1117** **Approximants de Padé**

**HIS1 N1763** **Fraction rationnelle**

$$\frac{6}{(1 - z)\,(1 - 2\,z)\,(1 - 3\,z)}$$

**6, 36, 150, 540, 1806, 5796, 18150, 55980, 171006, 519156, 1569750, 4733820, 14250606, 42850116, 128746950, 386634060, 1160688606, 3483638676**

## Differences of 0

**Réf.** VO11 31. DA63 2 212. R1 33.

**HIS2 A1118** **Approximants de Padé**

**HIS1 N2334** **Fraction rationnelle**

$$\frac{120}{(1 - z)\,(1 - 2z)\,(1 - 3z)\,(1 - 4z)\,(1 - 5z)}$$

**120, 1800, 16800, 126000, 834120, 5103000, 29607600, 165528000, 901020120, 4809004200, 25292030400, 131542866000, 678330198120, 3474971465400**

## Double factorials

**Réf.** AMM 55 425 48. MOC 24 231 70.

**HIS2 A1147** **Hypergéométrique** **Suite P-récurrente**

**HIS1 N1217** **exponentielle (algébrique)**

**Inverse fonctionnel de A1710**

**Inverse de A0698**

$$\frac{2z}{1 + (1 - 2z)^{1/2}}$$

**1, 1, 3, 15, 105, 945, 10395, 135135, 2027025, 34459425, 654729075, 13749310575, 316234143225, 7905853580625, 213458046676875, 6190283353629375**

## Partitions of n into squares

**Réf.** BIT 19 298 79.

**HIS2 A1156** **Euler**

**HIS1 N0079** **Produit infini**

$$\prod_{n\geq 1} \frac{1}{(1 - z^{c(n)})}$$

**c(n) = 1,4,9,16,..., les carrés parfaits.**

**1, 1, 1, 1, 2, 2, 2, 2, 3, 4, 4, 4, 5, 6, 6, 6, 8, 9, 10, 10, 12, 13, 14, 14, 16, 19, 20, 21, 23, 26, 27, 28, 31, 34, 37, 38, 43, 46, 49, 50, 55, 60, 63, 66, 71, 78, 81, 84, 90, 98, 104, 107, 116**

## Board-pile polyominoes with n cells

**Réf.** JCT 6 103 69. AB71 363. JSP 58 477 90.

**HIS2 A1169** **Approximants de Padé**

**HIS1 N0639** **Fraction rationnelle**

$$\frac{(1 - z)^3}{1 - 5z + 7z^2 - 4z^3}$$

**1, 2, 6, 19, 61, 196, 629, 2017, 6466, 20727, 66441, 212980, 682721, 2188509, 7015418, 22488411, 72088165, 231083620, 740754589, 2374540265, 7611753682**

## Baxter permutations of length 2n-1

**Réf.** MAL 2 25 67. JCT A24 393 78. FQ 27 166 89.

**HIS2 A1181** **P-récurrences** **Suite P-récurrente**

**HIS1 N0652**

$$(n + 3) (n + 2) a(n) = (7 n^2 + 7 n - 2) a(n - 1) + (8 n^2 - 24 n + 16) a(n - 2)$$

**1, 2, 6, 22, 92, 422, 2074, 10754, 58202, 326240, 1882960, 11140560, 67329992, 414499438, 2593341586, 16458756586, 105791986682, 687782586844, 4517543071924**

## Degree n permutations of order exactly 2

**Réf.** CJM 7 159 55.

**HIS2 A1189** **P-récurrences** **Suite P-récurrente**

**HIS1 N1127** **exponentielle**

$$a(n) = 3 a(n - 1) + (n - 3) a(n - 2) + (- 2 n + 3) a(n - 3) + (n - 2) a(n - 4)$$

$$\exp(1/2 z (2 + z)) - \exp(z)$$

**0, 1, 3, 9, 25, 75, 231, 763, 2619, 9495, 35695, 140151, 568503, 2390479, 10349535, 46206735, 211799311, 997313823, 4809701439, 23758664095, 119952692895**

## Expansion of an integral

**Réf.** C1 167.

**HIS2 A1193** **Hypergéométrique** **Suite P-récurrente**

**HIS1 N0770** **exponentielle (algébrique)**

**(n - 1) a(n) = (2 n - 3) n a(n - 1)**

$$\frac{z}{(1 - 2z)^{1/2}}$$

**1, 2, 9, 60, 525, 5670, 72765, 1081080, 18243225**

## Expansion of an integral

**Réf.** C1 167.

**HIS2 A1194** **Hypergéométrique** **Suite P-récurrente.**

**HIS1 N1139** **exponentielle (algébrique)** **double exponentielle**

$$\frac{z (2 - 3z)}{(1 - 2z)^{3/2}}$$

**3, 9, 54, 450, 4725, 59535, 873180, 14594580**

## Clouds with n points

**Réf.** C1 276.

**HIS2 A1205** **Dérivée logarithmique** **Suite P-récurrente.**

**HIS1 N1181** **exponentielle (algébrique)**

**2 a(n) = (n - 2) (n - 3) a(n - 3) + (2 n - 4) a(n - 1)**

$$\frac{\exp(-1/4\, z\, (z + 2))}{(1 - z)^{1/2}}$$

**1, 0, 0, 1, 3, 12, 70, 465, 3507, 30016, 286884, 3026655, 34944085, 438263364, 5933502822, 86248951243, 1339751921865, 22148051088480, 388246725873208**

## Packing a box with n dominoes

**Réf.** AMM 69 61 62.

**HIS2 A1224** **Approximants de Padé**

**HIS1 N0117** **Fraction rationnelle**

$$\frac{1 + z - 2 z^2 - z^3 - z^4 - z^5}{(z^4 + z^2 - 1)(z^2 + z - 1)}$$

**1, 2, 2, 4, 5, 9, 12, 21, 30, 51, 76, 127, 195, 322, 504, 826, 1309, 2135, 3410, 5545, 8900, 14445, 23256, 37701, 60813, 98514, 159094, 257608, 416325, 673933, 1089648**

## Stirling numbers of first kind

**Réf.** AS1 833. DKB 226.

**HIS2 A1233** **Tableaux généralisés** **Suite P-récurrente**

**HIS1 N2216** **exponentielle (log)**

$$\frac{- \ln(1 - z)^{5}}{120\ (1 - z)}$$

**1, 21, 322, 4536, 63273, 902055, 13339535, 206070150, 3336118786, 56663366760, 1009672107080, 18861567058880, 369012649234384**

## Stirling numbers of first kind

**Réf.** AS1 834. DKB 226.

**HIS2 A1234** **Tableaux généralisés** **Suite P-récurrente**

**HIS1 N2264** **exponentielle (log)**

$$\frac{\ln(1 - z)^{6}}{720\ (1 - z)}$$

**1, 28, 546, 9450, 157773, 2637558, 44990231, 790943153, 14409322928, 272803210680, 5374523477960, 110228466184200, 2353125040549984**

## Differences of reciprocals of unity

**Réf.** DKB 228.

**HIS2 A1240** **Approximants de Padé**

**HIS1 N2049** **Fraction rationnelle**

$$\frac{1}{(1 - 2z)(1 - 3z)(1 - 6z)}$$

**1, 11, 85, 575, 3661, 22631, 137845, 833375, 5019421, 30174551**

## Differences of reciprocals of unity

**Réf.** DKB 228.

**HIS2 A1241** **Approximants de Padé**

**HIS1 N2305** **Fraction rationnelle**

$$\frac{1}{(1 - 6z)(1 - 8z)(1 - 12z)(1 - 24z)}$$

**1, 50, 1660, 46760, 1217776, 30480800, 747497920, 18139003520, 437786795776**

## Permutations of length n by length of runs

**Réf.** AMM 65 534 58. DKB 262. C1 261.

**HIS2 A1250** **Inverse fonctionnel** **Relié aux nombres tangents**
**HIS1 N0472** **exponentielle (complexe)**

$$2 \tan(1/4 \pi + 1/2 z)$$

**2, 4, 10, 32, 122, 544, 2770, 15872, 101042, 707584, 5405530, 44736512, 398721962, 3807514624, 38783024290, 419730685952, 4809759350882**

## Permutations of length n by rises

**Réf.** DKB 263.

**HIS2 A1260** **P-récurrences** **Suite P-récurrente**
**HIS1 N1657**

$$a(n)(1 - n) = -(n + 3)(n + 2)a(n - 2) - (n + 3)(n - 1)a(n - 1)$$

**1, 5, 45, 385, 3710, 38934, 444990, 5506710, 73422855, 1049946755, 16035550531, 260577696015**

## Lah numbers

**Réf.** R1 44. C1 156.

**HIS2 A1286** **Dérivée logarithmique** **f.g. exponentielle**

**HIS1 N1766** **Fraction rationnelle**

$$\frac{2z+1}{(1-z)^{4}}$$

**1, 6, 36, 240, 1800, 15120, 141120, 1451520, 16329600, 199584000, 2634508800, 37362124800, 566658892800, 9153720576000, 156920924160000**

## Binomial coefficients C(n,10)

**Réf.** D1 2 7. RS3. B1 196. AS1 828.

**HIS2 A1287** **Approximants de Padé**

**HIS1 N2046** **Fraction rationnelle**

$$\frac{1}{(1-z)^{11}}$$

**1, 11, 66, 286, 1001, 3003, 8008, 19448, 43758, 92378, 184756, 352716, 646646, 1144066, 1961256, 3268760, 5311735, 8436285, 13123110, 20030010, 30045015**

## Binomial coefficients C(n,11)

**Réf.** D1 2 7. RS3. B1 196. AS1 828.

**HIS2 A1288 Approximants de Padé**

**HIS1 N2073 Fraction rationnelle**

$$\frac{1}{(1 - z)^{12}}$$

**1, 12, 78, 364, 1365, 4368, 12376, 31824, 75582, 167960, 352716, 705432, 1352078, 2496144, 4457400, 7726160, 13037895, 21474180, 34597290, 54627300**

## Stirling numbers of second kind

**Réf.** AS1 835. DKB 223.

**HIS2 A1296 Approximants de Padé**

**HIS1 N1845 Fraction rationnelle**

$$\frac{1 + 2 z}{(1 - z)^{5}}$$

**1, 7, 25, 65, 140, 266, 462, 750, 1155, 1705, 2431, 3367, 4550, 6020, 7820, 9996, 12597, 15675, 19285, 23485, 28336, 33902, 40250, 47450, 55575, 64701, 74907, 86275**

## Stirling numbers of second kind

**Réf.** AS1 835. DKB 223.

**HIS2 A1297** **Approximants de Padé**

**HIS1 N2136** **Fraction rationnelle**

$$\frac{1 + 8z + 6z^2}{(1 - z)^7}$$

**1, 15, 90, 350, 1050, 2646, 5880, 11880, 22275, 39325, 66066, 106470, 165620, 249900, 367200, 527136, 741285, 1023435, 1389850, 1859550, 2454606, 3200450**

---

## Stirling numbers of second kind

**Réf.** AS1 835. DKB 223.

**HIS2 A1298** **Approximants de Padé**

**HIS1 N2272** **Fraction rationnelle**

$$\frac{1 + 22z + 58z^2 + 24z^3}{(1 - z)^9}$$

**1, 31, 301, 1701, 6951, 22827, 63987, 159027, 359502, 752752, 1479478, 2757118, 4910178, 8408778, 13916778, 22350954, 34952799, 53374629, 79781779**

## Stirling numbers of first kind

**Réf.** AS1 833. DKB 226.

**HIS2 A1303** **Approximants de Padé**

**HIS1 N1779** **Fraction rationnelle**

$$\frac{6 + 8z + z^2}{(1 - z)^7}$$

**6, 50, 225, 735, 1960, 4536, 9450, 18150, 32670, 55770, 91091, 143325, 218400, 323680, 468180, 662796, 920550, 1256850, 1689765, 2240315, 2932776**

## Generalized pentagonal numbers

**Réf.** NZ66 231. AMM 76 884 69. HO70 119.

**HIS2 A1318** **Approximants de Padé**

**HIS1 N0511** **Fraction rationnelle**

$$\frac{z^2 + z + 1}{(1 + z)^2 (1 - z)^3}$$

**1, 2, 5, 7, 12, 15, 22, 26, 35, 40, 51, 57, 70, 77, 92, 100, 117, 126, 145, 155, 176, 187, 210, 222, 247, 260, 287, 301, 330, 345, 376, 392, 425, 442, 477, 495, 532, 551, 590**

**Réf.** MQET 1 9 16. AMM 56 445 49.
**HIS2 A1333** **Approximants de Padé**
**HIS1 N1064** **Fraction rationnelle**

$$\frac{1 + z}{1 - 2z - z^2}$$

**1, 3, 7, 17, 41, 99, 239, 577, 1393, 3363, 8119, 19601, 47321, 114243, 275807, 665857, 1607521, 3880899, 9369319, 22619537, 54608393, 131836323, 318281039**

---

## Binomial coefficient sums

**Réf.** CJM 22 26 70.
**HIS2 A1338** **Recoupements**
**HIS1 N0697** **exponentielle**

$$- \exp(z) (\ln(1 - z) + 1) + 2$$

**1, 0, 2, 7, 23, 88, 414, 2371, 16071, 125672, 1112082**

**Réf.** CJM 22 26 70. AD74 70.
**HIS2 A1339** **Dérivée logarithmique exponentielle** **Suite P-récurrente**
**HIS1 N1164**

**a(n) = (n + 1) a(n - 1) + (- n + 2) a(n - 2)**

**∑(n+1)! C(n,k), k=0...n**

$$\frac{\exp(z)}{(1 - z)^2}$$

**1, 3, 11, 49, 261, 1631, 11743, 95901, 876809, 8877691, 98641011, 1193556233, 15624736141, 220048367319, 3317652307271, 53319412081141, 909984632851473**

---

**Réf.** CJM 22 26 70.
**HIS2 A1340** **Dérivée logarithmique exponentielle** **Suite P-récurrente**
**HIS1 N0736**

$$\frac{2 \exp(z)}{(1 - z)^3}$$

**2, 8, 38, 212, 1370, 10112, 84158, 780908, 8000882**

**Réf.** CJM 22 26 70.

**HIS2 A1341** **Dérivée logarithmique exponentielle** **Suite P-récurrente**

**HIS1 N1755**

$$\frac{6\ \exp(z)}{(1 - z)^4}$$

**6, 30, 174, 1158, 8742, 74046, 696750, 7219974**

---

**Réf.** CJM 22 26 70.

**HIS2 A1342** **Dérivée logarithmique exponentielle** **Suite P-récurrente**

**HIS1 N2233**

$$\frac{24\ \exp(z)}{(1 - z)^5}$$

**24, 144, 984, 7584, 65304, 622704, 6523224**

**Réf.** CJM 22 26 70.

**HIS2 A1344** **Dérivée exponentielle** **Suite P-récurrente**

**HIS1 N0548**

**a(n) = (n - 3) a(n - 2) + (n - 1) a(n - 1)**

$$\frac{1}{(z-1)^2} - \frac{2}{z-1} - \ln(z-1)$$

**2, 5, 11, 38, 174, 984, 6600, 51120, 448560, 4394880, 47537280, 562464000, 7224940800, 100111334400, 1488257971200, 23625316915200, 398840682240000, 7134671351808000**

---

**Réf.** EUR 11 22 49.

**HIS2 A1350** **Approximants de Padé**

**HIS1 N1311** **Fraction rationnelle**

$$\frac{1+z^2}{(1-z)(1+z)(1-z-z^2)}$$

**1, 1, 4, 5, 11, 16, 29, 45, 76, 121, 199, 320, 521, 841, 1364, 2205, 3571, 5776, 9349, 15125, 24476, 39601, 64079, 103680, 167761, 271441, 439204, 710645, 1149851**

## Associated Mersenne numbers

**Réf.** EUR 11 22 49.

**HIS2 A1351** **Approximants de Padé** **expression factorisée**

**HIS1 N0879** **Fraction rationnelle**

$$\frac{z\,(1 - z + z^2)\,(z^2 + 3\,z + 1)}{(1 - z - z^3)\,(1 - z^2 - z^3)}$$

**0, 1, 3, 1, 3, 11, 9, 8, 27, 37, 33, 67, 117, 131, 192, 341, 459, 613, 999, 1483, 2013, 3032, 4623, 6533, 9477, 14311, 20829, 30007, 44544, 65657, 95139, 139625, 206091**

---

**Réf.** MOC 24 180 70.

**HIS2 A1352** **Approximants de Padé**

**HIS1 N1731** **Fraction rationnelle**

$$\frac{(1 + z)^2}{1 - 4\,z + z^2}$$

**1, 6, 24, 90, 336, 1254, 4680, 17466, 65184, 243270, 907896, 3388314, 12645360, 47193126, 176127144, 657315450, 2453134656, 9155223174, 34167758040**

**Réf.** MMAG 40 78 67. MOC 24 180 70; 25 799 71.

**HIS2 A1353** **Approximants de Padé**

**HIS1 N1420** **Fraction rationnelle**

$$\frac{1}{1 - 4z + z^2}$$

**1, 4, 15, 56, 209, 780, 2911, 10864, 40545, 151316, 564719, 2107560, 7865521, 29354524, 109552575, 408855776, 1525870529, 5694626340, 21252634831**

---

**n-node trees of height at most 3**

**Réf.** IBMJ 4 475 60. KU64.

**HIS2 A1383** **Euler**

**HIS1 N0422** **Produit infini**

$$\prod_{n\geq 1} \frac{1}{(1 - z^n)^{c(n)}}$$

**c(n) = partages de n**

**1, 1, 2, 4, 8, 15, 29, 53, 98, 177, 319, 565, 1001, 1749, 3047, 5264, 9054, 15467, 26320, 44532, 75054, 125904, 210413, 350215, 580901, 960035, 158153**

## n-node trees of height at most 4

**Réf.** IBMJ 4 475 60. KU64.

**HIS2 A1384** **Euler** **a(n) = suite précédente**

**HIS1 N0449** **Produit infini**

$$\prod_{n\geq 1} \frac{1}{(1-z^n)^{c(n)}}$$

**c(n) = arbres de hauteur au plus 3**

**1, 1, 2, 4, 9, 19, 42, 89, 191, 402, 847, 1763, 3667, 7564, 15564, 31851, 64987, 132031, 267471, 539949, 1087004, 2181796, 4367927, 8721533, 17372967, 34524291**

---

## n-node trees of height at most 5

**Réf.** IBMJ 4 475 60. KU64.

**HIS2 A1385** **Euler** **a(n) = suite précédente**

**HIS1 N0453** **Produit infini**

$$\prod_{n\geq 1} \frac{1}{(1-z^n)^{c(n)}}$$

**c(n) = arbres de hauteur au plus 4**

**1, 1, 2, 4, 9, 20, 47, 108, 252, 582, 1345, 3086, 7072, 16121, 36667, 83099, 187885, 423610, 953033, 2139158, 4792126, 10714105, 23911794, 53273599, 118497834**

**Réf.** QAM 14 407 56. MOC 29 216 75.
**HIS2 A1392** **Hypergéométrique** **Suite P-récurrente**
**HIS1 N1981** **algébrique**
**$_2F_1([5, 9/2], [10], 4z)$**

$$\frac{512 z^4}{\left(1 + (1 - 4z)^{1/2}\right)^9}$$

**1, 9, 54, 273, 1260, 5508, 23256, 95931, 389367, 1562275, 6216210, 24582285, 96768360, 379629720, 1485507600, 5801732460, 22626756594, 88152205554**

## Partitions into at most 3 parts

**Réf.** RS4 2. AMM 86 687 79.
**HIS2 A1399** **Approximants de Padé**
**HIS1 N0186** **Fraction rationnelle**

$$\frac{1}{(1 - z)(1 - z^2)(1 - z^3)}$$

**1, 1, 2, 3, 4, 5, 7, 8, 10, 12, 14, 16, 19, 21, 24, 27, 30, 33, 37, 40, 44, 48, 52, 56, 61, 65, 70, 75, 80, 85, 91, 96, 102, 108, 114, 120, 127, 133, 140, 147, 154, 161, 169, 176, 184**

## Partitions into at most 4 parts

**Réf.** RS4 2.

**HIS2 A1400** **Approximants de Padé**

**HIS1 N0229** **Fraction rationnelle**

$$\frac{1}{(1-z)(1-z^2)(1-z^3)(1-z^4)}$$

**1, 2, 3, 5, 6, 9, 11, 15, 18, 23, 27, 34, 39, 47, 54, 64, 72, 84, 94, 108, 120, 136, 150, 169, 185, 206, 225, 249, 270, 297, 321, 351, 378, 411, 441, 478, 511, 551, 588, 632, 672**

## Partitions of n into at most 5 parts

**Réf.** RS4 2.

**HIS2 A1401** **Recoupement**

**HIS1 N0237** **Fraction rationnelle**

$$\frac{1}{(1-z)(1-z^2)(1-z^3)(1-z^4)(1-z^5)}$$

**1, 2, 3, 5, 7, 10, 13, 18, 23, 30, 37, 47, 57, 70, 84, 101, 119, 141, 164, 192, 221, 255, 291, 333, 377, 427, 480, 540, 603, 674, 748, 831, 918, 1014, 1115, 1226, 1342, 1469**

## Partitions of n into at most 6 parts

**Réf.** CAY 10 415. RS4 2.

**HIS2 A1402** **Euler**

**HIS1 N0243** **Fraction rationnelle**

$$\frac{1}{(1 - z)(1 - z^2)(1 - z^3)(1 - z^4)(1 - z^5)(1 - z^6)}$$

**1, 1, 2, 3, 5, 7, 11, 14, 20, 26, 35, 44, 58, 71, 90, 110, 136, 163, 199, 235, 282, 331, 391, 454, 532, 612, 709, 811, 931, 1057, 1206, 1360, 1540, 1729, 1945, 2172, 2432**

## Central binomial coefficients

**Réf.** RS3. AS1 828. JCT 1 299 66.

**HIS2 A1405** **LLL** **Suite P-récurrente**

**HIS1 N0294** **algébrique**

**C(n,[n/2])**

$$\frac{1 - 4 z^2 - (1 - 4 z^2)^{1/2}}{2 (2 z^3 - z^2)}$$

**1, 2, 3, 6, 10, 20, 35, 70, 126, 252, 462, 924, 1716, 3432, 6435, 12870, 24310, 48620, 92378, 184756, 352716, 705432, 1352078, 2704156, 5200300, 10400600**

## Catalan numbers -1

**Réf.** MOC 22 390 68.

**HIS2 A1453** **LLL** **Suite P-récurrente**

**HIS1 N1409** **algébrique**

**(n + 2) a(n) = (6 n + 4) a(n - 1) + (- 9 n + 4) a(n - 2) + (4 n - 6) a(n - 3)**

$$\frac{1 - 4z + 3z^2 - (-(4z - 1)(z - 1)^4)^{1/2}}{2(z^3 - 2z^4 + z^5)}$$

**1, 4, 13, 41, 131, 428, 1429, 4861, 16795, 58785, 208011, 742899, 2674439, 9694844, 35357669, 129644789, 477638699, 1767263189, 6564120419, 24466267019**

## Degree n permutations of order dividing 3

**Réf.** CJM 7 159 55.

**HIS2 A1470** **Dérivée logarithmique** **Suite P-récurrente**

**HIS1 N1118** **exponentielle**

**a(n) = a(n - 1) + (n^2 - 3 n + 2) a(n - 3)**

$$(1 + z^2) \exp(1/3\, z\, (3 + z^2))$$

**1, 1, 3, 9, 21, 81, 351, 1233, 5769, 31041, 142011, 776601, 4874013, 27027729, 168369111, 1191911841, 7678566801, 53474964993, 418199988339**

## Degree n permutations of order dividing 4

**Réf.** CJM 7 159 55.

**HIS2 A1472** **Dérivée logarithmique** **Suite P-récurrente**

**HIS1 N0495** **exponentielle**

**a(n) = a(n - 1) + (n^3 - 6 n^2 + 11 n - 6) a(n - 4) + (n - 1) a(n - 2)**

$$(1 + z + z^3) \exp(1/4\, z\, (4 + z^3 + 2\, z))$$

**1, 2, 4, 16, 56, 256, 1072, 6224, 33616, 218656, 1326656, 9893632, 70186624, 574017536, 4454046976, 40073925376, 347165733632, 3370414011904**

---

**Réf.** R1 86 (divided by 2).

**HIS2 A1475** **Dérivée logarithmique** **Suite P-récurrente**

**HIS1 N0573** **exponentielle**

**a(n) = a(n - 1) + n a(n - 2)**

$$\exp(1/2\, z^2 + z + \ln(2 + 2\, z + z^2))$$

**1, 2, 5, 13, 38, 116, 382, 1310, 4748, 17848, 70076, 284252, 1195240, 5174768, 23103368, 105899656, 498656912, 2404850720, 11879332048, 59976346448**

## Stochastic matrices of integers

**Réf.** DMJ 35 659 68.

**HIS2 A1495** **Recoupements** **Suite P-récurrente**

**HIS1 N1188** **exponentielle:algébrique**

$$\frac{\exp(z\,(z^2 + 3\,z - 2)/(1-z))}{(1 - z)^{3/2}}$$

**0, 1, 1, 1, 3, 13, 70, 462, 3592, 32056, 322626, 3611890, 44491654, 597714474, 8693651092, 136059119332, 2279212812480, 40681707637888, 770631412413148**

## 4 x 4 stochastic matrices of integers

**Réf.** SS70. CJN 13 283 70. SIAC 4 477 75. ANS 4 1179 76.

**HIS2 A1496** **Dérivée logarithmique**

**HIS1 N2240** **Fraction rationnelle**

$$\frac{(z^4 + 12\,z^3 + 62\,z^2 + 12\,z + 1)\,(z + 1)^2}{(z - 1)^{10}}$$

**1, 24, 282, 2008, 10147, 40176, 132724, 381424, 981541, 2309384, 5045326, 10356424, 20158151, 37478624, 66952936, 115479776, 193077449, 313981688, 498033282, 772409528**

## Stochastic matrices of integers

**Réf.** SS70. DMJ 33 763 66.

**HIS2 A1499** **équations différentielles** **Formule de B. Salvy**

**HIS1 N1792** **exponentielle**

$$\frac{(z^2 - 2z + 4)\exp(-1/2\, z)}{(1 - z)^{5/2}}$$

**0, 1, 6, 90, 2040, 67950, 3110940, 187530840, 14398171200, 1371785398200, 158815387962000, 21959547410077200, 3574340599104475200**

## Bessel polynomial $y_n(1)$

**Réf.** RCI 77.

**HIS2 A1514** **P-récurrences** **Suite P-récurrente**

**HIS1 N1993**

$$a(n) = (2n + 4)\, a(n-1) + a(n-4) + (-6n + 9)\, a(n-2) + (2n - 10)\, a(n-3)$$

**0, 1, 9, 81, 835, 9990, 137466, 2148139, 37662381, 733015845, 15693217705, 366695853876, 9289111077324, 253623142901401, 7425873460633005**

**Réf.** RCI 77.

**HIS2 A1515** **équations différentielles** **Suite P-récurrente**

**HIS1 N0713** **exponentielle:algébrique** **Formule de B. Salvy**

**a(n) = (2n-1) a(n-1) + a(n-2)**

$$\frac{\exp(1 - (1 - 2 z)^{1/2})}{(1 - 2 z)^{1/2}}$$

**1, 2, 7, 37, 266, 2431, 27007, 353522, 5329837, 90960751, 1733584106, 36496226977, 841146804577, 21065166341402, 569600638022431**

---

**Denominators of convergents to e = exp(1)**

**Réf.** BAT 17 1871. MOC 2 69 46.

**HIS2 A1517** **équations différentielles** **Suite P-récurrente**

**HIS1 N1240** **exponentielle** **Voir A2119**

**a(n) = (4 n - 6) a(n - 1) + a(n - 2)**

$$\frac{\exp(1/2 - 1/2 (1 - 4 z)^{1/2})}{(1 - 4 z)^{1/2}}$$

**1, 3, 19, 193, 2721, 49171, 1084483, 28245729, 848456353, 28875761731, 1098127402131, 46150226651233, 2124008553358849, 106246577894593683**

## Bessel polynomial $y_n$ (3)

**Réf.** RCI 77.

**HIS2 A1518** **équations différentielles** **Suite P-récurrente**

**HIS1 N1495** **exponentielle** **Formule de B. Salvy**

**a(n) = (6 n - 9) a(n - 1) + a(n - 2)**

$$\frac{\exp\left(1/3 - 1/3\,(1 - 6z)^{1/2}\right)}{(1 - 6z)^{1/2}}$$

**1, 4, 37, 559, 11776, 318511, 10522639, 410701432, 18492087079, 943507142461, 53798399207356, 3390242657205889, 233980541746413697**

## Bisection of Fibonacci sequence

**Réf.** R1 39. FQ 9 283 71.

**HIS2 A1519** **Approximants de Padé**

**HIS1 N0569** **Fraction rationnelle**

$$\frac{1 - z}{1 - 3z + z^2}$$

**1, 2, 5, 13, 34, 89, 233, 610, 1597, 4181, 10946, 28657, 75025, 196418, 514229, 1346269, 3524578, 9227465, 24157817, 63245986, 165580141, 433494437**

## Stacks, or planar partitions of n

**Réf.** PCPS 47 686 51. QJMO 23 153 72.

**HIS2 A1522** **Approximants de Padé** **Conjecture**

**HIS1 N0238** **Fraction rationnelle**

$$\frac{z^{10} + z^{8} - 2z^{7} - z^{6} + 2z^{5} + z^{3} - z^{2} - z + 1}{(z + 1)(z^{4} + z^{3} - 1)(z - 1)^{3}}$$

**1, 1, 1, 2, 3, 5, 7, 10, 14, 19, 26, 35, 47, 62, 82, 107, 139, 179, 230, 293**

---

## Transpositions needed to generate permutations of length n

**Réf.** CJN 13 155 70.

**HIS2 A1540** **Inverse fonctionnel** **Suite P-récurrente**

**HIS1 N0734** **exponentielle**

**a(n) = - n a(n - 3) + (n + 2) a(n - 1) + (- n + 1) a(n - 2) + (n - 2) a(n - 4)**

**[cosh(1)*n!] -1**

$$\frac{(2z + 3 + 2z^{3} - 5z^{2})\exp(z)}{2(z - 1)^{3}} + \frac{1 - z^{2}}{(z - 1)^{3}\, 2\exp(z)}$$

**0, 2, 8, 36, 184, 1110, 7776, 62216, 559952, 5599530, 61594840, 739138092, 9608795208, 134523132926, 2017846993904, 32285551902480**

**Réf.** NCM 4 166 1878. QJM 45 14 14. ANN 36 644 35. AMM 75 683 68.

**HIS2 A1541 Approximants de Padé**

**HIS1 N1231 Fraction rationnelle**

$$\frac{1 - 3z}{1 - 6z + z^2}$$

**1, 3, 17, 99, 577, 3363, 19601, 114243, 665857, 3880899, 22619537, 131836323, 768398401, 4478554083, 26102926097, 152139002499, 886731088897**

---

**Réf.** NCM 4 166 1878. ANN 30 72 28. AMM 75 683 68.

**HIS2 A1542 Approximants de Padé**

**HIS1 N0802 Fraction rationnelle**

$$\frac{2z}{z^2 - 6z + 1}$$

**0, 2, 12, 70, 408, 2378, 13860, 80782, 470832, 2744210, 15994428, 93222358, 543339720, 3166815962, 18457556052, 107578520350, 627013566048**

**1^n + 2^n + 3^n**

**Réf.** AS1 813.

**HIS2 A1550** **Approximants de Padé**

**HIS1 N1020** **Fraction rationnelle**

$$\frac{3 - 12 z + 11 z^2}{(1 - z) (1 - 2 z) (1 - 3 z)}$$

**3, 6, 14, 36, 98, 276, 794, 2316, 6818, 20196, 60074, 179196, 535538, 1602516, 4799354, 14381676, 43112258, 129271236, 387682634, 1162785756, 3487832978**

---

**1^n + 2^n + 3^n + 4^n**

**Réf.** AS1 813.

**HIS2 A1551** **Approximants de Padé**

**HIS1 N1375** **Fraction rationnelle**

$$\frac{2 (5 z - 2) (5 z^2 - 5 z + 1)}{(1 - z) (1 - 2 z) (1 - 3 z) (1 - 4 z)}$$

**4, 10, 30, 100, 354, 1300, 4890, 18700, 72354, 282340, 1108650, 4373500, 17312754, 68711380, 273234810, 1088123500, 4338079554, 17309140420**

## 1^n + 2^n + 3^n + 4^n + 5^n

**Réf.** AS1 813.

**HIS2 A1552** **Approximants de Padé**

**HIS1 N1584** **Fraction rationnelle**

$$\frac{5 - 60 z + 255 z^2 - 450 z^3 + 274 z^4}{(1 - z)(1 - 2 z)(1 - 3 z)(1 - 4 z)(1 - 5 z)}$$

**5, 15, 55, 225, 979, 4425, 20515, 96825, 462979, 2235465, 10874275, 53201625, 261453379, 1289414505, 6376750435, 31605701625, 156925970179**

## 1^n + 2^n + 3^n + 4^n + 5^n + 6^n

**Réf.** AS1 813.

**HIS2 A1553** **Approximants de Padé**

**HIS1 N1723** **Fraction rationnelle**

$$\frac{(2 - 7 z)(252 z^4 - 392 z^3 + 203 z^2 - 42 z + 3)}{(1 - z)(1 - 2 z)(1 - 3 z)(1 - 4 z)(1 - 5 z)(1 - 6 z)}$$

**6, 21, 91, 441, 2275, 12201, 67171, 376761, 2142595, 12313161, 71340451, 415998681, 2438235715, 14350108521, 84740914531, 501790686201**

## 1^n + 2^n + 3^n + 4^n + 5^n + 6^n + 7^n

**Réf.** AS1 813.

**HIS2 A1554** **Approximants de Padé**

**HIS1 N1850** **Fraction rationnelle**

$$\frac{8028 z^7 - 13196 z^6 + 7175 z^5 - 1071 z^4 - 350 z^3 + 154 z^2 - 21 z + 1}{(1 - z)(1 - 2 z)(1 - 3 z)(1 - 4 z)(1 - 5 z)(1 - 6 z)(1 - 7 z)}$$

**1, 7, 28, 140, 784, 4676, 29008, 184820, 1200304, 7907396, 52666768, 353815700, 2393325424, 16279522916, 111239118928, 762963987380, 5249352196144**

## 1^n + 2^n + 3^n + 4^n + 5^n + 6^n + 7^n + 8^n

**Réf.** AS1 813.

**HIS2 A1555** **Recoupements**

**HIS1 N1914** **Fraction rationnelle**

$$\frac{8 - 252 z + 3276 z^2 - 22680 z^3 + 89796 z^4 - 201852 z^5 + 236248 z^6 - 109584 z^7}{(1 - z)(1 - 2 z)(1 - 3 z)(1 - 4 z)(1 - 5 z)(1 - 6 z)(1 - 7 z)(1 - 8 z)}$$

**8, 36, 204, 1296, 8772, 61776, 446964, 3297456, 24684612, 186884496, 1427557524, 10983260016, 84998999652, 660994932816, 5161010498484**

## A simple recurrence

**Réf.** IC 16 351 70.

**HIS2 A1558** **LLL**

**HIS1 N1143** **algébrique**

**(n + 3) a(n) = (- 11/2 n + 21/2) a(n - 3) + (9/2 n + 11/2) a(n - 1) + (- 1/2 n + 9/2) a(n - 2) + (- 2 n + 5) a(n - 4)**

$$\frac{1 - 3z - z^2 - (-(-1 + 4z)(-1 + z + z^2)^2)^{1/2}}{2(2z^4 + z^5)}$$

**1, 3, 10, 33, 111, 379, 1312, 4596, 16266, 58082, 209010, 757259, 2760123, 10114131, 37239072, 137698584, 511140558, 1904038986, 7115422212, 26668376994**

## A simple recurrence

**Réf.** IC 16 351 70.

**HIS2 A1559** **LLL** **Suite P-récurrente**

**HIS1 N1418** **algébrique**

**(n + 4) a(n) = (- 15/2 n + 4) a(n - 3) + (11/2 n + 12) a(n - 1) + (- 4 n + 3) a(n - 2) + (- 2 n + 3) a(n - 4)**

$$\frac{1 - 4z + z^2 + 2z^3 - (-(-1 + 4z)(z^2 + 2z - 1)^2)^{1/2}}{2(2z^5 + z^6)}$$

**1, 4, 15, 54, 193, 690, 2476, 8928, 32358, 117866, 431381, 1585842, 5853849, 21690378, 80650536, 300845232, 1125555054, 4222603968, 15881652606**

**Réf.** JRAM 198 61 57.
**HIS2 A1563** **Hypergéométrique**
**HIS1 N1436** **exponentielle**
**a(n) = (n + 2) a(n-1) + (n - 1) a(n-2)**
**$_3F_2$([1, 1, 1/2], [2, 2], 4 z)**

$$\frac{1 + z}{(1 - z)^3}$$

**1, 4, 18, 96, 600, 4320, 35280, 322560, 3265920, 36288000, 439084800, 5748019200, 80951270400, 1220496076800, 19615115520000, 334764638208000**

---

## 2nd differences of factorial numbers

**Réf.** JRAM 198 61 57.
**HIS2 A1564** **Dérivée logarithmique** **Suite P-récurrente**
**HIS1 N1202** **Fraction rationnelle** **f.g. exponentielle**
**a(n) = (n + 2) a(n - 1) + (- n + 2) a(n - 2)**

$$\frac{(1 + z)^2}{(1 - z)^3}$$

**1, 3, 14, 78, 504, 3720, 30960, 287280, 2943360, 33022080, 402796800, 5308934400, 75203251200, 1139544806400, 18394619443200, 315149522688000**

## 3rd differences of factorial numbers

**Réf.** JRAM 198 61 57.

**HIS2 A1565** **Dérivée logarithmique** **Suite P-récurrente**

**HIS1 N0793** **exponentielle** **f.g. exponentielle**

**a(n) = (3 - n) a(n - 2) + (2 + n) a(n - 1)**

$$-\frac{2}{(z-1)^3}-\frac{3}{(z-1)^2}-\frac{3}{z-1}+\ln(z-1)-1$$

**1, 2, 11, 64, 426, 3216, 27240, 256320, 2656080, 30078720, 369774720, 4906137600, 69894316800, 1064341555200, 17255074636800, 296754903244800**

## From the solution to a Pellian

**Réf.** AMM 56 174 49.

**HIS2 A1570** **Approximants de Padé**

**HIS1 N2108** **Fraction rationnelle**

$$\frac{1-z}{1-14z+z^2}$$

**1, 13, 181, 2521, 35113, 489061, 6811741, 94875313, 1321442641, 18405321661, 256353060613, 3570537526921, 49731172316281, 692665874901013**

## From the solution to a Pellian

**Réf.** AMM 56 175 49.

**HIS2 A1571** **Approximants de Padé**

**HIS1 N0762** **Fraction rationnelle**

$$\frac{z\,(2 - z)}{(1 - z)\,(1 - 4\,z + z^2)}$$

**0, 2, 9, 35, 132, 494, 1845, 6887, 25704, 95930, 358017, 1336139, 4986540, 18610022, 69453549, 259204175, 967363152, 3610248434, 13473630585, 50284273907**

## Winning moves in Fibonacci nim

**Réf.** FQ 3 62 65.

**HIS2 A1581** **Approximants de Padé**

**HIS1 N1359** **Fraction rationnelle**

$$\frac{(1 + z)\,(3\,z^5 + 2\,z^3 + z^2 + z + 2)}{(z^6 + z^5 + z^4 + z^3 + z^2 + z + 1)\,(z - 1)^2}$$

**4, 10, 14, 20, 24, 30, 36, 40, 46, 50, 56, 60, 66, 72, 76, 82, 86, 92, 96, 102, 108, 112, 118, 122, 128, 132, 138, 150, 160, 169, 176, 186, 192, 196, 202, 206, 212, 218, 222**

## Product of Fibonacci and Pell numbers

**Réf.** FQ 3 213 65.

**HIS2 A1582** **Approximants de Padé**

**HIS1 N0779** **Fraction rationnelle**

$$\frac{(1 - z)\,(1 + z)}{1 - 2z - 7z^{2} - 2z^{3} + z^{4}}$$

**1, 2, 10, 36, 145, 560, 2197, 8568, 33490, 130790, 510949, 1995840, 7796413, 30454814, 118965250, 464711184, 1815292333, 7091038640, 27699580729**

## A generalized Fibonacci sequence

**Réf.** FQ 4 244 66.

**HIS2 A1584** **Approximants de Padé**

**HIS1 N0080** **Fraction rationnelle**

$$\frac{(z - 1)^{2}\,(z^{2} + z + 1)}{(z^{4} - z^{3} + 1)\,(z^{4} + z^{3} - 1)}$$

**1, 1, 1, 1, 1, 1, 1, 1, 2, 2, 2, 4, 4, 4, 7, 7, 8, 12, 12, 16, 21, 21, 31, 37, 38, 58, 65, 71, 106, 114, 135, 191, 201, 257, 341, 359, 485, 605, 652, 904, 1070, 1202, 1664, 1894, 2237, 3029, 3370**

**Réf.** FQ 5 288 67.
**HIS2 A1588** **Approximants de Padé**
**HIS1 N0901** **Fraction rationnelle**

$$\frac{1 + z - 3 z^2}{(1 - z)(1 - z - z^2)}$$

**1, 3, 3, 5, 7, 11, 17, 27, 43, 69, 111, 179, 289, 467, 755, 1221, 1975, 3195, 5169, 8363, 13531, 21893, 35423, 57315, 92737, 150051, 242787, 392837, 635623, 1028459**

---

## Tribonacci numbers

**Réf.** FQ 5 211 67.
**HIS2 A1590** **Approximants de Padé**
**HIS1 N0296** **Fraction rationnelle**

$$\frac{249 z^{14} + 249 z^{13} + 249 z^{12} - 249 z^{11} + z - 1}{z^3 + z^2 + z - 1}$$

**1, 0, 1, 2, 3, 6, 11, 20, 37, 68, 125, 479, 423, 778, 1431, 2632, 4841, 8904, 16377, 30122, 55403, 101902, 187427, 344732, 634061, 1166220, 2145013, 3945294, 7256527**

## Pentanacci numbers

**Réf.** FQ 5 260 67.

**HIS2 A1591** **Approximants de Padé**

**HIS1 N0429** **Fraction rationnelle**

$$\frac{1}{1 - z - z^2 - z^3 - z^4 - z^5}$$

**1, 1, 2, 4, 8, 16, 31, 61, 120, 236, 464, 912, 1793, 3525, 6930, 13624, 26784, 52656, 103519, 203513, 400096, 786568, 1546352, 3040048, 5976577, 11749641**

## Hexanacci numbers

**Réf.** FQ 5 260 67.

**HIS2 A1592** **Approximants de Padé**

**HIS1 N0431** **Fraction rationnelle**

$$\frac{1}{1 - z - z^2 - z^3 - z^4 - z^5 - z^6}$$

**1, 1, 2, 4, 8, 16, 32, 63, 125, 248, 492, 976, 1936, 3840, 7617, 15109, 29970, 59448, 117920, 233904, 463968, 920319, 1825529, 3621088, 7182728, 14247536**

**Réf.** FQ 8 267 70.

**HIS2 A1595** **Approximants de Padé**

**HIS1 N0974** **Fraction rationnelle**

$$\frac{1 - z + z^2}{(1 - z)(1 - z - z^2)}$$

**1, 1, 3, 5, 9, 15, 25, 41, 67, 109, 177, 287, 465, 753, 1219, 1973, 3193, 5167, 8361, 13529, 21891, 35421, 57313, 92735, 150049, 242785, 392835, 635621, 1028457**

---

**Related to factors of Fibonacci numbers**

**Réf.** JA66 20.

**HIS2 A1603** **Approximants de Padé**

**HIS1 N2051** **Fraction rationnelle**

$$\frac{1 + 13 z^2 + z^4}{(1 - z)(1 - 3 z + z^2)(z^2 - 7 z + 1)}$$

**1, 11, 101, 781, 5611, 39161, 270281, 1857451, 12744061, 87382901, 599019851, 4105974961, 28143378001, 192899171531, 1322154751061, 9062194370461**

## Related to factors of Fibonacci numbers

**Réf.** JA66 20.

**HIS2 A1604** **Approximants de Padé**

**HIS1 N2042** **Fraction rationnelle**

$$\frac{11 - 90 z + 173 z^2 - 90 z^3 + 11 z^4}{(1 - z)(1 - 3 z + z^2)(z^2 - 7 z + 1)}$$

**11, 31, 151, 911, 5951, 40051, 272611, 1863551, 12760031, 87424711, 599129311, 4106261531, 28144128251, 192901135711, 1322159893351**

---

**Réf.** AMM 15 209 08. JA66 90. FQ 6(3) 68 68.

**HIS2 A1608** **Approximants de Padé**

**HIS1 N0163** **Fraction rationnelle**

$$\frac{z (2 + 3 z)}{1 - z^2 - z^3}$$

**0, 2, 3, 2, 5, 5, 7, 10, 12, 17, 22, 29, 39, 51, 68, 90, 119, 158, 209, 277, 367, 486, 644, 853, 1130, 1497, 1983, 2627, 3480, 4610, 6107, 8090, 10717, 14197, 18807, 24914**

**Réf.** JA66 91. FQ 6(3) 68 68.

**HIS2 A1609** **Approximants de Padé**

**HIS1 N1308** **Fraction rationnelle**

$$\frac{1 + 3 z^2}{1 - z - z^3}$$

**1, 1, 4, 5, 6, 10, 15, 21, 31, 46, 67, 98, 144, 211, 309, 453, 664, 973, 1426, 2090, 3063, 4489, 6579, 9642, 14131, 20710, 30352, 44483, 65193, 95545, 140028, 205221**

---

**Réf.** JA66 96. MOC 15 397 71.

**HIS2 A1610** **Approximants de Padé**

**HIS1 N0291** **Fraction rationnelle**

$$\frac{z (z - 2)}{(z - 1) (1 - z - z^2)}$$

**0, 2, 3, 6, 10, 17, 28, 46, 75, 122, 198, 321, 520, 842, 1363, 2206, 3570, 5777, 9348, 15126, 24475, 39602, 64078, 103681, 167760, 271442, 439203, 710646, 1149850**

## Fibonacci numbers + 1

**Réf.** JA66 97.

**HIS2 A1611** **Approximants de Padé**

**HIS1 N0103** **Fraction rationnelle**

$$\frac{1 - 2z^2}{(z - 1)(1 - z - z^2)}$$

**1, 2, 2, 3, 4, 6, 9, 14, 22, 35, 56, 90, 145, 234, 378, 611, 988, 1598, 2585, 4182, 6766, 10947, 17712, 28658, 46369, 75026, 121394, 196419, 317812, 514230, 832041**

---

**Réf.** JA66 97.

**HIS2 A1612** **Approximants de Padé**

**HIS1 N0364** **Fraction rationnelle**

$$\frac{3z^2 - 2}{(z - 1)(1 - z - z^2)}$$

**2, 4, 5, 8, 12, 19, 30, 48, 77, 124, 200, 323, 522, 844, 1365, 2208, 3572, 5779, 9350, 15128, 24477, 39604, 64080, 103683, 167762, 271444, 439205, 710648, 1149852**

## Convolved Fibonacci numbers

**Réf.** RCI 101. FQ 15 118 77.

**HIS2 A1628** **Approximants de Padé**

**HIS1 N1124** **Fraction rationnelle**

$$\frac{1}{(1 - z - z^2)^3}$$

**1, 3, 9, 22, 51, 111, 233, 474, 942, 1836, 3522, 6666, 12473, 23109, 42447, 77378, 140109, 252177, 451441, 804228, 1426380, 2519640, 4434420, 7777860**

---

## Convolved Fibonacci numbers

**Réf.** RCI 101. FQ 15 118 77.

**HIS2 A1629** **Approximants de Padé**

**HIS1 N0537** **Fraction rationnelle**

$$\frac{1}{(1 - z - z^2)^2}$$

**1, 2, 5, 10, 20, 38, 71, 130, 235, 420, 744, 1308, 2285, 3970, 6865, 11822, 20284, 34690, 59155, 100610, 170711, 289032, 488400, 823800, 1387225, 2332418, 3916061**

## Tetranacci numbers

**Réf.** FQ 8 7 70.

**HIS2 A1630** **Approximants de Padé**

**HIS1 N0301** **Fraction rationnelle**

$$\frac{z\ (1 + z)}{1 - z - z^2 - z^3 - z^4}$$

**0, 0, 1, 2, 3, 6, 12, 23, 44, 85, 164, 316, 609, 1174, 2263, 4362, 8408, 16207, 31240, 60217, 116072, 223736, 431265, 831290, 1602363, 3088654, 5953572, 11475879**

---

## Tetranacci numbers

**Réf.** FQ 8 7 70.

**HIS2 A1631** **Approximants de Padé**

**HIS1 N0410** **Fraction rationnelle**

$$\frac{1 - z}{1 - z - z^2 - z^3 - z^4}$$

**1, 0, 1, 2, 4, 7, 14, 27, 52, 100, 193, 372, 717, 1382, 2664, 5135, 9898, 19079, 36776, 70888, 136641, 263384, 507689, 978602, 1886316, 3635991, 7008598, 13509507**

**Réf.** IDM 8 64 01. FQ 6(3) 68 68.

**HIS2 A1634** **Approximants de Padé**

**HIS1 N0281** **Fraction rationnelle**

$$\frac{z\,(2 + 3\,z + 4\,z^2)}{(1 + z)\,(1 - z - z^3)}$$

**0, 2, 3, 6, 5, 11, 14, 22, 30, 47, 66, 99, 143, 212, 308, 454, 663, 974, 1425, 2091, 3062, 4490, 6578, 9643, 14130, 20711, 30351, 44484, 65192, 95546, 140027, 205222**

---

## A Fielder sequence

**Réf.** FQ 6(3) 68 68.

**HIS2 A1635** **Approximants de Padé**

**HIS1 N0289** **Fraction rationnelle**

$$\frac{z\,(2 + 3\,z + 4\,z^2 + 5\,z^3)}{1 - z^2 - z^3 - z^4 - z^5}$$

**0, 2, 3, 6, 10, 11, 21, 30, 48, 72, 110, 171, 260, 401, 613, 942, 1445, 2216, 3401, 5216, 8004, 12278, 18837, 28899, 44335, 68018, 104349, 160089, 245601, 376791**

## A Fielder sequence

**Réf.** FQ 6(3) 68 68.

**HIS2 A1636** **Approximants de Padé**

**HIS1 N0290** **Fraction rationnelle**

$$\frac{z\,(2 + 3\,z + 4\,z^{2} + 5\,z^{3} + 6\,z^{4})}{(z - 1)\,(z^{5} + z^{3} + z - 1)}$$

**0, 2, 3, 6, 10, 17, 21, 38, 57, 92, 143, 225, 351, 555, 868, 1366, 2142, 3365, 5282, 8296, 13023, 20451, 32108, 50417, 79160, 124295, 195159, 306431, 481139, 755462**

## A Fielder sequence

**Réf.** FQ 6(3) 68 68.

**HIS2 A1638** **Approximants de Padé**

**HIS1 N1348** **Fraction rationnelle**

$$\frac{(1 + z)\,(4\,z^{2} - z + 1)}{(1 - z - z^{2})\,(1 + z^{2})}$$

**1, 1, 4, 9, 11, 16, 29, 49, 76, 121, 199, 324, 521, 841, 1364, 2209, 3571, 5776, 9349, 15129, 24476, 39601, 64079, 103684, 167761, 271441, 439204, 710649, 1149851**

## A Fielder sequence

**Réf.** FQ 6(3) 68 68.

**HIS2 A1639** **Approximants de Padé**

**HIS1 N1349** **Fraction rationnelle**

$$\frac{1 + 3z^2 + 4z^3 + 5z^4}{1 - z - z^3 - z^4 - z^5}$$

**1, 1, 4, 9, 16, 22, 36, 65, 112, 186, 309, 522, 885, 1492, 2509, 4225, 7124, 12010, 20236, 34094, 57453, 96823, 163163, 274946, 463316, 780755, 1315687, 2217112**

## A Fielder sequence

**Réf.** FQ 6(3) 68 68.

**HIS2 A1640** **Approximants de Padé**

**HIS1 N1352** **Fraction rationnelle**

$$\frac{1 + 3z^2 + 4z^3 + 5z^4 + 6z^5}{1 - z - z^3 - z^4 - z^5 - z^6}$$

**1, 1, 4, 9, 16, 28, 43, 73, 130, 226, 386, 660, 1132, 1947, 3349, 5753, 9878, 16966, 29147, 50074, 86020, 147764, 253829, 436036, 749041, 1286728, 2210377, 3797047**

## A Fielder sequence

**Réf.** FQ 6(3) 69 68.

**HIS2 A1641** **Approximants de Padé**

**HIS1 N0935** **Fraction rationnelle**

$$\frac{1 + 2 z + 4 z^{3}}{(1 + z) (z^{3} - z^{2} + 2 z - 1)}$$

**1, 3, 4, 11, 16, 30, 50, 91, 157, 278, 485, 854, 1496, 2628, 4609, 8091, 14196, 24915, 43720, 76726, 134642, 236283, 414645, 727654, 1276941, 2240878, 3932464**

## A Fielder sequence

**Réf.** FQ 6(3) 69 68.

**HIS2 A1642** **Approximants de Padé**

**HIS1 N0937** **Fraction rationnelle**

$$\frac{(1 + z) (5 z^{3} - z^{2} + z + 1)}{1 - z - z^{2} - z^{4} - z^{5}}$$

**1, 3, 4, 11, 21, 36, 64, 115, 211, 383, 694, 1256, 2276, 4126, 7479, 13555, 24566, 44523, 80694, 146251, 265066, 480406, 870689, 1578040, 2860046, 5183558, 9394699**

## A Fielder sequence

**Réf.** FQ 6(3) 69 68.

**HIS2 A1643** **Approximants de Padé**

**HIS1 N0938** **Fraction rationnelle**

$$\frac{1 + 2z + 4z^3 + 5z^4 + 6z^5}{(1 + z)(1 - z - z^2 - z^3)(1 - z + z^2)}$$

**1, 3, 4, 11, 21, 42, 71, 131, 238, 443, 815, 1502, 2757, 5071, 9324, 17155, 31553, 58038, 106743, 196331, 361106, 664183, 1221623, 2246918, 4132721, 7601259**

## A Fielder sequence

**Réf.** FQ 6(3) 69 68.

**HIS2 A1644** **Approximants de Padé**

**HIS1 N1040** **Fraction rationnelle**

$$\frac{1 + 2z + 3z^2}{1 - z - z^2 - z^3}$$

**1, 3, 7, 11, 21, 39, 71, 131, 241, 443, 815, 1499, 2757, 5071, 9327, 17155, 31553, 58035, 106743, 196331, 361109, 664183, 1221623, 2246915, 4132721, 7601259**

## A Fielder sequence

**Réf.** FQ 6(3) 69 68.

**HIS2 A1645** **Approximants de Padé**

**HIS1 N1041** **Fraction rationnelle**

$$\frac{1 + 2z + 3z^{2} + 5z^{4}}{1 - z - z^{2} - z^{3} - z^{5}}$$

**1, 3, 7, 11, 26, 45, 85, 163, 304, 578, 1090, 2057, 3888, 7339, 13862, 26179, 49437, 93366, 176321, 332986, 628852, 1187596, 2242800, 4235569, 7998951**

## A Fielder sequence

**Réf.** FQ 6(3) 70 68.

**HIS2 A1648** **Approximants de Padé**

**HIS1 N1055** **Fraction rationnelle**

$$\frac{1 + 2z + 3z^{2} + 4z^{3}}{1 - z - z^{2} - z^{3} - z^{4}}$$

**1, 3, 7, 15, 26, 51, 99, 191, 367, 708, 1365, 2631, 5071, 9775, 18842, 36319, 70007, 134943, 260111, 501380, 966441, 1862875, 3590807, 6921503, 13341626**

## A Fielder sequence

**Réf.** FQ 6(3) 70 68.

**HIS2 A1649** **Approximants de Padé**

**HIS1 N1056** **Fraction rationnelle**

$$\frac{1 + 2z + 3z^2 + 4z^3 + 6z^5}{1 - z - z^2 - z^3 - z^4 - z^6}$$

**1, 3, 7, 15, 26, 57, 106, 207, 403, 788, 1530, 2985, 5812, 11322, 22052, 42959, 83675, 162993, 317491, 618440, 1204651, 2346534, 4570791, 8903409, 17342876**

---

**Réf.** FQ 6(3) 261 68.

**HIS2 A1651** **Approximants de Padé**

**HIS1 N0357** **Fraction rationnelle**

$$\frac{z^2 + z + 1}{(1 + z)\,(z - 1)^2}$$

**1, 2, 4, 5, 7, 8, 10, 11, 13, 14, 16, 17, 19, 20, 22, 23, 25, 26, 28, 29, 31, 32, 34, 35, 37, 38, 40, 41, 43, 44, 46, 47, 49, 50, 52, 53, 55, 56, 58, 59, 61, 62, 64, 65, 67, 68, 70, 71**

## Pythagorean triangles

**Réf.** MLG 2 322 10. FQ 6(3) 104 68.

**HIS2 A1652** **Approximants de Padé**

**HIS1 N1247** **Fraction rationnelle**

$$\frac{z\,(z-3)}{(z-1)\,(z^2-6z+1)}$$

**0, 3, 20, 119, 696, 4059, 23660, 137903, 803760, 4684659, 27304196, 159140519, 927538920, 5406093003, 31509019100, 183648021599, 1070379110496**

---

**Réf.** AMM 4 25 1897. MLG 2 322 10. FQ 6(3) 104 68.

**HIS2 A1653** **Approximants de Padé**

**HIS1 N1630** **Fraction rationnelle**

$$\frac{1-5z}{z^2-6z+1}$$

**1, 1, 5, 29, 169, 985, 5741, 33461, 195025, 1136689, 6625109, 38613965, 225058681, 1311738121, 7645370045, 44560482149, 259717522849, 1513744654945**

## Product of successive Fibonacci numbers

**Réf.** FQ 6 82 68. BR72 17.

**HIS2 A1654** **Approximants de Padé**

**HIS1 N0628** **Fraction rationnelle**

$$\frac{1}{(1 + z)\,(1 - 3z + z^2)}$$

**1, 2, 6, 15, 40, 104, 273, 714, 1870, 4895, 12816, 33552, 87841, 229970, 602070, 1576239, 4126648, 10803704, 28284465, 74049690, 193864606, 507544127**

## Fibonomial coefficients

**Réf.** FQ 6 82 68. BR72 74.

**HIS2 A1655** **Approximants de Padé**

**HIS1 N1208** **Fraction rationnelle**

$$\frac{1}{(z^2 - z - 1)\,(- 1 + 4z + z^2)}$$

**1, 3, 15, 60, 260, 1092, 4641, 19635, 83215, 352440, 1493064, 6324552, 26791505, 113490195, 480752895, 2036500788, 8626757644, 36543528780**

## Fibonomial coefficients

**Réf.** FQ 6 82 68. BR72 74.

**HIS2 A1656** **Approximants de Padé**

**HIS1 N1653** **Fraction rationnelle**

$$\frac{1}{(1 - z)(z^2 - 7z + 1)(z^2 + 3z + 1)}$$

**1, 5, 40, 260, 1820, 12376, 85085, 582505, 3994320, 27372840, 187628376, 1285992240, 8814405145, 60414613805, 41408893560, 2838203264876, 19453338487220**

## Fibonomial coefficients

**Réf.** FQ 6 82 68. BR72 74.

**HIS2 A1657** **Approximants de Padé**

**HIS1 N1945** **Fraction rationnelle**

$$\frac{1}{(z^2 + 11z - 1)(z^2 - 4z - 1)(1 - z - z^2)}$$

**1, 8, 104, 1092, 12376, 136136, 1514513, 16776144, 186135312, 2063912136, 22890661872, 253854868176, 2815321003313, 31222272414424, 34620798314872**

## Fibonomial coefficients

**Réf.** FQ 6 82 68. BR72 74.

**HIS2 A1658** **Approximants de Padé**

**HIS1 N2112** **Fraction rationnelle**

$$\frac{1}{(z + 1)\,(z^2 - 18\,z + 1)\,(z^2 - 3\,z + 1)\,(z^2 + 7\,z + 1)}$$

**1, 13, 273, 4641, 85085, 1514513, 27261234, 488605194, 8771626578, 157373300370, 2824135408458, 50675778059634, 909348684070099**

## Coefficients of iterated exponentials

**Réf.** SMA 11 353 45. PRV A32 2342 85.

**HIS2 A1669** **Recoupements**

**HIS1 N1879** **exponentielle**

$$\exp(\exp(\exp(\exp(\exp(\exp(\exp(z) - 1) - 1) - 1) - 1) - 1) - 1)$$

**1, 1, 7, 70, 910, 14532, 274778, 5995892, 148154860, 4085619622, 124304629050, 4133867297490, 149114120602860, 5796433459664946, 241482353893283349**

## The partition function G(n,3)

**Réf.** CMB 1 87 58.

**HIS2 A1680** **Dérivée logarithmique** **Suite P-récurrente**

**HIS1 N0579** **exponentielle**

**2 a(n) = (n^2 - 5 n + 6) a(n - 3) + 2 a(n - 1) + (2 n - 4) a(n - 2)**

$$\exp(z + 1/2\, z^{2} + 1/6\, z^{3})$$

**1, 1, 2, 5, 14, 46, 166, 652, 2780, 12644, 61136, 312676, 1680592, 9467680, 55704104, 341185496, 2170853456, 14314313872, 97620050080, 687418278544**

## The partition function G(n,4)

**Réf.** CMB 1 87 58.

**HIS2 A1681** **Dérivée logarithmique** **Suite P-récurrente**

**HIS1 N0584** **exponentielle**

**6 a(n) = (6 n - 12) a(n - 2) + 6 a(n - 1) + (3 n^2 - 15 n + 18) a(n - 3)**
**+ (n^3 - 9 n^2 + 26 n - 24) a(n - 4)**

$$\exp(z + 1/2\, z^{2} + 1/6\, z^{3} + 1/24\, z^{4})$$

**1, 1, 2, 5, 15, 51, 196, 827, 3795, 18755, 99146, 556711, 3305017, 20655285, 135399720, 927973061, 6631556521, 49294051497, 380306658250, 3039453750685**

**Réf.** MMAG 41 17 68.

**HIS2 A1687** **Approximants de Padé**

**HIS1 N0338** **Fraction rationnelle**

$$\frac{z}{1 - z^2 - z^5}$$

**0, 1, 0, 1, 0, 1, 1, 1, 2, 1, 3, 2, 4, 4, 5, 7, 7, 11, 11, 16, 18, 23, 29, 34, 45, 52, 68, 81, 102, 126, 154, 194, 235, 296, 361, 450, 555, 685, 851, 1046, 1301, 1601, 1986, 2452, 3032, 3753, 4633**

---

## 4th differences of factorial numbers

**Réf.** JRAM 198 61 57.

**HIS2 A1688** **Dérivée exponentielle** **Suite P-récurrente**

**HIS1 N1980**

**a(n) = (3 + n) a(n - 1) + (3 - n) a(n - 2)**

$$\frac{2 z (2 z^2 + 3 z - 4)}{(1 - z)^4} - \ln(- z + 1) + 1$$

**1, 9, 53, 362, 2790, 24024, 229080, 2399760, 27422640, 339696000, 4536362880, 64988179200, 994447238400, 16190733081600, 279499828608000**

## 5th differences of factorial numbers

**Réf.** JRAM 198 61 57.

**HIS2 A1689** **Dérivée** **Suite P-récurrente**

**HIS1 N1920** **exponentielle**

**a(n) = (4 + n) a(n - 1) + (3 - n) a(n - 2)**

$$\ln(1 - z) + \frac{5 z^4 - 10 z^3 + 20 z^2 + 9}{(1 - z)^5} - 1$$

**8, 44, 309, 2428, 21234, 205056, 2170680, 25022880, 312273360, 4196666880, 60451816320, 929459059200, 15196285843200, 263309095526400**

---

**Réf.** RS3.

**HIS2 A1700** **Hypergéométrique** **Suite P-récurrente**

**HIS1 N1144** **algébrique**

$_2F_1([1, 3/2], [2], 4 z)$

$$\frac{- 1 + 4 z + (1 - 4 z)^{1/2}}{2 (1 - 4 z)}$$

**1, 3, 10, 35, 126, 462, 1716, 6435, 24310, 92378, 352716, 1352078, 5200300, 20058300, 77558760, 300540195, 1166803110, 4537567650, 17672631900**

## Generalized Stirling numbers

**Réf.** PEF 77 7 62.

**HIS2 A1701** **Approximants de Padé**

**HIS1 N1735** **Fraction rationnelle**

$$\frac{1 - z - 6 z^{2} + 9 z^{3} - 5 z^{4} + z^{5}}{(1 - z)^{5}}$$

**1, 6, 26, 71, 155, 295, 511, 826, 1266, 1860, 2640, 3641, 4901, 6461, 8365, 10660, 13396, 16626, 20406, 24795, 29855, 35651, 42251, 49726, 58150, 67600, 78156**

## Generalized Stirling numbers

**Réf.** PEF 77 7 62.

**HIS2 A1702** **Approximants de Padé**

**HIS1 N2234** **Fraction rationnelle**

$$\frac{1 - 17 z - 7 z^{2} + 29 z^{3} - 34 z^{4} + 21 z^{5} - 7 z^{6} + z^{7}}{(1 - z)^{7}}$$

**1, 24, 154, 580, 1665, 4025, 8624, 16884, 30810, 53130, 87450, 138424, 211939, 315315, 457520, 649400, 903924, 1236444, 1664970, 2210460, 2897125, 3752749**

## Generalized Stirling numbers

**Réf.** PEF 77 7 62.

**HIS2 A1705** **Tableaux généralisés** **Suite P-récurrente**

**HIS1 N1625** **exponentielle (log)**

**a(n) = (1 + 2 n) a(n - 1) - n^2 a(n - 2)**

$$\frac{-\ln(-z+1)}{(1-z)^2}$$

**1, 5, 26, 154, 1044, 8028, 69264, 663696, 6999840, 80627040, 1007441280, 13575738240, 196287356160, 3031488633600, 49811492505600**

## Generalized Stirling numbers

**Réf.** PEF 77 7 62.

**HIS2 A1706** **Tableaux généralisés** **Suite P-récurrente**

**HIS1 N1988** **exponentielle (log)**

**a(n) = (3 n^2 + 3 n^3 ) a(n - 1) + (- 3 n^2 - 3 n^3 - n ) a(n - 2) + a(n - 3)**

$$\frac{\ln(1-z)^2}{(1-z)^2}$$

**1, 9, 71, 580, 5104, 48860, 509004, 5753736, 70290936, 924118272, 13020978816, 195869441664, 3134328981120, 53180752331520, 953884282141440**

## Generalized Stirling numbers

**Réf.** PEF 77 7 62.

**HIS2 A1707** **Tableaux généralisés exponentielle (log)** **Suite P-récurrente**

**HIS1 N2119**

$$\frac{\ln(1 - z)^3}{6\ (z - 1)^2}$$

**1, 14, 155, 1665, 18424, 214676, 2655764, 34967140, 489896616, 7292774280, 115119818736, 1922666722704, 33896996544384, 629429693586048**

## Generalized Stirling numbers

**Réf.** PEF 77 7 62.

**HIS2 A1708** **Tableaux généralisés exponentielle (log)** **Suite P-récurrente**

**HIS1 N2206**

$$\frac{\ln(1 - z)^4}{24\ (1 - z)^2}$$

**1, 20, 295, 4025, 54649, 761166, 11028590, 167310220, 2664929476, 44601786944, 784146622896, 14469012689040, 279870212258064, 5667093514231200**

## Generalized Stirling numbers

**Réf.** PEF 77 7 62.

**HIS2 A1709** **Tableaux généralisés exponentielle (log)** **Suite P-récurrente**

**HIS1 N2259**

$$\frac{\ln(1 - z)^{5}}{120\ (z - 1)^{2}}$$

**1, 27, 511, 8624, 140889, 2310945, 38759930, 671189310, 12061579816, 225525484184, 4392554369840, 89142436976320, 1884434077831824**

---

**Réf.** PEF 77 26 62.

**HIS2 A1710** **Dérivée logarithmique Fraction rationnelle** **f.g. exponentielle**

**HIS1 N1179**

$$\frac{1}{(1 - z)^{3}}$$

**1, 3, 12, 60, 360, 2520, 20160, 181440, 1814400, 19958400, 239500800, 3113510400, 43589145600, 653837184000, 10461394944000, 177843714048000**

## Generalized Stirling numbers

**Réf.** PEF 77 26 62.

**HIS2 A1711** **Tableaux généralisés** **Suite P-récurrente**

**HIS1 N1873** **exponentielle**

**a(n) = -(n^2 + 2 n + 1) a(n - 2) + (2 n + 3) a(n - 1)**

$$\frac{-\ln(1-z)}{(1-z)^3}$$

**1, 7, 47, 342, 2754, 24552, 241128, 2592720, 30334320, 383970240, 5231113920, 76349105280, 1188825724800, 19675048780800, 344937224217600**

## Generalized Stirling numbers

**Réf.** PEF 77 26 62.

**HIS2 A1712** **Tableaux généralisés** **Suite P-récurrente**

**HIS1 N2077** **exponentielle (log)**

**a(n) = (3 n^2 + 6 n^3 ) a(n - 1) - ( 3 n + 9 n^2 + 7 n^3 ) a(n - 2)**
**+ (1 + 3 n + 3 n^2 + n^3 ) a(n - 3)**

$$\frac{\ln(1-z)^2}{2(1-z)^3}$$

**1, 12, 119, 1175, 12154, 133938, 1580508, 19978308, 270074016, 3894932448, 59760168192, 972751628160, 16752851775360, 304473528961920**

## Generalized Stirling numbers

**Réf.** PEF 77 26 62.

**HIS2 A1713** **Tableaux généralisés exponentielle** **Suite P-récurrente**
**HIS1 N2190**

$$\frac{\ln(1 - z)^3}{6\ (z - 1)^3}$$

**1, 18, 245, 3135, 40369, 537628, 7494416, 109911300, 1698920916, 27679825272, 474957547272, 8572072384512, 162478082312064, 3229079010579072**

---

**Réf.** PEF 77 26 62.

**HIS2 A1714** **Tableaux généralisés exponentielle** **Suite P-récurrente**
**HIS1 N2252**

$$\frac{\ln(1 - z)^4}{24\ (1 - z)^3}$$

**1, 25, 445, 7140, 111769, 1767087, 28699460, 483004280, 8460980836, 154594537812, 2948470152264, 58696064973000, 1219007251826064**

**Réf.** PEF 77 44 62.
**HIS2 A1715** **Dérivée logarithmique** **f.g. exponentielle**
**HIS1 N1445** **Fraction rationnelle**

$$\frac{1}{(z - 1)^4}$$

**1, 4, 20, 120, 840, 6720, 60480, 604800, 6652800, 79833600, 1037836800, 14529715200, 217945728000, 3487131648000, 59281238016000**

---

## Generalized Stirling numbers

**Réf.** PEF 77 44 62.
**HIS2 A1716** **Tableaux généralisés** **Suite P-récurrente**
**HIS1 N1990** **exponentielle**
**a(n) = - ( n^2 + 4 n + 4) a(n - 2) + (2 n + 5) a(n - 1)**

$$\frac{4 \ln(1 - z) - 1}{(1 - z)^5}$$

**1, 9, 74, 638, 5944, 60216, 662640, 7893840, 101378880, 1397759040, 20606463360, 323626665600, 5395972377600, 95218662067200, 1773217155225600**

## Generalized Stirling numbers

**Réf.** PEF 77 44 62.

**HIS2 A1717** **Tableaux généralisés** **Suite P-récurrente**

**HIS1 N2143** **exponentielle** **Formule de B. Salvy**

**a(n) = - ( 9 n^3 + 3 n^2 ) a(n - 1) + (19 n^3 + 15 n^2 + 3 n) a(n - 2)**
**- ( 8 n^3 + 12 n^2 + 6 n + 1) a(n - 3)**

$$\frac{10 \ln(1 - z)^2 - 9 \ln(1 - z) + 1}{(1 - z)^6}$$

**1, 15, 179, 2070, 24574, 305956, 4028156, 56231712, 832391136, 13051234944, 216374987520, 3785626465920, 69751622298240, 1350747863435520**

---

**Réf.** PEF 77 61 62.

**HIS2 A1720** **Approximants de Padé** **f.g. exponentielle**

**HIS1 N1634** **Fraction rationnelle**

$$\frac{1}{(1 - z)^5}$$

**1, 5, 30, 210, 1680, 15120, 151200, 1663200, 19958400, 259459200, 3632428800, 54486432000, 871782912000, 14820309504000, 266765571072000**

## Generalized Stirling numbers

**Réf.** PEF 77 61 62.

**HIS2 A1721** **Tableaux généralisés** **Suite P-récurrente**

**HIS1 N2052** **exponentielle**

**a(n) = (2 n + 7) a(n - 1) - (n^2 + 6 n + 9) a(n - 2)**

$$\frac{1 - 5 \ln(1 - z)}{(z - 1)^{6}}$$

**1, 11, 107, 1066, 11274, 127860, 1557660, 20355120, 284574960, 4243508640, 67285058400, 1131047366400, 20099588140800, 376612896038400**

## Generalized Stirling numbers

**Réf.** PEF 77 61 62.

**HIS2 A1722** **Tableaux généralisés** **Suite P-récurrente**

**HIS1 N2191** **exponentielle:log**

**a(n) = (3n + 12) a(n-1) - (3 n^2 - 21 n - 37) a(n-2)**
**+ (n^3 + 9 n^2 + 27 n + 27) a(n-3)**

$$\frac{1 + 15 \ln(1 - z)^{2} - 11 \ln(1 - z)}{(1 - z)^{7}}$$

**1, 18, 251, 3325, 44524, 617624, 8969148, 136954044, 2201931576, 37272482280, 663644774880, 12413008539360, 243533741849280, 5003753991174720**

**Réf.** PEF 107 5 63.

**HIS2 A1725** **Dérivée logarithmique** **f.g. exponentielle**

**HIS1 N1772** **Fraction rationnelle**

$$\frac{1}{(1 - z)^{6}}$$

**1, 6, 42, 336, 3024, 30240, 332640, 3991680, 51891840, 726485760, 10897286400, 174356582400, 2964061900800, 53353114214400, 1013709170073600**

---

**Réf.** PEF 107 19 63.

**HIS2 A1730** **Dérivée logarithmique** **f.g. exponentielle**

**HIS1 N1876** **Fraction rationnelle**

$$\frac{1}{(1 - z)^{7}}$$

**1, 7, 56, 504, 5040, 55440, 665280, 8648640, 121080960, 1816214400, 29059430400, 494010316800, 8892185702400, 168951528345600, 3379030566912000**

## Lah numbers

**Réf.** R1 44. C1 156.

**HIS2 A1754** **Dérivée logarithmique**

**HIS1 N2079** **Fraction rationnelle**

$$\frac{3 z^2 + 6 z + 1}{(z - 1)^6}$$

**1, 12, 120, 1200, 12600, 141120, 1693440, 21772800, 299376000, 4390848000, 68497228800, 1133317785600, 19833061248000, 366148823040000**

## Lah numbers

**Réf.** R1 44. C1 156.

**HIS2 A1755** **Dérivée logarithmique**

**HIS1 N2207** **Fraction rationnelle**

$$\frac{4 z^3 + 18 z^2 + 12 z + 1}{(z - 1)^8}$$

**1, 20, 300, 4200, 58800, 846720, 12700800, 199584000, 3293136000, 57081024000, 1038874636800, 19833061248000, 396661224960000**

## Expansion of an integral

**Réf.** C1 167.

**HIS2 A1756** **Hypergéométrique** **Suite P-récurrente**
**HIS1 N2131** **algébrique** **f.g. exponentielle**

$$\frac{15 z (2 - 6 z + 5 z^2)}{2 (1 - 2 z)^{5/2}}$$

**15, 60, 450, 4500, 55125, 793800, 13097700**

## Dissections of a disk

**Réf.** CMA 2 25 70. MAN 191 98 71.

**HIS2 A1761** **Hypergéométrique** **Inverse fonctionnel de A1561**
**HIS1 N1478** **algébrique** **Suite P-récurrente.**

**3F2([1, 1, 1/2],[2, 2],4 z)**
**n a(n) = 2 (n - 1) (2 n - 3) a(n - 1)**

$$\frac{1 - (1 - 4 z)^{1/2}}{2 z}$$

**1, 1, 4, 30, 336, 5040, 95040, 2162160, 57657600**

## Dissections of a ball

**Réf.** CMA 2 25 70. MAN 191 98 71.

**HIS2 A1763** **Inverse fonctionnel** **Suite P-récurrente**

**HIS1 N1788** **algébrique 3è degré**

**S(z) est l'inverse de**

$$\frac{z}{(1+z)^3}$$

**1, 1, 6, 72, 1320, 32760, 1028160, 39070080**

---

## Binomial coefficients C(3n,n-1)/n

**Réf.** CMA 2 25 70. MAN 191 98 71. FQ 11 125 73. DM 9 355 74.

**HIS2 A1764** **Hypergéométrique** **Suite P-récurrente**

**HIS1 N1174** **algébrique 3è degré** **f.g. exponentielle**

**$_3F_2$ ([1, 5/3, 4/3],[2, 5/2], 27 z /4)**

**S(z) est racine**

**de**

$$1 - S(z) + 3\,S(z)\,z + 3\,S(z)^2\,z^2 + S(z)^3\,z^3$$

**1, 3, 12, 55, 273, 1428, 7752, 43263, 246675, 1430715, 8414640, 50067108, 300830572, 1822766520, 11124755664, 68328754959, 422030545335, 2619631042665**

## Coefficients of iterated exponentials

**Réf.** SMA 11 353 45. PRV A32 2342 85.

**HIS2 A1765** **Recoupements**

**HIS1 N1882** **exponentielle**

$$-\ln(1+\ln(1+\ln(1+\ln(1+\ln(1+\ln(1+\ln(1-z)))))))+1$$

**1, 1, 7, 77, 1155, 21973, 506989, 13761937, 429853851, 15192078027, 599551077881, 26140497946017, 1248134313062231, 64783855286002573**

## Number of comparisons for merge sort of n elements

**Réf.** AMM 66 389 59. WE71 207. KN1 3 187.

**HIS2 A1768** **Approximants de Padé**

**HIS1 N0954** **Fraction rationnelle**

$$\frac{(z+1)(z^6-z^3+z+1)(z^2-z+1)}{(z-1)^2}$$

**0, 1, 3, 5, 7, 10, 13, 16, 19, 22, 26, 30, 34, 38, 42, 46, 50, 54, 58, 62, 66, 71, 76, 81, 86, 91, 96, 101, 106, 111, 116, 121, 126**

## Lah numbers

**Réf.** R1 44. C1 156.

**HIS2 A1777** **Dérivée logarithmique**

**HIS1 N2267** **exponentielle**

$$\frac{5 z^4 + 40 z^3 + 60 z^2 + 20 z + 1}{(z - 1)^{10}}$$

**1, 30, 630, 11760, 211680, 3810240, 69854400, 1317254400, 25686460800, 519437318400, 10908183686400, 237996734976000, 5394592659456000**

## Lah numbers

**Réf.** R1 44. C1 156.

**HIS2 A1778** **Dérivée logarithmique**

**HIS1 N2297** **exponentielle**

$$\frac{6 z^5 + 75 z^4 + 200 z^3 + 150 z^2 + 30 z + 1}{(z - 1)^{12}}$$

**1, 42, 1176, 28224, 635040, 13970880, 307359360, 6849722880, 155831195520, 3636061228800, 87265469491200, 2157837063782400, 55024845126451200**

**Réf.** PRSE 62 190 46. BIO 46 422 59. AS1 796.
**HIS2 A1787** **Approximants de Padé**
**HIS1 N1398** **Fraction rationnelle**

$$\frac{1}{(1 - 2z)^2}$$

**1, 4, 12, 32, 80, 192, 448, 1024, 2304, 5120, 11264, 24576, 53248, 114688, 245760, 524288, 1114112, 2359296, 4980736, 10485760, 22020096, 46137344**

**Réf.** PRSE 62 190 46. AS1 796. MFM 74 62 70.
**HIS2 A1788** **Approximants de Padé**
**HIS1 N1729** **Fraction rationnelle**

$$\frac{1}{(1 - 2z)^3}$$

**1, 6, 24, 80, 240, 672, 1792, 4608, 11520, 28160, 67584, 159744, 372736, 860160, 1966080, 4456448, 10027008, 22413312, 49807360, 110100480, 242221056**

**Réf.** PRSE 62 190 46. AS1 796. MFM 74 62 70.

**HIS2 A1789** **Approximants de Padé**

**HIS1 N1916** **Fraction rationnelle**

$$\frac{1}{(1 - 2z)^4}$$

**1, 8, 40, 160, 560, 1792, 5376, 15360, 42240, 112640, 292864, 745472, 1863680, 4587520, 11141120, 26738688, 63504384, 149422080, 348651520, 807403520**

---

## Binomial coefficients C(2n,n-1)

**Réf.** LA56 517. AS1 828. PLC 1 292 70.

**HIS2 A1791** **Hypergéométrique** **Suite P-récurrente**

**HIS1 N1421** **algébrique**

$$\frac{4z}{(1 - 4z)^{1/2}\,(1 + (1 - 4z)^{1/2})^2}$$

**1, 4, 15, 56, 210, 792, 3003, 11440, 43758, 167960, 646646, 2496144, 9657700, 37442160, 145422675, 565722720, 2203961430, 8597496600, 33578000610**

**Réf.** PRSE 62 190 46. AS1 795.
**HIS2 A1792** **Approximants de Padé**
**HIS1 N1100** **Fraction rationnelle**

$$\frac{4z - 3}{(1 - 2z)^2}$$

**3, 8, 20, 48, 112, 256, 576, 1280, 2816, 6144, 13312, 28672, 61440, 131072, 278528, 589824, 1245184, 2621440, 5505024, 11534336, 24117248, 50331648**

---

## Coefficients of Chebyshev polynomials

**Réf.** PRSE 62 190 46. AS1 795.
**HIS2 A1793** **Approximants de Padé**
**HIS1 N1591** **Fraction rationnelle**

$$\frac{1 - z}{(1 - 2z)^3}$$

**1, 5, 18, 56, 160, 432, 1120, 2816, 6912, 16640, 39424, 92160, 212992, 487424, 1105920, 2490368, 5570560**

## Coefficients of Chebyshev polynomials

**Réf.** PRSE 62 190 46. AS1 795.

**HIS2 A1794** **Approximants de Padé**

**HIS1 N1859** **Fraction rationnelle**

$$\frac{1 - z}{(1 - 2z)^4}$$

**1, 7, 32, 120, 400, 1232, 3584, 9984, 26880, 70400, 180224, 452608, 1118208, 2723840, 6553600**

---

**Réf.** AS1 799.

**HIS2 A1804** **Dérivée logarithmique** **Suite P-récurrente**

**HIS1 N0834** **exponentielle**

**a(n) = (n + 7) a(n-1) - (4 n + 6) a(n-2) + (2 n - 2) a(n-3)**

$$\frac{z (z + 2)}{(1 - z)^4}$$

**2, 18, 144, 1200, 10800, 105840, 1128960, 13063680, 163296000, 2195424000, 31614105600, 485707622400, 7933224499200, 137305808640000, 2510734786560000**

## Coefficients of Laguerre polynomials

**Réf.** AS1 799.

**HIS2 A1805** **Hypergéométrique** **f.g. exponentielle**

**HIS1 N1794** **Fraction rationnelle**

$$\frac{2z(z^2+6z+3)}{(z-1)^6}$$

**6, 96, 1200, 14400, 176400, 2257920, 30481920, 435456000, 6586272000, 105380352000**

## Coefficients of Laguerre polynomials

**Réf.** AS1 799.

**HIS2 A1806** **Hypergéométrique** **f.g. exponentielle**

**HIS1 N2242** **Fraction rationnelle**

$$\frac{6z(4+18z+12z^2+z^3)}{(z-1)^8}$$

**24, 600, 10800, 176400, 2822400, 45722880, 762048000, 13172544000, 237105792000**

## Coefficients of Laguerre polynomials

**Réf.** AS1 799.

**HIS2 A1807** **Hypergéométrique** **f.g. exponentielle**

**HIS1 N2337** **Fraction rationnelle**

$$\frac{24\,(5 + 40\,z + 60\,z^2 + 20\,z^3 + z^4)\,z}{(z - 1)^{10}}$$

**120, 4320, 105840, 2257920, 45722880, 914457600, 18441561600, 379369267200**

## Coefficients of Laguerre polynomials

**Réf.** LA56 519. AS1 799.

**HIS2 A1809** **Hypergéométrique** **f.g. exponentielle**

**HIS1 N1989** **Fraction rationnelle**

$$\frac{z\,(2 + z)}{2\,(z - 1)^4}$$

**1, 9, 72, 600, 5400, 52920, 564480, 6531840, 81648000, 1097712000, 15807052800**

## Coefficients of Laguerre polynomials

**Réf.** LA56 519. AS1 799.

**HIS2 A1810** **Hypergéométrique** **f.g. exponentielle**

**HIS1 N2163** **Fraction rationnelle**

$$\frac{(z^2 + 6z + 3)\, z}{3\,(z-1)^6}$$

**1, 16, 200, 2400, 29400, 376320, 5080320, 72576000, 1097712000, 17563392000**

## Coefficients of Laguerre polynomials

**Réf.** LA56 519. AS1 799.

**HIS2 A1811** **Hypergéométrique** **f.g. exponentielle**

**HIS1 N2253** **Fraction rationnelle**

$$\frac{z\,(18z + 4 + 12z^2 + z^3)}{4\,(z-1)^8}$$

**1, 25, 450, 7350, 117600, 1905120, 31752000, 548856000, 9879408000**

## Coefficients of Laguerre polynomials

**Réf.** LA56 519. AS1 799.

**HIS2 A1812** **Hypergéométrique** **f.g. exponentielle**

**HIS1 N2289** **Fraction rationnelle**

$$\frac{(40 z + z^4 + 60 z^2 + 20 z^3 + 5) z}{5 (z - 1)^{10}}$$

**1, 36, 882, 18816, 381024, 7620480, 153679680, 3161410560**

## Produit des nombres impairs : 1.3.5.7. ... x (2^n)

**Réf.** MOC 3 168 48.

**HIS2 A1813** **Hypergéométrique** **Suite P-récurrente**

**HIS1 N0808** **algébrique** **f.g. exponentielle**

$$\frac{2 z}{1 + (1 - 4 z)^{1/2}}$$

**1, 2, 12, 120, 1680, 30240, 665280, 17297280, 518918400, 17643225600, 670442572800, 28158588057600, 1295295050649600, 64764752532480000**

## Coefficients of Hermite polynomials

**Réf.** MOC 3 168 48.

**HIS2 A1814** **Hypergéométrique** **Suite P-récurrente**

**HIS1 N2088** **algébrique** **f.g. exponentielle**

$$\frac{(1 + 2z)}{(1 - 4z)^{5/2}}$$

**12, 180, 3360, 75600, 1995840, 60540480, 2075673600, 79394515200, 3352212864000, 154872234316800, 7771770303897600, 420970891461120000**

---

**Réf.** AS1 801.

**HIS2 A1815** **Approximants de Padé**

**HIS1 N0799** **Fraction rationnelle**

$$\frac{2z}{(1 - 2z)^{3}}$$

**0, 2, 12, 48, 160, 480, 1344, 3584, 9216, 23040, 56320, 135168, 319488, 745472, 1720320, 3932160, 8912896, 20054016, 44826624, 99614720, 220200960, 484442112, 1061158912**

## Coefficients of Hermite polynomials

**Réf.** AS1 801.

**HIS2 A1816** **Approximants de Padé**

**HIS1 N2078** **Fraction rationnelle**

$$\frac{12}{(1 - 2 z)^{5}}$$

**12, 120, 720, 3360, 13440, 48384, 161280, 506880, 1520640**

---

**Réf.** RCI 217.

**HIS2 A1818** **hypergéométrique** **Suite P-récurrente**

**HIS1 N1997** **intégrales elliptiques** **double exponentielle**

$${}_2F_1([1/2, 1/2], [1], 4 z) - 1$$

**1, 9, 225, 11025, 893025, 108056025, 18261468225, 4108830350625, 1187451971330625, 428670161650355625, 189043541287806830625**

## Central factorial numbers

**Réf.** RCI 217.

**HIS2 A1823** **Approximants de Padé**

**HIS1 N1998** **Fraction rationnelle**

$$\frac{9 + 196 z + 350 z^2 + 84 z^3 + z^4}{(1 - z)^7}$$

**9, 259, 1974, 8778, 28743, 77077, 179452, 375972, 725781, 1312311, 2249170, 3686670, 5818995, 8892009, 13211704, 19153288, 27170913, 37808043**

---

**Réf.** EUL (1) 1 375 11. MMAG 40 78 67.

**HIS2 A1834** **Approximants de Padé**

**HIS1 N1598** **Fraction rationnelle**

$$\frac{1 + z}{1 - 4 z + z^2}$$

**1, 5, 19, 71, 265, 989, 3691, 13775, 51409, 191861, 716035, 2672279, 9973081, 37220045, 138907099, 518408351, 1934726305, 7220496869, 26947261171**

**Réf.** EUL (1) 1 375 11. MMAG 40 78 67.

**HIS2 A1835** **Approximants de Padé**

**HIS1 N1160** **Fraction rationnelle**

$$\frac{1 - 3z}{1 - 4z + z^2}$$

**1, 1, 3, 11, 41, 153, 571, 2131, 7953, 29681, 110771, 413403, 1542841, 5757961, 21489003, 80198051, 299303201, 1117014753, 4168755811, 15558008491**

---

**Réf.** TI68 126 (divided by 2).

**HIS2 A1840** **Approximants de Padé**

**HIS1 N0233** **Fraction rationnelle**

$$\frac{1}{(z^2 + z + 1)(1 - z)^3}$$

**1, 2, 3, 5, 7, 9, 12, 15, 18, 22, 26, 30, 35, 40, 45, 51, 57, 63, 70, 77, 84, 92, 100, 108, 117, 126, 135, 145, 155, 165, 176, 187, 198, 210, 222, 234, 247, 260, 273, 287, 301**

## Related to Zarankiewicz's problem

**Réf.** TI68 126.

**HIS2 A1841** **Approximants de Padé** **Conjecture**

**HIS1 N0977** **Fraction rationnelle**

$$\frac{2z^4 + z^5 + 2z^3 + 2z^2 + 2z + 3}{(1 - z + z^2)(z^2 + z + 1)(1 + z)^2(1 - z)^3}$$

**3, 5, 10, 14, 21, 26, 36, 43, 55, 64, 78, 88, 105, 117, 136, 150, 171, 186, 210, 227, 253, 272, 300, 320, 351, 373, 406, 430, 465, 490, 528, 555, 595, 624, 666, 696, 741**

## Centered square numbers

**Réf.** MMAG 35 162 62. SIAR 12 277 70. INOC 24 4550 85.

**HIS2 A1844** **Approximants de Padé**

**HIS1 N1567** **Fraction rationnelle**

$$\frac{(1 + z)^2}{(1 - z)^3}$$

**1, 5, 13, 25, 41, 61, 85, 113, 145, 181, 221, 265, 313, 365, 421, 481, 545, 613, 685, 761, 841, 925, 1013, 1105, 1201, 1301, 1405, 1513, 1625, 1741, 1861, 1985, 2113, 2245**

**Réf.** SIAR 12 277 70. C1 81.
**HIS2 A1845** **Approximants de Padé**
**HIS1 N1844** **Fraction rationnelle**

$$\frac{(1 + z)^{3}}{(z - 1)^{4}}$$

**1, 7, 25, 63, 129, 231, 377, 575, 833, 1159, 1561, 2047, 2625, 3303, 4089, 4991, 6017, 7175, 8473, 9919, 11521, 13287, 15225, 17343, 19649, 22151, 24857, 27775**

---

**Réf.** SIAR 12 277 70. C1 81.
**HIS2 A1846** **Approximants de Padé**
**HIS1 N1974** **Fraction rationnelle**

$$\frac{(1 + z)^{4}}{(z - 1)^{5}}$$

**1, 9, 41, 129, 321, 681, 1289, 2241, 3649, 5641, 8361, 11969, 16641, 22569, 29961, 39041, 50049, 63241, 78889, 97281, 118721, 143529, 172041, 204609, 241601**

**Réf.** SIAR 12 277 70. C1 81.
**HIS2 A1847** **Approximants de Padé**
**HIS1 N2045** **Fraction rationnelle**

$$\frac{(1+z)^5}{(z-1)^6}$$

**1, 11, 61, 231, 681, 1683, 3653, 7183, 13073, 22363, 36365, 56695, 85305, 124515, 177045, 246047, 335137, 448427, 590557, 766727, 982729, 1244979, 1560549**

---

**Réf.** SIAR 12 277 70. C1 81.
**HIS2 A1848** **Approximants de Padé**
**HIS1 N2102** **Fraction rationnelle**

$$\frac{(1+z)^6}{(z-1)^7}$$

**1, 13, 85, 377, 1289, 3653, 8989, 19825, 40081, 75517, 134245, 227305, 369305, 579125, 880685, 1303777, 1884961, 2668525, 3707509, 5064793, 6814249**

**Réf.** SIAR 12 277 70. C1 81.
**HIS2 A1849** **Approximants de Padé**
**HIS1 N2139** **Fraction rationnelle**

$$\frac{(1 + z)^{7}}{(z - 1)^{8}}$$

**1, 15, 113, 575, 2241, 7183, 19825, 48639, 108545, 224143, 433905, 795455, 1392065, 2340495, 3800305, 5984767, 9173505, 13726991, 20103025, 28875327**

---

**Réf.** SIAR 12 277 70.
**HIS2 A1850** **Dérivée logarithmique**
**HIS1 N1184** **algébrique**
$\sum C(n,k).C(n+k,k)$, $k=0...n$

$$\frac{1}{(1 - 6 z + z^{2})^{1/2}}$$

**1, 3, 13, 63, 321, 1683, 8989, 48639, 265729, 1462563, 8097453, 45046719, 251595969, 1409933619, 7923848253, 44642381823, 252055236609, 1425834724419**

## Series-reduced planted trees with n nodes, n-3 endpoints

**Réf.** jr.

**HIS2 A1859** **Approximants de Padé**

**HIS1 N0531** **Fraction rationnelle**

$$\frac{1 + z^2 + 2z^3 - z^4}{(1 + z)(1 - z)^3}$$

1, 2, 5, 10, 16, 24, 33, 44, 56, 70, 85, 102, 120, 140, 161, 184, 208, 234, 261, 290, 320, 352, 385, 420, 456, 494, 533, 574, 616, 660, 705, 752, 800, 850, 901, 954, 1008, 1064, 1121, 1180

## Series-reduced planted trees with n nodes, n-4 endpoints

**Réf.** jr.

**HIS2 A1860** **Approximants de Padé**

**HIS1 N1171** **Fraction rationnelle**

$$\frac{3 + 3z + 2z^2}{(z^2 + z + 1)(z - 1)^4}$$

3, 12, 29, 57, 99, 157, 234, 333, 456, 606, 786, 998, 1245

## Values of Bell polynomials

**Réf.** jr. PSPM 19 173 71.

**HIS2 A1861** **équations différentielles** **Formule de B. Salvy**

**HIS1 N0653** **exponentielle**

$$\exp(2 \exp(z) - 2)$$

**1, 2, 6, 22, 94, 454, 2430, 14214, 89918, 610182, 4412798**

---

## Convolved Fibonacci numbers

**Réf.** RCI 101. FQ 15 118 77.

**HIS2 A1872** **Dérivée logarithmique**

**HIS1 N1413** **Fraction rationnelle**

$$\frac{1}{(1 - z - z^2)^4}$$

**1, 4, 14, 40, 105, 256, 594, 1324, 2860, 6020, 12402, 25088**

## Convolved Fibonacci numbers

**Réf.** RCI 101. FQ 15 118 77. DM 26 267 79.

**HIS2 A1873** **Dérivée logarithmique**

**HIS1 N1600** **Fraction rationnelle**

$$\frac{1}{(1 - z - z^2)^5}$$

**1, 5, 20, 65, 190, 511, 1295, 3130, 7285, 16435, 36122, 77645, 163730, 339535**

## Convolved Fibonacci numbers

**Réf.** RCI 101.

**HIS2 A1874** **Dérivée logarithmique** **erreurs dans la suite**

**HIS1 N1738** **Fraction rationnelle** **corrigées par la formule**

$$\frac{1}{(1 - z - z^2)^6}$$

**1, 6, 27, 98, 315, 924, 2534, 6588, 16407, 39430, 91959, 209034, 464723, 1013292, 2171850, 4584620, 9546570, 19635840, 39940460, 80421600, 160437690, 317354740, 622844730, 1213580820**

## Convolved Fibonacci numbers

**Réf.** RCI 101. DM 26 267 79.

**HIS2 A1875** **Dérivée logarithmique**

**HIS1 N1865** **Fraction rationnelle**

$$\frac{1}{(1 - z - z^2)^7}$$

**1, 7, 35, 140, 490, 1554, 4578, 12720, 33705, 85855, 211519**

---

**Réf.** RCI 77.

**HIS2 A1879** **Hypergéométrique** **Suite P-récurrente**

**HIS1 N1775** **algébrique** **f.g. exponentielle**

**a(n) = (2 n + 2) a(n-1) + (-2 n + 3) a(n-2)**

$$\frac{z}{(1 - 2 z)^{3/2}}$$

**1, 6, 45, 420, 4725, 62370, 945945, 16216200, 310134825, 6547290750, 151242416325, 3794809718700, 102776096548125, 2988412653476250, 92854250304440625**

## Coefficients of Bessel polynomials $y_n(x)$

**Réf.** RCI 77.

**HIS2 A1880** **Tableaux généralisés** **f.g. exponentielle**

**HIS1 N2146** **algébrique**

$$\frac{z\,(2+z)}{2\,(1-2z)^{7/2}}$$

**1, 15, 210, 3150, 51975, 945945, 18918900**

## Coefficients of Bessel polynomials $y_n(x)$

**Réf.** RCI 77.

**HIS2 A1881** **Tableaux généralisés** **f.g. exponentielle**

**HIS1 N2217** **algébrique**

$$\frac{z\,(2+3z)}{2\,(1-2z)^{9/2}}$$

**1, 21, 378, 6930, 135135, 2837835**

**Réf.** AMM 72 1024 65.

**HIS2 A1882** **Approximants de Padé**

**HIS1 N0273** **Fraction rationnelle**

$$\frac{2 + 3z - 3z^2 - z^3}{1 - 4z^2 + 2z^4}$$

**2, 3, 5, 11, 16, 38, 54, 130, 184, 444, 628, 1516, 2144, 5176, 7320, 17672, 24992, 60336, 85328, 206000, 291328, 703328, 994656, 2401312, 3395968, 8198592**

---

## Hit polynomials

**Réf.** RI63.

**HIS2 A1891** **Approximants de Padé**

**HIS1 N1365** **Fraction rationnelle**

$$\frac{z(1 + z)}{(1 - z - z^2)(z - 1)^2}$$

**0, 1, 4, 10, 21, 40, 72, 125, 212**

## Bisection of Fibonacci sequence

**Réf.** IDM 22 23 15. PLMS 21 729 70. FQ 9 283 71.

**HIS2 A1906** **Approximants de Padé**

**HIS1 N1101** **Fraction rationnelle**

$$\frac{1}{1 - 3z + z^2}$$

**1, 3, 8, 21, 55, 144, 377, 987, 2584, 6765, 17711, 46368, 121393, 317811, 832040, 2178309, 5702887, 14930352, 39088169, 102334155, 267914296, 701408733**

---

## Permutations with no cycles of length 4

**Réf.** R1 83.

**HIS2 A1907** **Dérivée logarithmique**

**HIS1 N1261** **exponentielle**

**a(n) = (4 n - 5) a(n-1) + (4 n - 8) a(n-2)**

$$\frac{1}{(1 - 4z)\exp(z)}$$

**1, 3, 25, 299, 4785, 95699, 2296777, 64309755, 2057912161, 74084837795, 2963393511801, 130389314519243, 6258687096923665, 325451729040030579**

**Réf.** R1 83.
**HIS2 A1908** **Dérivée logarithmique** **Suite P-récurrente**
**HIS1 N1500** **exponentielle**
**a(n) = (5 n - 6) a(n-1) + (5 n - 10) a(n-2)**

$$\frac{1}{(1 - 5z)\exp(z)}$$

**1, 4, 41, 614, 12281, 307024, 9210721, 322375234, 12895009361, 580275421244, 29013771062201, 1595757408421054, 95745444505263241**

---

**Réf.** R1 188.
**HIS2 A1909** **Dérivée logarithmique**
**HIS1 N1450** **exponentielle**
**a(n) = (n + 2) a(n-1) + (n - 2) a(n-2)**

$$\frac{1}{(1 - z)^5 \exp(z)}$$

**0, 1, 4, 21, 134, 1001, 8544, 81901, 870274, 10146321, 128718044, 1764651461, 25992300894, 409295679481, 6860638482424, 121951698034461**

**Réf.** R1 188.
**HIS2 A1910** **Dérivée logarithmique**
**HIS1 N1637** **exponentielle**
**a(n) = (n + 3) a(n-1) + (n-2) a(n-2)**

$$\frac{1}{(1 - z)^{6} \exp(z)}$$

**0, 1, 5, 31, 227, 1909, 18089, 190435, 2203319, 27772873, 378673901, 5551390471, 87057596075, 1453986832381, 25762467303377, 482626240281739**

---

**Réf.** R1 233. LNM 748 151 79.
**HIS2 A1911** **Approximants de Padé**
**HIS1 N1007** **Fraction rationnelle**

$$\frac{1 + z}{(1 - z)\,(1 - z - z^{2})}$$

**1, 3, 6, 11, 19, 32, 53, 87, 142, 231, 375, 608, 985, 1595, 2582, 4179, 6763, 10944, 17709, 28655, 46366, 75023, 121391, 196416, 317809, 514227, 832038, 1346267**

## Quadrinomial coefficients

**Réf.** JCT 1 372 66. C1 78.

**HIS2 A1919** **Approximants de Padé**

**HIS1 N1769** **Fraction rationnelle**

$$\frac{3 z^{2} - 8 z + 6}{(z - 1)^{8}}$$

**6, 40, 155, 456, 1128, 2472, 4950, 9240, 16302, 27456, 44473, 69680, 106080, 157488, 228684, 325584, 455430, 627000, 850839, 1139512, 1507880, 1973400, 2556450, 3280680**

---

**Réf.** AMM 53 465 46.

**HIS2 A1921** **Approximants de Padé**

**HIS1 N1885** **Fraction rationnelle**

$$\frac{z (z - 7)}{(z - 1) (1 - 14 z + z^{2})}$$

**0, 7, 104, 1455, 20272, 282359, 3932760, 54776287, 762935264, 10626317415, 148005508552, 2061450802319, 28712305723920, 399910829332567**

**Réf.** AMM 53 465 46.

**HIS2 A1922** **Approximants de Padé**

**HIS1 N1946** **Fraction rationnelle**

$$\frac{7z - 1}{(z - 1)(1 - 14z + z^2)}$$

**1, 8, 105, 1456, 20273, 282360, 3932761, 54776288, 762935265, 10626317416, 148005508553, 2061450802320, 28712305723921, 399910829332568**

---

## From rook polynomials

**Réf.** SMA 20 18 54.

**HIS2 A1924** **Approximants de Padé**

**HIS1 N1053** **Fraction rationnelle**

$$\frac{1}{(1 - z - z^2)(z - 1)^2}$$

**1, 3, 7, 14, 26, 46, 79, 133, 221, 364, 596, 972, 1581, 2567, 4163, 6746, 10926, 17690, 28635, 46345, 75001, 121368, 196392, 317784, 514201, 832011, 1346239**

## From rook polynomials

**Réf.** SMA 20 18 54.

**HIS2 A1925** **Approximants de Padé**

**HIS1 N1724** **Fraction rationnelle**

$$\frac{1 + z}{(1 - z - z^2)^2 (z - 1)^3}$$

**1, 6, 22, 64, 162, 374, 809, 1668, 3316, 6408, 12108, 22468, 41081, 74202, 132666, 235160, 413790, 723530, 1258225, 2177640, 3753096, 6444336, 11028792**

## From rook polynomials

**Réf.** SMA 20 18 54.

**HIS2 A1926** **Approximants de Padé**

**HIS1 N1978** **Fraction rationnelle**

$$\frac{(1 + z)^2}{(1 - z - z^2)^3 (z - 1)^4}$$

**1, 9, 46, 177, 571, 1632, 4270, 10446, 24244, 53942, 115954, 242240, 494087, 987503, 1939634, 3753007, 7167461, 13532608, 25293964, 46856332, 86110792**

## Sum of Fibonacci and Pell numbers

**Réf.**

**HIS2 A1932** **Approximants de Padé**

**HIS1 N0319** **Fraction rationnelle**

$$\frac{(2 + z)(1 - 2z)}{(1 - z - z^2)(1 - 2z - z^2)}$$

**2, 3, 7, 15, 34, 78, 182, 429, 1019, 2433, 5830, 14004, 33694, 81159, 195635, 471819, 1138286, 2746794, 6629290, 16001193, 38624911, 93240069, 225087338**

---

## Coefficients of an elliptic function

**Réf.** CAY 9 128.

**HIS2 A1934** **Euler**

**HIS1 N1397** **Produit infini**

$$\prod_{n\geq 1} \frac{1}{(1 - z^n)^{c(n)}}$$

$$c(n) = 4,2,4,2,4,2,4,2,4,2,...$$

**1, 4, 12, 32, 76, 168, 352, 704**

## Coefficients of an elliptic function

**Réf.** CAY 9 128.

**HIS2 A1935** **Euler**

**HIS1 N0204** **Produit infini**

$$\prod_{n\geq 1} \frac{1}{(1 - z^n)^{c(n)}}$$

**c(n) = 1,2,3 mod 4**

**1, 1, 2, 3, 4, 6, 9, 12, 16, 22, 29, 38, 50, 64, 82, 105, 132, 166, 208, 258, 320, 395, 484, 592, 722, 876, 1060**

## Coefficients of an elliptic function

**Réf.** CAY 9 128. MOC 29 852 75.

**HIS2 A1936** **Euler**

**HIS1 N0532** **Produit infini**

$$\prod_{n\geq 1} \frac{1}{(1 - z^n)^{c(n)}}$$

**c(n) = 2,2,2,0,2,2,2,0,...**

**1, 2, 5, 10, 18, 32, 55, 90, 144, 226, 346, 522, 777, 1138, 1648, 2362, 3348, 4704, 6554, 9056, 12425, 16932, 22922, 30848, 41282, 54946, 72768, 95914, 125842, 164402**

## Coefficients of an elliptic function

**Réf.** CAY 9 128.

**HIS2 A1937** **Euler** **erreurs dans la suite corrigées avec la formule.**
**HIS1 N1120** **Produit infini**

$$\prod_{n\geq 1} \frac{1}{(1 - z^n)^{c(n)}}$$

$$c(n) = 3,3,3,0,3,3,3,0,...$$

**1, 3, 9, 22, 48, 99, 194, 363, 657, 1155, 1977, 3312, 5443, 8787, 13968, 21894, 33873, 51795, 78345, 117412, 174033, 255945**

## Coefficients of an elliptic function

**Réf.** CAY 9 128.

**HIS2 A1938** **Euler**
**HIS1 N1412** **Produit infini**

$$\prod_{n\geq 1} \frac{1}{(1 - z^n)^{c(n)}}$$

$$c(n) = 4,4,4,0,4,4,4,0,...$$

**1, 4, 14, 40, 101, 236, 518, 1080, 2162, 4180, 7840, 14328, 25591, 44776, 76918, 129952, 216240, 354864, 574958**

## Coefficients of an elliptic function

**Réf.** CAY 9 128.

**HIS2 A1939** **Euler**

**HIS1 N1599** **Produit infini**

$$\prod_{n\geq 1} \frac{1}{(1-z^n)^{c(n)}}$$

$$c(n) = 5,5,5,0,5,5,5,0,...$$

**1, 5, 20, 65, 185, 481, 1165, 2665, 5820, 12220, 24802, 48880, 93865, 176125, 323685, 583798, 1035060, 1806600, 3108085**

---

## Coefficients of an elliptic function

**Réf.** CAY 9 128.

**HIS2 A1940** **Euler** **erreurs dans la suite corrigées avec**

**HIS1 N1737** **Produit infini** **la formule.**

$$\prod_{n\geq 1} \frac{1}{(1-z^n)^{c(n)}}$$

$$c(n) = 6,6,6,0,6,6,6,0,...$$

**1, 6, 27, 98, 309, 882, 2330, 5784, 13644, 30826, 67107, 141444, 289746, 578646, 1129527, 2159774, 4052721, 7474806, 15063859**

## Coefficients of an elliptic function

**Réf.** CAY 9 128.

**HIS2 A1941** **Euler**

**HIS1 N1864** **Produit infini**

$$\prod_{n\geq 1} \frac{1}{(1 - z^n)^{c(n)}}$$

$$c(n) = 7,7,7,0,7,7,7,0,...$$

**1, 7, 35, 140, 483, 1498, 4277, 11425, 28889, 69734, 161735, 362271, 786877, 1662927, 3428770, 6913760, 13660346, 26492361, 50504755**

---

**Réf.** JLMS 8 166 33.

**HIS2 A1945** **Approximants de Padé**

**HIS1 N1525** **Fraction rationnelle**

$$\frac{z\,(1 + 2\,z + z^2 + 2\,z^3 + z^4)}{(z^3 - z - 1)\,(-\,1 + z^2 + z^3)}$$

**0, 1, 1, 1, 5, 1, 7, 8, 5, 19, 11, 23, 35, 27, 64, 61, 85, 137, 133, 229, 275, 344, 529, 599, 875, 1151, 1431, 2071, 2560, 3481, 4697, 5953, 8245, 10649, 14111, 19048, 24605**

**Réf.** RCI 139.
**HIS2 A1946** **Approximants de Padé**
**HIS1 N0794** **Fraction rationnelle**

$$\frac{11z - 2}{z^2 + 11z - 1}$$

**2, 11, 123, 1364, 15127, 167761, 1860498, 20633239, 228826127, 2537720636, 28143753123, 312119004989, 3461452808002, 38388099893011**

**Related to Bernoulli numbers**

**Réf.** RCI 141.
**HIS2 A1947** **Approximants de Padé**
**HIS1 N1265** **Fraction rationnelle**

$$\frac{4z - 3}{z^2 + 11z - 1}$$

**3, 29, 322, 3571, 39603, 439204, 4870847, 54018521, 599074578, 6643838879, 73681302247, 817138163596, 9062201101803, 100501350283429**

## A probability difference equation

**Réf.** AMM 32 369 25.

**HIS2 A1949** **Approximants de Padé**

**HIS1 N0430** **Fraction rationnelle**

$$\frac{1}{(1-z)(1-z-z^2-z^3-z^4-z^5)}$$

**1, 2, 4, 8, 16, 32, 63, 124, 244, 480, 944, 1856, 3649, 7174, 14104, 27728, 54512, 107168, 210687, 414200, 814296, 1600864, 3147216, 6187264, 12163841**

## Restricted partitions

**Réf.** CAY 2 277.

**HIS2 A1971** **Approximants de Padé**

**HIS1 N0227** **Fraction rationnelle**

$$\frac{1-z^6}{(1-z)(1-z^2)(1-z^3)(1-z^4)}$$

**1, 1, 2, 3, 5, 6, 8, 10, 13, 15, 18, 21, 25, 28, 32, 36, 41, 45, 50**

## Restricted partitions

**Réf.** CAY 2 277.

**HIS2 A1972** **Approximants de Padé**

**HIS1 N0199** **Fraction rationnelle**

$$\frac{2 - z + z^{3} - 2 z^{4} + z^{5}}{(1 + z)(1 + z^{2})(z - 1)^{3}}$$

**2, 3, 4, 6, 8, 10, 12, 15, 18, 21, 24, 28, 32, 36, 40, 45, 50**

---

**Réf.** CAY 2 278.

**HIS2 A1973** **Approximants de Padé**

**HIS1 N0969** **Fraction rationnelle**

$$\frac{1 - z + z^{2}}{(1 + z)(z^{2} + z + 1)(z - 1)^{4}}$$

**1, 1, 3, 5, 8, 12, 18, 24, 33, 43, 55, 69, 86, 104, 126, 150, 177, 207, 241, 277, 318, 362, 410, 462, 519, 579, 645, 715, 790, 870, 956, 1046, 1143, 1245, 1353, 1467, 1588, 1714, 1848, 1988**

## Expansion of a generating function

**Réf.** CAY 10 414.

**HIS2 A1993** **Euler**

**HIS1 N0973** **Fraction rationnelle**

$$\frac{1}{(1-z)(1-z^2)^2(1-z^3)^2(1-z^4)}$$

**1, 1, 3, 5, 9, 13, 22, 30, 45, 61, 85, 111**

## Expansion of a generating function

**Réf.** CAY 10 415.

**HIS2 A1994** **Euler**

**HIS1 N0927** **Fraction rationnelle**

$$\frac{1}{(1-z)(1-z^2)^2(1-z^3)(1-z^4)(1-z^5)}$$

**1, 1, 3, 4, 8, 11, 18, 24, 36, 47, 66, 84, 113, 141, 183, 225, 284, 344, 425, 508, 617, 729, 872, 1020, 1205, 1397, 1632, 1877, 2172, 2480, 2846, 3228, 3677**

## Expansion of a generating function

**Réf.** CAY 10 415.

**HIS2 A1996** **Euler**

**HIS1 N0112** **Fraction rationnelle**

$$\frac{1}{(1-z^2)(1-z^3)(1-z^4)(1-z^5)(1-z^6)(1-z^7)}$$

**1, 0, 1, 1, 2, 2, 4, 4, 6, 7, 10, 11, 16, 17, 23, 26, 33, 37, 47, 52, 64, 72, 86, 96, 115, 127, 149, 166, 192, 212, 245, 269, 307, 338, 382, 419, 472, 515, 576, 629, 699, 760, 843, 913**

## Folding a piece of wire of length n

**Réf.** AMM 44 51 37. GMJ 15 146 74.

**HIS2 A1998** **Approximants de Padé**

**HIS1 N0468** **Fraction rationnelle**

$$\frac{3z^4 - 8z^3 + 2z^2 + 3z - 1}{(z-1)(3z-1)(3z^2-1)}$$

**1, 1, 2, 4, 10, 25, 70, 196, 574, 1681, 5002, 14884, 44530, 133225, 399310, 1196836, 3589414, 10764961, 32291602, 96864964, 290585050, 871725625, 2615147350**

**Réf.** AMM 43 29 36.
**HIS2 A2002** **LLL** **suite P-récurrente**
**HIS1 N1621** **algébrique**
**n a(n) = (7 n - 5) a(n - 1) + (- 7 n + 16) a(n - 2) + (n - 3) a(n - 3)**
**a(n) = ∑C(n,k+1).C(n+k,k), k=0..n-1**

$$\frac{z + (1 - 6 z + z^2)^{1/2} - 1}{- 2 (1 - 6 z + z^2)^{1/2} z}$$

**1, 5, 25, 129, 681, 3653, 19825, 108545, 598417, 3317445, 18474633, 103274625, 579168825, 3256957317**

---

**Réf.** AMM 43 29 36.
**HIS2 A2003** **LLL** **Suite P-récurrente**
**HIS1 N0735** **algébrique**
**n a(n) = (5 n - 1) a(n - 1) + (5 n - 14) a(n - 2) + (- n + 3) a(n - 3)**
**a(n) = 2 ∑ C(n-1,k) C(n+k,k) , k = 0 ..n-1**

$$\frac{z + 1 + (1 - 6 z + z^2)^{1/2}}{- 2 (1 - 6 z + z^2)^{1/2} z}$$

**2, 8, 38, 192, 1002, 5336, 28814, 157184, 864146, 4780008, 26572086, 148321344, 830764794, 4666890936**

## Almost trivalent maps

**Réf.** PLC 1 292 70.

**HIS2 A2011** **Hypergéométrique**

**HIS1 N1458** **algébrique**

$$\frac{4}{(1 - 4z)^{3/2}}$$

**4, 24, 120, 560, 2520, 11088, 48048**

## n appears n times

**Réf.** MMAG 38 186 65. KN1 1 43.

**HIS2 A2024** **Euler**

**HIS1 N0089** **Produit infini**

**a(n)=[(1+[√(8n-7)])/2]**

$$\prod_{n \geq 1} \frac{1}{(1 - z^n)^{c(n)}}$$

$$c(n) = 2,-1,1,-1,1,-1,1, \ ...$$

**1, 2, 2, 3, 3, 3, 4, 4, 4, 4, 5, 5, 5, 5, 5, 6, 6, 6, 6, 6, 6, 7, 7, 7, 7, 7, 7, 7, 8, 8, 8, 8, 8, 8, 8, 8, 9, 9, 9, 9, 9, 9, 9, 9, 9, 10, 10, 10, 10, 10, 10, 10, 10, 10, 10, 11, 11, 11, 11, 11, 11**

## Related to partitions

**Réf.** AMM 76 1036 69.

**HIS2 A2040** **Approximants de Padé**

**HIS1 N0442** **Fraction rationnelle**

$$\frac{1}{1 - 2z - 5z^4 - 7z^6}$$

**1, 2, 4, 8, 21, 52, 131, 316, 765, 1846, 4494**

---

**Réf.** AMM 3 244 1896.

**HIS2 A2041** **Approximants de Padé**

**HIS1 N1759** **Fraction rationnelle**

$$\frac{1}{(z - 1)(1 + 2z)(1 - 2z)(5z - 1)}$$

**1, 6, 35, 180, 921, 4626, 23215, 116160, 581141, 2906046, 14531595, 72659340, 363302161, 1816516266, 9082603175, 45413037720, 227065275981, 1135326467286**

## Simplices in barycentric subdivisions of n-simplex

**Réf.** SKA 11 95 28. MMAG 37 132 64.

**HIS2 A2050** **Recoupements**

**HIS1 N1622** **exponentielle**

$$\frac{\exp(z)\,(1-\exp(z))}{\exp(z)-2}$$

**1, 5, 25, 149, 1081, 9365, 94585, 1091669, 14174521, 204495125, 3245265145, 56183135189, 1053716696761, 21282685940885, 460566381955705**

## Binomial coefficients C(2n+1,n-1)

**Réf.** CAY 13 95. AS1 828.

**HIS2 A2054** **Hypergéométrique** **Suite P-récurrente**

**HIS1 N1607** **algébrique**

**2F1([2, 5/2], [4], 4 z)**

$$\frac{8z}{(1-4z)^{1/2}\,\left(1+(1-4z)^{1/2}\right)^3}$$

**1, 5, 21, 84, 330, 1287, 5005, 19448, 75582, 293930, 1144066, 4457400, 17383860, 67863915, 265182525, 1037158320, 4059928950, 15905368710**

## Dissections of a polygon by number of parts

**Réf.** CAY 13 95. AEQ 18 385 78.

**HIS2 A2055** **Hypergéométrique** **Suite P-récurrente**

**HIS1 N1982** **algébrique**

$$\frac{\left(z - (1 - 4z)^{1/2}\right) z}{\left(1 + (1 - 4z)^{1/2}\right)^4 (1 - 4z)^{3/2}}$$

**1, 9, 56, 300, 1485, 7007, 32032, 143208, 629850, 2735810, 11767536, 50220040, 212952285**

## Dissections of a polygon by number of parts

**Réf.** CAY 13 95. AEQ 18 385 78.

**HIS2 A2056** **Hypergéométrique** **simplifiée avec LLL**

**HIS1 N2115** **algébrique 2è degré**

$$\frac{1/2\,(1 - 21z + 180z^2 - 800z^3 + 1920z^4 - 2304z^5 + 1024z^6}{(z^5 (4z-1)^5)}$$

$$- \frac{- (- (10z^4 - 50z^3 + 40z^2 - 11z + 1)^2 (4z-1)^5)^{1/2})}{(z^5 (4z-1)^5)}$$

**1, 14, 120, 825, 5005, 28028, 148512, 755820, 3730650, 17978180, 84987760, 395482815**

## 4 C(2n+1,n-1)/(n+3)

**Réf.** CAY 13 95. FQ 14 397 76. DM 14 84 76.

**HIS2 A2057** **Hypergéométrique**

**HIS1 N1415** **algébrique**

**$_2F_1([2, 5/2], [5], 4z)$**

$$\frac{16 z}{\left(1 + (1 - 4 z)^{1/2}\right)^4}$$

**1, 4, 14, 48, 165, 572, 2002, 7072, 25194, 90440, 326876, 1188640, 4345965, 15967980, 58929450, 218349120, 811985790, 3029594040, 11338026180, 42550029600**

---

## Partitions of a polygon by number of parts

**Réf.** CAY 13 95.

**HIS2 A2059** **Hypergéométrique**

**HIS1 N1269** **algébrique**

$$\frac{\left(2 z - 3 (1 - 4 z)^{1/2}\right) z}{\left(1 + (1 - 4 z)^{1/2}\right)^6 (1 - 4 z)^{3/2}}$$

**3, 32, 225, 1320, 7007, 34944, 167076, 775200, 3517470, 15690048**

## Central polygonal numbers

**Réf.** HO50 22. HO70 87.

**HIS2 A2061** **Approximants de Padé**

**HIS1 N1049** **Fraction rationnelle**

$$\frac{1 - 2 z + 3 z^{2}}{(1 - z)^{3}}$$

**1, 1, 3, 7, 13, 21, 31, 43, 57, 73, 91, 111, 133, 157, 183, 211, 241, 273, 307, 343, 381, 421, 463, 507, 553, 601, 651, 703, 757, 813, 871, 931, 993, 1057, 1123, 1191, 1261**

## n'th Fibonacci number + n

**Réf.** HO70 96.

**HIS2 A2062** **Approximants de Padé**

**HIS1 N0240** **Fraction rationnelle**

$$\frac{z (3 z - 2)}{(1 - z - z^{2}) (1 - z)^{2}}$$

**0, 2, 3, 5, 7, 10, 14, 20, 29, 43, 65, 100, 156, 246, 391, 625, 1003, 1614, 2602, 4200, 6785, 10967, 17733, 28680, 46392, 75050, 121419, 196445, 317839, 514258**

## Cullen numbers

**Réf.** SI64a 346. UPNT B20.

**HIS2 A2064** **Approximants de Padé**

**HIS1 N1125** **Fraction rationnelle**

$$\frac{1 - 2 z + 2 z^{2}}{(1 - z)^{2} (2 z - 1)}$$

**1, 3, 9, 25, 65, 161, 385, 897, 2049, 4609, 10241, 22529, 49153, 106497, 229377, 491521, 1048577, 2228225, 4718593, 9961473, 20971521, 44040193, 92274689**

## First differences are periodic

**Réf.** TCPS 2 219 1827.

**HIS2 A2081** **Approximants de Padé**

**HIS1 N0426** **Fraction rationnelle**

$$\frac{2 (1 + 2 z^{2} + 2 z^{3})}{(1 + z^{2}) (z - 1)^{2}}$$

**2, 4, 8, 16, 22, 24, 28, 36, 42, 44, 48, 56, 62, 64, 68, 76, 82, 84, 88, 96, 102, 104, 108, 116, 122, 124, 128, 136, 142, 144, 148, 156, 162, 164, 168, 176, 182, 184, 188, 196, 202, 204, 208, 216**

## Partitions of n into non-prime parts

**Réf.** JNSM 9 91 69.

**HIS2 A2095** **Euler**

**HIS1 N0094** **Produit infini**

$$\prod_{n\geq 1} \frac{1}{(1 - z^{n})^{c(n)}}$$

**c(n) = Les nombres non-premiers**

**1, 1, 1, 1, 2, 2, 3, 3, 5, 6, 8, 8, 12, 13, 17, 19, 26, 28, 37, 40, 52, 58, 73, 79, 102, 113, 139, 154, 191, 210, 258, 284, 345, 384, 462, 509, 614, 679, 805, 893, 1060, 1171, 1382**

---

## Logarithmic numbers

**Réf.** MAS 31 78 63. CACM 13 726 70.

**HIS2 A2104** **équations différentielles** **Suite P-récurrente**

**HIS1 N1105** **exponentielle** **Formule de B. Salvy**

**a(n) = (n + 1) a(n-1) + (-2 n + 2) a(n-2) + (n - 2) a(n-3)**

$$- \exp(z) \ln(1 - z)$$

**1, 3, 8, 24, 89, 415, 2372, 16072, 125673, 1112083, 10976184, 119481296, 1421542641, 18348340127, 255323504932, 3809950977008, 60683990530225**

## The square of Euler's product

**Réf.** PLMS 21 190 1889.

**HIS2 A2107** **Recoupements**

**HIS1 N0028** **Produit infini**

$$\prod_{n\geq 1} \frac{1}{(1 - z^n)^{c(n)}}$$

**c(n) = -2,-2,-2,-2,-2-,2,...**

**1, 2, 1, 2, 1, 2, 2, 0, 2, 2, 1, 0, 0, 2, 3, 2, 2, 0, 0, 2, 2, 0, 0, 2, 1, 0, 2, 2, 2, 2, 1, 2, 0, 2, 2, 2, 2, 0, 2, 0, 4, 0, 0, 0, 1, 2, 0, 0, 2, 0, 2, 2, 1, 2, 0, 2, 2, 0, 0, 2, 0, 2, 0, 2, 2, 0, 4, 0, 0**

## Numerators of convergents to exp(1)

**Réf.** BAT 17 1871. MOC 2 69 46.

**HIS2 A2119** **équations différentielles** **formule de B. Salvy**

**HIS1 N1880** **exponentielle**

**a(n) = (4 n - 6) a(n - 1) + a(n - 2)**

$$\frac{\exp(1/2\,(1 - 4z)^{1/2} - 1/2)}{(1 - 4z)^{1/2}}$$

**1, 1, 7, 71, 1001, 18089, 398959, 10391023, 312129649, 10622799089, 403978495031, 16977719590391, 781379079653017, 39085931702241241**

## From symmetric functions

**Réf.** PLMS 23 314 23.

**HIS2 A2124** **Approximants de Padé**

**HIS1 N0062** **Fraction rationnelle**

$$\frac{1 - z^{6}}{1 - z^{3} - z^{5} - z^{6} - z^{7} + z^{9}}$$

**1, 0, 0, 1, 0, 1, 1, 1, 2, 1, 3, 4, 3, 7, 7, 8, 14, 15, 21, 28, 33, 47, 58, 76, 103, 125, 169, 220, 277, 373**

## From symmetric functions

**Réf.** PLMS 23 315 23.

**HIS2 A2125** **Approximants de Padé**

**HIS1 N0006** **Fraction rationnelle**

$$\frac{(1 - z^{6})^{2}}{(1 - z^{3} - z^{5} - z^{6} - z^{7} + z^{9})^{2}}$$

**1, 0, 0, 2, 0, 2, 3, 2, 6, 4, 9, 14, 11, 26, 29, 34, 62, 68, 99, 140, 169, 252, 322, 430, 607, 764, 1059, 1424, 1845, 2546**

**Réf.** CAY 9 190. PLMS 17 29 17. EMN 34 1 44. AMM 79 519 72.
**HIS2 A2135** **Dérivée logarithmique**
**HIS1 N0594** **exponentielle**
**a(n) = (n - 1) a(n - 1) + (- 1/2 n^2 + 5/2 n - 3) a(n - 3)**

$$\frac{\exp(1/4\, z\, (z + 2))}{(1 - z)^{1/2}}$$

**1, 1, 2, 5, 17, 73, 388, 2461, 18155, 152531, 1436714, 14986879, 171453343, 2134070335, 28708008128, 415017867707, 6416208498137, 105630583492969**

---

**Matrices with 2 rows**

**Réf.** PLMS 17 29 17.
**HIS2 A2136** **Dérivée logarithmique** **Suite P-récurrente**
**HIS1 N0656** **exponentielle**
**a(n) = n a(n - 1) + (- 1/2 n^2 + 5/2 n - 3) a(n - 3)**

$$\frac{\exp(1/4\, z\, (z + 2))}{(1 - z)^{3/2}}$$

**1, 2, 6, 23, 109, 618, 4096, 31133, 267219, 2557502**

## Pell numbers

**Réf.** AJM 1 187 1878. FQ 4 373 66. RI89 43.

**HIS2 A2203** **Approximants de Padé**

**HIS1 N0136** **Fraction rationnelle**

$$\frac{2\,(1 - z)}{1 - 2\,z - z^{2}}$$

**2, 2, 6, 14, 34, 82, 198, 478, 1154, 2786, 6726, 16238, 39202, 94642, 228486, 551614, 1331714, 3215042, 7761798, 18738638, 45239074, 109216786, 263672646**

## Restricted hexagonal polyominoes with n cells

**Réf.** EMS 17 11 70. rcr.

**HIS2 A2212** **Inverse fonctionnel** **Suite P-récurrente**

**HIS1 N1145** **algébrique**

**(n + 1) a(n) = (6 n - 3) a(n - 1) + (- 5 n + 10) a(n - 2)**

$$\frac{-1 + 3\,z + (1 - 6\,z + 5\,z^{2})^{1/2}}{2\,z}$$

**1, 3, 10, 36, 137, 543, 2219, 9285, 39587, 171369, 751236, 3328218, 14878455, 67030785, 304036170, 1387247580, 6363044315, 29323149825, 135700543190**

## Dissections of a polygon

**Réf.** DM 11 388 75.

**HIS2 A2293** **Inverse fonctionnel algébrique** **Suite P-récurrente**

**HIS1 N1454**

**1/9 (n - 1) (3 n - 4) (3 n - 2) a(n) = 8/27 (4 n - 5) (4 n - 7) (2 n - 3) a(n - 1)**

$$ {}_4F_3([1, 3/2, 5/4, 7/4], [2, 5/3, 7/3], 256 z / 27) $$

**1, 1, 4, 22, 140, 969, 7084, 53820, 420732, 3362260, 27343888, 225568798, 1882933364, 15875338990, 134993766600, 1156393243320, 9969937491420**

---

## C(5n,n)/(4n+1)

**Réf.** DM 11 388 75.

**HIS2 A2294** **Hypergéométrique algébrique** **Suite P-récurrente**

**HIS1 N1646**

**1/32 (4 n - 5) (n - 1) (4 n - 3) (2 n - 3) a(n) = 5/256 (5 n - 9) (5 n - 8) (5 n - 7) (5 n - 6) a(n - 1)**

$$ {}_5F_4([1, 9/5, 7/5, 8/5, 6/5], [2, 3/2, 9/4, 7/4], 3125 z / 256) $$

**1, 1, 5, 35, 285, 2530, 23751, 231880, 2330445, 23950355, 250543370, 2658968130, 28558343775, 309831575760, 3390416787880, 37377257159280, 414741863546285**

## Dissections of a polygon

**Réf.** DM 11 388 75.

**HIS2 A2295** **Hypergéométrique** **Suite P-récurrente**

**HIS1 N1780** **algébrique**

**1/625 (n - 1) (5 n - 4) (5 n - 8) (5 n - 7) (5 n - 6) a(n) =**
**72 / 3125 (3 n - 5) (6 n - 11) (6 n - 7) (3 n - 4) (2 n - 3) a(n - 1)**

**6F5([1, 3/2, 5/3, 4/3, 7/6, 11/6],**

**[2, 11/5, 9/5, 7/5, 8/5],46656 z / 3125)**

**1, 1, 6, 51, 506, 5481, 62832, 749398, 9203634, 115607310, 1478314266, 19180049928, 251857119696, 3340843549855, 44700485049720, 602574657427116**

## Dissections of a polygon

**Réf.** DM 11 389 75.

**HIS2 A2296** **Hypergéométrique** **Suite P-récurrente**

**HIS1 N1878** **algébrique**

**1/648 (n - 1) (6 n - 7) (3 n - 4) (2 n - 3) (3 n - 5) (6 n - 5) a(n) =**
**7 / 46656 (7 n - 11) (7 n - 10) (7 n - 13) (7 n - 9) (7 n - 12) (7 n - 8) a(n - 1)**

**7F6([1,8/7,9/7,11/7,10/7,13/7,12/7],**

**[2,3/2,5/3,13/6,4/3,11/6],823543z/46656)**

**1, 1, 7, 70, 819, 10472, 141778, 1997688, 28989675, 430321633, 6503352856, 99726673130, 1547847846090, 24269405074740, 383846168712104**

**Réf.** TOH 42 152 36.
**HIS2 A2301** **Dérivée logarithmique** **f.g. exponentielle**
**HIS1 N0737** **Fraction rationnelle**

$$\frac{2}{(z - 1)^{4}}$$

**2, 8, 40, 240, 1680, 13440, 120960, 1209600, 13305600, 159667200, 2075673600, 29059430400, 435891456000, 6974263296000, 118562476032000**

---

## Sums of fourth powers of odd numbers

**Réf.** AMS 2 358 31 (divided by 2). CC55 742.
**HIS2 A2309** **Approximants de Padé**
**HIS1 N2327** **Fraction rationnelle**

$$\frac{1 + 76 z + 230 z^{2} + 76 z^{3} + z^{4}}{(z - 1)^{6}}$$

**1, 82, 707, 3108, 9669, 24310, 52871, 103496, 187017, 317338, 511819, 791660, 1182285, 1713726, 2421007, 3344528, 4530449, 6031074, 7905235, 10218676**

## NSW numbers

**Réf.** AMM 4 25 1897. IDM 10 236 03. ANN 36 644 35. RI89 288.

**HIS2 A2315** **Approximants de Padé**

**HIS1 N1869** **Fraction rationnelle**

$$\frac{1 + z}{z^2 - 6z + 1}$$

**1, 7, 41, 239, 1393, 8119, 47321, 275807, 1607521, 9369319, 54608393, 318281039, 1855077841, 10812186007, 63018038201, 367296043199, 2140758220993**

## The pronic numbers

**Réf.** D1 2 232.

**HIS2 A2378** **Approximants de Padé**

**HIS1 N0616** **Fraction rationnelle**

$$\frac{2z}{(1 - z)^3}$$

**0, 2, 6, 12, 20, 30, 42, 56, 72, 90, 110, 132, 156, 182, 210, 240, 272, 306, 342, 380, 420, 462, 506, 552, 600, 650, 702, 756, 812, 870, 930, 992, 1056, 1122, 1190, 1260**

**Réf.** MFM 74 62 70 (divided by 5).

**HIS2 A2409** **Approximants de Padé**

**HIS1 N1668** **Fraction rationnelle**

$$\frac{1}{(1 - 2z)^7}$$

**1, 14, 112, 672, 3360, 14784, 59136, 219648, 768768, 2562560, 8200192, 25346048, 76038144, 222265344, 635043840, 1778122752, 4889837568, 13231325184, 35283533824**

---

## Pentagonal pyramidal numbers

**Réf.** D1 2 2. B1 194.

**HIS2 A2411** **Approximants de Padé**

**HIS1 N1709** **Fraction rationnelle**

$$\frac{1 + 2z}{(z - 1)^4}$$

**1, 6, 18, 40, 75, 126, 196, 288, 405, 550, 726, 936, 1183, 1470, 1800, 2176, 2601, 3078, 3610, 4200, 4851, 5566, 6348, 7200, 8125, 9126, 10206, 11368, 12615, 13950**

## Hexagonal pyramidal numbers

**Réf.** D1 2 2. B1 194.

**HIS2 A2412** **Approximants de Padé**

**HIS1 N1839** **Fraction rationnelle**

$$\frac{1 + 3z}{(z - 1)^4}$$

**1, 7, 22, 50, 95, 161, 252, 372, 525, 715, 946, 1222, 1547, 1925, 2360, 2856, 3417, 4047, 4750, 5530, 6391, 7337, 8372, 9500, 10725, 12051, 13482, 15022, 16675, 18445**

## Heptagonal pyramidal numbers

**Réf.** D1 2 2. B1 194.

**HIS2 A2413** **Approximants de Padé**

**HIS1 N1904** **Fraction rationnelle**

$$\frac{1 + 4z}{(z - 1)^4}$$

**1, 8, 26, 60, 115, 196, 308, 456, 645, 880, 1166, 1508, 1911, 2380, 2920, 3536, 4233, 5016, 5890, 6860, 7931, 9108, 10396, 11800, 13325, 14976, 16758, 18676, 20735**

## Octagonal pyramidal numbers

**Réf.** D1 2 2. B1 194.

**HIS2 A2414** **Approximants de Padé**

**HIS1 N1966** **Fraction rationnelle**

$$\frac{1 + 5z}{(z - 1)^4}$$

**1, 9, 30, 70, 135, 231, 364, 540, 765, 1045, 1386, 1794, 2275, 2835, 3480, 4216, 5049, 5985, 7030, 8190, 9471, 10879, 12420, 14100, 15925, 17901, 20034, 22330, 24795**

## 4-dimensional pyramidal numbers

**Réf.** B1 195.

**HIS2 A2415** **Approximants de Padé**

**HIS1 N1714** **Fraction rationnelle**

$$\frac{1 + z}{(1 - z)^5}$$

**1, 6, 20, 50, 105, 196, 336, 540, 825, 1210, 1716, 2366, 3185, 4200, 5440, 6936, 8721, 10830, 13300, 16170, 19481, 23276, 27600, 32500, 38025, 44226, 51156, 58870**

## 4-dimensional figurate numbers

**Réf.** B1 195.

**HIS2 A2417** **Approximants de Padé**

**HIS1 N1907** **Fraction rationnelle**

$$\frac{1 + 3z}{(1 - z)^5}$$

**1, 8, 30, 80, 175, 336, 588, 960, 1485, 2200, 3146, 4368, 5915, 7840, 10200, 13056, 16473, 20520, 25270, 30800, 37191, 44528, 52900, 62400, 73125, 85176, 98658**

---

## 4-dimensional figurate numbers

**Réf.** B1 195.

**HIS2 A2418** **Approximants de Padé**

**HIS1 N1970** **Fraction rationnelle**

$$\frac{1 + 4z}{(1 - z)^5}$$

**1, 9, 35, 95, 210, 406, 714, 1170, 1815, 2695, 3861, 5369, 7280, 9660, 12580, 16116, 20349, 25365, 31255, 38115, 46046, 55154, 65550, 77350, 90675, 105651**

## 4-dimensional figurate numbers

**Réf.** B1 195.

**HIS2 A2419** **Approximants de Padé**

**HIS1 N2008** **Fraction rationnelle**

$$\frac{1 + 5 z}{(1 - z)^5}$$

**1, 10, 40, 110, 245, 476, 840, 1380, 2145, 3190, 4576, 6370, 8645, 11480, 14960, 19176, 24225, 30210, 37240, 45430, 54901, 65780, 78200, 92300, 108225, 126126**

---

**Réf.** TH09 164. FMR 1 55.

**HIS2 A2420** **Recoupements** **Suite P-récurrente**

**HIS1 N0128** **algébrique**

**a(n) (n - 1) (n - 2) = 2 a(n - 1) (n - 2) (2 n - 5)**

$$(1 - 4 z)^{1/2}$$

**1, 2, 2, 4, 10, 28, 84, 264, 858, 2860, 9724, 33592, 117572, 416024, 1485800, 5348880, 19389690, 70715340, 259289580, 955277400, 3534526380, 13128240840, 48932534040**

**Réf.** TH09 164. FMR 1 55.

**HIS2 A2421** **Recoupements** **Inverse de A2457**

**HIS1 N1683** **algébrique**

$$(1 - 4 z)^{3/2}$$

**1, 6, 6, 4, 6, 12, 28, 72, 198, 572, 1716, 5304, 16796, 54264, 178296, 594320, 2005830, 6843420, 23571780, 81880920, 286583220, 1009864680, 3580429320, 12765008880**

---

**Réf.** TH09 164. FMR 1 55.

**HIS2 A2422** **Recoupements** **Inverse de A2802**

**HIS1 N2003** **algébrique**

$$(1 - 4 z)^{5/2}$$

**1, 10, 30, 20, 10, 12, 20, 40, 90, 220, 572, 1560, 4420, 12920, 38760, 118864, 371450, 1179900, 3801900, 12406200, 40940460, 136468200, 459029400, 1556708400, 5318753700**

**Réf.** TH09 164. FMR 1 55.
**HIS2 A2423** **Recoupements**
**HIS1 N2114** **algébrique**

$$(1 - 4 z)^{7/2}$$

**1, 14, 70, 140, 70, 28, 28, 40, 70, 140, 308, 728, 1820, 4760, 12920, 36176, 104006, 305900, 917700, 2801400, 8684340, 27293640, 86843400, 279409200, 908079900, 2978502072**

**Réf.** TH09 164. FMR 1 55.
**HIS2 A2424** **Recoupements**
**HIS1 N2188** **algébrique**

$$(1 - 4 z)^{9/2}$$

**1, 18, 126, 420, 630, 252, 84, 72, 90, 140, 252, 504, 1092, 2520, 6120, 15504, 40698, 110124, 305900, 869400, 2521260, 7443720, 22331160, 67964400, 209556900, 653817528**

**From expansion of (1+x+x^2) ^n**

**Réf.** EUL (1) 15 59 27. FQ 7 341 69. HE74 1 42.

**HIS2 A2426** **Hypergéométrique**

**HIS1 N1070** **algébrique**

$$\frac{1}{(1+z)^{1/2}\,(3z-1)^{1/2}}$$

**1, 1, 3, 7, 19, 51, 141, 393, 1107, 3139, 8953, 25653, 73789, 212941, 616227, 1787607, 5196627, 15134931, 44152809, 128996853, 377379369**

---

**Réf.** QJM 47 110 16. FMR 1 112. DA63 2 283.

**HIS2 A2446** **Approximants de Padé**

**HIS1 N1748** **Fraction rationnelle**

$$\frac{6z}{(1-4z)(1-z)}$$

**0, 6, 30, 126, 510, 2046, 8190, 32766, 131070, 524286, 2097150, 8388606, 33554430, 134217726, 536870910, 2147483646, 8589934590, 34359738366**

**Réf.** TH09 35. FMR 1 112. RCI 217.
**HIS2 A2450** **Approximants de Padé**
**HIS1 N1608** **Fraction rationnelle**

$$\frac{1}{(1 - 4z)(1 - z)}$$

**1, 5, 21, 85, 341, 1365, 5461, 21845, 87381, 349525, 1398101, 5592405, 22369621, 89478485, 357913941, 1431655765, 5726623061, 22906492245**

**Réf.** TH09 35. FMR 1 112. RCI 217.
**HIS2 A2451** **Approximants de Padé**
**HIS1 N2118** **Fraction rationnelle**

$$\frac{1}{(1 - z)(1 - 4z)(1 - 9z)}$$

**1, 14, 147, 1408, 13013, 118482, 1071799, 9668036, 87099705, 784246870, 7059619931, 63542171784, 571901915677, 5147206719578, 46325218390143, 416928397167052**

## Central factorial numbers

**Réf.** TH09 36. FMR 1 112. RCI 217.

**HIS2 A2452** **Approximants de Padé**

**HIS1 N2025** **Fraction rationnelle**

$$\frac{1}{(1 - z)(1 - 9z)}$$

**1, 10, 91, 820, 7381, 66430, 597871, 5380840, 48427561, 435848050, 3922632451, 35303692060, 317733228541, 2859599056870, 25736391511831**

## Central factorial numbers

**Réf.** TH09 36. FMR 1 112. RCI 217.

**HIS2 A2453** **Approximants de Padé**

**HIS1 N2283** **Fraction rationnelle**

$$\frac{1}{(1 - z)(1 - 9z)(1 - 25z)}$$

**1, 35, 966, 24970, 631631, 15857205, 397027996**

## Central factorial numbers

**Réf.** OP80 7. FMR 1 110. RCI 217.

**HIS2 A2454** **Hypergéométrique** **Suite P-récurrente**

**HIS1 N1510** **Fraction rationnelle** **f.g. exponentielle double**

**a(n) = 4 (n - 1)^2 a(n - 1)**

$$ {}_3F_2([1, 1, 1], [2, 2], 4z) $$

**1, 4, 64, 2304, 147456, 14745600, 2123366400, 416179814400, 106542032486400, 34519618525593600**

## Central differences of 0

**Réf.** QJM 47 110 16. FMR 1 112. DA63 2 283.

**HIS2 A2456** **Hypergéométrique** **Suite P-récurrente**

**HIS1 N2270** **algébrique** **f.g. exponentielle double**

$$ \frac{z (2 + z)}{2 (1 - 2 z)^{7/2}} $$

**1, 30, 1260, 75600, 6237000, 681080400, 95351256000, 16672848192000, 3563821301040000, 914714133933600000, 277707211062240960000**

**Réf.** OP80 21. SE33 92. JO39 449. SAM 22 120 43. LA56 514.

**HIS2 A2457** **Hypergéométrique** **Suite P-récurrente**

**HIS1 N1752** **algébrique**

$$\frac{1}{(1 - 4z)^{3/2}}$$

**1, 6, 30, 140, 630, 2772, 12012, 51480, 218790, 923780, 3879876, 16224936, 67603900, 280816200, 1163381400, 4808643120, 19835652870, 81676217700, 335780006100**

---

**The game of Mousetrap with n cards**

**Réf.** QJM 15 241 1878. jos.

**HIS2 A2467** **Recoupements** **A0166 - 1**

**HIS1 N1423** **exponentielle**

$$\frac{1 - \exp(z)}{(z - 1)\exp(z)}$$

**1, 1, 4, 15, 76, 455, 3186, 25487, 229384, 2293839, 25232230, 302786759, 3936227868, 55107190151, 826607852266, 13225725636255, 224837335816336, 4047072044694047**

## Wonderful Demlo numbers

**Réf.** MAS 6 68 38.

**HIS2 A2477** **Approximants de Padé** **Demlo est une ville aux E.U.**

**HIS1 N2339** **Fraction rationnelle**

**a(n) = 1,11*11, 111*111, 1111*1111,...**

$$\frac{1 + 10z}{(1 - z)(1 - 10z)(1 - 100z)}$$

**1, 121, 12321, 1234321, 123454321, 12345654321, 1234567654321, 123456787654321, 12345678987654321, 1234567900987654321**

## Bisection of A0930

**Réf.** EUL (1) 1 322 11.

**HIS2 A2478** **Approximants de Padé**

**HIS1 N1017** **Fraction rationnelle**

$$\frac{1}{1 - z - 2z^2 - z^3}$$

**1, 1, 3, 6, 13, 28, 60, 129, 277, 595, 1278, 2745, 5896, 12664, 27201, 58425, 125491, 269542, 578949, 1243524, 2670964, 5736961, 12322413, 26467299, 56849086**

**Réf.** ELM 2 95 47. WW 114.
**HIS2 A2487** **Euler**
**HIS1 N0056** **Produit infini**
**a(2n+1) = a(n) et a(2n) = a(n) + a(n-1)**

$$\prod_{n\geq 0} \left(1 + z^{2^n} + z^{2^{(n+1)}}\right)$$

**1, 1, 2, 1, 3, 2, 3, 1, 4, 3, 5, 2, 5, 3, 4, 1, 5, 4, 7, 3, 8, 5, 7, 2, 7, 5, 8, 3, 7, 4, 5, 1, 6, 5, 9, 4, 11, 7, 10, 3, 11, 8, 13, 5, 12, 7, 9, 2, 9, 7, 12, 5, 13, 8, 11, 3, 10, 7, 11, 4, 9, 5, 6, 1, 7**

---

**Réf.** MOC 4 23 50.
**HIS2 A2492** **Approximants de Padé**
**HIS1 N1444** **Fraction rationnelle**

$$\frac{4\,(1 + z)}{(z - 1)^4}$$

**4, 20, 56, 120, 220, 364, 560, 816, 1140, 1540, 2024, 2600, 3276, 4060, 4960, 5984, 7140, 8436, 9880, 11480, 13244, 15180, 17296, 19600, 22100, 24804, 27720**

## Expansion of a modular function

**Réf.** PLMS 9 386 59.

**HIS2 A2512** **Euler**

**HIS1 N0539** **Produit infini**

**Conjecture : erreurs dans la suite à partie du 12è terme ?**

$$\prod_{n\geq 1} \frac{1}{(1 - z^{n})^{c(n)}}$$

**c(n) = 2,2,2,4,2,2,2,4,...**

**1, 2, 5, 10, 22, 40, 75, 130, 230, 382, 636, 1016, 1633, 2540, 3942, 5978, 9057**

---

## Expansion of a modular function

**Réf.** PLMS 9 387 59.

**HIS2 A2513** **Euler** **erreur probable à partir du 13è**

**HIS1 N0931** **Produit infini** **terme**

*** Le motif [1,2] est périodique**

$$\prod_{n\geq 1} \frac{1}{(1 - z^{c(n)})}$$

**c(n) = 1,2,...***

**1, 1, 3, 4, 9, 12, 23, 31, 54, 73, 118, 159, 246, 340, 500, 684, 984, 1341, 1883**

## Permutations of length n within distance 2

**Réf.** AENS 79 207 62.

**HIS2 A2524** **Approximants de Padé**

**HIS1 N0626** **Fraction rationnelle**

$$\frac{1 - z}{1 - 2z - 2z^3 + z^5}$$

1, 1, 2, 6, 14, 31, 73, 172, 400, 932, 2177, 5081, 11854, 27662, 64554

---

## Permutations according to distance

**Réf.** AENS 79 207 62.

**HIS2 A2525** **Approximants de Padé**

**HIS1 N0463** **Fraction rationnelle**

$$\frac{z}{1 - 2z - 2z^3 + z^5}$$

0, 1, 2, 4, 10, 24, 55, 128, 300, 700, 1632, 3809, 8890, 20744, 48406

**Réf.** MQET 1 10 16. NZ66 181.
**HIS2 A2530** **Approximants de Padé**
**HIS1 N0934** **Fraction rationnelle**

$$\frac{1 - z - z^2}{1 - 4 z^2 + z^4}$$

**1, 1, 3, 4, 11, 15, 41, 56, 153, 209, 571, 780, 2131, 2911, 7953, 10864, 29681, 40545, 110771, 151316, 413403, 564719, 1542841, 2107560, 5757961, 7865521**

---

**Réf.** MQET 1 10 16. NZ66 181.
**HIS2 A2531** **Approximants de Padé**
**HIS1 N0513** **Fraction rationnelle**

$$\frac{1 + z - 2 z^2 + z^3}{1 - 4 z^2 + z^4}$$

**1, 1, 2, 5, 7, 19, 26, 71, 97, 265, 362, 989, 1351, 3691, 5042, 13775, 18817, 51409, 70226, 191861, 262087, 716035, 978122, 2672279, 3650401, 9973081, 13623482**

**Réf.** MQET 1 11 16.
**HIS2 A2532** **Approximants de Padé**
**HIS1 N0758** **Fraction rationnelle**

$$\frac{z}{1 - 2 z - 5 z^{2}}$$

**0, 1, 2, 9, 28, 101, 342, 1189, 4088, 14121, 48682, 167969, 579348, 1998541, 6893822, 23780349, 82029808, 282961361, 976071762, 3366950329, 11614259468**

---

**Réf.** MQET 1 11 16.
**HIS2 A2533** **Approximants de Padé**
**HIS1 N1834** **Fraction rationnelle**

$$\frac{1 - z}{1 - 2 z - 5 z^{2}}$$

**1, 1, 7, 19, 73, 241, 847, 2899, 10033, 34561, 119287, 411379, 1419193, 4895281, 16886527, 58249459, 200931553, 693110401, 2390878567, 8247309139**

**Réf.** MQET 1 11 16.

**HIS2 A2534** **Approximants de Padé**

**HIS1 N0814** **Fraction rationnelle**

$$\frac{z}{1 - 2z - 9z^2}$$

**0, 1, 2, 13, 44, 205, 806, 3457, 14168, 59449, 246410, 1027861, 4273412, 17797573, 74055854, 308289865, 1283082416, 5340773617, 22229288978, 92525540509**

---

**Réf.** MQET 1 11 16.

**HIS2 A2535** **Approximants de Padé**

**HIS1 N2043** **Fraction rationnelle**

$$\frac{1 - z}{1 - 2z - 9z^2}$$

**1, 1, 11, 31, 161, 601, 2651, 10711, 45281, 186961, 781451, 3245551, 13524161, 56258281, 234234011, 974792551, 4057691201, 16888515361, 70296251531**

**Réf.** MQET 1 12 16.

**HIS2 A2536** **Approximants de Padé**

**HIS1 N1540** **Fraction rationnelle**

$$\frac{z\,(1 + z - 3\,z^2)}{1 - 8\,z^2 + 9\,z^4}$$

**0, 1, 1, 5, 8, 31, 55, 203, 368, 1345, 2449, 8933, 16280, 59359, 108199**

---

**Réf.** MQET 1 12 16.

**HIS2 A2537** **Approximants de Padé**

**HIS1 N1379** **Fraction rationnelle**

$$\frac{1 + z - 4\,z^2 + 3\,z^3}{1 - 8\,z^2 + 9\,z^4}$$

**1, 1, 4, 11, 23, 79, 148, 533, 977, 3553, 6484, 23627, 43079, 157039, 286276, 1043669, 1902497, 6936001, 12643492, 46094987, 84025463, 306335887, 558412276, 2035832213**

## Coefficients for numerical differentiation

**Réf.** OP80 21. SE33 92. SAM 22 120 43. LA56 514.

**HIS2 A2544** **Hypergéométrique** **Suite P-récurrente**

**HIS1 N2075** **algébrique**

**$_2F_1$ ([2, 3/2], [1], 4 z)**

$$\frac{1 + 2z}{(1 - 4z)^{5/2}}$$

**1, 12, 90, 560, 3150, 16632, 84084, 411840, 1969110, 9237800, 42678636, 194699232, 878850700, 3931426800, 17450721000**

---

## From a definite integral

**Réf.** EMS 10 184 57.

**HIS2 A2570** **Approximants de Padé**

**HIS1 N1698** **Fraction rationnelle**

$$\frac{1}{(1 - z)\,(1 - 3z + z^2)\,(1 + z)^3}$$

**1, 1, 6, 11, 36, 85, 235, 600, 1590, 4140, 10866, 28416, 74431, 194821, 510096, 1335395, 3496170, 9153025, 23963005, 62735880**

## From a definite integral

**Réf.** EMS 10 184 57.

**HIS2 A2571** **Approximants de Padé**

**HIS1 N1553** **Fraction rationnelle**

$$\frac{1 + 4z + z^2 - z^3}{(1 - 3z + z^2)(1 + z)^2}$$

**1, 5, 10, 30, 74, 199, 515, 1355, 3540, 9276, 24276, 63565, 166405, 435665, 1140574, 2986074, 7817630, 20466835, 53582855, 140281751**

---

**Réf.** CC55 742. JO61 7.

**HIS2 A2593** **Approximants de Padé**

**HIS1 N2262** **Fraction rationnelle**

$$\frac{z(1 + z)(z^2 + 22z + 1)}{(z - 1)^5}$$

**0, 1, 28, 153, 496, 1225, 2556, 4753, 8128, 13041, 19900, 29161, 41328, 56953, 76636, 101025, 130816, 166753, 209628, 260281, 319600, 388521, 468028, 559153**

## Sums of 5th powers of odd numbers

**Réf.** CC55 742.

**HIS2 A2594** **Approximants de Padé**

**HIS1 N2354** **Fraction rationnelle**

$$\frac{(1 + z)\,(z^4 + 236\,z^3 + 1446\,z^2 + 236\,z + 1)}{(1 - z)^7}$$

**1, 244, 3369, 20176, 79225, 240276, 611569, 1370944, 2790801, 5266900, 9351001, 15787344, 25552969, 39901876, 60413025, 89042176, 128177569, 180699444**

## A generalized partition function

**Réf.** PNISI 17 237 51.

**HIS2 A2597** **LLL**

**HIS1 N1000** **Fraction rationnelle**

$$\frac{1}{(z + 1)^2\,(z^2 + z + 1)^3\,(z - 1)^6\,z}$$

**1, 3, 6, 9, 15, 25, 34, 51, 73, 97, 132, 178, 226, 294, 376, 466, 582, 722, 872, 1062, 1282, 1522, 1812, 2147, 2507, 2937, 3422, 3947, 4557, 5243, 5978, 6825, 7763, 8771**

**Réf.** AMS 26 304 55.

**HIS2 A2620** **Approximants de Padé**

**HIS1 N0374** **Fraction rationnelle**

$$\frac{1}{(1 + z)\ (z - 1)^{3}}$$

**1, 2, 4, 6, 9, 12, 16, 20, 25, 30, 36, 42, 49, 56, 64, 72, 81, 90, 100, 110, 121, 132, 144, 156, 169, 182, 196, 210, 225, 240, 256, 272, 289, 306, 324, 342, 361, 380, 400, 420**

---

**Réf.** AMS 26 304 55.

**HIS2 A2621** **Approximants de Padé**

**HIS1 N0394** **Fraction rationnelle**

$$\frac{1}{(1 + z^{2})\ (z^{2} + z + 1)\ (1 + z)^{2}\ (z - 1)^{5}}$$

**1, 2, 4, 7, 12, 18, 27, 38, 53, 71, 94, 121, 155, 194, 241, 295, 359, 431, 515, 609, 717, 837, 973, 1123, 1292, 1477, 1683, 1908, 2157, 2427, 2724, 3045, 3396, 3774, 4185**

## A partition function

**Réf.** AMS 26 304 55.

**HIS2 A2622** **Approximants de Padé**

**HIS1 N0395** **Fraction rationnelle**

$$\frac{1}{(1-z)^2 (1-z^2)(1-z^3)(1-z^4)(1-z^5)}$$

**1, 2, 4, 7, 12, 19, 29, 42, 60, 83, 113, 150, 197, 254, 324, 408, 509, 628, 769, 933, 1125, 1346, 1601, 1892, 2225, 2602, 3029, 3509, 4049, 4652, 5326, 6074, 6905, 7823**

---

**Réf.** AMS 26 308 55. PGEC 22 1050 73.

**HIS2 A2623** **Approximants de Padé**

**HIS1 N1050** **Fraction rationnelle**

$$\frac{1}{(1+z)(z-1)^4}$$

**1, 3, 7, 13, 22, 34, 50, 70, 95, 125, 161, 203, 252, 308, 372, 444, 525, 615, 715, 825, 946, 1078, 1222, 1378, 1547, 1729, 1925, 2135, 2360, 2600, 2856, 3128, 3417, 3723**

## A partition function

**Réf.** AMS 26 308 55.

**HIS2 A2624** **Approximants de Padé**

**HIS1 N1091** **Fraction rationnelle**

$$\frac{1}{(1 + z)^2 (1 - z)^5}$$

**1, 3, 8, 16, 30, 50, 80, 120, 175, 245, 336, 448, 588, 756, 960, 1200, 1485, 1815, 2200, 2640, 3146, 3718, 4368, 5096, 5915, 6825, 7840, 8960, 10200, 11560, 13056**

---

**Réf.** AMS 26 308 55.

**HIS2 A2625** **Approximants de Padé**

**HIS1 N1093** **Fraction rationnelle**

$$\frac{1}{(z^2 + z + 1)(1 + z)^2 (z - 1)^6}$$

**1, 3, 8, 17, 33, 58, 97, 153, 233, 342, 489, 681, 930, 1245, 1641, 2130, 2730, 3456, 4330, 5370, 6602, 8048, 9738, 11698, 13963, 16563, 19538, 22923, 26763, 31098, 35979**

**Réf.** AMS 26 308 55.
**HIS2 A2626** **Approximants de Padé**
**HIS1 N1094** **Fraction rationnelle**

$$\frac{1}{(z^2 + 1)(z^2 + z + 1)(z + 1)^3 (1 - z)^7}$$

**1, 3, 8, 17, 34, 61, 105, 170, 267, 403, 594, 851, 1197, 1648, 2235, 2981, 3927, 5104, 6565, 8351, 10529, 13152, 16303, 20049, 24492, 29715, 35841, 42972, 51255**

---

**Réf.** MFM 73 18 69.
**HIS2 A2662** **Approximants de Padé**
**HIS1 N1585** **Fraction rationnelle**

$$\frac{z^2}{(2z - 1)(z - 1)^3}$$

**0, 0, 1, 5, 16, 42, 99, 219, 466, 968, 1981, 4017, 8100, 16278, 32647, 65399, 130918, 261972, 524097, 1048365, 2096920, 4194050, 8388331, 16776915, 33554106, 67108512**

**Réf.** MFM 73 18 69.

**HIS2 A2663** **Approximants de Padé**

**HIS1 N1725** **Fraction rationnelle**

$$\frac{1}{(2z-1)(1-z)^4}$$

**1, 6, 22, 64, 163, 382, 848, 1816, 3797, 7814, 15914, 32192, 64839, 130238, 261156, 523128, 1047225, 2095590, 4192510, 8386560, 16774891, 33551806, 67105912, 134214424**

**Réf.** MFM 73 18 69.

**HIS2 A2664** **Approximants de Padé**

**HIS1 N1851** **Fraction rationnelle**

$$\frac{1}{(2z-1)(1-z)^5}$$

**1, 7, 29, 93, 256, 638, 1486, 3302, 7099, 14913, 30827, 63019, 127858, 258096, 519252, 1042380, 2089605, 4185195, 8377705, 16764265, 33539156, 67090962, 134196874, 268411298**

## Coefficients for central differences

**Réf.** SAM 42 162 63.

**HIS2 A2671** **Hypergéométrique**

**HIS1 N2246** **algébrique**

$$\frac{1}{(1 - 16 z)^{3/2}}$$

**1, 24, 1920, 322560, 92897280, 40874803200, 25505877196800, 21424936845312000, 23310331287699456000, 31888533201572855808000**

## Coefficients for central differences

**Réf.** SAM 42 162 63.

**HIS2 A2674** **Hypergéométrique** **f.g. exponentielle double**

**HIS1 N2092** **algébrique**

$$\frac{1}{2 (1 - 4 z)^{1/2}}$$

**1, 12, 360, 20160, 1814400, 239500800, 43589145600, 10461394944000, 3201186852864000, 1216451004088320000, 562000363888803840000**

## Coefficients of orthogonal polynomials

**Réf.** MOC 9 174 55.

**HIS2 A2690** **Dérivée logarithmique** **Suite P-récurrente**

**HIS1 N1491** **exponentielle:algébrique**

**a(n) = (4 n-4) a(n - 1) + (8 n - 20) a(n - 2)**

$$\frac{1 - 2 z}{(1 - 4 z)^{3/2}}$$

**1, 4, 36, 480, 8400, 181440, 4656960, 138378240, 4670265600, 176432256000, 7374868300800, 337903056691200**

## Coefficients of orthogonal polynomials

**Réf.** MOC 9 174 55.

**HIS2 A2691** **Dérivée logarithmique** **Suite P-récurrente**

**HIS1 N1996** **exponentielle**

**n a(n) = 2 (n + 1) (2 n - 1) a(n - 1)**

$$\frac{1 - z}{(1 - 4 z)^{5/2}}$$

**1, 9, 120, 2100, 45360, 1164240, 34594560, 1167566400, 44108064000, 1843717075200, 84475764172800**

## Binomial coefficients C(2 n , n - 2)

**Réf.** LA56 517. AS1 828.

**HIS2 A2694** **Hypergéométrique**

**HIS1 N1741** **algébrique**

$$\frac{16}{(1-4z)^{1/2}\,\left(1+(1-4z)^{1/2}\right)^{4}}$$

**1, 6, 28, 120, 495, 2002, 8008, 31824, 125970, 497420, 1961256, 7726160, 30421755, 119759850, 471435600, 1855967520, 7307872110, 28781143380**

---

## Spheroidal harmonics

**Réf.** MES 52 75 24.

**HIS2 A2695** **LLL** **Suite P-récurrente**

**HIS1 N1985** **algébrique**

**(n - 2) a(n) = (6 n - 9) a(n - 1) + (- n + 1) a(n - 2)**

$$\frac{z}{(z^{2}-6z+1)^{3/2}}$$

**0, 1, 9, 66, 450, 2955, 18963, 119812, 748548, 4637205, 28537245**

**Réf.** LA56 517. AS1 828.
**HIS2 A2696** **Hypergéométrique**
**HIS1 N1921** **algébrique**
$_2F_1([7/2, 4], [7], 4 z)$

$$\frac{64}{(1 - 4 z)^{1/2} \left(1 + (1 - 4 z)^{1/2}\right)^6}$$

**1, 8, 45, 220, 1001, 4368, 18564, 77520, 319770, 1307504, 5311735, 21474180, 86493225, 347373600, 1391975640, 5567902560, 22239974430, 88732378800**

## Coefficients of Chebyshev polynomials

**Réf.** LA56 516.
**HIS2 A2697** **Approximants de Padé**
**HIS1 N1923** **Fraction rationnelle**

$$\frac{1}{(4 z - 1)^2}$$

**1, 8, 48, 256, 1280, 6144, 28672, 131072, 589824, 2621440, 11534336, 50331648**

## Coefficients of Chebyshev polynomials

**Réf.** LA56 516.

**HIS2 A2698** **Approximants de Padé**

**HIS1 N2189** **Fraction rationnelle**

$$\frac{1 + 6 z - 8 z^2}{(1 - 4 z)^3}$$

**1, 18, 160, 1120, 6912, 39424, 212992, 1105920, 5570560, 27394048, 132120576**

---

**Réf.** LA56 518.

**HIS2 A2699** **Approximants de Padé**

**HIS1 N0825** **Fraction rationnelle**

$$\frac{2 z}{(4 z - 1)^2}$$

**0, 2, 16, 96, 512, 2560, 12288, 57344, 262144, 1179648, 5242880, 23068672, 100663296, 436207616, 1879048192, 8053063680, 34359738368, 146028888064, 618475290624, 2611340115968**

## Coefficients of Chebyshev polynomials

**Réf.** LA56 518.

**HIS2 A2700** **Approximants de Padé**

**HIS1 N1275** **Fraction rationnelle**

$$\frac{4z - 3}{(4z - 1)^3}$$

**3, 40, 336, 2304, 14080, 79872, 430080, 2228224, 11206656, 55050240, 265289728, 1258291200**

## Keys

**Réf.** MAG 53 11 69.

**HIS2 A2714** **Approximants de Padé**

**HIS1 N1832** **Fraction rationnelle**

$$\frac{7 - 9z - 9z^2 + 3z^3}{1 - 4z + 2z^2 + 4z^3 - z^4}$$

**7, 19, 53, 149, 421, 1193, 3387, 9627, 27383, 77923**

**Réf.** MAG 46 55 62; 55 440 71. MMAG 47 290 74.
**HIS2 A2717** **Approximants de Padé**
**HIS1 N1569** **Fraction rationnelle**

$$\frac{1 + 2 z}{(1 + z) (z - 1)^4}$$

**1, 5, 13, 27, 48, 78, 118, 170, 235, 315, 411, 525, 658, 812, 988, 1188, 1413, 1665, 1945, 2255, 2596, 2970, 3378, 3822, 4303, 4823, 5383, 5985, 6630, 7320, 8056, 8840**

---

**Réf.** SE33 78.
**HIS2 A2720** **Dérivée logarithmique** **Suite P-récurrente**
**HIS1 N0708** **exponentielle (rationnel)**
**a(n) = (2 n - 2) a(n-1) + (-n^2 + 4 n-4) a(n-2)**

$$\frac{1}{(1 - z) \exp(z/(z-1))}$$

**1, 2, 7, 34, 209, 1546, 13327, 130922, 1441729, 17572114, 234662231, 3405357682, 53334454417, 896324308634, 16083557845279, 306827170866106, 6199668952527617**

## Apéry numbers

**Réf.** SE33 93. MI 1 195 78.

**HIS2 A2736** **Hypergéométrique**

**HIS1 N0848** **algébrique**

$$\frac{1 + 2 z}{(1 - 4 z)^{5/2}}$$

**0, 2, 24, 180, 1120, 6300, 33264, 168168, 823680, 3938220, 18475600, 85357272, 389398464, 1757701400, 7862853600, 34901442000, 153876579840, 674412197580, 2940343837200**

## Coefficients for extrapolation

**Réf.** SE33 97.

**HIS2 A2740** **Hypergéométrique** **Suite P-récurrente**

**HIS1 N0821** **algébrique**

$$\frac{6 z^2 - 6 z + 1 + (1 - 4 z)^{3/2}}{- 2 (1 - 4 z)^{3/2} z^3}$$

**0, 2, 15, 84, 420, 1980, 9009, 40040**

## Logarithmic numbers

**Réf.** MAS 31 77 63. jos.

**HIS2 A2741** **Recoupements** **Suite P-récurrente**

**HIS1 N0010** **exponentielle**

**a(n) = (n - 3) a(n - 1) + (n - 2) a(n - 3) + (2 n - 4) a(n - 2)**

$$\frac{\ln(1 - z)}{\exp(z)}$$

**1, 1, 2, 0, 9, 35, 230, 1624, 13209, 120287, 1214674, 13469896, 162744945, 2128047987, 29943053062, 451123462672, 7245940789073, 123604151490591**

---

## Logarithmic numbers

**Réf.** MAS 31 78 63. jos.

**HIS2 A2747** **Dérivée logarithmique** **Suite P-récurrente**

**HIS1 N0759** **exponentielle**

**a(n) = 2 a(n - 1) + (n^2 - n - 1) a(n - 2) + (- 2 n^2 + 6 n - 4) a(n - 3)**
**+ (n^2 - 5 n + 6) a(n - 4)**

$$\frac{\exp(z)\,(z^3 - z^2 - z - 1)}{(1 - z)^2\,(z + 1)^2}$$

**1, 2, 9, 28, 185, 846, 7777, 47384, 559953, 4264570, 61594841, 562923252, 9608795209, 102452031878, 2017846993905, 24588487650736, 548854382342177**

## Terms in certain determinants

**Réf.** PLMS 10 122 1879.

**HIS2 A2775** **Dérivée logarithmique**

**HIS1 N1927** **Fraction rationnelle**

$$\frac{z^2 + 4z + 1}{(z - 1)^4}$$

**0, 1, 8, 54, 384, 3000, 25920, 246960, 2580480**

---

**Réf.** IJ1 11 162 69.

**HIS2 A2783** **Approximants de Padé**

**HIS1 N1159** **Fraction rationnelle**

$$\frac{1 - 3z + 4z^2}{(1 - z)(1 - 2z)(1 - 3z)}$$

**1, 3, 11, 39, 131, 423, 1331, 4119, 12611, 38343, 116051, 350199, 1054691, 3172263, 9533171, 28632279, 85962371, 258018183, 774316691, 2323474359, 6971471651, 20916512103**

**Réf.** JRAM 227 49 67.

**HIS2 A2798** **Approximants de Padé**

**HIS1 N2186** **Fraction rationnelle**

$$\frac{3\,(6 + 9\,z + 2\,z^{2})}{(1 + z)\,(z - 1)^{2}}$$

**18, 45, 69, 96, 120, 147, 171**

---

**Réf.** AJM 2 94 1879. LU91 1 223.

**HIS2 A2801** **équations différentielles** **Suite P-récurrente**

**HIS1 N0744** **exponentielle (algébrique)** **Formule de B. Salvy**

**a(n) = (2 n - 3) a(n - 1) + (- n + 2) a(n - 2)**

$$\frac{\exp(1/2\,z)\,2^{3/4}}{(-1 + 2\,z)^{1/4}}$$

**1, 1, 2, 8, 50, 418, 4348, 54016, 779804, 12824540, 236648024, 4841363104, 108748223128, 2660609220952, 70422722065040, 2005010410792832**

**Réf.** JO39 449. JCT 13 215 72.
**HIS2 A2802** **Hypergéométrique**
**HIS1 N2019** **algébrique**
**$_2F_1$([5/2], [ ], 4 z)**

$$\frac{1}{(1 - 4z)^{5/2}}$$

**1, 10, 70, 420, 2310, 12012, 60060, 291720, 1385670, 6466460, 29745716, 135207800, 608435100, 2714556600, 12021607800, 52895074320, 231415950150, 1007340018300**

---

**Réf.** JO39 449. JCT B18 258 75.
**HIS2 A2803** **Hypergéométrique** **Suite P-récurrente**
**HIS1 N2140** **algébrique**
**$_2F_1$([5/2], [ ], 4 z)**

$$\frac{1 + z}{(1 - 4z)^{7/2}}$$

**1, 15, 140, 1050, 6930, 42042, 240240, 1312740, 6928350, 35565530, 178474296, 878850700, 4259045700, 20359174500, 96172862400, 449608131720, 2082743551350**

**Réf.** PIEE 115 763 68. DM 55 272 85.

**HIS2 A2807** **P-récurrences** **Suite P-récurrente**

**HIS1 N1867**

$$a(n) = n\, a(n-5) + (6n+1)\, a(n-3) - (4n+7)\, a(n-2) + (n+5)\, a(n-1) - 2\, a(n-5) + (-4n+4)\, a(n-4)$$

**0, 0, 1, 7, 37, 197, 1172, 8018, 62814, 556014, 5488059, 59740609, 710771275, 9174170011, 127661752406, 1904975488436, 30341995265036, 513771331467372, 9215499383109573**

---

## Doubly triangular numbers

**Réf.** TCPS 9 477 1856. SIAC 4 477 75. ANS 4 1178 76.

**HIS2 A2817** **Approximants de Padé**

**HIS1 N1718** **Fraction rationnelle**

$$\frac{1 + z + z^2}{(1 - z)^5}$$

**1, 6, 21, 55, 120, 231, 406, 666, 1035, 1540, 2211, 3081, 4186, 5565, 7260, 9316, 11781, 14706, 18145, 22155, 26796, 32131, 38226, 45150, 52975, 61776, 71631, 82621**

## Partitions of n into parts 1/2, 3/4, 7/8, etc

**Réf.** EMS 11 224 59.

**HIS2 A2843** **Approximants de Padé** **Conjecture**

**HIS1 N0405** **Fraction rationnelle**

$$\frac{(z^2 + z + 1)(z - 1)^2}{1 - 2z - z^3 + 3z^4}$$

**1, 1, 2, 4, 7, 13, 24, 43, 78, 141, 253, 456**

---

## Partitions of n with no part of size 1

**Réf.** TAIT 1 334. AS1 836.

**HIS2 A2865** **Euler**

**HIS1 N0113** **Produit infini**

$$\prod_{n \geq 1} \frac{1}{(1 - z^n)^{c(n)}}$$

$$c(n) = 0,1,1,1,1,1,...$$

**1, 0, 1, 1, 2, 2, 4, 4, 7, 8, 12, 14, 21, 24, 34, 41, 55, 66, 88, 105, 137, 165, 210, 253, 320, 383, 478, 574, 708, 847, 1039, 1238, 1507, 1794, 2167, 2573, 3094, 3660, 4378, 5170**

**Réf.** PSPM 19 172 71.
**HIS2 A2866** **Dérivée logarithmique** **f.g. exponentielle**
**HIS1 N1463** **Fraction rationnelle**
**a(n) = 2 ^ (n-1) Γ (n+1)**

$$\frac{1}{(1 - 2 z)^{2}}$$

**1, 4, 24, 192, 1920, 23040, 322560, 5160960, 92897280, 1857945600, 40874803200, 980995276800, 25505877196800, 714164561510400, 21424936845312000, 685597979049984000**

---

**Réf.** PSPM 19 172 71.
**HIS2 A2867** **Dérivée logarithmique** **Suite P-récurrente**
**HIS1 N0806** **algébrique** **f.g. exponentielle**
**a(n) = 2 a(n - 1) + (4 n^2 - 12 n + 8) a(n - 2)**

$$\frac{1}{(1 - 2 z)^{3/2} (2 z + 1)^{1/2}}$$

**1, 2, 12, 72, 720, 7200, 100800, 1411200, 25401600, 457228800, 10059033600, 221298739200, 5753767219200, 149597947699200, 4487938430976000, 134638152929280000**

## Sorting numbers

**Réf.** PSPM 19 173 71.

**HIS2 A2871** **équations différentielles** **Formule de B. Salvy**

**HIS1 N0483** **exponentielle**

$$\exp(1/2 \exp(2 z) + \exp(z) - 3/2)$$

**1, 2, 4, 12, 48, 200, 1040, 5600, 33600**

---

## Sorting numbers

**Réf.** PSPM 19 173 71.

**HIS2 A2874** **équations différentielles** **Formule de B. Salvy**

**HIS1 N0738** **exponentielle**

$$\exp(1/3 \exp(3 z) + \exp(z) - 4/3)$$

**1, 2, 8, 42, 268, 1994, 16852**

## Bisection of Lucas sequence

**Réf.** FQ 9 284 71.

**HIS2 A2878** **Approximants de Padé**

**HIS1 N1384** **Fraction rationnelle**

$$\frac{1 + z}{1 - 3z + z^2}$$

**1, 4, 11, 29, 76, 199, 521, 1364, 3571, 9349, 24476, 64079, 167761, 439204, 1149851, 3010349, 7881196, 20633239, 54018521, 141422324, 370248451, 969323029**

---

**Réf.** AIP 9 345 60. SIAR 17 168 75.

**HIS2 A2893** **P-récurrences** **Suite P-récurrente**

**HIS1 N1214**

**a(n) = ∑ C(n,k) ^ 2 .C(2k,k), k=0..n**

$$(n - 1)^2 a(n) = (10 n^2 - 30 n + 23) a(n - 1) + (- 9 n^2 + 36 n - 36) a(n - 2)$$

**1, 3, 15, 93, 639, 4653, 35169, 272835, 2157759, 17319837, 140668065, 1153462995, 9533639025, 79326566595, 663835030335, 5582724468093, 47152425626559, 399769750195965**

## 2n-step polygons on square lattice

**Réf.** AIP 9 345 60.

**HIS2 A2894** **hypergéométrique** **Suite P-récurrente**

**HIS1 N1490** **Intégrales elliptiques**

$$2F1\ ([1/2,\ 1/2],\ [1],\ 16\ z)$$

**1, 4, 36, 400, 4900, 63504, 853776**

---

## 2n-step polygons on b.c.c. lattice

**Réf.** AIP 9 345 60.

**HIS2 A2897** **hypergéométrique** **Suite P-récurrente**

**HIS1 N1952** **Intégrales elliptiques**

$$3F2\ ([1/2,\ 1/2,\ 1/2],\ [1,\ 1],\ 64\ z)$$

**1, 8, 216, 8000, 343000, 16003008, 788889024**

**Réf.** JALG 20 173 72.

**HIS2 A2965** **Approximants de Padé**

**HIS1** **Fraction rationnelle**

$$\frac{1 + 2z + z^2 + z^3}{1 - 2z^2 - z^4}$$

**1, 2, 3, 5, 7, 12, 17, 29, 41, 70, 99, 169, 239, 408, 577, 985, 1393, 2378, 3363, 5741, 8119, 13860, 19601, 33461, 47321, 80782, 114243, 195025, 275807, 470832**

---

## Problimes (second definition)

**Réf.** AMM 80 677 73.

**HIS2 A3067** **Approximants de Padé**

**HIS1** **Fraction rationnelle**

**Conjecture seulement , le dernier terme aurait dû être : 89**

$$\frac{z^9 + z^5 + z^2 + 2}{(z - 1)^2}$$

**2, 4, 7, 10, 13, 17, 21, 25, 29, 34, 39, 44, 49, 54, 59, 64, 69, 74, 79, 84, 90**

## Partitions of n into parts 6n+1 or 6n-1

**Réf.**

**HIS2 A3105** **Euler**

**HIS1** **Produit infini**

$$\prod_{n\geq 1} \frac{1}{(1 - z^{c(n)})}$$

**c(n) = n congru à 1,5 mod 6**

**1, 1, 1, 1, 1, 2, 2, 3, 3, 3, 4, 5, 6, 7, 8, 9, 10, 12, 14, 16, 18, 20, 23, 26, 30, 34, 38, 42, 47, 53, 60, 67, 74, 82, 91, 102, 114, 126, 139, 153, 169, 187, 207, 228, 250, 274, 301, 331, 364**

---

## Partitions of n into parts 5n+2 or 5n+3

**Réf.** AN76 238. AMM 95 711 88; 96 403 89.

**HIS2 A3106** **Euler**

**HIS1** **Produit infini**

$$\prod_{n\geq 1} \frac{1}{(1 - z^{c(n)})}$$

**c(n) = n congru à 2,3 mod 5**

**1, 0, 1, 1, 1, 1, 2, 2, 3, 3, 4, 4, 6, 6, 8, 9, 11, 12, 15, 16, 20, 22, 26, 29, 35, 38, 45, 50, 58, 64, 75, 82, 95, 105, 120, 133, 152, 167, 190, 210, 237, 261, 295, 324, 364, 401, 448, 493, 551**

## Partitions of n into Fibonacci parts

Réf.
HIS2 A3107 Euler
HIS1 Produit infini

$$\prod_{n\geq 1} \frac{1}{(1 - z^{c(n)})}$$

**c(n) = Nombres de Fibonacci.**

1, 1, 2, 3, 4, 6, 8, 10, 14, 17, 22, 27, 33, 41, 49, 59, 71, 83, 99, 115, 134, 157, 180, 208, 239, 272, 312, 353, 400, 453, 509, 573, 642, 717, 803, 892, 993, 1102, 1219, 1350

## Partitions of n into cubes

Réf.
HIS2 A3108 Euler
HIS1 Produit infini

$$\prod_{n\geq 1} \frac{1}{(1 - z^{c(n)})}$$

**c(n) = 1,8,27,64,... Cubes**

1, 1, 1, 1, 1, 1, 1, 1, 2, 2, 2, 2, 2, 2, 2, 2, 3, 3, 3, 3, 3, 3, 3, 3, 4, 4, 4, 5, 5, 5, 5, 5, 6, 6, 6, 7, 7, 7, 7, 7, 8, 8, 8, 9, 9, 9, 9, 9, 10, 10, 10, 11, 11, 11, 12, 12, 13, 13, 13, 14, 14, 14, 15, 15, 17, 17

## Partitions of n into parts 5n+1 and 5n-1

**Réf.** AN76 238. AMM 95 711 88; 96 403 89.

**HIS2 A3114** **Euler**

**HIS1** **Produit infini**

$$\prod_{n\geq 1} \frac{1}{(1 - z^{c(n)})}$$

$$c(n) = n \text{ congru à } 1,4 \bmod 5$$

**1, 1, 1, 1, 2, 2, 3, 3, 4, 5, 6, 7, 9, 10, 12, 14, 17, 19, 23, 26, 31, 35, 41, 46, 54, 61, 70, 79, 91, 102, 117, 131, 149, 167, 189, 211, 239, 266, 299, 333, 374, 415, 465, 515, 575, 637**

---

## Arborescences of type (n,1)

**Réf.** DM 5 197 73.

**HIS2 A3120** **Approximants de Padé** **Conjecture**

**HIS1** **Fraction rationnelle**

$$\frac{(z - 1)(3z^2 + z - 1)}{1 - 3z - z^2 + 7z^3 - 3z^4}$$

**1, 1, 2, 3, 7, 13, 31, 66, 159**

**Réf.** KN1 3 207.
**HIS2 A3143** **Approximants de Padé**
**HIS1** **Fraction rationnelle**

$$\frac{1 + z^3 - z^4 + z^5 - z^6 + z^7}{(z - 1)(1 - z + z^2)(z^2 + z + 1)(-1 + 2 z^2)}$$

**1, 1, 2, 3, 4, 6, 9, 13, 19, 27, 38, 54, 77, 109, 155, 219, 310, 438, 621, 877, 1243, 1755, 2486, 3510, 4973, 7021, 9947, 14043, 19894, 28086, 39789, 56173, 79579, 112347**

---

**Réf.** FQ 10 171 72.
**HIS2 A3148** **Dérivée logarithmique** **Suite P-récurrente**
**HIS1** **algébrique** **f.g. exponentielle**
**a(n) = a(n - 1) + (4 n^2 - 14 n + 12) a(n - 2)**

$$\frac{1}{(1 - 2 z)(1 + 2 z)^{1/2}}$$

**1, 1, 7, 27, 321, 2265, 37575, 390915, 8281665, 114610545, 2946939975, 51083368875, 1542234996225, 32192256321225, 1114841223671175**

## Star numbers

**Réf.** GA88 20.

**HIS2 A3154** **Approximants de Padé**

**HIS1** **Fraction rationnelle**

$$\frac{z^2 + 10z + 1}{(1 - z)^3}$$

**1, 13, 37, 73, 121, 181, 253, 337, 433, 541, 661, 793, 937, 1093, 1261, 1441, 1633, 1837, 2053, 2281, 2521, 2773, 3037, 3313, 3601, 3901, 4213, 4537, 4873, 5221, 5581**

---

## If n appears, 2n doesn't

**Réf.** FQ 10 501 72. AMM 87 671 80.

**HIS2 A3159** **Euler** **Suite reliée à la suite de**

**HIS1** **Produit infini** **Thue-Morse.**

*** Voir [AABBJPS]**

$$\frac{(1 + Z) \prod_{n \geq 0} (1 + Z^{c(n)})}{(1 - Z)}$$

$$c(n) = 1,3,5,11,21,43,85,171,...*$$

**1, 3, 4, 5, 7, 9, 11, 12, 13, 15, 16, 17, 19, 20, 21, 23, 25, 27, 28, 29, 31, 33, 35, 36, 37, 39, 41, 43, 44, 45, 47, 48, 49, 51, 52, 53, 55, 57, 59, 60, 61, 63, 64, 65, 67, 68, 69, 71**

## ∑ C(n,k).C(2n+k,k-1)/n, k=1...n

**Réf.** FQ 11 123 73.

**HIS2 A3168** **Inverse fonctionnel algébrique** **Suite p-récurrente**

**HIS1** **Inverse ordinaire de A3169**

**L'inverse fonctionnel est rationnel.**

**Solution de** $\left(\dfrac{z}{(1+2z)(z+1)^2}\right)^{<-1>}$

**1, 1, 4, 21, 126, 818, 5594, 39693, 289510, 2157150, 16348960, 125642146, 976789620, 7668465964, 60708178054, 484093913917, 3884724864390**

---

## 2-line arrays

**Réf.** FQ 11 124 73; 14 232 76.

**HIS2 A3169** **Inverse fonctionnel algébrique** **Suite p-récurrente**

**HIS1** **Inverse ordinaire de A3168**

**Solution de** $\left(\dfrac{1+z}{3-2z^2+z^3}\right)^{<-1>}$

**1, 3, 14, 79, 494, 3294, 22952, 165127, 1217270, 9146746, 69799476, 539464358, 4214095612, 33218794236, 263908187100, 2110912146295, 16985386737830**

## Hex numbers

**Réf.** INOC 24 4550 85. AMM 95 701 88. GA88 18.

**HIS2 A3215** **Approximants de Padé**

**HIS1** **Fraction rationnelle**

$$\frac{1 + 4 z + z^{2}}{(1 - z)^{3}}$$

**1, 7, 19, 37, 61, 91, 127, 169, 217, 271, 331, 397, 469, 547, 631, 721, 817, 919, 1027, 1141, 1261, 1387, 1519, 1657, 1801, 1951, 2107, 2269, 2437, 2611, 2791, 2977**

---

## Even permutations of length n with no fixed points

**Réf.** AMM 79 394 72.

**HIS2 A3221** **Dérivée logarithmique** **Suite P-récurrente**

**HIS1** **exponentielle**

**a(n) = 3 n a(n - 2) + (n - 1) a(n - 1) + (3 n - 1) a(n - 3) + (n - 1) a(n - 4)**

$$\frac{4 - 6 z + 16 z^{2} - 13 z^{3} + 6 z^{4} - z^{5}}{2 (z - 1)^{4} \exp(z)}$$

**0, 0, 2, 3, 24, 130, 930, 7413, 66752, 667476, 7342290, 88107415, 1145396472, 16035550518, 240533257874, 3848532125865, 65425046139840, 1177650830516968**

**Réf.** DT76.
**HIS2 A3229** **Approximants de Padé**
**HIS1** **Fraction rationnelle**

$$\frac{1 + 2 z^2}{1 - z - 2 z^3}$$

**1, 1, 3, 5, 7, 13, 23, 37, 63, 109, 183, 309, 527, 893, 1511, 2565, 4351, 7373, 12503, 21205, 35951, 60957, 103367, 175269, 297183, 503917, 854455, 1448821, 2456655**

---

**Réf.** DT76.
**HIS2 A3230** **Approximants de Padé**
**HIS1** **Fraction rationnelle**

$$\frac{1}{(z - 1) (2 z - 1) (1 - z - 2 z^3)}$$

**1, 4, 11, 28, 67, 152, 335, 724, 1539, 3232, 6727, 13900, 28555, 58392, 118959, 241604, 489459, 989520, 1997015, 4024508, 8100699, 16289032, 32726655, 65705268, 131837763**

## Partially achiral planted trees

**Réf.** JRAM 278 334 75.

**HIS2 A3237** **Approximants de Padé** **conjecture faible**

**HIS1** **Fraction rationnelle**

$$\frac{z\,(1 - z^{2} - z^{3} - z^{4} + z^{5})}{1 - z - 2\,z^{2} + 3\,z^{5}}$$

**0, 1, 1, 2, 3, 6, 10, 19, 33, 62, 110, 204**

## Partially achiral trees

**Réf.** JRAM 278 334 75.

**HIS2 A3243** **Approximants de Padé** **conjecture faible**

**HIS1** **Fraction rationnelle**

$$\frac{1 - z^{2} - 2\,z^{3} - 8\,z^{4} + 7\,z^{5} + 4\,z^{6}}{1 - z - z^{2} - 2\,z^{3} - 6\,z^{4} + 14\,z^{5}}$$

**1, 1, 1, 2, 3, 6, 9, 19, 30, 61, 99, 208**

## Related to Fibonacci representations

**Réf.** FQ 11 386 73.

**HIS2 A3253** **Approximants de Padé** **conjecture seulement**

**HIS1** **Fraction rationnelle**

$$\frac{1 + z + z^{2} + z^{15} - z^{16}}{1 - z - z^{2} + z^{3}}$$

**1, 2, 4, 5, 7, 8, 10, 11, 13, 14, 16, 17, 19, 20, 22, 24, 25, 27, 28, 30, 31, 33, 34, 36, 37, 39, 40, 42, 43, 45, 46, 48, 49, 51, 52, 54, 55, 57, 58, 60, 62, 63, 65, 66, 68, 69, 71, 72**

## Woodall numbers

**Réf.** BR73 159.

**HIS2 A3261** **Approximants de Padé**

**HIS1** **Fraction rationnelle**

$$\frac{1 + 2 z - 4 z^{2}}{(1 - z) (2 z - 1)^{2}}$$

**1, 7, 23, 63, 159, 383, 895, 2047, 4607, 10239, 22527, 49151, 106495, 229375, 491519, 1048575, 2228223, 4718591, 9961471, 20971519, 44040191, 92274687**

**Réf.** BR72 120.

**HIS2 A3269** **Approximants de Padé**

**HIS1** **Fraction rationnelle**

$$\frac{1}{1 - z - z^4}$$

**1, 1, 1, 1, 2, 3, 4, 5, 7, 10, 14, 19, 26, 36, 50, 69, 95, 131, 181, 250, 345, 476, 657, 907, 1252, 1728, 2385, 3292, 4544, 6272, 8657, 11949, 16493, 22765, 31422, 43371, 59864**

---

## Key permutations of length n

**Réf.** CJN 14 152 71.

**HIS2 A3274** **Approximants de Padé**

**HIS1** **Fraction rationnelle**

$$\frac{1 - z + 3 z^2 - 2 z^3 + z^5}{(1 - z - z^3)(z - 1)^2}$$

**1, 2, 6, 12, 20, 34, 56, 88, 136, 208, 314, 470, 700, 1038, 1534, 2262, 3330, 4896, 7192, 10558, 15492, 22724, 33324, 48860, 71630, 105002, 153912, 225594, 330650**

## 4-line partitions of n decreasing across rows

**Réf.** MOC 26 1004 72.

**HIS2 A3292** **Euler**

**HIS1** **Produit infini**

$$\prod_{n\geq 1} \frac{1}{(1 - Z^n)^{c(n)}}$$

$$c(n) = 1,1,2,2,2,2,...$$

**1, 2, 4, 7, 11, 19, 29, 46, 70, 106, 156, 232, 334, 482, 686, 971, 1357, 1894, 2612, 3592, 4900, 6656, 8980, 12077, 16137, 21490, 28476, 37600, 49422, 64763, 84511**

## Planar partitions of n decreasing across rows

**Réf.** MOC 26 1004 72.

**HIS2 A3293** **Euler**

**HIS1** **Produit infini**

$$\prod_{n\geq 1} \frac{1}{(1 - Z^n)^{c(n)}}$$

$$c(n) = 1,1,2,2,3,3,4,4,5,5,6,6,...$$

**1, 2, 4, 7, 12, 21, 34, 56, 90, 143, 223, 348, 532, 811, 1224, 1834, 2725, 4031, 5914, 8638, 12540, 18116, 26035, 37262, 53070, 75292, 106377, 149738, 209980**

## Certain triangular arrays of integers

**Réf.** P4BC 112.

**HIS2 A3402** **Euler**

**HIS1** **Fraction rationnelle**

$$\frac{1}{(1-z)(1-z^2)(1-z^3)^2(1-z^4)(1-z^5)}$$

**1, 1, 2, 4, 6, 9, 14, 19, 27, 37, 49, 64, 84, 106, 134, 168, 207, 253, 309, 371, 445, 530, 626, 736, 863, 1003, 1163, 1343, 1543, 1766, 2017, 2291, 2597, 2935, 3305, 3712, 4161**

## Certain triangular arrays of integers

**Réf.** P4BC 118.

**HIS2 A3403** **Euler**

**HIS1** **Fraction rationnelle**

*** c(n) : suite finie.**

$$\prod_{n\geq 1} \frac{1}{(1-z^n)^{c(n)}}$$

$$c(n) = 1,1,2,2,2,1,1.*$$

**1, 1, 2, 4, 7, 11, 18, 27, 41, 60, 87, 122, 172, 235, 320, 430, 572, 751, 982, 1268, 1629, 2074, 2625, 3297, 4123, 5118, 6324, 7771, 9506, 11567, 14023, 16917, 20335**

## Connected ladder graphs with n nodes

**Réf.** DM 9 355 74.

**HIS2 A3409** **Recoupements algébrique** **Suite P-récurrente**

**HIS1**

$$\frac{6}{(1 - 4 z)^{1/2} \left(1 + (1 - 4 z)^{1/2}\right)}$$

**3, 9, 30, 105, 378, 1386, 5148, 19305**

---

**Réf.** rkg.

**HIS2 A3410** **Approximants de Padé**

**HIS1** **Fraction rationnelle**

$$\frac{(1 + z)(1 + z^2)}{1 + z + z^3}$$

**1, 2, 3, 5, 7, 10, 15, 22, 32, 47, 69, 101, 148, 217, 318, 466, 683, 1001, 1467, 2150, 3151, 4618, 6768, 9919, 14537, 21305, 31224, 45761, 67066, 98290, 144051, 211117**

**Réf.** rkg.
**HIS2 A3411** **Approximants de Padé**
**HIS1** **Fraction rationnelle**

$$\frac{z^4 + z^3 + z^2 + z + 1}{1 + z + z^4}$$

**1, 2, 3, 4, 6, 8, 11, 15, 21, 29, 40, 55, 76, 105, 145, 200, 276, 381, 526, 726, 1002, 1383, 1909, 2635, 3637, 5020, 6929, 9564, 13201, 18221, 25150, 34714, 47915, 66136**

---

**From a nim-like game**

**Réf.** rkg.
**HIS2 A3413** **Approximants de Padé**
**HIS1** **Fraction rationnelle**

$$\frac{(z^5 + z^3 + 1)\,(z^2 + z + 1)}{z^6 + z - 1}$$

**1, 2, 3, 4, 5, 7, 9, 12, 15, 19, 24, 31, 40, 52, 67, 86, 110, 141, 181, 233, 300, 386, 496, 637, 818, 1051, 1351, 1737, 2233, 2870, 3688, 4739, 6090, 7827, 10060, 12930**

## Continued fraction expansion of e = exp(1)

**Réf.** PE29 134.

**HIS2 A3417** **Approximants de Padé**

**HIS1** **Fraction rationnelle**

$$\frac{2 + z + 2 z^2 - 3 z^3 - z^4 + z^6}{(z - 1)^2 (z^2 + z + 1)^2}$$

**2, 1, 2, 1, 1, 4, 1, 1, 6, 1, 1, 8, 1, 1, 10, 1, 1, 12, 1, 1, 14, 1, 1, 16, 1, 1, 18, 1, 1, 20, 1, 1, 22, 1, 1, 24, 1, 1, 26, 1, 1, 28, 1, 1, 30, 1, 1, 32, 1, 1, 34, 1, 1, 36, 1, 1, 38, 1, 1, 40, 1, 1, 42**

---

## Hamiltonian circuits on n-octahedron

**Réf.** JCT B19 2 75.

**HIS2 A3436** **P-récurrences** **Suite P-récurrente**

**HIS1** **exponentielle (algébrique)**

**Une relation élémentaire existe avec A0806.**

$$a(n) = (2 n + 2) a(n - 1) - a(n - 3) + (- 2 n + 4) a(n - 2)$$

**1, 4, 31, 293, 3326, 44189, 673471, 11588884, 222304897, 4704612119, 108897613826, 2737023412199, 74236203425281, 2161288643251828**

## Dissections of a polygon

**Réf.** AEQ 18 387 78.

**HIS2 A3451** **Approximants de Padé**

**HIS1** **Fraction rationnelle**

$$\frac{z^2 - 2z - 1}{(z - 1)^4 (z + 1)^2}$$

**1, 4, 8, 16, 25, 40, 56, 80, 105, 140, 176, 224**

## Dissections of a polygon

**Réf.** AEQ 18 388 78.

**HIS2 A3453** **Approximants de Padé**

**HIS1** **Fraction rationnelle**

$$\frac{z^2 - z - 1}{(z - 1)^4 (z + 1)^2}$$

**1, 3, 6, 11, 17, 26, 36, 50, 65, 85, 106, 133**

## Bode numbers

**Réf.** SKY 43 281 72. MCL1.

**HIS2 A3461** **Approximants de Padé**

**HIS1** **Fraction rationnelle**

$$\frac{4 - 5z - 3z^2}{(2z - 1)(z - 1)}$$

**4, 7, 10, 16, 28, 52, 100, 196, 388, 772, 1540, 3076, 6148, 12292, 24580, 49156, 98308, 196612, 393220, 786436, 1572868, 3145732, 6291460, 12582916, 25165828**

---

**Réf.** RI89 60.

**HIS2 A3462** **Approximants de Padé**

**HIS1** **Fraction rationnelle**

$$\frac{1}{(1 - z)(1 - 3z)}$$

**1, 4, 13, 40, 121, 364, 1093, 3280, 9841, 29524, 88573, 265720, 797161, 2391484, 7174453, 21523360, 64570081, 193710244, 581130733, 1743392200, 5230176601**

## Minimal covers of an n-set

**Réf.** DM 5 249 73.

**HIS2 A3467** **P-récurrences** **Suite P-récurrente**

**HIS1** **Fraction rationnelle** **Formule de B. Salvy**

**(n - 1) (n - 2) a(n) = (n + 2) (5 n - 10) a(n - 1) + (n + 2) (- 4 n - 4) a(n - 2)**

$$1 + \frac{1}{(4 z - 1)^4} + \frac{3}{(z - 1)^4}$$

**5, 28, 190, 1340, 9065, 57512, 344316, 1966440, 10813935, 57672340, 299893594, 1526727748, 7633634645, 37580965520, 182536112120, 876173330832**

## Minimal covers of an n-set

**Réf.** DM 5 249 73.

**HIS2 A3468** **Approximants de Padé**

**HIS1** **Fraction rationnelle**

$$\frac{1}{(1 - 4 z) (1 - 5 z) (1 - 6 z) (1 - 7 z)}$$

**1, 22, 305, 3410, 33621, 305382, 2619625, 21554170, 171870941, 1337764142, 10216988145, 76862115330, 571247591461, 4203844925302, 30687029023865**

## Minimal covers of an n-set

**Réf.** DM 5 249 73.

**HIS2 A3469** **Approximants de Padé**

**HIS1** **Fraction rationnelle**

$$\frac{1 - z - z^2}{(2z - 1)^2 (1 - z)^3}$$

**1, 6, 22, 65, 171, 420, 988, 2259, 5065, 11198, 24498, 53157, 114583, 245640, 524152, 1113959, 2359125, 4980546, 10485550, 22019865, 46137091, 96468716**

---

**Réf.** PRSE 62 190 46. AS1 796. MFM 74 62 70 (divided by 2).

**HIS2 A3472** **Approximants de Padé**

**HIS1** **Fraction rationnelle**

$$\frac{1}{(1 - 2z)^5}$$

**1, 10, 60, 280, 1120, 4032, 13440, 42240, 126720, 366080, 1025024, 2795520, 7454720, 19496960, 50135040, 127008768, 317521920, 784465920, 1917583360**

**Réf.** DT76.
**HIS2 A3476** **Approximants de Padé**
**HIS1** **Fraction rationnelle**

$$\frac{1 + z + z^2}{1 - z - 2z^3}$$

**1, 2, 3, 5, 9, 15, 25, 43, 73, 123, 209, 355, 601, 1019, 1729, 2931, 4969, 8427, 14289, 24227, 41081, 69659, 118113, 200275, 339593, 575819, 976369, 1655555**

---

**Réf.** DT76.
**HIS2 A3477** **Approximants de Padé**
**HIS1** **Fraction rationnelle**

$$\frac{1}{(1 - 2z)(1 - z - 2z^3)(1 + z^2)}$$

**1, 3, 6, 14, 33, 71, 150, 318, 665, 1375, 2830, 5798, 11825, 24039, 48742, 98606, 199113, 401455, 808382, 1626038, 3267809, 6562295, 13169814, 26416318, 52962681**

**Réf.** DT76.
**HIS2 A3478** **Approximants de Padé**
**HIS1** **Fraction rationnelle**

$$\frac{1}{(1 - 2z)(1 - z - 2z^3)}$$

**1, 3, 7, 17, 39, 85, 183, 389, 815, 1693, 3495, 7173, 14655, 29837, 60567, 122645, 247855, 500061, 1007495, 2027493, 4076191, 8188333, 16437623, 32978613, 66132495**

**Réf.** DT76.
**HIS2 A3479** **Approximants de Padé**
**HIS1** **Fraction rationnelle**

$$\frac{1}{(1 - z)(1 - z - 2z^3)}$$

**1, 2, 3, 6, 11, 18, 31, 54, 91, 154, 263, 446, 755, 1282, 2175, 3686, 6251, 10602, 17975, 30478, 51683, 87634, 148591, 251958, 427227, 724410, 1228327, 2082782**

**Réf.** MOC 29 220 75. DM 75 95 89.
**HIS2 A3480** **Approximants de Padé**
**HIS1** **Fraction rationnelle**

$$\frac{(z-1)^2}{1-4z+2z^2}$$

**1, 2, 7, 24, 82, 280, 956, 3264, 11144, 38048, 129904, 443520, 1514272, 5170048, 17651648, 60266496, 205762688, 702517760, 2398545664, 8189147136, 27959497216**

---

**Réf.** DM 9 89 74.
**HIS2 A3481** **Approximants de Padé**
**HIS1** **Fraction rationnelle**

$$\frac{2+4z-z^2}{1-8z+8z^2-z^3}$$

**2, 20, 143, 986, 6764, 46367, 317810, 2178308, 14930351, 102334154, 701408732, 4807526975, 32951280098, 225851433716, 1548008755919**

**Réf.** DM 9 89 74.
**HIS2 A3482** **Approximants de Padé**
**HIS1** **Fraction rationnelle**

$$\frac{5 - z}{1 - 8z + 8z^2 - z^3}$$

**0, 5, 39, 272, 1869, 12815, 87840, 602069, 4126647, 28284464, 193864605, 1328767775, 9107509824, 62423800997, 427859097159, 2932589879120**

---

## Hurwitz-Radon function at powers of 2

**Réf.** LA73a 131.
**HIS2 A3485** **Approximants de Padé**
**HIS1** **Fraction rationnelle**

$$\frac{1 + z + 2z^2 + 4z^3}{(1 - z)(1 - z^4)}$$

**1, 2, 4, 8, 9, 10, 12, 16, 17, 18, 20, 24, 25, 26, 28, 32, 33, 34, 36, 40, 41, 42, 44, 48, 49, 50, 52, 56, 57, 58, 60, 64, 65, 66, 68, 72, 73, 74, 76, 80, 81, 82, 84, 88, 89, 90, 92, 96**

**Réf.** B1 198. MMAG 48 209 75.
**HIS2 A3499** **Approximants de Padé**
**HIS1** **Fraction rationnelle**

$$\frac{2 - 6z}{1 - 6z + z^2}$$

**2, 6, 34, 198, 1154, 6726, 39202, 228486, 1331714, 7761798, 45239074, 263672646, 1536796802, 8957108166, 52205852194, 304278004998, 1773462177794**

---

**Réf.** FQ 11 29 73. MMAG 48 209 75.
**HIS2 A3500** **Approximants de Padé**
**HIS1** **Fraction rationnelle**

$$\frac{2 - 4z}{1 - 4z + z^2}$$

**2, 4, 14, 52, 194, 724, 2702, 10084, 37634, 140452, 524174, 1956244, 7300802, 27246964, 101687054, 379501252, 1416317954, 5285770564, 19726764302**

**Réf.** MMAG 48 209 75.

**HIS2 A3501** **Approximants de Padé**

**HIS1** **Fraction rationnelle**

$$\frac{2 - 5 z}{1 - 5 z + z^2}$$

**2, 5, 23, 110, 527, 2525, 12098, 57965, 277727, 1330670, 6375623, 30547445, 146361602, 701260565, 3359941223, 16098445550, 77132286527, 369562987085**

---

**Binomial coefficients C (2n + 1, n - 2)**

**Réf.** AS1 828.

**HIS2 A3516** **Hypergéométrique** **Suite P-récurrente**

**HIS1** **algébrique**

$_2F_1([3, 7/2], [6], 4 z)$

$$\frac{32}{(1 - 4 z)^{1/2} (1 + (1 - 4 z)^{1/2})^5}$$

**1, 7, 36, 165, 715, 3003, 12376, 50388, 203490, 817190, 3268760, 13037895, 51895935, 206253075, 818809200, 3247943160, 12875774670, 51021117810**

## Binomial coefficients 6C(2n+1,n-2)/(n+4)

**Réf.** FQ 14 397 76. DM 14 84 76.

**HIS2 A3517** **Hypergéométrique** **Suite P-récurrente**

**HIS1** **algébrique**

**2F1([3, 7/2], [7], 4 z)**

$$\frac{64}{\left(1 + (1 - 4z)^{1/2}\right)^6}$$

**1, 6, 27, 110, 429, 1638, 6188, 23256, 87210, 326876, 1225785, 4601610, 17298645, 65132550, 245642760, 927983760, 3511574910, 13309856820, 50528160150**

## Binomial coefficients 8C(2n+1,n-3)/(n+5)

**Réf.** FQ 14 397 76. DM 14 84 76.

**HIS2 A3518** **Hypergéométrique** **Suite P-récurrente**

**HIS1** **algébrique**

**2F1([9/2 , 4] , [9], 4 z)**

$$\frac{256\, z}{\left(1 + (1 - 4z)^{1/2}\right)^8}$$

**1, 8, 44, 208, 910, 3808, 15504, 62016, 245157, 961400, 3749460, 14567280, 56448210, 218349120, 843621600, 3257112960, 12570420330, 48507033744**

## Binomial coefficients 10C(2n+1,n-4)/(n+6)

**Réf.** FQ 14 397 76.

**HIS2 A3519** **Hypergéométrique** **Suite P-récurrente**

**HIS1** **algébrique**

$_2F_1([11/2,5],[11],4z)$

$$\frac{1024}{\left(1 + (1 - 4z)^{1/2}\right)^{10}}$$

**1, 10, 65, 350, 1700, 7752, 33915, 144210, 600875, 2466750, 10015005, 40320150, 161280600, 641886000, 2544619500, 10056336264, 39645171810**

---

**Réf.** BR72 119. FQ 14 38 76.

**HIS2 A3520** **Approximants de Padé**

**HIS1** **Fraction rationnelle**

$$\frac{1}{(1 - z^2 - z^3)(1 - z + z^2)}$$

**1, 1, 1, 1, 1, 2, 3, 4, 5, 6, 8, 11, 15, 20, 26, 34, 45, 60, 80, 106, 140, 185, 245, 325, 431, 571, 756, 1001, 1326, 1757, 2328, 3084, 4085, 5411, 7168, 9496, 12580, 16665, 22076, 29244**

**Réf.** BR72 113.
**HIS2 A3522** **Approximants de Padé**
**HIS1** **Fraction rationnelle**

$$\frac{(z-1)^2}{1 - 3z + 3z^2 - z^3 - z^4}$$

**1, 1, 1, 1, 2, 5, 11, 21, 37, 64, 113, 205, 377, 693, 1266, 2301, 4175, 7581, 13785, 25088, 45665, 83097, 151169, 274969, 500162, 909845, 1655187, 3011157, 5477917, 9965312**

---

**Réf.** JCT A29 122 80. MOC 37 479 81.
**HIS2 A4004** **Approximants de Padé**
**HIS1** **Fraction rationnelle**

$$\frac{z(1+3z)}{(1-9z)(z-1)^2}$$

**0, 1, 14, 135, 1228, 11069, 99642, 896803, 8071256, 72641337, 653772070, 5883948671, 52955538084, 476599842805, 4289398585298, 38604587267739, 347441285409712, 3126971568687473**

## Coefficients of elliptic function sn

**Réf.** CA95 56. TM93 4 92. JCT A29 122 80. MOC 37 480 81.

**HIS2 A4005** **Approximants de Padé**

**HIS1** **Fraction rationnelle**

$$\frac{1 + 89 z - 69 z^{2} - 405 z^{3}}{(1 - z)^{3} (1 - 9 z)^{2} (1 - 25 z)}$$

**1, 135, 5478, 165826, 4494351, 116294673, 2949965020, 74197080276, 1859539731885, 46535238000235, 1163848723925346, 29100851707716150, 727566807977891803**

---

## Theta series of square lattice

**Réf.** SPLAG 106.

**HIS2 A4018** **Euler**

**HIS1** **Produit infini**

*** Le motif [4, -6, 4, -2] est périodique**

$$\prod_{n \geq 1} \frac{1}{(1 - z^{n})^{c(n)}}$$

$$c(n) = 4,-6,4,-2,...*$$

**1, 4, 4, 0, 4, 8, 0, 0, 4, 4, 8, 0, 0, 8, 0, 0, 4, 8, 4, 0, 8, 0, 0, 0, 0, 12, 8, 0, 0, 8, 0, 0, 4, 0, 8, 0, 4, 8, 0, 0, 8, 8, 0, 0, 0, 8, 0, 0, 0, 4, 12, 0, 8, 8, 0, 0, 0, 0, 8, 0, 0, 8, 0, 0, 4, 16, 0, 0, 8, 0**

## Theta series of square lattice w.r.t. edge.

**Réf.** SPLAG 106.

**HIS2 A4020** **Euler**

**HIS1** **Produit infini**

*** Le motif [2, -2] est périodique**

$$\prod_{n \geq 1} \frac{1}{(1 - z^n)^{c(n)}}$$

$$c(n) = 2, -2, ...*$$

**2, 4, 2, 4, 4, 0, 6, 4, 0, 4, 4, 4, 2, 4, 0, 4, 8, 0, 4, 0, 2, 8, 4, 0, 4, 4, 0, 4, 4, 4, 2, 8, 0, 0, 4, 0, 8, 4, 4, 4, 0, 0, 6, 4, 0, 4, 8, 0, 4, 4, 0, 8, 0, 0, 0, 8, 6, 4, 4, 0, 4, 4, 0, 0, 4, 4, 8, 4**

## Theta series of b.c.c. lattice w.r.t. deep hole

**Réf.** JCP 83 6532 85.

**HIS2 A4024** **Euler**

**HIS1** **Produit infini**

*** Le motif [1, 1, 1, -3] est périodique**

$$\prod_{n \geq 1} \frac{1}{(1 - z^n)^{c(n)}}$$

$$c(n) = 1, 1, 1, -3, ...*$$

**4, 4, 8, 12, 4, 12, 12, 12, 16, 16, 8, 8, 28, 12, 20, 24, 8, 16, 28, 12, 16, 28, 20, 32, 20, 16, 16, 32, 20, 24, 28, 8, 36, 44, 12, 32, 36, 16, 24, 20, 28, 20, 56, 28, 16, 40, 20, 40, 44, 12**

## Theta series of b.c.c. lattice w.r.t. long edge

**Réf.** JCP 6532 85.

**HIS2 A4025** **Euler**

**HIS1** **Produit infini**

*** Le motif [2, -3, 2, 1, 2, -3, 2, -3] est périodique**

$$\prod_{n\geq 1} \frac{1}{(1-z^n)^{c(n)}}$$

$$c(n) = 2,-3,2,1,2,-3,2,-3,...*$$

**2, 4, 0, 0, 8, 8, 0, 0, 10, 8, 0, 0, 8, 16, 0, 0, 16, 12, 0, 0, 16, 8, 0, 0, 10, 24, 0, 0, 24, 16, 0, 0, 16, 16, 0, 0, 8, 24, 0, 0, 32, 16, 0, 0, 24, 16, 0, 0, 18, 28, 0, 0, 24, 32, 0, 0, 16, 8, 0**

---

**Réf.** AMM 87 206 80.

**HIS2 A4116** **Approximants de Padé**

**HIS1** **Fraction rationnelle**

$$\frac{z^3 - z - 1}{(1+z)(z-1)^3}$$

**1, 3, 6, 9, 13, 17, 22, 27, 33, 39, 46, 53, 61, 69, 78, 87, 97, 107, 118, 129, 141, 153, 166, 179, 193, 207, 222, 237, 253, 269, 286, 303, 321, 339, 358, 377, 397, 417, 438, 459**

**Réf.** MOC 30 660 76.
**HIS2 A4119** **Approximants de Padé**
**HIS1** **Fraction rationnelle**

$$\frac{1 + z - 3z^2}{(2z - 1)(z - 1)}$$

**1, 4, 7, 13, 25, 49, 97, 193, 385, 769, 1537, 3073, 6145, 12289, 24577, 49153, 98305, 196609, 393217, 786433, 1572865, 3145729, 6291457, 12582913, 25165825**

---

**Réf.** SIAR 12 296 70.
**HIS2 A4120** **Approximants de Padé**
**HIS1** **Fraction rationnelle**

$$\frac{1 + z - z^5}{(1 - z)^3}$$

**1, 4, 9, 16, 25, 35, 46, 58, 71, 85, 100, 116, 133, 151, 170, 190, 211**

## Postage stamp problem

**Réf.** SIAA 1 383 80.

**HIS2 A4129** **Approximants de Padé** **Conjecture**

**HIS1** **Fraction rationnelle**

$$\frac{(z^4 + z^3 + 2z^2 + 2z + 1)(z^2 + z + 1)}{(z - 1)(z^5 + z^4 + z^3 - z - 1)}$$

**1, 3, 6, 9, 13, 17, 22, 27, 33, 40, 47, 56, 65**

## A counter moving problem

**Réf.** BA62 38.

**HIS2 A4138** **Approximants de Padé**

**HIS1** **Fraction rationnelle**

$$\frac{1 - z^2 + 4z^3 - 2z^4}{(z - 1)(2z^4 - z^3 + z^2 + z - 1)}$$

**1, 2, 3, 8, 13, 24, 37, 66, 107, 186, 303, 516, 849, 1436, 2377, 3998, 6639, 11134, 18531, 31024, 51701, 86464, 144205, 241018, 402163, 671906, 1121463, 1873244**

## Alternate Lucas numbers - 2

**Réf.** FQ 13 51 75.

**HIS2 A4146** **Approximants de Padé** **Suite P-récurrente**

**HIS1** **fraction rationnelle** **Suite corrigée au 12è terme.**

$$\frac{1 + z}{1 - 4z + 4z^2 - z^3}$$

**1, 5, 16, 45, 121, 320, 841, 2205, 5776, 15125, 39601, 103680*, 271441, 710645, 1860496, 4870845, 12752041, 33385280, 87403801, 228826125, 599074576**

## Generalized Catalan numbers

**Réf.** DM 26 264 79. JCT B29 89 80.

**HIS2 A4148** **LLL** **Suite P-récurrente**

**HIS1** **algébrique**

**(n + 2) a(n) = (4 - n) a(n - 4) + (2 n + 1) a(n - 1)**
**+ (n - 1) a(n - 2) + (2 n - 5) a(n - 3)**

$$\frac{1 - z - z^2 - (1 - 2z - z^2 - 2z^3 + z^4)^{1/2}}{2z^3}$$

**1, 1, 2, 4, 8, 17, 37, 82, 185, 423, 978, 2283, 5373, 12735, 30372, 72832, 175502, 424748, 1032004, 2516347**

## Related to symmetric groups

**Réf.** DM 21 320 78.

**HIS2 A4211** **équations différentielles** **Formule de B. Salvy**

**HIS1** **exponentielle**

$$\exp(1/2 \exp(2 z) + 2 z - 1/2)$$

**1, 3, 11, 49, 257, 1539, 10299, 75905**

## Related to symmetric groups

**Réf.** DM 21 320 78.

**HIS2 A4212** **équations différentielles** **Formule de B. Salvy**

**HIS1** **exponentielle**

$$\exp(1/3 \exp(3 z) + 3 z - 1/3)$$

**1, 4, 19, 109, 742, 5815, 51193, 498118**

## Related to symmetric groups

**Réf.** DM 21 320 78.

**HIS2 A4213** **équations différentielles** **Formule de B. Salvy**

**HIS1** **exponentielle**

$$\exp(1/4 \exp(4 z) + 4 z - 1/4)$$

**1, 5, 29, 201, 1657, 15821, 170389, 2032785**

## Pythagoras theorem generalized

**Réf.** BU71 75.

**HIS2 A4253** **Approximants de Padé**

**HIS1** **Fraction rationnelle**

$$\frac{1 - z}{1 - 5 z + z^2}$$

**1, 4, 19, 91, 436, 2089, 10009, 47956, 229771, 1100899, 5274724, 25272721, 121088881, 580171684, 2779769539, 13318676011, 63813610516, 305749376569**

## Pythagoras theorem generalized

**Réf.** BU71 75.

**HIS2 A4254** **Approximants de Padé**

**HIS1** **Fraction rationnelle**

$$\frac{1}{1 - 5z + z^2}$$

**1, 5, 24, 115, 551, 2640, 12649, 60605, 290376, 1391275, 6665999, 31938720, 153027601, 733199285, 3512968824, 16831644835, 80645255351, 386394631920**

---

**Réf.** dsk.

**HIS2 A4255** **Approximants de Padé**

**HIS1** **Fraction rationnelle**

$$\frac{1 - 2z + 4z^2}{(1 - z)^5}$$

**1, 3, 9, 25, 60, 126, 238, 414, 675, 1045, 1551, 2223, 3094, 4200, 5580, 7276, 9333, 11799, 14725, 18165, 22176, 26818, 32154, 38250, 45175, 53001, 61803, 71659**

**Réf.** JCT B21 75 76.

**HIS2 A4303** **LLL** **Suite P-récurrente**

**HIS1** **algébrique**

**(n + 1) a(n) = 68 n a(n - 5) - 16 n a(n - 6) + (11 n - 2) a(n - 1)**
**+ (- 47 n + 61) a(n - 2) + (101 n - 240) a(n - 3)**
**+ (- 116 n + 398) a(n - 4) - 304 a(n - 5) + 88 a(n - 6)**

$$\frac{-1/2\,(-1 + 10 z - 42 z^2 + 98 z^3 - 137 z^4 + 112 z^5 - 48 z^6 + 8 z^7)}{(z^2 (2 z - 1)^2 (z - 1)^4)} + \frac{(-(-1 + 4 z)(2 z - 1)^4 (z - 1)^8)^{1/2})}{(z^2 (2 z - 1)^2 (z - 1)^4)}$$

**1, 1, 1, 3, 16, 75, 309, 1183, 4360, 15783, 56750, 203929, 734722, 2658071, 9662093, 35292151, 129513736, 477376575, 1766738922, 6563071865, 24464169890**

---

## Davenport-Schinzel numbers

**Réf.** ARS 1 47 76. UPNT E20.

**HIS2 A5004** **Approximants de Padé** **Conjecture**

**HIS1** **Fraction rationnelle**

$$\frac{(z^3 - z^2 + z + 1)(z^2 + z + 1)}{(1 + z)(z - 1)^2}$$

**1, 3, 5, 8, 10, 14, 16, 20, 22, 26**

**Related to symmetric groups**

**Réf.** DM 21 320 78.

**HIS2 A5011** **équations différentielles exponentielle** **Formule de B. Salvy**
**HIS1**

$$\exp(1/5 \exp(5 z) + 5 z - 1/5)$$

**1, 6, 41, 331, 3176, 35451, 447981, 6282416**

---

**Related to symmetric groups**

**Réf.** DM 21 320 78.

**HIS2 A5012** **équations différentielles exponentielle** **Formule de B. Salvy**
**HIS1**

$$\exp(1/6 \exp(6 z) + 6 z - 1/6)$$

**1, 7, 55, 505, 5497, 69823, 1007407, 16157905**

**Réf.** LNM 748 57 79.

**HIS2 A5013** **Approximants de Padé**

**HIS1** **Fraction rationnelle**

$$\frac{z^2 + z + 1}{(z^2 - z - 1)(z^2 + z - 1)}$$

**0, 1, 1, 4, 3, 11, 8, 29, 21, 76, 55, 199, 144, 521, 377, 1364, 987, 3571, 2584, 9349, 6765, 24476, 17711, 64079, 46368, 167761, 121393, 439204, 317811, 1149851, 832040**

---

## Random walks

**Réf.** DM 17 44 77.

**HIS2 A5021** **Approximants de Padé**

**HIS1** **Fraction rationnelle**

$$\frac{(1 - z)(z - 5)}{1 - 5z + 6z^2 - z^3}$$

**5, 19, 66, 221, 728, 2380, 7753, 25213, 81927, 266110, 864201, 2806272, 9112264, 29587889, 96072133, 311945595, 1012883066, 3288813893, 10678716664**

## Random walks

**Réf.** DM 17 44 77. TCS 9 105 79.

**HIS2 A5022** **Approximants de Padé**

**HIS1** **Fraction rationnelle**

$$\frac{1}{(1 - 2z)(1 - 4z + 2z^2)}$$

**1, 6, 26, 100, 364, 1288, 4488, 15504, 53296, 182688, 625184, 2137408, 7303360, 24946816, 85196928, 290926848, 993379072, 3391793664, 11580678656, 39539651584**

## Random walks

**Réf.** DM 17 44 77.

**HIS2 A5023** **Approximants de Padé**

**HIS1** **Fraction rationnelle**

$$\frac{7 - 15z + 10z^2 - z^3}{(1 - z)(z^3 - 9z^2 + 6z - 1)}$$

**7, 34, 143, 560, 2108, 7752, 28101, 100947, 360526, 1282735, 4552624, 16131656, 57099056, 201962057, 714012495, 2523515514, 8916942687, 31504028992**

## Random walks

**Réf.** DM 17 44 77.

**HIS2 A5024** **Approximants de Padé**

**HIS1** **Fraction rationnelle**

$$\frac{8 - 21 z + 20 z^{2} - 5 z^{3}}{(5 z^{2} - 5 z + 1) (1 - 3 z + z^{2})}$$

**8, 43, 196, 820, 3264, 12597, 47652, 177859, 657800, 2417416, 8844448, 32256553, 117378336, 426440955, 1547491404, 5610955132, 20332248992**

## Random walks

**Réf.** DM 17 44 77.

**HIS2 A5025** **Approximants de Padé**

**HIS1** **Fraction rationnelle**

$$\frac{9 - 28 z + 35 z^{2} - 15 z^{3} + z^{4}}{1 - 9 z + 28 z^{2} - 35 z^{3} + 15 z^{4} - z^{5}}$$

**9, 53, 260, 1156, 4845, 19551, 76912, 297275, 1134705, 4292145, 16128061, 60304951, 224660626, 834641671, 3094322026, 11453607152, 42344301686**

**Réf.** JCT A23 293 77. JCP 67 5027 77. TAMS 272 406 82.
**HIS2 A5043** **LLL** **Suite P-récurrente**
**HIS1** **algébrique**
**(n + 2) a(n) = 2 n a(n - 1) + 3 n a(n - 2)**

$$\frac{1 - z - 2 z^2 - (1 - 2 z - 3 z^2)^{1/2}}{2 (z^3 + z^4)}$$

**0, 1, 1, 3, 6, 15, 36, 91, 232, 603, 1585, 4213, 11298, 30537, 83097, 227475, 625992, 1730787, 4805595, 13393689, 37458330, 105089229, 295673994, 834086421**

---

**Réf.** AMM 86 477 79; 86 687 79.
**HIS2 A5044** **Approximants de Padé**
**HIS1** **Fraction rationnelle**

$$\frac{1}{(1 + z^2) (z^2 + z + 1) (1 + z)^2 (z - 1)^3}$$

**1, 0, 1, 1, 2, 1, 3, 2, 4, 3, 5, 4, 7, 5, 8, 7, 10, 8, 12, 10, 14, 12, 16, 14, 19, 16, 21, 19, 24, 21, 27, 24, 30, 27, 33, 30, 37, 33, 40, 37, 44, 40, 48, 44, 52, 48, 56, 52, 61, 56, 65, 61, 70, 65**

## 3 times 3 matrices with row and column sums n

**Réf.** MO78. NAMS 26 A-27 (763-05-13) 79.

**HIS2 A5045** **Approximants de Padé**

**HIS1** **Fraction rationnelle**

$$\frac{z^6 - z^5 + z^3 - z - 1}{(1 + z^2)(z^2 + z + 1)(1 + z)^2 (z - 1)^5}$$

**1, 3, 6, 10, 17, 25, 37, 51, 70, 92, 121, 153, 194, 240, 296, 358, 433, 515, 612, 718, 841, 975, 1129, 1295, 1484, 1688, 1917, 2163, 2438, 2732, 3058, 3406, 3789, 4197, 4644**

## Minimal determinant of n-dimensional norm 3 lattice

**Réf.** SPLAG 180.

**HIS2 A5103** **Approximants de Padé** **Conjecture**

**HIS1** **Fraction rationnelle**

$$\frac{1 + z + 2 z^2 + 2 z^3 + 6 z^4}{1 - 2 z + 2 z^3}$$

**1, 3, 8, 16, 32, 48, 64, 64**

**Réf.** clm.
**HIS2 A5126** **Approximants de Padé**
**HIS1** **Fraction rationnelle**

$$\frac{2 - 4 z + z^2}{(1 - 2 z) (z - 1)^2}$$

**2, 4, 7, 12, 21, 38, 71, 136, 265, 522, 1035, 2060, 4109, 8206, 16399, 32784, 65553, 131090, 262163, 524308, 1048597, 2097174, 4194327, 8388632, 16777241, 33554458, 67108891**

---

**Réf.** CACM 23 704 76. LNM 829 122 80. MBIO 54 8 81.
**HIS2 A5172** **équations différentielles** **Formule de B. Salvy**
**HIS1** **exponentielle**

$$- 1/2 - W(- 1/2 \exp(z - 1/2))$$

**1, 4, 32, 416, 7552, 176128, 5018624, 168968192, 6563282944, 288909131776, 14212910809088, 772776684683264, 46017323176296448, 2978458881388183550**

## Trees of subsets of an n-set

**Réf.** MBIO 54 9 81.

**HIS2 A5173** **Approximants de Padé**

**HIS1** **Fraction rationnelle**

$$\frac{z\,(1 + 6\,z)}{(1 - z)\,(1 + 2\,z)\,(1 + 3\,z)}$$

**0, 1, 12, 61, 240, 841, 2772, 8821, 27480, 84481, 257532, 780781, 2358720, 7108921, 21392292, 64307941, 193185960, 580082161, 1741295052, 5225982301, 15682141200**

## Trees of subsets of an n-set

**Réf.** MBIO 54 9 81.

**HIS2 A5174** **Approximants de Padé**

**HIS1** **Fraction rationnelle**

$$\frac{2\,z^{2}\,(5 + 12\,z)}{(1 - z)\,(1 + 2\,z)\,(1 + 3\,z)\,(1 - 4\,z)}$$

**0, 0, 10, 124, 890, 5060, 25410, 118524, 527530, 2276020, 9613010, 40001324, 164698170, 672961380, 2734531810, 11066546524, 44652164810, 179768037140**

## Trees of subsets of an n-set

**Réf.** MBIO 54 9 81.

**HIS2 A5175** **Approximants de Padé**

**HIS1** **Fraction rationnelle**

$$\frac{z^2 (3 + 86 z + 120 z^2)}{(1 - z) (1 + 2 z) (1 + 3 z) (1 - 4 z) (1 - 5 z)}$$

**0, 0, 3, 131, 1830, 16990, 127953, 851361, 5231460, 30459980, 170761503, 931484191, 4979773890, 26223530970, 136522672653, 704553794621, 3611494269120, 18415268221960**

---

**Réf.** MMAG 63 15 90.

**HIS2 A5183** **Approximants de Padé**

**HIS1** **Fraction rationnelle**

$$\frac{1 - 3 z + 3 z^2}{(z - 1) (2 z - 1)^2}$$

**1, 2, 5, 13, 33, 81, 193, 449, 1025, 2305, 5121, 11265, 24577, 53249, 114689, 245761, 524289, 1114113, 2359297, 4980737, 10485761, 22020097, 46137345, 96468993, 201326593**

## (F(2n)+F(n+1))/2, where F(n) is a Fibonacci number

**Réf.** CJN 25 391 82.

**HIS2 A5207** **Approximants de Padé**

**HIS1** **Fraction rationnelle**

$$\frac{z^3 - z^2 - 2z + 1}{(1 - 3z + z^2)(1 - z - z^2)}$$

**1, 2, 4, 9, 21, 51, 127, 322, 826, 2135, 5545, 14445, 37701, 98514, 257608, 673933, 1763581, 4615823, 12082291, 31628466, 82798926, 216761547, 567474769, 1485645049**

## n-bead necklaces with 4 red beads

**Réf.** JAuMS 33 12 82. AJMG 22 5231 85.

**HIS2 A5232** **Approximants de Padé**

**HIS1** **Fraction rationnelle**

$$\frac{z^7 - 2z^6 + 2z^4 - 2z^3 + 2z^2 - z - 1}{(z^2 + 1)(z + 1)^2(1 - z)^4}$$

**1, 3, 4, 8, 10, 16, 20, 29, 35, 47, 56, 72, 84, 104, 120, 145, 165, 195, 220, 256, 286, 328, 364, 413, 455, 511, 560, 624, 680, 752, 816, 897, 969, 1059, 1140, 1240, 1330, 1440, 1540, 1661**

**Réf.** MAG 69 263 85.
**HIS2 A5246** **Approximants de Padé**
**HIS1** **Fraction rationnelle**

$$\frac{1 + z - 2 z^2 - z^3}{1 - 4 z^2 + z^4}$$

**1, 1, 2, 3, 7, 11, 26, 41, 97, 153, 362, 571, 1351, 2131, 5042, 7953, 18817, 29681, 70226, 110771, 262087, 413403, 978122, 1542841, 3650401, 5757961, 13623482, 21489003, 50843527**

---

**Réf.** MAG 69 264 85.
**HIS2 A5247** **Approximants de Padé**
**HIS1** **Fraction rationnelle**

$$\frac{(1 + z) (1 + z - 3 z^2)}{(z^2 - z - 1) (1 - z - z^2)}$$

**1, 2, 1, 3, 2, 7, 5, 18, 13, 47, 34, 123, 89, 322, 233, 843, 610, 2207, 1597, 5778, 4181, 15127, 10946, 39603, 28657, 103682, 75025, 271443, 196418, 710647, 514229, 1860498, 1346269**

**Réf.** FQ 9 284 71. MMAG 48 209 75. MAG 69 264 85.
**HIS2 A5248** **Approximants de Padé**
**HIS1** **Fraction rationnelle**

$$\frac{2 - 3z}{1 - 3z + z^2}$$

**2, 3, 7, 18, 47, 123, 322, 843, 2207, 5778, 15127, 39603, 103682, 271443, 710647, 1860498, 4870847, 12752043, 33385282, 87403803, 228826127, 599074578, 1568397607, 4106118243**

---

**Réf.** BR72 112. FQ 16 85 78. LAA 62 113 84.
**HIS2 A5251** **Approximants de Padé**
**HIS1** **Fraction rationnelle**

$$\frac{z - 1}{z^3 - z^2 + 2z - 1}$$

**1, 1, 1, 2, 4, 7, 12, 21, 37, 65, 114, 200, 351, 616, 1081, 1897, 3329, 5842, 10252, 17991, 31572, 55405, 97229, 170625, 299426, 525456, 922111, 1618192, 2839729, 4983377, 8745217**

**Réf.** FQ 7 341 69; 16 85 78.
**HIS2 A5252** **Approximants de Padé**
**HIS1** **Fraction rationnelle**
**$\sum$C(n-2k,2k), k=0...n**

$$\frac{z - 1}{(1 - z + z^2)(-1 + z + z^2)}$$

**1, 1, 1, 1, 2, 4, 7, 11, 17, 27, 44, 72, 117, 189, 305, 493, 798, 1292, 2091, 3383, 5473, 8855, 14328, 23184, 37513, 60697, 98209, 158905, 257114, 416020, 673135, 1089155, 1762289**

---

**Binary words not containing ..01110...**

**Réf.** FQ 16 85 78.
**HIS2 A5253** **Approximants de Padé**
**HIS1** **Fraction rationnelle**

$$\frac{1 - z + z^4}{1 - 2z + z^2 - z^5}$$

**1, 1, 1, 1, 2, 4, 7, 11, 16, 23, 34, 52, 81, 126, 194, 296, 450, 685, 1046, 1601, 2452, 3753, 5739, 8771, 13404, 20489, 31327, 47904, 73252, 112004, 171245, 261813, 400285**

## Apéry numbers

**Réf.** AST 61 12 79. JNT 25 201 87.

**HIS2 A5258** **P-récurrences** **Suite P-récurrente**

**HIS1**

$$(n-1)^2 a(n) = (n^2 - 4n + 4)\, a(n-2) + (11 n^2 - 33 n + 25)\, a(n-1)$$

**1, 3, 19, 147, 1251, 11253, 104959, 1004307, 9793891, 96918753, 970336269, 9807518757, 99912156111, 1024622952993, 10567623342519, 109527728400147**

## Apéry numbers

**Réf.** AST 61 13 79. JNT 25 201 87.

**HIS2 A5259** **P-récurrences** **Suite P-récurrente**

**HIS1**

$$(n-1)^3 a(n) = (-n^3 + 6 n^2 - 12 n + 8)\, a(n-2) + (34 n^3 - 153 n^2 + 231 n - 117)\, a(n-1)$$

**1, 5, 73, 1445, 33001, 819005, 21460825, 584307365, 16367912425, 468690849005, 13657436403073, 403676083788125, 12073365010564729, 364713572395983725**

**Réf.** JNT 25 201 87.

**HIS2 A5260** **P-récurrences** **Suite P-récurrente**

**HIS1**

**∑C(n,k)^4 , k=0...n**

$$(n - 1)^3 a(n) =$$

$$+ (12 n^3 - 54 n^2 + 82 n - 42) a(n - 1)$$

$$(64 n^3 - 384 n^2 + 764 n - 504) a(n - 2)$$

**1, 2, 18, 164, 1810, 21252, 263844, 3395016, 44916498, 607041380, 8345319268, 116335834056, 1640651321764, 23365271704712, 335556407724360, 4854133484555664**

---

**Réf.** CRUX 13 331 87.

**HIS2 A5262** **Approximants de Padé**

**HIS1** **Fraction rationnelle**

$$\frac{1 + z^2 + 4 z^3}{(1 + z) (2 z - 1) (1 - z)^2}$$

**1, 3, 9, 25, 59, 131, 277, 573, 1167, 2359, 4745, 9521, 19075, 38187, 76413, 152869, 305783, 611615, 1223281, 2446617, 4893291, 9786643, 19573349, 39146765, 78293599**

## Greg trees

**Réf.** MANU 34 127 90.

**HIS2 A5263** **équations différentielles** **Formule de B. Salvy**

**HIS1** **exponentielle**

$$1/4 - 1/4\,(2 + 2\,W(-\exp(-1/2)\,(1/2 + 1/2\,z)))^2$$

**1, 1, 4, 32, 396, 6692, 143816, 3756104, 115553024, 4093236352, 164098040448, 7345463787136**

## From Euclid's proof

**Réf.** SZ 27 31 78. LNM 829 122 80. MANU 34 127 90.

**HIS2 A5264** **Inverse fonctionnel**

**HIS1** **exponentielle** **f.g. exponentielle**

**L'inverse est (1+2 z-exp(z))/exp(z)**

$$-W(-\exp(-1/2)\,(1/2 + 1/2\,z)) - 1/2$$

**1, 3, 22, 262, 4336, 91984, 2381408, 72800928, 2566606784, 102515201984, 4575271116032, 225649908491264, 12187240730230208, 715392567595384832**

**Réf.** NET 96. MMAG 61 28 88. rkg.

**HIS2 A5286** **Approximants de Padé**

**HIS1** **Fraction rationnelle**

$$\frac{1 + 2z - 3z^2 + z^3}{(z - 1)^4}$$

**1, 6, 15, 29, 49, 76, 111, 155, 209, 274, 351, 441, 545, 664, 799, 951, 1121, 1310, 1519, 1749, 2001, 2276, 2575, 2899, 3249, 3626, 4031, 4465, 4929, 5424, 5951, 6511, 7105, 7734**

## Permutations by inversions

**Réf.** NET 96. DKB 241. MMAG 61 28 88. rkg.

**HIS2 A5287** **Approximants de Padé**

**HIS1** **Fraction rationnelle**

$$\frac{5 - 5z - z^2 - 3z^3 - z^4}{(1 - z)^5}$$

**5, 20, 49, 98, 174, 285, 440, 649, 923, 1274, 1715, 2260, 2924, 3723, 4674, 5795, 7105**

## Permutations by inversions

**Réf.** NET 96. DKB 241. MMAG 61 28 88. rkg.

**HIS2 A5288** **Approximants de Padé**

**HIS1** **Fraction rationnelle**

$$\frac{3 + 4z - 16z^2 + 13z^3 - z^4 - 3z^5 + z^6}{(z-1)^6}$$

**3, 22, 71, 169, 343, 628, 1068, 1717, 2640, 3914, 5629, 7889, 10813, 14536, 19210, 25005, 32110**

## Graphs on n nodes with 3 cliques

**Réf.** AMM 80 1124 73; 82 997 75. JLMS 8 97 74. rkg.

**HIS2 A5289** **Approximants de Padé**

**HIS1** **Fraction rationnelle**

$$\frac{z^2 (3z^3 + z^2 + z + 1)}{(z^2 + z + 1)(1 + z)^2 (z-1)^6}$$

**0, 0, 1, 4, 12, 31, 67, 132, 239, 407, 657, 1019, 1523, 2211, 3126, 4323, 5859, 7806, 10236, 13239, 16906, 21346, 26670, 33010, 40498, 49290, 59543, 71438, 85158, 100913**

## Representation degeneracies for Raymond strings

**Réf.** NUPH B274 544 86.

**HIS2 A5303** **Euler**

**HIS1** **Produit infini**

*** Le motif [4, 2] est périodique**

$$\prod_{n\geq 1} \frac{1}{(1 - z^n)^{c(n)}}$$

**c(n) = 0,2,4,3,4,2,4,2,...***

**1, 0, 2, 4, 6, 12, 22, 36, 62, 104, 166, 268, 426, 660, 1022, 1564, 2358, 3540, 5266, 7756, 11362, 16524, 23854, 34252, 48890, 69368, 97942, 137588, 192314, 267628, 370798, 511524, 702886**

---

## Representation degeneracies for Raymond strings

**Réf.** NUPH B274 548 86.

**HIS2 A5304** **Euler**

**HIS1** **Produit infini**

*** Le motif [4, 2] est périodique**

$$\prod_{n\geq 1} \frac{1}{(1 - z^n)^{c(n)}}$$

**c(n) = 1,1,3,3,4,3,4,2,...***

**2, 2, 4, 10, 18, 32, 58, 98, 164, 274, 442, 704, 1114, 1730, 2660, 4058, 6114, 9136, 13554, 19930**

## Representation degeneracies for Raymond strings

**Réf.** NUPH B274 548 86.

**HIS2 A5305** **Euler**

**HIS1** **Produit infini**

*** Le motif [4, 2] est périodique**

$$\prod_{n\geq 1} \frac{1}{(1 - Z^n)^{c(n)}}$$

**c(n) = 2,1,2,2,4,3,4,3,4,2,4,2,...***

**2, 4, 8, 16, 30, 56, 100, 172, 290, 480, 780, 1248, 1970, 3068, 4724, 7200, 10862, 16240, 24080**

---

## Representation degeneracies for Raymond strings

**Réf.** NUPH B274 548 86.

**HIS2 A5306** **Euler**

**HIS1** **Produit infini**

*** Le motif [4, 2] est périodique**

$$\prod_{n\geq 1} \frac{1}{(1 - Z^n)^{c(n)}}$$

**c(n) = 2,2,3,0,3,3,4,3,4,3,4,2,4,2,...***

**2, 4, 10, 22, 40, 76, 138, 238, 408, 682, 1112, 1792, 2844, 4444, 6872, 10510, 15896, 23834**

## Bosonic string states

**Réf.** CU86.

**HIS2 A5308** **Euler**

**HIS1** **Produit infini**

$$\prod_{n\geq 1} \frac{1}{(1-z^n)^{c(n)}}$$

$$c(n) = 0,0,0,1,1,2,2,3,3,4,4,...$$

**1, 0, 0, 0, 1, 1, 2, 2, 4, 4, 7, 8, 14, 16, 25, 31**

---

## Fermionic string states

**Réf.** CU86.

**HIS2 A5309** **Approximants de Padé** **conjecture**

**HIS1** **Fraction rationnelle**

$$\frac{1 - 2z + 2z^2}{1 - 2z}$$

**1, 0, 2, 4, 8, 16, 32, 60, 114, 212**

## Fermionic string states

**Réf.** CU86.

**HIS2 A5310** **Approximants de Padé**

**HIS1** **Fraction rationnelle**

$$\frac{2\,(1 - 2\,z + 2\,z^{2})}{(2\,z - 1)\,(z - 1)}$$

**2, 2, 6, 14, 30, 62, 126, 246, 472**

## Triangular anti-Hadamard matrices of order n

**Réf.** LAA 62 117 84.

**HIS2 A5313** **Approximants de Padé**

**HIS1** **Fraction rationnelle**

$$\frac{1 - z - 3\,z^{2} + z^{3}}{(1 + z)\,(1 - 3\,z + z^{2})\,(z - 1)^{2}}$$

**1, 3, 6, 13, 29, 70, 175, 449, 1164, 3035, 7931, 20748, 54301, 142143, 372114, 974185, 2550425, 6677074, 17480779, 45765245, 119814936, 313679543, 821223671, 2149991448**

**Réf.** LAA 62 130 84.
**HIS2 A5314** **Approximants de Padé**
**HIS1** **Fraction rationnelle**

$$\frac{(z - 1)(1 + z)^2}{z^3 - z^2 + 2z - 1}$$

**1, 1, 2, 3, 5, 9, 16, 28, 49, 86, 151, 265, 465, 816, 1432, 2513, 4410, 7739, 13581, 23833, 41824, 73396, 128801, 226030, 396655, 696081, 1221537, 2143648, 3761840, 6601569**

---

**(2 ^ n + C(2n,n))/2**

**Réf.** pcf.
**HIS2 A5317** **LLL** **Suite P-récurrente**
**HIS1** **algébrique**

$$\frac{4z + 2(-4z + 1)^{1/2} z - (-4z + 1)^{1/2} - 1}{2(1 - 4z)(1 - 2z)}$$

**1, 2, 5, 14, 43, 142, 494, 1780, 6563, 24566, 92890, 353740, 1354126, 5204396, 20066492, 77575144, 300572963, 1166868646, 4537698722, 17672894044, 68923788698**

## Column of Motzkin triangle

**Réf.** JCT A23 293 77.

**HIS2 A5322** **LLL** **Suite P-récurrente**

**HIS1** **algébrique**

**a(n) (5 + n) = (13 + 4 n) a(n - 1) - n a(n - 2) - 6 n a(n - 3)**

$$\frac{1 - 3z + 2z^3 - \left(-(3z^2 + 2z - 1)(-1 + 2z)^2\right)^{1/2}}{2z^6}$$

**1, 3, 9, 25, 69, 189, 518, 1422, 3915, 10813, 29964, 83304, 232323, 649845, 1822824, 5126520, 14453451, 40843521, 115668105, 328233969, 933206967, 2657946907, 7583013474**

## Column of Motzkin triangle

**Réf.** JCT A23 293 77.

**HIS2 A5323** **LLL** **Suite P-récurrente**

**HIS1** **algébrique**

**(n + 7) (n - 1) a(n) = (n + 2) (2 n + 5) a(n - 1) + (n + 2) (3 n + 3) a(n - 2)**

$$\frac{1 - 4z + 2z^2 + 4z^3 - z^4 - \left(-(-1 + 2z + 3z^2)(1 - 3z + z^2 + z^3)^2\right)^{1/2}}{z^8}$$

**1, 4, 14, 44, 133, 392, 1140, 3288, 9438, 27016, 77220, 220584, 630084, 1800384, 5147328, 14727168, 42171849, 120870324, 346757334, 995742748, 2862099185**

## Column of Motzkin triangle

**Réf.** JCT A23 293 77.

**HIS2 A5324** **LLL** **Suite P-récurrente**

**HIS1** **algébrique**

**a(n) (n + 9) (n - 1) = (n + 3) (3 n + 6) a(n - 2) + (n + 3) (2 n + 7) a(n - 1)**

$$\frac{-1/2\,(-1 + 5z - 5z^2 - 5z^3 + 5z^4 + z^5}{z^{10}} + \frac{(-(z+1)(3z-1)(z^2+z-1)^2 (z^2 - 3z + 1)^2)^{1/2})}{z^{10}}$$

**1, 5, 20, 70, 230, 726, 2235, 6765, 20240, 60060, 177177, 520455, 1524120, 4453320, 12991230, 37854954, 110218905, 320751445, 933149470, 2714401580, 7895719634**

## Column of Motzkin triangle

**Réf.** JCT A23 293 77.

**HIS2 A5325** **LLL** **Suite P-récurrente**

**HIS1** **algébrique**

**a(n) (n + 11) (n - 1) = (n + 4) (3 n + 9) a(n - 2) + (n + 4) (2 n + 9) a(n - 1)**

$$\frac{1/2\,(1 - 6z + 9z^2 + 4z^3 - 12z^4 + 2z^6}{z^{12}} - \frac{(-(z+1)(3z-1)(z-1)^2 (2z-1)^2 (2z^2 + 2z - 1)^2)^{1/2})}{z^{12}}$$

**1, 6, 27, 104, 369, 1242, 4037, 12804, 39897, 122694, 373581, 1128816, 3390582, 10136556, 30192102, 89662216, 265640691, 785509362, 2319218869, 6839057544**

## Putting balls into 4 boxes

**Réf.** SIAR 12 296 70.

**HIS2 A5337** **Approximants de Padé**

**HIS1** **Fraction rationnelle**

$$\frac{15 - 20 z + 6 z^{2}}{(z - 1)^{4}}$$

**15, 40, 76, 124, 185, 260, 350, 456, 579, 720, 880, 1060, 1211**

---

## Low discrepancy sequences in base 3

**Réf.** JNT 30 68 88.

**HIS2 A5357** **Approximants de Padé**

**HIS1** **Fraction rationnelle**

$$\frac{1 + z^{3} + z^{11}}{(z - 1)^{2}}$$

**0, 0, 0, 1, 2, 3, 5, 7, 9, 11, 13, 15, 17, 19, 22, 25, 28, 31, 34, 37, 40, 43, 46, 49, 52, 55, 58, 61, 64, 67**

## Hoggatt sequence

**Réf.** FQ 27 167 89. FA90.

**HIS2 A5362** **P-récurrences** **Suite P-récurrente**

**HIS1**

$$(n + 5) (n + 4) (n + 3) (n + 2) a(n) =$$

$$(12 n^4 + 78 n^3 + 162 n^2 + 108 n) a(n - 1)$$

$$+ (64 n^4 - 64 n^3 - 196 n^2 + 76 n + 120) a(n - 2)$$

**1, 2, 7, 32, 177, 1122, 7898, 60398, 494078, 4274228, 38763298, 366039104, 3579512809, 36091415154, 373853631974, 3966563630394, 42997859838010, 47519**

---

**Réf.** FA90.

**HIS2 A5367** **Approximants de Padé**

**HIS1** **Fraction rationnelle**

$$\frac{1 - z + z^3}{(1 + z) (z - 1)^3}$$

**1, 1, 2, 3, 5, 7, 10, 13, 17, 21, 26, 31, 37, 43, 50, 57, 65, 73, 82, 91, 101, 111, 122, 133, 145, 157, 170, 183, 197, 211, 226, 241, 257, 273, 290, 307, 325, 343, 362, 381, 401, 421, 442, 463**

## Low discrepancy sequences in base 4

**Réf.** JNT 30 69 88.

**HIS2 A5377** **Approximants de Padé** **Conjecture**

**HIS1** **Fraction rationnelle**

$$\frac{z^4 (1 + z^2)(z^4 - z^2 + 1)}{(z - 1)^2}$$

**0, 0, 0, 0, 1, 2, 3, 4, 5, 6, 8, 10, 12, 14, 16, 18, 20, 22, 24, 26, 28, 30, 32, 34, 36, 38, 40, 42, 44, 46**

---

**Réf.** SAM 273 71. DM 75 94 89.

**HIS2 A5380** **Euler**

**HIS1** **Produit infini**

$$\prod_{n \geq 1} \frac{1}{(1 - z^n)^{c(n)}}$$

$$c(n) = 2, 3, 4, 5, ...$$

**1, 2, 6, 14, 33, 70, 149, 298, 591, 1132, 2139, 3948, 7199, 12894, 22836, 39894, 68982, 117948, 199852, 335426, 558429, 922112, 1511610, 2460208, 3977963, 6390942, 10206862, 16207444, 25596941, 40214896**

## Area of nth triple of squares around a triangle

**Réf.** PYTH 14 81 75.

**HIS2 A5386** **Approximants de Padé**

**HIS1** **Fraction rationnelle**

$$\frac{1 - z}{(1 + z)(1 - 5z + z^2)}$$

**1, 3, 16, 75, 361, 1728, 8281**

## Partitional matroids on n elements

**Réf.** SMH 9 249 74.

**HIS2 A5387** **Dérivée logarithmique**

**HIS1** **exponentielle**

$$\exp(\exp(z)\, z - \exp(z) + 2z + 1)$$

**1, 2, 5, 16, 62, 276, 1377, 7596, 45789, 298626, 2090910, 15621640, 123897413, 1038535174, 9165475893, 84886111212, 822648571314, 8321077557124, 87648445601429**

## Hamiltonian circuits on 2n 4 rectangle

**Réf.** JPA 17 445 84.

**HIS2 A5389** **Approximants de Padé**

**HIS1** **Fraction rationnelle**

$$\frac{1 - 2z - z^2}{1 - 8z + 10z^2 + z^4}$$

**1, 6, 37, 236, 1517, 9770, 62953, 405688, 2614457, 16849006, 108584525, 699780452, 4509783909, 29063617746, 187302518353, 1207084188912, 7779138543857, 50133202843990**

## The odd numbers

**Réf.**

**HIS2 A5408** **Approximants de Padé**

**HIS1** **Fraction rationnelle**

$$\frac{1 + z}{(z - 1)^2}$$

**1, 3, 5, 7, 9, 11, 13, 15, 17, 19, 21, 23, 25, 27, 29, 31, 33, 35, 37, 39, 41, 43, 45, 47, 49, 51, 53, 55, 57, 59, 61, 63, 65, 67, 69, 71, 73, 75, 77, 79, 81, 83, 85, 87, 89, 91, 93, 95, 97, 99, 101**

## Polynomials of height n

**Réf.** CR41 103. smd.

**HIS2 A5409** **Approximants de Padé**

**HIS1** **Fraction rationnelle**

$$\frac{1 - 2z + 2z^2 + z^3}{(1 - z)(1 - 2z - z^2)}$$

**1, 1, 4, 11, 28, 69, 168, 407, 984**

## Binary grids

**Réf.** TYCM 9 267 78.

**HIS2 A5418** **Approximants de Padé**

**HIS1** **Fraction rationnelle**

$$\frac{3z^2 - 1}{(1 - 2z)(2z^2 - 1)}$$

**1, 2, 3, 6, 10, 20, 36, 72, 136, 272, 528, 1056, 2080, 4160, 8256, 16512, 32896, 65792, 131328, 262656, 524800, 1049600, 2098176, 4196352, 8390656, 16781312, 33558528, 67117056**

## States of telephone exchange with n subscribers

**Réf.** JCT A21 162 1976.

**HIS2 A5425** **Dérivée logarithmique** **Suite P-récurrente**

**HIS1** **exponentielle**

**a(n) = 2 a(n - 1) + (n - 2) a(n - 2)**

$$\exp(2 z + 1/2 z^2)$$

**1, 2, 5, 14, 43, 142, 499, 1850, 7193, 29186, 123109, 538078, 2430355, 11317646, 54229907, 266906858, 1347262321, 6965034370, 36833528197, 199037675054, 1097912385851**

---

## Apéry numbers

**Réf.** MI 1 195 78. JNT 20 92 85.

**HIS2 A5429** **Hypergéométrique** **Suite P-récurrente.**

**HIS1** **algébrique**

$$\frac{4 z^2 + 10 z + 1}{(1 - 4 z)^{7/2}}$$

**0, 2, 48, 540, 4480, 31500, 199584, 1177176, 6589440, 35443980, 184756000, 938929992, 4672781568, 22850118200, 110079950400, 523521630000, 2462025277440, 11465007358860**

## Apéry numbers

**Réf.** MI 1 195 78. JNT 20 92 85.

**HIS2 A5430** **Hypergéométrique** **Suite P-récurrente**

**HIS1** **algébrique**

$$\frac{2z}{(1-4z)^{3/2}}$$

**0, 2, 12, 60, 280, 1260, 5544, 24024, 102960, 437580, 1847560, 7759752, 32449872, 135207800, 561632400, 2326762800, 9617286240, 39671305740, 163352435400**

---

## Convex polygons of length 2n on square lattice

**Réf.** TCS 34 179 84. JPA 21 L472 88.

**HIS2 A5436** **LLL** **Suite P-récurrente**

**HIS1** **algébrique**

**(n - 3) a(n) = (12 n - 42) a(n - 1) + (- 48 n + 192) a(n - 2) + (64 n - 288) a(n - 3)**

$$\frac{-4z^3-4z^2(1-4z)^{1/2}+11z^2-6z+1}{(4z-1)^2}$$

**1, 2, 7, 28, 120, 528, 2344, 10416, 46160, 203680, 894312, 3907056, 16986352, 73512288, 316786960, 1359763168, 5815457184, 24788842304, 105340982248, 446389242480**

## From a Fibonacci-like differential equation

**Réf.** FQ 27 306 89.

**HIS2 A5442** **Approximants de Padé** **f.g. exponentielle**

**HIS1** **Fraction rationnelle**

$$\frac{1}{1 - z - z^2}$$

**1, 1, 4, 18, 120, 960, 9360, 105840, 1370880, 19958400**

## From a Fibonacci-like differential equation

**Réf.** FQ 27 306 89.

**HIS2 A5443** **Dérivée logarithmique** **f.g. exponentielle**

**HIS1** **Fraction rationnelle**

$$\frac{1 - z^2}{1 - z - z^2}$$

**0, 1, 2, 12, 72, 600, 5760, 65520, 846720, 12337920**

## Centered triangular numbers

**Réf.** INOC 24 4550 85.

**HIS2 A5448** **Approximants de Padé**

**HIS1** **Fraction rationnelle**

$$\frac{z^2 + z + 1}{(1 - z)^3}$$

**1, 4, 10, 19, 31, 46, 64, 85, 109, 136, 166, 199, 235, 274, 316, 361, 409, 460, 514, 571, 631, 694, 760, 829, 901, 976, 1054, 1135, 1219, 1306, 1396, 1489, 1585, 1684, 1786, 1891, 1999**

---

**Réf.** rkg.

**HIS2 A5460** **Dérivée logarithmique**

**HIS1** **exponentielle**

$$\frac{2 z + 1}{(1 - z)^5}$$

**1, 7, 50, 390, 3360, 31920, 332640, 3780000, 46569600, 618710400, 8821612800, 134399865600, 2179457280000, 37486665216000, 681734237184000, 13071512982528000**

## Simplices in barycentric subdivision of n-simplex

**Réf.** rkg.

**HIS2 A5461** **Approximants de Padé** **Suite P-récurrente**

**HIS1** **Fraction rationnelle**

**a(n) = (n + 13) a(n - 1) + (- 8 n - 36) a(n - 2) + (12 n + 12) a(n - 3)**

$$\frac{6 z^2 + 8 z + 1}{(1 - z)^7}$$

**1, 15, 180, 2100, 25200, 317520, 4233600, 59875200, 898128000, 14270256000, 239740300800, 4249941696000, 79332244992000, 1556132497920000**

## Simplices in barycentric subdivision of n-simplex

**Réf.** rkg.

**HIS2 A5462** **Dérivée logarithmique** **f.g. exponentielle**

**HIS1** **Fraction rationnelle**

$$\frac{24 z^3 + 58 z^2 + 22 z + 1}{(1 - z)^9}$$

**1, 31, 602, 10206, 166824, 2739240, 46070640, 801496080, 14495120640, 273158645760, 5368729766400, 110055327782400, 2351983118284800**

**Réf.** JCT A24 316 78.
**HIS2 A5491** **Approximants de Padé**
**HIS1** **Fraction rationnelle**

$$\frac{3z^3 + z^2 + z + 1}{(z - 1)^4}$$

**1, 5, 15, 37, 77, 141, 235, 365, 537, 757, 1031, 1365, 1765, 2237, 2787, 3421, 4145, 4965, 5887, 6917, 8061, 9325, 10715, 12237, 13897, 15701, 17655, 19765, 22037**

---

**From expansion of falling factorials**

**Réf.** JCT A24 316 78.
**HIS2 A5492** **Approximants de Padé**
**HIS1** **Fraction rationnelle**

$$\frac{15 - 23z + 41z^2 - 13z^3 + 4z^4}{(1 - z)^5}$$

**15, 52, 151, 372, 799, 1540, 2727, 4516, 7087, 10644, 15415, 21652, 29631, 39652, 52039, 67140, 85327, 106996, 132567, 162484, 197215, 237252, 283111**

## From sum of 1/F(n)

**Réf.** FQ 15 46 77.

**HIS2 A5522** **Approximants de Padé** **Conjecture**

**HIS1** **Fraction rationnelle**

**F(n) : Nombres de Fibonacci**

$$\frac{3 - 9z + z^2 + 10z^3 - 4z^4}{(1 - z)(1 - 3z + z^2)(1 - z - z^2)}$$

**3, 6, 10, 21, 46, 108, 263, 658, 1674, 4305, 11146, 28980**

---

## Sums of successive Motzkin numbers

**Réf.** JCT B29 82 80.

**HIS2 A5554** **LLL** **Suite P-récurrente**

**HIS1** **algébrique**

**(n + 1) a(n) = 2 n a(n - 1) + (3 n - 9) a(n - 2)**

$$\frac{1 - z^2 - (-(3z - 1)(z + 1)^3)^{1/2}}{2z^2}$$

**1, 2, 3, 6, 13, 30, 72, 178, 450, 1158, 3023, 7986, 21309, 57346, 155469, 424206, 1164039, 3210246, 8893161, 24735666, 69051303, 193399578, 543310782, 1530523638**

## Walks on square lattice

**Réf.** GU90.

**HIS2 A5555** **Approximants de Padé**

**HIS1** **Fraction rationnelle**

$$\frac{5 - 6 z + 2 z^2}{(z - 1)^4}$$

**5, 14, 28, 48, 75, 110, 154, 208, 273, 350, 440, 544, 663, 798, 950, 1120, 1309, 1518, 1748, 2000, 2275, 2574, 2898, 3248, 3625, 4030, 4464, 4928, 5423, 5950, 6510, 7104, 7733**

## Walks on square lattice

**Réf.** GU90.

**HIS2 A5556** **Approximants de Padé**

**HIS1** **Fraction rationnelle**

$$\frac{14 - 28 z + 20 z^2 - 5 z^3}{(1 - z)^5}$$

**14, 42, 90, 165, 275, 429, 637, 910, 1260, 1700, 2244, 2907, 3705, 4655, 5775, 7084, 8602, 10350, 12350, 14625, 17199, 20097, 23345, 26970, 31000, 35464, 40392, 45815, 51765**

## Walks on square lattice

**Réf.** GU90.

**HIS2 A5557** **Approximants de Padé**

**HIS1** **Fraction rationnelle**

$$\frac{42 - 120 z + 135 z^2 - 70 z^3 + 14 z^4}{(z - 1)^6}$$

**42, 132, 297, 572, 1001, 1638, 2548, 3808, 5508, 7752, 10659, 14364, 19019, 24794, 31878, 40480, 50830, 63180, 77805, 95004, 115101, 138446, 165416, 196416, 231880**

## Walks on square lattice

**Réf.** GU90.

**HIS2 A5558** **P-récurrences** **Suite P-récurrente**

**HIS1**

$$(n + 2) (n + 1) a(n) = (- 64 n^2 + 320 n - 384) a(n - 3)$$

$$+ (16 n^2 - 48 n + 16) a(n - 2) + (4 n^2 + 4 n - 4) a(n - 1)$$

**1, 1, 3, 6, 20, 50, 175, 490, 1764, 5292, 19404, 60984, 226512, 736164, 2760615, 9202050, 34763300, 118195220, 449141836, 1551580888, 5924217936, 20734762776**

## Walks on square lattice

**Réf.** GU90.

**HIS2 A5559** **P-récurrences** **Suite P-récurrente**

**HIS1**

$$(n-1)(n+4)(n+3)\,a(n) = (64/5\,n^3 - 192/5\,n^2 + 128/5\,n)\,a(n-3) + (16\,n^3 + 96/5\,n^2 - 128/5\,n)\,a(n-2) + (-4/5\,n^3 + 12/5\,n^2 + 76/5\,n + 132/5)\,a(n-1)$$

**1, 2, 8, 20, 75, 210, 784, 2352, 8820, 27720, 104544, 339768, 1288287, 4294290, 16359200, 55621280, 212751396, 734959368, 2821056160, 9873696560, 38013731756**

## Walks on square lattice

**Réf.** GU90.

**HIS2 A5563** **Approximants de Padé**

**HIS1** **Fraction rationnelle**

$$\frac{z-3}{(z-1)^3}$$

**3, 8, 15, 24, 35, 48, 63, 80, 99, 120, 143, 168, 195, 224, 255, 288, 323, 360, 399, 440, 483, 528, 575, 624, 675, 728, 783, 840, 899, 960, 1023, 1088**

## Walks on square lattice

**Réf.** GU90.

**HIS2 A5564** **Approximants de Padé**

**HIS1** **Fraction rationnelle**

$$\frac{6 - 4 z + z^2}{(z - 1)^4}$$

**6, 20, 45, 84, 140, 216, 315, 440, 594, 780, 1001, 1260, 1560, 1904, 2295, 2736, 3230, 3780, 4389, 5060, 5796, 6600, 7475, 8424, 9450, 10556**

## Walks on square lattice

**Réf.** GU90.

**HIS2 A5565** **Approximants de Padé**

**HIS1** **Fraction rationnelle**

$$\frac{20 - 25 z + 14 z^2 - 3 z^3}{(1 - z)^5}$$

**20, 75, 189, 392, 720, 1215, 1925, 2904, 4212, 5915, 8085, 10800, 14144, 18207, 23085, 28880, 35700, 43659, 52877, 63480, 75600, 89375, 104949, 122472, 142100**

## Walks on square lattice

**Réf.** GU90.

**HIS2 A5566** **P-récurrences** **Suite P-récurrente**

**HIS1**

$$(n + 1) n a(n) = (16 n^2 - 48 n + 32) a(n - 2) + (8 n - 4) a(n - 1)$$

**1, 2, 6, 18, 60, 200, 700, 2450, 8820, 31752, 116424, 426888, 1585584, 5889312, 22084920, 82818450, 312869700, 1181952200, 4491418360, 17067389768**

## Walks on square lattice

**Réf.** GU90.

**HIS2 A5567** **Approximants de Padé**

**HIS1** **Fraction rationnelle**

$$\frac{2 (5 - 10 z + 4 z^2)}{(2 z - 1)^3 (z - 1)^3}$$

**10, 70, 308, 1092, 3414, 9834, 26752, 69784, 176306, 434382, 1048812, 2490636, 5833006, 13500754, 30933368, 70255008, 158335434, 354419190, 788529700**

## Product of successive Catalan numbers

**Réf.** JCT A43 1 86.

**HIS2 A5568** **Hypergéométrique**

**HIS1** **Intégrales elliptiques**

$$\frac{{}_2F_1([1/2, -1/2],[2],16 z) + 1/2 z}{2 z}$$

**1, 2, 10, 70, 588, 5544, 56628, 613470, 6952660, 81662152, 987369656, 12228193432, 154532114800, 1986841476000, 25928281261800, 342787130211150, 4583937702039300**

## Walks on square lattice

**Réf.** GU90.

**HIS2 A5569** **Hypergéométrique** **Suite P-récurrente**

**HIS1**

$1/5 (n - 1) (5 n + 2) (n + 3) (n + 2) a(n) = 4/5 (5 n + 7) (2 n + 1) (2 n - 1) n a(n - 1)$

$$4 ({}_4F_3([2, 17/5, 5/2, 3/2], [4, 5, 12/5],16 z))$$

**4, 34, 308, 3024, 31680, 349206, 4008004, 47530912, 579058896, 7215393640, 91644262864, 1183274479040, 15497363512800, 205519758825150**

## Walks on cubic lattice

**Réf.** GU90.

**HIS2 A5570** **Approximants de Padé**

**HIS1** **Fraction rationnelle**

$$\frac{z - 17}{(z - 1)^3}$$

**17, 50, 99, 164, 245, 342, 455, 584, 729, 890, 1067, 1260, 1469, 1694, 1935, 2192, 2465, 2754, 3059, 3380, 3717, 4070, 4439, 4824, 5225, 5642, 6075, 6524, 6989, 7470**

## Walks on cubic lattice

**Réf.** GU90.

**HIS2 A5571** **Approximants de Padé**

**HIS1** **Fraction rationnelle**

$$\frac{4\,(19 - 4\,z + z^2)}{(z - 1)^4}$$

**76, 288, 700, 1376, 2380, 3776, 5628, 8000, 10956, 14560, 18876, 23968, 29900, 36736, 44540, 53376, 63308, 74400, 86716**

## Walks on cubic lattice

**Réf.** GU90.

**HIS2 A5572** **inverse fonctionnel** **Suite P-récurrente**

**HIS1** **algébrique**

**(n + 1) a(n) = (- 12 n + 24) a(n - 2) + (8 n - 4) a(n - 1)**

$$\frac{1 - 4 z - (1 - 8 z + 12 z^2)^{1/2}}{2 z}$$

**1, 4, 17, 76, 354, 1704, 8421, 42508, 218318, 1137400, 5996938, 31940792, 171605956, 928931280, 5061593709**

---

## Walks on cubic lattice

**Réf.** GU90.

**HIS2 A5573** **inverse fonctionnel** **Suite P-récurrente**

**HIS1** **algébrique**

**n a(n) = (- 12 n + 24) a(n - 2) + (8 n - 6) a(n - 1)**

$$\frac{1 - 6 z - (1 - 8 z + 12 z^2)^{1/2}}{2 z}$$

**1, 5, 26, 139, 758, 4194, 23460, 132339, 751526, 4290838, 24607628, 141648830, 817952188, 4736107172, 27487711752, 159864676803**

**Réf.** GTA91 603.

**HIS2 A5578** **Approximants de Padé**

**HIS1** **Fraction rationnelle**

$$\frac{1 - z - z^2}{(z - 1)(2z - 1)(1 + z)}$$

**1, 1, 2, 3, 6, 11, 22, 43, 86, 171, 342, 683, 1366, 2731, 5462, 10923, 21846, 43691, 87382, 174763, 349526, 699051, 1398102, 2796203, 5592406, 11184811, 22369622**

---

**Réf.** AS1 797.

**HIS2 A5581** **Approximants de Padé**

**HIS1** **Fraction rationnelle**

$$\frac{2 - z}{(z - 1)^4}$$

**2, 7, 16, 30, 50, 77, 112, 156, 210, 275, 352, 442, 546, 665, 800, 952, 1122, 1311, 1520, 1750, 2002, 2277, 2576, 2900, 3250, 3627, 4032, 4466, 4930, 5425, 5952, 6512, 7106, 7735, 8400**

**Réf.** AS1 797.

**HIS2 A5582** **Approximants de Padé**

**HIS1** **Fraction rationnelle**

$$\frac{2 - z}{(z - 1)^5}$$

**2, 9, 25, 55, 105, 182, 294, 450, 660, 935, 1287, 1729, 2275, 2940, 3740, 4692, 5814, 7125, 8645, 10395, 12397, 14674, 17250, 20150, 23400, 27027, 31059, 35525, 40455, 45880, 51832**

## Coefficients of Chebyshev polynomials

**Réf.** AS1 797.

**HIS2 A5583** **Approximants de Padé**

**HIS1** **Fraction rationnelle**

$$\frac{2 - z}{(z - 1)^6}$$

**2, 11, 36, 91, 196, 378, 672, 1122, 1782, 2717, 4004, 5733, 8008, 10948, 14688, 19380, 25194, 32319, 40964, 51359, 63756, 78430, 95680, 115830, 139230, 166257, 197316, 232841**

## Coefficients of Chebyshev polynomials

**Réf.** AS1 797.

**HIS2 A5584** **Approximants de Padé**

**HIS1** **Fraction rationnelle**

$$\frac{2 - z}{(z - 1)^{7}}$$

**2, 13, 49, 140, 336, 714, 1386, 2508, 4290, 7007, 11011, 16744, 24752, 35700, 50388, 69768, 94962, 127281, 168245, 219604, 283360, 361790, 457470, 573300, 712530, 878787**

## 5-dimensional pyramidal numbers

**Réf.** AS1 797.

**HIS2 A5585** **Approximants de Padé**

**HIS1** **Fraction rationnelle**

$$\frac{1 + z}{(z - 1)^{6}}$$

**1, 7, 27, 77, 182, 378, 714, 1254, 2079, 3289, 5005, 7371, 10556, 14756, 20196, 27132, 35853, 46683, 59983, 76153, 95634, 118910, 146510, 179010, 217035, 261261, 312417, 371287**

**Réf.** AS1 796.
**HIS2 A5586** **Approximants de Padé**
**HIS1** **Fraction rationnelle**

$$\frac{z\,(5 - 6\,z + 2\,z^2)}{(z - 1)^4}$$

**0, 5, 14, 28, 48, 75, 110, 154, 208, 273, 350, 440, 544, 663, 798, 950, 1120, 1309, 1518, 1748, 2000, 2275, 2574, 2898, 3248, 3625, 4030, 4464, 4928, 5423, 5950, 6510, 7104, 7733, 8398**

---

**Réf.** AS1 796.
**HIS2 A5587** **Approximants de Padé**
**HIS1** **Fraction rationnelle**

$$\frac{z\,(-\,14 + 28\,z - 20\,z^2 + 5\,z^3)}{(z - 1)^5}$$

**0, 14, 42, 90, 165, 275, 429, 637, 910, 1260, 1700, 2244, 2907, 3705, 4655, 5775, 7084, 8602, 10350, 12350, 14625, 17199, 20097, 23345, 26970, 31000, 35464, 40392, 45815, 51765**

**Réf.** CJN 25 391 82.
**HIS2 A5592** **Approximants de Padé**
**HIS1** **Fraction rationnelle**

$$\frac{2 - 2z + z^2}{(1 - z)(1 - 3z + z^2)}$$

**2, 6, 17, 46, 122, 321, 842, 2206, 5777, 15126, 39602, 103681, 271442, 710646, 1860497, 4870846, 12752042, 33385281, 87403802, 228826126, 599074577, 1568397606, 4106118242**

---

**Réf.** CJN 25 391 82.
**HIS2 A5593** **Approximants de Padé**
**HIS1** **Fraction rationnelle**

$$\frac{2 - 5z + z^2 + 2z^3 - z^4}{(1 - z)(1 - z - z^2)(1 - 3z + z^2)}$$

**2, 5, 12, 29, 71, 177, 448, 1147, 2960, 7679, 19989, 52145, 136214, 356121, 931540, 2437513, 6379403, 16698113, 43710756, 114427391, 299560472, 784236315, 2053119817, 5375076769**

## Functions realized by cascades of n gates

**Réf.** BU77.

**HIS2 A5609** **Approximants de Padé**

**HIS1** **Fraction rationnelle**

$$\frac{16\ (7\ z\ -\ 4)}{(28\ z\ -\ 1)\ (1\ -\ z)}$$

**64, 1744, 48784, 1365904, 38245264, 1070867344, 29984285584, 839559996304**

## Functions realized by cascades of n gates

**Réf.** BU77.

**HIS2 A5610** **Approximants de Padé**

**HIS1** **Fraction rationnelle**

$$\frac{2\ (7\ -\ 6\ z)}{(1\ -\ 6\ z)\ (1\ -\ z)}$$

**14, 86, 518, 3110, 18662, 111974, 671846, 4031078**

**Disjunctively-realizable functions of n variables**

**Réf.** PGEC 24 687 75.

**HIS2 A5616** **Inverse fonctionnel** **f.g. exponentielle**

**HIS1** **exponentielle**

**L'inverse de S(z) est**

$$\ln(z + 1) - z + \ln(z + 2) - \ln(2)$$

**2, 10, 114, 2154, 56946, 1935210, 80371122, 3944568042, 223374129138, 14335569726570, 1028242536825906, 81514988432370666, 7077578056972377714**

---

**Réf.** PGEC 11 140 62.

**HIS2 A5618** **Approximants de Padé**

**HIS1** **Fraction rationnelle**

$$\frac{3 z - 1}{(1 - 6 z) (z - 1)}$$

**4, 16, 88, 520, 3112, 18664, 111976, 671848, 4031080, 24186472, 145118824, 870712936, 5224277608, 31345665640, 188073993832, 1128443962984, 6770663777896**

## Functions realized by n-input cascades

**Réf.** PGEC 27 790 78.

**HIS2 A5619** **Approximants de Padé**

**HIS1** **Fraction rationnelle**

$$\frac{16\ (1 - 18\ z + 20\ z^{2})}{(z - 1)\ (80\ z^{2} - 32\ z + 1)}$$

**16, 240, 6448, 187184, 5474096, 160196400, 4688357168, 137211717424, 4015706384176**

---

**Réf.** JACM 23 705 76. PGEC 27 315 78. LNM 829 122 80.

**HIS2 A5640** **Inverse fonctionnel**

**HIS1** **exponentielle**

$$- 2\ W(- 1/2\ \exp(z - 1/2))$$

**1, 2, 8, 64, 832, 15104, 352256, 10037248, 337936384, 13126565888**

## From sum of inverse binomial coefficients

**Réf.** C1 294.

**HIS2 A5649** **Recoupements**

**HIS1** **exponentielle**

$$\frac{1}{(\exp(z) - 2)^2}$$

**1, 2, 8, 44, 308, 2612, 25988, 296564, 3816548, 54667412, 862440068, 14857100084, 277474957988, 5584100659412, 120462266974148, 2772968936479604, 67843210855558628**

---

## Tower of Hanoi with cyclic moves only

**Réf.** IPL 13 118 81. GKP 18.

**HIS2 A5665** **Approximants de Padé**

**HIS1** **Fraction rationnelle**

$$\frac{z\ (1 + 2\ z)}{(z - 1)\ (2\ z^2 + 2\ z - 1)}$$

**0, 1, 5, 15, 43, 119, 327, 895, 2447, 6687, 18271, 49919, 136383, 372607, 1017983, 2781183, 7598335, 20759039, 56714751, 154947583, 423324671, 1156544511, 3159738367**

## Tower of Hanoi with cyclic moves only

**Réf.** IPL 13 118 81. GKP 18.

**HIS2 A5666** **Approximants de Padé**

**HIS1** **Fraction rationnelle**

$$\frac{z\,(2 + z)}{(z - 1)\,(2\,z^2 + 2\,z - 1)}$$

**0, 2, 7, 21, 59, 163, 447, 1223, 3343, 9135, 24959, 68191, 186303, 508991, 1390591, 3799167, 10379519, 28357375, 77473791, 211662335, 578272255, 1579869183, 4316282879**

---

**Réf.** rkg.

**HIS2 A5667** **Approximants de Padé**

**HIS1** **Fraction rationnelle**

$$\frac{1 - 3\,z}{1 - 6\,z - z^2}$$

**1, 3, 19, 117, 721, 4443, 27379, 168717, 1039681, 6406803, 39480499, 243289797, 1499219281, 9238605483, 56930852179, 350823718557, 2161873163521, 13322062699683**

## Convergents to square root of 10

**Réf.** rkg.

**HIS2 A5668** **Approximants de Padé**

**HIS1** **Fraction rationnelle**

$$\frac{z}{1 - 6z - z^2}$$

**0, 1, 6, 37, 228, 1405, 8658, 53353, 328776, 2026009, 12484830, 76934989, 474094764, 2921503573**

## F(n) - 2 ^ [n/2]

**Réf.** rkg.

**HIS2 A5672** **Approximants de Padé**

**HIS1** **Fraction rationnelle**

$$\frac{z^3}{(1 - z - z^2)(1 - 2z^2)}$$

**0, 0, 0, 1, 1, 4, 5, 13, 18, 39, 57, 112, 169, 313, 482, 859, 1341, 2328, 3669, 6253, 9922, 16687, 26609, 44320, 70929, 117297, 188226, 309619, 497845, 815656, 1313501, 2145541**

**Réf.** rkg.
**HIS2 A5673** **Approximants de Padé**
**HIS1** **Fraction rationnelle**

$$\frac{z^4}{(1 - z)(2 z^2 - 1)(z^2 + z - 1)}$$

**0, 0, 0, 0, 1, 2, 6, 11, 24, 42, 81, 138, 250, 419, 732, 1214, 2073, 3414, 5742, 9411, 15664, 25586, 42273, 68882, 113202, 184131, 301428, 489654, 799273, 1297118, 2112774**

---

**Réf.** rkg.
**HIS2 A5674** **Approximants de Padé**
**HIS1** **Fraction rationnelle**

$$\frac{z^4}{(1 - 2 z)(2 z^2 - 1)(z^2 + z - 1)}$$

**0, 0, 0, 0, 1, 3, 10, 25, 63, 144, 327, 711, 1534, 3237, 6787, 14056, 28971, 59283, 120894, 245457, 497167, 1004256, 2025199, 4077007, 8198334, 16467597, 33052491, 66293208**

## $\sum C(n-k,4k)$, k=0...n

**Réf.**

**HIS2 A5676** **Approximants de Padé**

**HIS1** **Fraction rationnelle**

$$\frac{(1 - z)^3}{1 - 4z + 6z^2 - 4z^3 + z^4 - z^5}$$

**1, 1, 1, 1, 1, 2, 6, 16, 36, 71, 128, 220, 376, 661, 1211, 2290, 4382, 8347, 15706, 29191, 53824, 99009, 182497, 337745, 627401, 1167937, 2174834, 4046070, 7517368, 13951852, 25880583**

---

## Twopins positions

**Réf.** GU81.

**HIS2 A5682** **Approximants de Padé**

**HIS1** **Fraction rationnelle**

$$\frac{1}{(z^3 - z^2 + 2z - 1)(-1 + z^2 + z^3)}$$

**1, 2, 4, 8, 15, 28, 51, 92, 165, 294, 522, 924, 1632, 2878, 5069, 8920, 15686, 27570, 48439, 85080, 149405, 262320, 460515, 808380, 1418916, 2490432**

## Numbers of Twopins positions

**Réf.** GU81.

**HIS2 A5683** **Approximants de Padé**

**HIS1** **Fraction rationnelle**

$$\frac{1 - z^2 - z^3 - z^4 - z^5}{(z^3 - z^2 + 2 z - 1) (1 - z^2 - z^3)}$$

**1, 2, 3, 5, 8, 13, 22, 37, 63, 108, 186, 322, 559, 973, 1697, 2964, 5183, 9071, 15886, 27835, 48790, 85545, 150021, 263136, 461596, 809812, 1420813, 2492945**

## Twopins positions

**Réf.** GU81.

**HIS2 A5684** **Approximants de Padé**

**HIS1** **Fraction rationnelle**

$$\frac{1}{(1 - z + z^2) (1 - z - z^2) (1 - z^2 - z^4)}$$

**1, 2, 4, 6, 11, 18, 32, 52, 88, 142, 236, 382, 629, 1018, 1664, 2692, 4383, 7092, 11520, 18640, 30232, 48916, 79264, 128252, 207705, 336074, 544084**

## Twopins positions

**Réf.** GU81.

**HIS2 A5685** **Approximants de Padé**

**HIS1** **Fraction rationnelle**

$$\frac{1 - z^2 + z^3 - 2z^4 - z^5 - z^6 - z^7}{(1 - z + z^2)(1 - z - z^3)(1 - z^2 - z^4)}$$

**1, 2, 3, 5, 7, 11, 16, 26, 40, 65, 101, 163, 257, 416, 663, 1073, 1719, 2781, 4472, 7236, 11664, 18873, 30465, 49293, 79641, 128862, 208315, 337061, 545071**

---

## Twopins positions

**Réf.** GU81.

**HIS2 A5686** **Approximants de Padé**

**HIS1** **Fraction rationnelle**

$$\frac{(1 + z)(z^3 + z + 1)}{1 + z^2 + z^5}$$

**1, 2, 2, 3, 3, 4, 5, 6, 8, 9, 12, 14, 18, 22, 27, 34, 41, 52, 63, 79, 97, 120, 149, 183, 228, 280, 348, 429, 531, 657, 811**

## Twopins positions

**Réf.** GU81.

**HIS2 A5687** **Approximants de Padé**

**HIS1** **Fraction rationnelle**

$$\frac{1}{(1 - 2z + z^{2} - z^{5})(1 - z^{2} - z^{5})}$$

**1, 2, 4, 6, 9, 14, 22, 36, 57, 90, 139, 214, 329, 506, 780, 1200, 1845, 2830, 4337, 6642, 10170, 15572, 23838, 36486, 55828, 85408, 130641, 199814, 305599**

## Twopins positions

**Réf.** FQ 16 85 78. GU81.

**HIS2 A5689** **Approximants de Padé**

**HIS1** **Fraction rationnelle**

$$\frac{1 + z^{2} + z^{3} + z^{4} + z^{5}}{(1 - z - z^{3})(z^{3} - z + 1)}$$

**1, 2, 4, 7, 11, 16, 22, 30, 42, 61, 91, 137, 205, 303, 443, 644, 936, 1365, 1999, 2936, 4316, 6340, 9300, 13625, 19949, 29209, 42785, 62701, 91917, 134758, 197548, 289547**

## Twopins positions

**Réf.** GU81.

**HIS2 A5690** **Approximants de Padé**

**HIS1** **Fraction rationnelle**

$$\frac{1}{(1 - z - z^3)(1 - z + z^3)(1 - z^2 - z^6)}$$

**1, 2, 4, 6, 9, 12, 18, 26, 41, 62, 96, 142, 212, 308, 454, 662, 979, 1438, 2128, 3126, 4606, 6748, 9910, 14510, 21298, 31212, 45820, 67176, 98571, 144476**

---

## Dyck paths

**Réf.** LNM 1234 118 86.

**HIS2 A5700** **hypergéométrique** **Suite P-récurrente**

**HIS1** **Intégrales elliptiques**

$$_3F_2([1, 1/2, 3/2],[3 , 4],16 z)$$

**1, 1, 3, 14, 84, 594, 4719, 40898, 379236, 3711916, 37975756, 403127256**

**Réf.** R1 150. rkg.

**HIS2 A5704** **Euler**

**HIS1** **Produit infini**

$$\frac{1}{(1-z)^2\,(1-z^3)\,(1-z^9)\,(1-z^{27})\ldots}$$

**1, 1, 2, 4, 8, 19, 44, 112, 287, 763**

---

**Réf.** AMM 95 555 88.

**HIS2 A5708** **Approximants de Padé**

**HIS1** **Fraction rationnelle**

$$\frac{1}{1-z-z^6}$$

**1, 1, 1, 1, 1, 1, 2, 3, 4, 5, 6, 7, 9, 12, 16, 21, 27, 34, 43, 55, 71, 92, 119, 153, 196, 251, 322, 414, 533, 686, 882, 1133, 1455, 1869, 2402, 3088, 3970, 5103**

**Réf.** AMM 95 555 88.
**HIS2 A5709** **Approximants de Padé**
**HIS1** **Fraction rationnelle**

$$\frac{1}{1 - z - z^7}$$

**1, 1, 1, 1, 1, 1, 1, 2, 3, 4, 5, 6, 7, 8, 10, 13, 17, 22, 28, 35, 43, 53, 66, 83, 105, 133, 168, 213, 266, 332, 415, 520, 653, 821, 1034, 1300, 1632, 2047, 2567, 3220, 4041**

---

**Réf.** AMM 95 555 88.
**HIS2 A5710** **Approximants de Padé**
**HIS1** **Fraction rationnelle**

$$\frac{1}{1 - z - z^8}$$

**1, 1, 1, 1, 1, 1, 1, 1, 2, 3, 4, 5, 6, 7, 8, 9, 11, 14, 18, 23, 29, 36, 44, 53, 64, 78, 96, 119, 148, 184, 228, 281, 345, 423, 519, 638, 786, 970, 1198, 1479, 1824, 2247, 2766, 3404**

**Réf.** AMM 95 555 88.
**HIS2 A5711** **Approximants de Padé**
**HIS1** **Fraction rationnelle**

$$\frac{1 + z^{8}}{1 - z - z^{9}}$$

**1, 1, 1, 1, 1, 1, 1, 1, 2, 3, 4, 5, 6, 7, 8, 9, 10, 12, 15, 19, 24, 30, 37, 45, 54, 64, 76, 91, 110, 134, 164, 201, 246, 300, 364, 440, 531, 641, 775, 939, 1140, 1386, 1686, 2050, 2490, 3021**

---

**From expansion of (1 + x + x ^ 2) ^ n**
**Réf.** C1 78.
**HIS2 A5712** **Approximants de Padé**
**HIS1** **Fraction rationnelle**

$$\frac{z^{2} - z - 1}{(z - 1)^{5}}$$

**1, 6, 19, 45, 90, 161, 266, 414, 615, 880, 1221, 1651, 2184, 2835, 3620, 4556, 5661, 6954, 8455, 10185, 12166, 14421, 16974, 19850, 23075, 26676, 30681, 35119, 40020, 45415**

## From expansion of (1 + x + x ^ 2) ^ n

**Réf.** C1 78.

**HIS2 A5714** **Approximants de Padé**

**HIS1** **Fraction rationnelle**

$$\frac{1 + 3z - 4z^2 + z^3}{(1 - z)^7}$$

**1, 10, 45, 141, 357, 784, 1554, 2850, 4917, 8074, 12727, 19383, 28665, 41328, 58276, 80580, 109497, 146490, 193249, 251713, 324093, 412896, 520950, 651430, 807885**

## From expansion of (1 + x + x ^ 2) ^ n

**Réf.** C1 78.

**HIS2 A5715** **Approximants de Padé**

**HIS1** **Fraction rationnelle**

$$\frac{(2 - z)(z^2 - 2)}{(1 - z)^8}$$

**4, 30, 126, 393, 1016, 2304, 4740, 9042, 16236, 27742, 45474, 71955, 110448, 165104, 241128, 344964, 484500, 669294, 910822, 1222749, 1621224, 2125200, 2756780**

## From expansion of (1 + x + x ^ 2) ^ n

**Réf.** C1 78.

**HIS2 A5716** **Approximants de Padé**

**HIS1** **Fraction rationnelle**

$$\frac{1 + 6 z - 9 z^{2} + 3 z^{3}}{(1 - z)^{9}}$$

**1, 15, 90, 357, 1107, 2907, 6765, 14355, 28314, 52624, 93093, 157950, 258570, 410346, 633726, 955434, 1409895, 2040885, 2903428, 4065963, 5612805, 7646925**

---

## From expansion of (1 + x + x ^ 2) ^ n

**Réf.** C1 78.

**HIS2 A5717** **LLL** **Suite P-récurrente**

**HIS1** **algébrique**

**(n + 1) a(n) = 3 n a(n - 1) + (- 3 n + 6) a(n - 3) + (n + 3) a(n - 2)**

$$\frac{z + (z + 1)^{1/2} (1 - 3 z)^{1/2} - 1}{2 (z^{2} (z + 1)^{1/2} (1 - 3 z)^{1/2} )}$$

**1, 2, 6, 16, 45, 126, 357, 1016, 2907, 8350, 24068, 69576, 201643, 585690, 1704510, 4969152, 14508939, 42422022, 124191258, 363985680, 1067892399, 3136046298, 9217554129**

## Quadrinomial coefficients

**Réf.** C1 78.

**HIS2 A5718** **Approximants de Padé**

**HIS1** **Fraction rationnelle**

$$\frac{z^2 - 3z + 3}{(1 - z)^5}$$

**3, 12, 31, 65, 120, 203, 322, 486, 705, 990, 1353, 1807, 2366, 3045, 3860, 4828, 5967, 7296, 8835, 10605, 12628, 14927, 17526, 20450, 23725, 27378, 31437, 35931, 40890, 46345**

## Quadrinomial coefficients

**Réf.** C1 78.

**HIS2 A5719** **Approximants de Padé**

**HIS1** **Fraction rationnelle**

$$\frac{2 - 2z^2 + z^3}{(z - 1)^6}$$

**2, 12, 40, 101, 216, 413, 728, 1206, 1902, 2882, 4224, 6019, 8372, 11403, 15248, 20060, 26010, 33288, 42104, 52689, 65296, 80201, 97704, 118130, 141830, 169182, 200592**

## Quadrinomial coefficients

**Réf.** C1 78.

**HIS2 A5720** **Approximants de Padé**

**HIS1** **Fraction rationnelle**

$$\frac{1 + 3 z - 5 z^2 + 2 z^3}{(1 - z)^7}$$

**1, 10, 44, 135, 336, 728, 1428, 2598, 4455, 7282, 11440, 17381, 25662, 36960, 52088, 72012, 97869, 130986, 172900, 225379, 290444, 370392, 467820, 585650, 727155, 895986**

---

## Quadrinomial coefficients

**Réf.** C1 78.

**HIS2 A5725** **P-récurrences** **Suite P-récurrente**

**HIS1** **algébrique**

**La méthode LLL permet de trouver l'expression algébrique du 3è degré.**

$$1/2 (n - 1) (2 n - 3) a(n) = (- 21/4 n^2 + 143/4 n - 50) a(n - 1)$$
$$+ (24 n^2 - 139 n + 200) a(n - 2) + (20 n^2 - 120 n + 180) a(n - 3)$$
$$+ (32 n^2 - 224 n + 384) a(n - 4)$$

**1, 1, 3, 10, 31, 101, 336, 1128, 3823, 13051, 44803, 154518, 534964, 1858156, 6472168, 22597760, 79067375, 277164295, 973184313, 3422117190, 12049586631, 42478745781**

**Réf.** LI68 20. MMAG 49 181 76.
**HIS2 A5732** **Approximants de Padé**
**HIS1** **Fraction rationnelle**

$$\frac{z^3 - z - 1}{(z - 1)^7}$$

**1, 8, 35, 111, 287, 644, 1302, 2430, 4257, 7084, 11297, 17381, 25935, 37688, 53516, 74460, 101745, 136800, 181279, 237083, 306383, 391644, 495650, 621530, 772785, 953316**

## Coefficients of a modular function

**Réf.** GMJ 8 29 67.
**HIS2 A5758** **Euler**
**HIS1** **Produit infini**
*** Le motif [12] est constant**

$$\prod_{n \geq 1} \frac{1}{(1 - z^n)^{c(n)}}$$

$$c(n) = 12,12,12,12,...*$$

**1, 12, 90, 520, 2535, 10908, 42614, 153960**

## Convex polygons of length 2n on square lattice

**Réf.** TCS 34 179 84.

**HIS2 A5770** **Approximants de Padé**

**HIS1** **Fraction rationnelle**

$$\frac{1 - 3z + 2z^2 + z^3}{(4z - 1)(2z - 1)(1 - 3z + z^2)^2}$$

**1, 9, 55, 286, 1362, 6143, 26729, 113471, 473471, 1951612, 7974660, 32384127, 130926391, 527657073, 2121795391, 8518575466, 34162154550, 136893468863, 548253828965**

## Directed animals of size n

**Réf.** AAM 9 340 88.

**HIS2 A5773** **Inverse fonctionnel algébrique** **Suite P-récurrente**

**HIS1** **Inverse des nombres de Motzkin**

$$\frac{-1 + 3z + (1 - 2z - 3z^2)^{1/2}}{2(1 - 3z)}$$

**1, 2, 5, 13, 35, 96, 267, 750, 2123, 6046, 17303, 49721, 143365, 414584, 1201917, 741365049, 2173243128, 6377181825, 18730782252, 3492117, 10165779, 29643870, 86574831, 253188111**

## Directed animals of size n

**Réf.** AAM 9 340 88.

**HIS2 A5774** **P-récurrences et LLL** **Suite P-récurrente**

**HIS1** **algébrique**

$$a(n)\,(2 + n) = (4 + 4n)\,a(n-1) - n\,a(n-2)$$
$$(12 - 6n)\,a(n-3)$$

$$\frac{1 - 3z - (-(3z^2 + 2z - 1)(-1 + 2z)^2)^{1/2}}{2(3z^4 - z^3)}$$

**1, 3, 9, 26, 75, 216, 623, 1800, 5211, 15115, 43923**

## 4-dimensional Catalan numbers

**Réf.** TS89. CN 75 124 90.

**HIS2 A5790** **Hypergéométrique** **Suite P-récurrente**

**HIS1**

$${}_4F_3([1, 5/4, 7/4, 3/2], [3, 4, 5], 256\,z)$$

**1, 14, 462, 24024, 1662804, 140229804, 13672405890, 1489877926680, 177295473274920**

## Permutations with subsequences of length <= 3

**Réf.** JCT A53 281 90.

**HIS2 A5802** **P-récurrences** **Suite P-récurrente**

**HIS1**

$$(n + 1)^2 \, a(n) = (10 n^2 - 18 n + 5) \, a(n - 1) + (- 9 n^2 + 36 n - 36) \, a(n - 2)$$

**1, 1, 2, 6, 23, 103, 513, 2761, 15767, 94359, 586590, 3763290, 24792705, 167078577, 1148208090, 8026793118, 56963722223, 409687815151, 2981863943718, 21937062144834**

## Second-order Eulerian numbers

**Réf.** JCT A24 28 78. GKP 256.

**HIS2 A5803** **Approximants de Padé**

**HIS1** **Fraction rationnelle**

$$\frac{2 z}{(1 - 2 z) (z - 1)^2}$$

**0, 2, 8, 22, 52, 114, 240, 494, 1004, 2026, 4072, 8166, 16356, 32738, 65504, 131038, 262108, 524250, 1048536, 2097110, 4194260, 8388562, 16777168, 33554382, 67108812, 134217674**

## Sums of adjacent Catalan numbers

**Réf.** dek.

**HIS2 A5807** **Hypergéométrique** **améliorée par**
**HIS1** **algébrique** **la méthode LLL**

$$\frac{1 - z - (-(4z - 1)(z + 1)^2)^{1/2}}{2 z^2}$$

**2, 3, 7, 19, 56, 174, 561, 1859, 6292, 21658, 75582, 266798, 950912, 3417340, 12369285, 45052515, 165002460, 607283490, 2244901890, 8331383610, 31030387440**

## Binomial coefficients

**Réf.** AS1 828.

**HIS2 A5809** **hypergéométrique-LLL** **suite P-récurrente**
**HIS1** **algébrique**

$$2F1([1/3, 2/3], [1/2], 27 z/4)$$

**1, 3, 15, 84, 495, 3003, 18564, 116280, 735471, 4686825, 30045015, 193536720, 1251677700, 8122425444, 52860229080, 344867425584, 2254848913647, 14771069086725**

## Binomial coefficients (4n,n)

**Réf.** AS1 828. dek.

**HIS2 A5810** **hypergéométrique-LLL** **suite P-récurrente**

**HIS1** **algébrique**

$$_3F_2([1/2,\ 3/4,\ 1/4],[2/3,\ 1/3],256\ z/27)$$

**1, 4, 28, 220, 1820, 15504, 134596, 1184040, 10518300, 94143280, 847660528, 7669339132, 69668534468, 635013559600, 5804731963800, 53194089192720, 488526937079580**

---

**Réf.** JCT A43 1 1986.

**HIS2 A5817** **P-récurrences** **Suite P-récurrente**

**HIS1**

$$(n + 4)\ (n + 3)\ a(n) =$$

$$(8\ n + 12)\ a(n - 1) + (16\ n^2 - 16\ n)\ a(n - 2)$$

**1, 2, 4, 10, 25, 70, 196, 588, 1764, 5544, 17424, 56628, 184041, 613470, 2044900, 6952660, 23639044, 81662152, 282105616, 987369656, 3455793796, 12228193432, 43268992144**

## Spanning trees in third power of cycle

**Réf.** FQ 23 258 85.

**HIS2 A5822** **Approximants de Padé**

**HIS1** **Fraction rationnelle**

$$\frac{(1 - z)(1 + z)(z^4 + z^3 - z^2 + z + 1)}{z^8 - 4z^6 - z^4 - 4z^2 + 1}$$

**1, 1, 2, 4, 11, 16, 49, 72, 214, 319, 947, 1408, 4187, 6223, 18502, 27504, 81769, 121552, 361379, 537196**

---

**Réf.** JSC 10 599 90.

**HIS2 A5824** **Approximants de Padé**

**HIS1** **Fraction rationnelle**

$$\frac{z(1 + 2z)(1 - z)}{1 - 5z^2 + 2z^4}$$

**0, 1, 1, 3, 5, 13, 23, 59, 105, 269, 479, 1227, 2185, 5597, 9967, 25531, 45465, 116461, 207391, 531243, 946025, 2423293, 4315343, 11053979, 19684665, 50423309, 89792639**

## Worst case of a Jacobi symbol algorithm

**Réf.** JSC 10 605 90.

**HIS2 A5825** **Approximants de Padé**

**HIS1** **Fraction rationnelle**

$$\frac{z\,(1 + 2z - 4z^2)}{(1 - 2z^2)\,(1 - 5z + 2z^2)}$$

**0, 1, 7, 31, 145, 659, 3013, 13739, 62685, 285931**

## Worst case of a Jacobi symbol algorithm

**Réf.** JSC 10 605 90.

**HIS2 A5826** **Approximants de Padé**

**HIS1** **Fraction rationnelle**

$$\frac{1 + 6z^2 - 4z^3}{(1 - 2z^2)\,(1 - 5z + 2z^2)}$$

**1, 5, 31, 141, 659, 3005, 13739, 62669, 285931, 1304285**

## Worst case of a Jacobi symbol algorithm

**Réf.** JSC 10 605 90.

**HIS2 A5827** **Approximants de Padé**

**HIS1** **Fraction rationnelle**

$$\frac{1 - 2z - 2z^2 + 2z^3}{(1 - 2z^2)(1 - 5z + 2z^2)}$$

**1, 3, 13, 57, 259, 1177, 5367, 24473, 111631, 509193**

---

**Réf.** ST89.

**HIS2 A5840** **Recoupements**

**HIS1** **exponentielle**

$$\frac{\exp(z)(1 - z)}{2 - \exp(z)}$$

**1, 1, 2, 8, 46, 332, 2874, 29024, 334982, 4349492, 62749906, 995818760, 17239953438, 323335939292, 6530652186218, 141326092842416, 3262247252671414, 80009274870905732**

## Packing a square with squares of sides 1...n

**Réf.** GA77 147. UPG D5.

**HIS2 A5842** **Euler** **Conjecture**

**HIS1** **Produit infini**

$$\frac{(1-z^{2})(1-z^{9})(1-z^{11})(1-z^{13})(1-z^{15})\ldots}{(1-z)^{3}(1-z^{8})(1-z^{10})(1-z^{12})(1-z^{14})(1-z^{16})\ldots}$$

**1, 3, 5, 7, 9, 11, 13, 15, 18, 21, 24, 27, 30, 33, 36, 39, 43**

## The even numbers

**Réf.**

**HIS2 A5843** **Approximants de Padé**

**HIS1** **Fraction rationnelle**

$$\frac{2}{(z-1)^{2}}$$

**2, 4, 6, 8, 10, 12, 14, 16, 18, 20, 22, 24, 26, 28, 30, 32, 34, 36, 38, 40, 42, 44, 46, 48, 50, 52, 54, 56, 58, 60, 62, 64, 66, 68, 70, 72, 74, 76, 78, 80, 82, 84, 86, 88, 90, 92, 94, 96, 98, 100, 102, 104**

## Theta series of b.c.c. lattice w.r.t. short edge

**Réf.** JCP 83 6526 85.

**HIS2 A5869** **Euler**

**HIS1** **Produit infini**

*** Le motif [3, -3] est périodique**

$$\prod_{n\geq 1} \frac{1}{(1-z^n)^{c(n)}}$$

$$c(n) = 3, -3, \ldots *$$

**2, 6, 6, 8, 12, 6, 12, 18, 6, 14, 18, 12, 18, 18, 12, 12, 30, 18, 14, 24, 6, 30, 30, 12, 24, 24, 18, 24, 30, 12, 26, 42, 24, 12, 30, 18, 24, 48, 18, 36, 24, 18, 36, 30, 24, 26, 48, 18, 30, 48, 12, 36, 54**

## Theta series of cubic lattice

**Réf.** SPLAG 107.

**HIS2 A5875** **Euler**

**HIS1** **Produit infini**

*** Le motif [6, -9, 6, -3] est périodique**

$$\prod_{n\geq 1} \frac{1}{(1-z^n)^{c(n)}}$$

$$c(n) = 6, -9, 6, -3, \ldots *$$

**1, 6, 12, 8, 6, 24, 24, 0, 12, 30, 24, 24, 8, 24, 48, 0, 6, 48, 36, 24, 24, 48, 24, 0, 24, 30, 72, 32, 0, 72, 48, 0, 12, 48, 48, 48, 30, 24, 72, 0, 24, 96, 48, 24, 24, 72, 48, 0, 8, 54, 84, 48, 24, 72, 96**

## Theta series of cubic lattice w.r.t. edge

**Réf.** SPLAG 107.

**HIS2 A5876** **Euler**

**HIS1** **Produit infini**

*** Le motif [4, -5, 4, -3] est périodique**

$$\prod_{n\geq 1} \frac{1}{(1 - z^n)^{c(n)}}$$

$$c(n) = 4,-5,4,-3,...*$$

**2, 8, 10, 8, 16, 16, 10, 24, 16, 8, 32, 24, 18, 24, 16, 24, 32, 32, 16, 32, 34, 16, 48, 16, 16, 56, 32, 24, 32, 40, 26, 48, 48, 16, 32, 32, 32, 56, 48, 24, 64, 32, 26, 56, 16, 40, 64, 64, 16, 40, 48, 32**

## Theta series of cubic lattice w.r.t. square

**Réf.** SPLAG 107.

**HIS2 A5877** **Euler**

**HIS1** **Produit infini**

*** Le motif [2, -1, 2, -3] est périodique**

$$\prod_{n\geq 1} \frac{1}{(1 - z^n)^{c(n)}}$$

$$c(n) = 2,-1,2,-3,...*$$

**4, 8, 8, 16, 12, 8, 24, 16, 16, 24, 16, 16, 28, 32, 8, 32, 32, 16, 40, 16, 16, 40, 40, 32, 36, 16, 24, 48, 32, 24, 40, 48, 16, 56, 32, 16, 64, 40, 32, 32, 36, 40, 48, 48, 32, 48, 48, 16, 80, 40, 24, 80**

## Theta series of $D_4$ lattice w.r.t. deep hole

**Réf.** SPLAG 118.

**HIS2 A5879** **Euler**

**HIS1** **Produit infini**

*** Le motif [4, -4] est périodique**

$$\prod_{n\geq 1} \frac{1}{(1-z^n)^{c(n)}}$$

$$c(n) = 4, -4, \ldots *$$

**8, 32, 48, 64, 104, 96, 112, 192, 144, 160, 256, 192, 248, 320, 240, 256, 384, 384, 304, 448, 336, 352, 624, 384, 456, 576, 432, 576, 640, 480, 496, 832, 672, 544, 768, 576, 592, 992, 768, 640**

## Theta series of $D_4$ lattice w.r.t. edge

**Réf.**

**HIS2 A5880** **Euler**

**HIS1** **Produit infini**

*** Le motif [4,-4] est périodique**

$$\prod_{n\geq 1} \frac{1}{(1-z^n)^{c(n)}}$$

$$c(n) = 4, -4, \ldots *$$

**2, 8, 12, 16, 26, 24, 28, 48, 36, 40, 64, 48, 62, 80, 60, 64, 96, 96, 76, 112, 84, 88, 156, 96, 114, 144, 108, 144, 160, 120, 124, 208, 168, 136, 192, 144, 148, 248, 192, 160, 242, 168, 216, 240**

## Theta series of planar hexagonal lattice with respect to edge

**Réf.** JCP 83 6523 85.

**HIS2 A5881** **Euler**

**HIS1** **Produit infini**

*** Le motif [1, -1, 2, -1, 1, -2] est périodique**

$$\prod_{n\geq 1} \frac{1}{(1-z^n)^{c(n)}}$$

$$c(n) = 1,-1,2,-1,1,-2,...^*$$

**2, 2, 0, 4, 2, 0, 4, 0, 0, 4, 4, 0, 2, 2, 0, 4, 0, 0, 4, 4, 0, 4, 0, 0, 6, 0, 0, 0, 4, 0, 4, 4, 0, 4, 0, 0, 4, 2, 0, 4, 2, 0, 0, 0, 0, 8, 4, 0, 4, 0, 0, 4, 0, 0, 4, 4, 0, 0, 4, 0, 2, 0, 0, 4, 4, 0, 8, 0, 0, 4, 0, 0, 0, 6**

## Theta series of planar hexagonal lattice w.r.t. deep hole

**Réf.** JCP 83 6524 85.

**HIS2 A5882** **Euler**

**HIS1** **Produit infini**

*** Le motif [1,1,-2] est périodique**

$$\prod_{n\geq 1} \frac{1}{(1-z^n)^{c(n)}}$$

$$c(n) = 1,1,-2,...^*$$

**3, 3, 6, 0, 6, 3, 6, 0, 3, 6, 6, 0, 6, 0, 6, 0, 9, 6, 0, 0, 6, 3, 6, 0, 6, 6, 6, 0, 0, 0, 12, 0, 6, 3, 6, 0, 6, 6, 0, 0, 3, 6, 6, 0, 12, 0, 6, 0, 0, 6, 6, 0, 6, 0, 6, 0, 9, 6, 6, 0, 6, 0, 0, 0, 6, 9, 6, 0, 0, 6, 6, 0, 12, 0, 6, 0, 6**

## Theta series of f.c.c. lattice w.r.t. edge

**Réf.** JCP 83 6526 85.

**HIS2 A5884** **Euler**

**HIS1** **Produit infini**

*** Le motif [2, -1, 2, -3] est périodique**

$$\prod_{n\geq 1} \frac{1}{(1-z^n)^{c(n)}}$$

$$c(n) = 2,-1,2,-3,...^*$$

**2, 4, 4, 8, 6, 4, 12, 8, 8, 12, 8, 8, 14, 16, 4, 16, 16, 8, 20, 8, 8, 20, 20, 16, 18, 8, 12, 24, 16, 12, 20, 24, 8, 28, 16, 8, 32, 20, 16, 16, 18, 20, 24, 24, 16, 24, 24, 8, 40, 20, 12, 40, 16, 12, 20**

---

## Theta series of f.c.c. lattice w.r.t. tetrahedral hole

**Réf.** JCP 83 6526 85.

**HIS2 A5886** **Euler**

**HIS1** **Produit infini**

*** Le motif [3,-3] est périodique**

$$\prod_{n\geq 1} \frac{1}{(1-z^n)^{c(n)}}$$

$$c(n) = 3,-3,...^*$$

**4, 12, 12, 16, 24, 12, 24, 36, 12, 28, 36, 24, 36, 36, 24, 24, 60, 36, 28, 48, 12, 60, 60, 24, 48, 48, 36, 48, 60, 24, 52, 84, 48, 24, 60, 36, 48, 96, 36, 72, 48, 36, 72, 60, 48, 52, 96, 36, 60, 96**

## Centered pentagonal numbers

**Réf.** INOC 24 4550 85.

**HIS2 A5891** **Approximants de Padé**

**HIS1** **Fraction rationnelle**

$$\frac{z^2 + 3z + 1}{(1 - z)^3}$$

**1, 6, 16, 31, 51, 76, 106, 141, 181, 226, 276, 331, 391, 456, 526, 601, 681, 766, 856, 951, 1051, 1156, 1266, 1381, 1501, 1626, 1756, 1891, 2031, 2176, 2326, 2481, 2641, 2806, 2976**

---

## Square octagonal numbers

**Réf.** INOC 24 4550 85.

**HIS2 A5892** **Approximants de Padé**

**HIS1** **Fraction rationnelle**

$$\frac{1 + 9z + 4z^2}{(1 - z)^3}$$

**1, 12, 37, 76, 129, 196, 277, 372, 481, 604, 741, 892, 1057, 1236, 1429, 1636, 1857, 2092, 2341, 2604, 2881, 3172, 3477, 3796, 4129, 4476, 4837, 5212, 5601, 6004, 6421, 6852, 7297**

## Points on surface of tetrahedron

**Réf.** MF73 46. CO74. INOC 24 4550 85.

**HIS2 A5893** **Approximants de Padé**

**HIS1** **Fraction rationnelle**

$$\frac{(1 + z)(1 + z^2)}{(1 - z)^3}$$

**1, 4, 10, 20, 34, 52, 74, 100, 130, 164, 202, 244, 290, 340, 394, 452, 514, 580, 650, 724, 802, 884, 970, 1060, 1154, 1252, 1354, 1460, 1570, 1684, 1802, 1924, 2050, 2180, 2314, 2452, 2594**

## Centered tetrahedral numbers

**Réf.** INOC 24 4550 85.

**HIS2 A5894** **Approximants de Padé**

**HIS1** **Fraction rationnelle**

$$\frac{(1 + z)(1 + z^2)}{(z - 1)^4}$$

**1, 5, 15, 35, 69, 121, 195, 295, 425, 589, 791, 1035, 1325, 1665, 2059, 2511, 3025, 3605, 4255, 4979, 5781, 6665, 7635, 8695, 9849, 11101, 12455, 13915, 15485, 17169, 18971, 20895**

## Points on surface of cube

**Réf.** MF73 46. CO74. INOC 24 4550 85.

**HIS2 A5897** **Approximants de Padé**

**HIS1** **Fraction rationnelle**

$$\frac{(1 + z)\,(1 + 4z + z^2)}{(1 - z)^3}$$

**1, 8, 26, 56, 98, 152, 218, 296, 386, 488, 602, 728, 866, 1016, 1178, 1352, 1538, 1736, 1946, 2168, 2402, 2648, 2906, 3176, 3458, 3752, 4058, 4376, 4706, 5048, 5402, 5768, 6146, 6536**

## Centered cube numbers

**Réf.** AMM 82 819 75. INOC 24 4550 85.

**HIS2 A5898** **Approximants de Padé**

**HIS1** **Fraction rationnelle**

$$\frac{(1 + z)\,(1 + 4z + z^2)}{(z - 1)^4}$$

**1, 9, 35, 91, 189, 341, 559, 855, 1241, 1729, 2331, 3059, 3925, 4941, 6119, 7471, 9009, 10745, 12691, 14859, 17261, 19909, 22815, 25991, 29449, 33201, 37259, 41635, 46341, 51389, 56791**

## Points on surface of octahedron

**Réf.** MF73 46. CO74. INOC 24 4550 85.

**HIS2 A5899** **Approximants de Padé**

**HIS1** **Fraction rationnelle**

$$\frac{(1 + z)^3}{(1 - z)^3}$$

**1, 6, 18, 38, 66, 102, 146, 198, 258, 326, 402, 486, 578, 678, 786, 902, 1026, 1158, 1298, 1446, 1602, 1766, 1938, 2118, 2306, 2502, 2706, 2918, 3138, 3366, 3602, 3846, 4098, 4358, 4626**

## Octahedral numbers

**Réf.** CO74. INOC 24 4550 85.

**HIS2 A5900** **Approximants de Padé**

**HIS1** **Fraction rationnelle**

$$\frac{(1 + z)^2}{(z - 1)^4}$$

**1, 6, 19, 44, 85, 146, 231, 344, 489, 670, 891, 1156, 1469, 1834, 2255, 2736, 3281, 3894, 4579, 5340, 6181, 7106, 8119, 9224, 10425, 11726, 13131, 14644, 16269, 18010, 19871, 21856**

## Points on surface of cuboctahedron (or icosahedron)

**Réf.** RO69 109. MF73 46. CO74. INOC 24 4550 85.

**HIS2 A5901** **Approximants de Padé**

**HIS1** **Fraction rationnelle**

$$\frac{(1 + z)\,(z^2 + 8\,z + 1)}{(1 - z)^3}$$

**1, 12, 42, 92, 162, 252, 362, 492, 642, 812, 1002, 1212, 1442, 1692, 1962, 2252, 2562, 2892, 3242, 3612, 4002, 4412, 4842, 5292, 5762, 6252, 6762, 7292, 7842, 8412, 9002, 9612, 10242, 10892**

## Centered icosahedral (or cuboctahedral) numbers

**Réf.** INOC 24 4550 85.

**HIS2 A5902** **Approximants de Padé**

**HIS1** **Fraction rationnelle**

$$\frac{(1 + z)\,(z^2 + 8\,z + 1)}{(z - 1)^4}$$

**1, 13, 55, 147, 309, 561, 923, 1415, 2057, 2869, 3871, 5083, 6525, 8217, 10179, 12431, 14993, 17885, 21127, 24739, 28741, 33153, 37995, 43287, 49049, 55301, 62063, 69355, 77197, 85609**

## Points on surface of dodecahedron

**Réf.** INOC 24 4550 85.

**HIS2 A5903** **Approximants de Padé**

**HIS1** **Fraction rationnelle**

$$\frac{(1 + z)\,(z^2 + 28\,z + 1)}{(1 - z)^3}$$

**1, 32, 122, 272, 482, 752, 1082, 1472, 1922, 2432, 3002, 3632, 4322, 5072, 5882, 6752, 7682, 8672, 9722, 10832, 12002, 13232, 14522, 15872, 17282, 18752, 20282, 21872, 23522, 25232**

## Centered dodecahedral numbers

**Réf.** INOC 24 4550 85.

**HIS2 A5904** **Approximants de Padé**

**HIS1** **Fraction rationnelle**

$$\frac{(1 + z)\,(z^2 + 28\,z + 1)}{(z - 1)^4}$$

**1, 33, 155, 427, 909, 1661, 2743, 4215, 6137, 8569, 11571, 15203, 19525, 24597, 30479, 37231, 44913, 53585, 63307, 74139, 86141, 99373, 113895, 129767, 147049, 165801, 186083**

## Points on surface of truncated tetrahedron

**Réf.** CO74. INOC 24 4552 85.

**HIS2 A5905** **Approximants de Padé**

**HIS1** **Fraction rationnelle**

$$\frac{(1 + z)\,(z^2 + 12\,z + 1)}{(1 - z)^3}$$

**1, 16, 58, 128, 226, 352, 506, 688, 898, 1136, 1402, 1696, 2018, 2368, 2746, 3152, 3586, 4048, 4538, 5056, 5602, 6176, 6778, 7408, 8066, 8752, 9466, 10208, 10978, 11776, 12602, 13456**

## Truncated tetrahedral numbers

**Réf.** CO74. INOC 24 4552 85.

**HIS2 A5906** **Approximants de Padé**

**HIS1** **Fraction rationnelle**

$$\frac{1 + 12\,z + 10\,z^2}{(z - 1)^4}$$

**1, 16, 68, 180, 375, 676, 1106, 1688, 2445, 3400, 4576, 5996, 7683, 9660, 11950, 14576, 17561, 20928, 24700, 28900, 33551, 38676, 44298, 50440, 57125, 64376, 72216, 80668, 89755**

## Truncated octahedral numbers

**Réf.** CO74. INOC 24 4552 85.

**HIS2 A5910** **Approximants de Padé**

**HIS1** **Fraction rationnelle**

$$\frac{1 + 34 z + 55 z^{2} + 6 z^{3}}{(z - 1)^{4}}$$

**1, 38, 201, 586, 1289, 2406, 4033, 6266, 9201, 12934, 17561, 23178, 29881, 37766, 46929, 57466, 69473, 83046, 98281, 115274, 134121, 154918, 177761, 202746, 229969, 259526, 291513**

## Points on surface of truncated cube

**Réf.** INOC 24 4552 85.

**HIS2 A5911** **Approximants de Padé**

**HIS1** **Fraction rationnelle**

$$\frac{(1 + z) (z^{2} + 44 z + 1)}{(1 - z)^{3}}$$

**1, 48, 186, 416, 738, 1152, 1658, 2256, 2946, 3728, 4602, 5568, 6626, 7776, 9018, 10352, 11778, 13296, 14906, 16608, 18402, 20288, 22266, 24336, 26498, 28752, 31098, 33536, 36066**

## Truncated cube numbers

**Réf.** INOC 24 4552 85.

**HIS2 A5912** **Approximants de Padé**

**HIS1** **Fraction rationnelle**

$$\frac{1 + 52 z + 93 z^2 + 8 z^3}{(z - 1)^4}$$

**1, 56, 311, 920, 2037, 3816, 6411, 9976, 14665, 20632, 28031, 37016, 47741, 60360, 75027, 91896, 111121, 132856, 157255, 184472, 214661, 247976, 284571, 324600, 368217, 415576, 466831**

## Points on surface of hexagonal prism

**Réf.** INOC 24 4552 85.

**HIS2 A5914** **Approximants de Padé**

**HIS1** **Fraction rationnelle**

$$\frac{(1 + z) (z^2 + 10 z + 1)}{(1 - z)^3}$$

**1, 14, 50, 110, 194, 302, 434, 590, 770, 974, 1202, 1454, 1730, 2030, 2354, 2702, 3074, 3470, 3890, 4334, 4802, 5294, 5810, 6350, 6914, 7502, 8114, 8750, 9410, 10094, 10802, 11534, 12290**

## Hexagonal prism numbers

**Réf.** INOC 24 4552 85.

**HIS2 A5915** **Approximants de Padé**

**HIS1** **Fraction rationnelle**

$$\frac{1 + 10 z + 7 z^{2}}{(z - 1)^{4}}$$

**1, 14, 57, 148, 305, 546, 889, 1352, 1953, 2710, 3641, 4764, 6097, 7658, 9465, 11536, 13889, 16542, 19513, 22820, 26481, 30514, 34937, 39768, 45025, 50726, 56889, 63532, 70673, 78330**

## Rhombic dodecahedral numbers

**Réf.** AMM 82 819 75. INOC 24 4552 85.

**HIS2 A5917** **Approximants de Padé**

**HIS1** **Fraction rationnelle**

$$\frac{(1 + z) (z^{2} + 10 z + 1)}{(z - 1)^{4}}$$

**1, 15, 65, 175, 369, 671, 1105, 1695, 2465, 3439, 4641, 6095, 7825, 9855, 12209, 14911, 17985, 21455, 25345, 29679, 34481, 39775, 45585, 51935, 58849, 66351, 74465, 83215, 92625**

## Points on surface of square pyramid

**Réf.** CO74. INOC 24 4552 85.

**HIS2 A5918** **Approximants de Padé**

**HIS1** **Fraction rationnelle**

$$\frac{(1 + z)\,(z^2 + z + 1)}{(1 - z)^3}$$

**1, 5, 14, 29, 50, 77, 110, 149, 194, 245, 302, 365, 434, 509, 590, 677, 770, 869, 974, 1085, 1202, 1325, 1454, 1589, 1730, 1877, 2030, 2189, 2354, 2525, 2702, 2885, 3074, 3269, 3470, 3677**

## Points on surface of tricapped prism

**Réf.** INOC 24 4552 85.

**HIS2 A5919** **Approximants de Padé**

**HIS1** **Fraction rationnelle**

$$\frac{(1 + z)\,(z^2 + 5\,z + 1)}{(1 - z)^3}$$

**1, 9, 30, 65, 114, 177, 254, 345, 450, 569, 702, 849, 1010, 1185, 1374, 1577, 1794, 2025, 2270, 2529, 2802, 3089, 3390, 3705, 4034, 4377, 4734, 5105, 5490, 5889, 6302, 6729, 7170, 7625**

## Tricapped prism numbers

**Réf.** INOC 24 4552 85.

**HIS2 A5920** **Approximants de Padé**

**HIS1** **Fraction rationnelle**

$$\frac{1 + 5z + 3z^2}{(z - 1)^4}$$

**1, 9, 33, 82, 165, 291, 469, 708, 1017, 1405, 1881, 2454, 3133, 3927, 4845, 5896, 7089, 8433, 9937, 11610, 13461, 15499, 17733, 20172, 22825, 25701, 28809, 32158, 35757, 39615, 43741**

## From solution to a difference equation

**Réf.** FQ 25 363 87.

**HIS2 A5921** **Dérivée logarithmique** **F.G. exponentielle**

**HIS1** **Fraction rationnelle**

$$\frac{(z + 1)^2}{z^2 - z + 1}$$

**1, 3, 10, 48, 312, 2520, 24480, 277200, 3588480, 52254720**

## n-step mappings with 4 inputs

**Réf.** PRV A32 2342 85.

**HIS2 A5945** **Approximants de Padé** **Conjecture**

**HIS1** **exponentielle**

$$\exp(z)\ (1 + 14 z + 31/2 z^2 + 3 z^3)$$

**1, 15, 60, 154, 315, 561, 910**

## Sum of cubes of Fibonacci numbers

**Réf.** BR72 18.

**HIS2 A5968** **Approximants de Padé**

**HIS1** **Fraction rationnelle**

$$\frac{1 - 2 z - z^2}{(z - 1)(1 - 4 z - z^2)(z^2 - z - 1)}$$

**1, 2, 10, 37, 162, 674, 2871, 12132, 51436, 217811, 922780, 3908764, 16558101, 70140734, 297121734, 1258626537, 5331629710, 22585142414, 95672204155, 405273951280**

## Sum of fourth powers of Fibonacci numbers

**Réf.** BR72 19.

**HIS2 A5969** **Approximants de Padé**

**HIS1** **Fraction rationnelle**

$$\frac{(1 + z)\,(1 - 5 z + z^2)}{(z^2 - 7 z + 1)\,(z^2 + 3 z + 1)\,(z - 1)^2}$$

**1, 2, 18, 99, 724, 4820, 33381, 227862, 1564198, 10714823, 73457064, 503438760, 3450734281, 23651386922, 162109796922, 1111115037483, 7615701104764, 52198777931900**

## Sum of squares of Lucas numbers

**Réf.** BR72 20.

**HIS2 A5970** **Approximants de Padé**

**HIS1** **Fraction rationnelle**

$$\frac{1 + 7 z - 4 z^2}{(1 - z)\,(1 + z)\,(1 - 3 z + z^2)}$$

**1, 10, 26, 75, 196, 520, 1361, 3570, 9346, 24475, 64076, 167760, 439201, 1149850, 3010346, 7881195, 20633236, 54018520, 141422321, 370248450, 969323026, 2537720635**

## Sum of cubes of Lucas numbers

**Réf.** BR72 21.

**HIS2 A5971** **Approximants de Padé**

**HIS1** **Fraction rationnelle**

$$\frac{1 + 24 z - 23 z^2 - 8 z^3}{(z - 1) (1 - 4 z - z^2) (z^2 - z - 1)}$$

**1, 28, 92, 435, 1766, 7598, 31987, 135810, 574786, 2435653, 10316252, 43702500, 185123261, 784200368, 3321916912, 14071880655, 59609419066, 252509590018, 1069647725567**

## Sum of fourth powers of Lucas numbers

**Réf.** BR72 21.

**HIS2 A5972** **Approximants de Padé**

**HIS1** **Fraction rationnelle**

$$\frac{1 + 76 z - 164 z^2 - 79 z^3 + 16 z^4}{(z^2 - 7 z + 1) (z^2 + 3 z + 1) (z - 1)^2}$$

**1, 82, 338, 2739, 17380, 122356, 829637, 5709318, 39071494, 267958135, 1836197336, 12586569192, 86266785673, 591288786874, 4052734152890, 27777904133691**

## Longest walk on edges of n-cube

**Réf.** clm.

**HIS2 A5985** **Approximants de Padé**

**HIS1** **Fraction rationnelle**

$$\frac{1 + 2z - 4z^2 + 4z^3}{(1 - z)(1 + 2z)(1 + z)(2z - 1)^2}$$

**1, 4, 9, 32, 65, 192, 385, 1024, 2049, 5120, 10241, 24576, 49153, 114688, 229377, 524288, 1048577, 2359296, 4718593, 10485760, 20971521, 46137344, 92274689, 201326592**

---

## Column-strict plane partitions of n

**Réf.** SAM 50 260 71.

**HIS2 A5986** **Euler**

**HIS1** **Produit infini**

$$\prod_{n\geq 1} \frac{1}{(1 - z^n)^{c(n)}}$$

$$c(n) = 2,2,3,3,4,4,5,5,...$$

**1, 2, 5, 11, 23, 45, 87, 160, 290, 512, 889, 1514, 2547, 4218, 6909, 11184, 17926, 28449, 44772, 69862, 108205, 166371, 254107, 385617, 581729, 872535, 1301722, 1932006, 2853530**

## Symmetric plane partitions of n

**Réf.** SAM 50 261 71.

**HIS2 A5987** **Euler**

**HIS1** **Produit infini**

*** c(n) = 1 si n est impair et [n/4] si n est pair.**

$$\prod_{n\geq 1} \frac{1}{(1 - z^n)^{c(n)}}$$

**c(n) = 1,0,1,1,1,1,1,2,1,2,1,...***

**1, 1, 1, 2, 3, 4, 6, 8, 12, 16, 22, 29, 41, 53, 71, 93, 125, 160, 211, 270, 354, 450, 581, 735, 948, 1191, 1517, 1902, 2414, 3008, 3791, 4709, 5909, 7311, 9119, 11246, 13981, 17178, 21249**

---

## Paraffins

**Réf.** BER 30 1919 1897.

**HIS2 A5993** **Euler**

**HIS1** **Fraction rationnelle**

$$\frac{1 - z^4}{(1 - z)^2 (1 - z^2)^3}$$

**1, 2, 6, 10, 19, 28, 44, 60, 85, 110**

## Paraffins

**Réf.** BER 30 1919 1897.

**HIS2 A5994** **Euler**

**HIS1** **Fraction rationnelle**

$$\frac{1 - z^{4}}{(1 - z)^{3} (1 - z^{2})^{3}}$$

**1, 3, 9, 19, 38, 66, 110, 170, 255, 365**

## Paraffins

**Réf.** BER 30 1919 1897.

**HIS2 A5995** **Euler**

**HIS1** **Fraction rationnelle**

$$\frac{(1 - z^{4})^{6} (1 - z^{8})^{18}}{(1 - z)^{3} (1 - z^{2})^{6} (1 - z^{6})^{8}}$$

**1, 3, 12, 28, 66, 126, 236, 396, 651, 1001**

## Paraffins

**Réf.** BER 30 1920 1897.

**HIS2 A5996** **Euler**

**HIS1** **Fraction rationnelle**

$$\frac{1 - z^3}{(1 - z)^3 (1 - z^2)^2}$$

**2, 6, 16, 30, 54, 84, 128, 180, 250, 330**

## Paraffins

**Réf.** BER 30 1922 1897.

**HIS2 A6000** **Approximants de Padé**

**HIS1** **Fraction rationnelle**

$$\frac{1 + 2 z^2}{(z - 1)^4}$$

**1, 4, 12, 28, 55, 96, 154, 232, 333**

## Paraffins

**Réf.** BER 30 1922 1897.

**HIS2 A6001** **Approximants de Padé**

**HIS1** **Fraction rationnelle**

$$\frac{1 + 2z^{3}}{(z - 1)^{4}}$$

**1, 4, 10, 22, 43, 76, 124, 190, 277**

## Paraffins

**Réf.** BER 30 1922 1897.

**HIS2 A6003** **Approximants de Padé**

**HIS1** **Fraction rationnelle**

$$\frac{1 - z^{3}}{(1 - z)^{5}}$$

**1, 5, 15, 34, 65, 111, 175, 260**

## Paraffins

**Réf.** BER 30 1922 1897.

**HIS2 A6004** **Euler**

**HIS1** **Fraction rationnelle**

$$\frac{(1 - z^4)(1 - z^5)(1 - z^6)}{(1 - z)^4 (1 - z^2)(1 - z^3)(1 - z^7)}$$

**1, 4, 11, 25, 49, 86, 139, 211**

## Paraffins

**Réf.** BER 30 1923 1897.

**HIS2 A6007** **Euler**

**HIS1** **Fraction rationnelle**

$$\frac{1 - z^4}{(1 - z)^5 (1 - z^2)}$$

**1, 5, 16, 40, 85, 161, 280, 456**

## Paraffins

**Réf.** BER 30 1923 1897. GA66 246.

**HIS2 A6008** **Approximants de Padé**

**HIS1** **Fraction rationnelle**

$$\frac{z\,(1+z)\,(1-z+z^2)}{(1-z)^5}$$

**0, 1, 5, 15, 36, 75, 141, 245, 400, 621, 925, 1331, 1860, 2535, 3381, 4425, 5696, 7225, 9045, 11191, 13700, 16611, 19965, 23805, 28176, 33125, 38701, 44955, 51940, 59711, 68325**

---

## Paraffins

**Réf.** BER 30 1923 1897.

**HIS2 A6011** **Approximants de Padé**

**HIS1** **Fraction rationnelle**

$$\frac{1+z}{(1-z)^5}$$

**3, 18, 60, 150, 315, 588, 1008, 1620**

**Réf.** GK90 86.

**HIS2 A6012** **Approximants de Padé**

**HIS1** **Fraction rationnelle**

$$\frac{1 - 2z}{1 - 4z + 2z^2}$$

**1, 2, 6, 20, 68, 232, 792, 2704, 9232, 31520, 107616, 367424, 1254464, 4283008, 14623104, 49926400, 170459392, 581984768, 1987020288, 6784111616, 23162405888, 79081400320**

---

**Réf.** dek.

**HIS2 A6013** **Inverse fonctionnel algébrique** **Suite P-récurrente**

**HIS1**

$${}_3F_2([1, 4/3, 2/3], [2, 3/2], 27z/4)$$

**1, 2, 7, 30, 143, 728, 3876, 21318, 120175, 690690, 4032015, 23841480, 142498692, 859515920, 5225264024, 31983672534, 196947587823, 1219199353190, 7583142491925, 47365474641870**

Réf. rkg.
HIS2 A6040 P-récurrences Suite P-récurrente
HIS1

$$a(n) = (- n^2 + 4 n - 4) a(n - 2) + (n^2 - 2 n + 2) a(n - 1)$$

1, 2, 9, 82, 1313, 32826, 1181737, 57905114, 3705927297, 300180111058, 30018011105801, 3632179343801922, 523033825507476769, 88392716510763573962

---

Réf. rkg.
HIS2 A6041 P-récurrences Suite P-récurrente
HIS1

$$(n - 1) a(n) = (n^2 - 3 n + 3) n a(n - 1) + (- n^2 + 4 n - 3) n a(n - 2)$$

0, 2, 9, 76, 1145, 27486, 962017, 46176824, 2909139921, 232731193690, 23040388175321, 2764846581038532, 395373061088510089, 66422674262869694966

## A traffic light problem

**Réf.** BIO 46 422 59.

**HIS2 A6043** **Hypergéométrique**

**HIS1** **Fraction rationnelle**

$$\frac{2}{(1 - 3z)^3}$$

**2, 18, 108, 540, 2430**

## Square hex numbers

**Réf.** GA88 19.

**HIS2 A6051** **Approximants de Padé**

**HIS1** **Fraction rationnelle**

$$\frac{1 - 26z + z^2}{(1 - z)(z^2 - 194z + 1)}$$

**1, 169, 32761, 6355441, 1232922769, 239180661721, 46399815451081, 9001325016847969, 1746210653453054881, 338755865444875798921, 65716891685652451935769**

## Triangular star numbers

**Réf.** GA88 20.

**HIS2 A6060** **Approximants de Padé**

**HIS1** **Fraction rationnelle**

$$\frac{1 + 58 z + z^2}{(1 - z) (z^2 - 194 z + 1)}$$

**1, 253, 49141, 9533161, 1849384153, 358770992581, 69599723176621, 13501987525271953, 2619315980179582321, 508133798167313698381, 98575337528478677903653**

## Square star numbers

**Réf.** GA88 22.

**HIS2 A6061** **Approximants de Padé**

**HIS1** **Fraction rationnelle**

$$\frac{z^2 + 22 z + 1}{(1 - z) (z^2 - 98 z + 1)}$$

**1, 121, 11881, 1164241, 114083761, 11179044361, 1095432263641, 107341182792481, 10518340481399521, 1030690025994360601, 100997104206965939401**

## Star-hex numbers

**Réf.** GA88 22. JRM 16 192 83.

**HIS2 A6062** **Approximants de Padé**

**HIS1** **Fraction rationnelle**

$$\frac{(1 + z)^2}{(1 - z)\,(z^2 - 34\,z + 1)}$$

**1, 37, 1261, 42841, 1455337, 49438621, 1679457781, 57052125937, 1938092824081, 65838103892821, 2236557439531837**

## Maximal length rook tour on n X n board

**Réf.** GA86 76.

**HIS2 A6071** **Approximants de Padé**

**HIS1** **Fraction rationnelle**

$$\frac{1 + z + 4\,z^2 + 6\,z^3 - 5\,z^4 + z^5}{(1 + z)\,(z - 1)^4}$$

**1, 4, 14, 38, 76, 136, 218, 330, 472, 652, 870, 1134**

## Gaussian binomial coefficient [n,2] for q=2

**Réf.** GJ83 99. ARS A17 328 84.

**HIS2 A6095** **Approximants de Padé**

**HIS1** **Fraction rationnelle**

$$\frac{1}{(1 - z)(1 - 2z)(1 - 4z)}$$

**1, 7, 35, 155, 651, 2667, 10795, 43435, 174251, 698027, 2794155, 11180715, 44731051, 178940587, 715795115, 2863245995, 11453115051, 45812722347, 183251413675**

## Gaussian binomial coefficient [n,3] for q=2

**Réf.** GJ83 99. ARS A17 328 84.

**HIS2 A6096** **Approximants de Padé**

**HIS1** **Fraction rationnelle**

$$\frac{1}{(1 - z)(1 - 2z)(1 - 4z)(1 - 8z)}$$

**1, 15, 155, 1395, 11811, 97155, 788035, 6347715, 50955971, 408345795, 3269560515, 26167664835, 209386049731, 1675267338435, 13402854502595, 107225699266755, 857817047249091**

## Gaussian binomial coefficient [n,4] for q=2

**Réf.** GJ83 99. ARS A17 328 84.

**HIS2 A6097** **Approximants de Padé**

**HIS1** **Fraction rationnelle**

$$\frac{1}{(1 - z)(1 - 2z)(1 - 4z)(1 - 8z)(1 - 16z)}$$

**1, 31, 651, 11811, 200787, 3309747, 53743987, 866251507, 13910980083, 222984027123, 3571013994483, 57162391576563, 914807651274739, 14638597687734259**

## Gaussian binomial coefficient [n,2] for q=3

**Réf.** GJ83 99. ARS A17 328 84.

**HIS2 A6100** **Approximants de Padé**

**HIS1** **Fraction rationnelle**

$$\frac{1}{(1 - z)(1 - 3z)(1 - 9z)}$$

**1, 13, 130, 1210, 11011, 99463, 896260, 8069620, 72636421, 653757313, 5883904390, 52955405230, 476599444231, 4289397389563, 38604583680520, 347441274648040, 3126971536402441**

## Gaussian binomial coefficient [n,3] for q=3

**Réf.** GJ83 99. ARS A17 328 84.

**HIS2 A6101** **Approximants de Padé**

**HIS1** **Fraction rationnelle**

$$\frac{1}{(1 - z)\ (1 - 3\ z)\ (1 - 9\ z)\ (1 - 27\ z)}$$

**1, 40, 1210, 33880, 925771, 25095280, 678468820, 18326727760, 494894285941, 13362799477720, 360801469802830, 9741692640081640, 263026177881648511, 7101711092201899360**

## Gaussian binomial coefficient [n,4] for q=3

**Réf.** GJ83 99. ARS A17 328 84.

**HIS2 A6102** **Approximants de Padé**

**HIS1** **Fraction rationnelle**

$$\frac{1}{(1 - z)\ (1 - 3\ z)\ (1 - 9\ z)\ (1 - 27\ z)\ (1 - 81\ z)}$$

**1, 121, 11011, 925771, 75913222, 6174066262, 500777836042, 40581331447162, 3287582741506063, 266307564861468823, 21571273555248777493, 1747282899667791058573**

## Gaussian binomial coefficient [n,2] for q=4

**Réf.** GJ83 99. ARS A17 328 84.

**HIS2 A6105** **Approximants de Padé**

**HIS1** **Fraction rationnelle**

$$\frac{1}{(1 - z)(1 - 4z)(1 - 16z)}$$

**1, 21, 357, 5797, 93093, 1490853, 23859109, 381767589, 6108368805, 97734250405, 1563749404581, 25019996065701, 400319959420837, 6405119440211877, 102481911401303973**

## Gaussian binomial coefficient [n,3] for q=4

**Réf.** GJ83 99. ARS A17 328 84.

**HIS2 A6106** **Approximants de Padé**

**HIS1** **Fraction rationnelle**

$$\frac{1}{(1 - z)(1 - 4z)(1 - 16z)(1 - 64z)}$$

**1, 85, 5797, 376805, 24208613, 1550842085, 99277752549, 6354157930725, 406672215935205, 26027119554103525, 1665737215212030181, 106607206793565997285**

## Gaussian binomial coefficient [n,5] for q=2

**Réf.** GJ83 99. ARS A17 328 84.

**HIS2 A6110** **Approximants de Padé**

**HIS1** **Fraction rationnelle**

$$\frac{1}{(1 - z)(1 - 2z)(1 - 4z)(1 - 8z)(1 - 16z)(1 - 32z)}$$

**1, 63, 2667, 97155, 3309747, 109221651, 3548836819, 114429029715, 3675639930963, 117843461817939, 3774561792168531, 120843139740969555, 3867895279362300499**

## Gaussian binomial coefficient [n,2] for q=5

**Réf.** GJ83 99. ARS A17 329 84.

**HIS2 A6111** **Approximants de Padé**

**HIS1** **Fraction rationnelle**

$$\frac{1}{(1 - z)(1 - 5z)(1 - 25z)}$$

**1, 31, 806, 20306, 508431, 12714681, 317886556, 7947261556, 198682027181, 4967053120931, 124176340230306, 3104408566792806, 77610214474995931, 1940255363400777181**

## Gaussian binomial coefficient [n,3] for q=5

**Réf.** GJ83 99. ARS A17 329 84.

**HIS2 A6112** **Approximants de Padé**

**HIS1** **Fraction rationnelle**

$$\frac{1}{(1 - z)(1 - 5z)(1 - 25z)(1 - 125z)}$$

**1, 156, 20306, 2558556, 320327931, 40053706056, 5007031143556, 625886840206056, 78236053707784181, 9779511680526143556, 1222439084242108174806**

---

**Réf.** FQ 15 24 77.

**HIS2 A6130** **Approximants de Padé**

**HIS1** **Fraction rationnelle**

$$\frac{1}{1 - z - 3z^2}$$

**1, 1, 4, 7, 19, 40, 97, 217, 508, 1159, 2683, 6160, 14209, 32689, 75316, 173383, 399331, 919480, 2117473, 4875913, 11228332, 25856071, 59541067, 137109280, 315732481**

**Réf.** FQ 15 24 77.
**HIS2 A6131** **Approximants de Padé**
**HIS1** **Fraction rationnelle**

$$\frac{1}{1 - z - 4z^2}$$

**1, 1, 5, 9, 29, 65, 181, 441, 1165, 2929, 7589, 19305, 49661, 126881, 325525, 833049, 2135149, 5467345, 14007941, 35877321, 91909085, 235418369, 603054709, 1544728185**

---

**Réf.** FQ 11 52 73.
**HIS2 A6138** **Approximants de Padé**
**HIS1** **Fraction rationnelle**

$$\frac{1 + z}{1 - z - 3z^2}$$

**1, 2, 5, 11, 26, 59, 137, 314, 725, 1667, 3842, 8843, 20369, 46898, 108005, 248699, 572714, 1318811, 3036953, 6993386, 16104245, 37084403, 85397138, 196650347, 452841761**

**Réf.** FQ 27 434 89.

**HIS2 A6139** **LLL** **Suite P-récurrente**

**HIS1** **algébrique**

**(n - 1) a(n) = (4 n - 6) a(n - 1) + (4 n - 8) a(n - 2)**

$$\frac{1}{(1 - 4 z - 4 z^2)^{1/2}}$$

**1, 2, 8, 32, 136, 592, 2624, 11776, 53344, 243392, 1116928, 5149696, 23835904, 110690816, 515483648, 2406449152, 11258054144, 52767312896, 247736643584**

---

**Dyck paths**

**Réf.** SC83.

**HIS2 A6149** **Hypergéométrique** **Suite P-récurrente**

**HIS1**

$$_4F_3 ([1, 1/2, 3/2, 5/2], [4, 5, 6], 64 z)$$

**1, 1, 4, 30, 330, 4719, 81796, 1643356, 37119160, 922268360, 24801924512, 713055329720**

## Dyck paths

**Réf.** SC83.

**HIS2 A6150** **Hypergéométrique** **Suite P-récurrente**

**HIS1**

$${}_5F_4\left([1, 1/2, 7/2, 5/2, 3/2], [5, 6, 7, 8], 256\,z\right)$$

**1, 1, 5, 55, 1001, 26026, 884884, 37119160, 1844536720, 105408179176, 6774025632340**

## Dyck paths

**Réf.** SC83.

**HIS2 A6151** **Recoupements** **Suite P-récurrente**

**HIS1**

$${}_6F_5\left([1, 1/2, 3/2, 5/2, 7/2, 9/2], [6, 7, 8, 9, 10], 1024\,z\right)$$

**1, 1, 6, 91, 2548, 111384, 6852768, 553361016, 55804330152, 6774025632340**

## Expansion of z exp(z/(1-z))

**Réf.** ARS 10 142 80.

**HIS2 A6152** **Dérivée logarithmique** **Suite P-récurrente**

**HIS1** **exponentielle**

**a(n) = (2 n - 2) a(n - 1) + (- n^2 + 5 n - 5) a(n - 2) + (- n^2 + 6 n - 8) a(n - 3)**

$$\frac{z^2 - z + 1}{\exp(1/(1-z))\,(z - 1)^2}$$

**1, 2, 9, 52, 365, 3006, 28357, 301064, 3549177, 45965530, 648352001, 9888877692, 162112109029, 2841669616982, 53025262866045, 1049180850990736, 21937381717388657**

---

**Réf.** RAIRO 12 58 78.

**HIS2 A6157** **Dérivée logarithmique** **f.g. exponentielle**

**HIS1** **Fraction rationnelle**

$$\frac{1 + z}{(1 - z)^4}$$

**1, 5, 28, 180, 1320, 10920, 100800, 1028160, 11491200, 139708800, 1836172800, 25945920000, 392302310400, 6320426112000, 108101081088000, 1956280854528000**

## From sum of 1/F(n)

**Réf.** FQ 16 169 78.

**HIS2 A6172** **Approximants de Padé**

**HIS1** **Fraction rationnelle**

**F(n) : Nombres de Fibonacci**

$$\frac{2 + 3 z - 19 z^2 + 17 z^3 - 4 z^4}{(z - 1) (z^2 - z - 1) (1 - 3 z + z^2)}$$

**2, 9, 10, 42, 79, 252, 582, 1645, 4106, 11070, 28459, 75348, 195898**

---

## $\sum$ (k+1)! C(n-2,k)/2 , k=0...n-2

**Réf.** DM 55 272 85.

**HIS2 A6183** **Dérivée logarithmique** **Suite P-récurrente**

**HIS1** **exponentielle**

**a(n) = (1 + n) a(n - 1) + (2 - n) a(n - 2)**

$$\frac{2 \exp(z)}{(1 - z)^2}$$

**2, 6, 22, 98, 522, 3262, 23486, 191802, 1753618, 17755382, 197282022, 2387112466, 31249472282, 440096734638, 6635304614542, 106638824162282, 1819969265702946**

**Réf.** FQ 15 292 77. ARS 6 168 78.

**HIS2 A6190** **Approximants de Padé**

**HIS1** **Fraction rationnelle**

$$\frac{1}{1 - 3 z - z^2}$$

**1, 3, 10, 33, 109, 360, 1189, 3927, 12970, 42837, 141481, 467280, 1543321, 5097243, 16835050, 55602393, 183642229, 606529080, 2003229469, 6616217487, 21851881930**

---

## Partitions into pairs

**Réf.** PLIS 23 65 78.

**HIS2 A6198** **équations différentielles** **Suite P-récurrente**

**HIS1** **exponentielle** **Formule de B. Salvy**

**$a(n) = (2 n - 2) a(n - 1) + (2 n - 4) a(n - 2) + a(n - 3)$**

$$\frac{2 - 2 z - (1 - 2 z)^{1/2}}{(1 - 2 z)^{3/2} \exp(1 - (1 - 2 z)^{1/2})}$$

**1, 1, 6, 41, 365, 3984, 51499, 769159, 13031514, 246925295, 5173842311, 118776068256, 2964697094281, 79937923931761, 2315462770608870, 71705109685449689**

## Partitions into pairs

**Réf.** PLIS 23 65 78.

**HIS2 A6199** **P-récurrences** **Suite P-récurrente**

**HIS1**

$$a(n) = 2\,n\,a(n-1) + (2\,n-6)\,a(n-3) + a(n-4) + (2\,n-3)\,a(n-2)$$

**1, 3, 21, 185, 2010, 25914, 386407, 6539679, 123823305, 2593076255, 59505341676, 1484818160748, 40025880386401, 1159156815431055, 35891098374564105**

## Partitions into pairs

**Réf.** PLIS 23 65 78.

**HIS2 A6200** **P-récurrences** **Suite P-récurrrente**

**HIS1**

$$a(n)\,(n-1) = (2 + 6\,n^2 - 2\,n)\,a(n-2) + (-6 + 2\,n + 2\,n^2)\,a(n-1) - n\,a(n-3)$$

**1, 6, 55, 610, 7980, 120274, 2052309, 39110490, 823324755, 18974858540, 475182478056, 12848667150956, 373081590628565, 11578264139795430, 382452947343624515**

## From continued fraction for Zeta(3)

**Réf.** LNM 751 68 79.

**HIS2 A6221** **Approximants de Padé**

**HIS1** **Fraction rationnelle**

$$\frac{z\ (1 + z)\ (5\ z^2 + 92\ z + 5)}{(z - 1)^4}$$

**0, 5, 117, 535, 1463, 3105, 5665, 9347, 14355, 20893, 29165, 39375, 51727, 66425, 83673, 103675, 126635, 152757, 182245, 215303, 252135, 292945, 337937, 387315, 441283, 500045**

---

**Réf.** LNM 751 68 79.

**HIS2 A6222** **Approximants de Padé**

**HIS1** **Fraction rationnelle**

$$\frac{3 + 16\ z + 3\ z^2}{(1 - z)^3}$$

**3, 25, 69, 135, 223, 333, 465, 619, 795, 993, 1213, 1455, 1719, 2005, 2313, 2643, 2995, 3369, 3765, 4183, 4623, 5085, 5569, 6075, 6603, 7153, 7725, 8319, 8935, 9573, 10233, 10915**

## Binary trees of height n requiring 3 registers

**Réf.** TCS 9 105 79.

**HIS2 A6223** **Approximants de Padé**

**HIS1** **Fraction rationnelle**

$$\frac{1}{(2z-1)(2z^4-16z^3+20z^2-8z+1)(1-4z+2z^2)}$$

**1, 14, 118, 780, 4466, 23276, 113620, 528840, 2375100, 10378056, 44381832, 186574864, 773564328, 3171317360, 12880883408, 51915526432, 207893871472, 827983736608**

---

**Réf.** AMM 28 114 21. JO61 150. jos.

**HIS2 A6228** **équations différentielles** **Suite P-récurrente**

**HIS1** **exponentielle**

$a(n) = (n^2 - 6n + 10)\, a(n-2)$

$$\exp(\arcsin(z))$$

**1, 1, 1, 2, 5, 20, 85, 520, 3145, 26000, 204425, 2132000, 20646925, 260104000, 2993804125, 44217680000, 589779412625, 9993195680000, 151573309044625, 2898026747200000**

## Bitriangular permutations

**Réf.** DMJ 13 267 46.

**HIS2 A6230** **Approximants de Padé**

**HIS1** **Fraction rationnelle**

$$\frac{(1 + z)\ (1 + 6\ z)}{(1 - z)\ (1 - 2\ z)\ (1 - 3\ z)}$$

**1, 13, 73, 301, 1081, 3613, 11593, 36301, 111961, 342013, 1038313, 3139501, 9467641, 28501213, 85700233, 257493901, 773268121, 2321377213, 6967277353, 20908123501**

---

## ∑ n(n-1) ... (n-k+1)/k, k=2..n

**Réf.** .rkg.

**HIS2 A6231** **P-récurrences** **Suite P-récurrente**

**HIS1** **exponentielle**

**Une solution de l'équation différentielle existe avec la fonction Ei(z), B. Salvy.**

$$\begin{aligned} a(n) = {} & (n + 3)\ a(n - 1) \\ & + (-\ 3\ n - 1)\ a(n - 2) \\ & + (3\ n - 3)\ a(n - 3) \\ & + (-\ n + 2)\ a(n - 4) \end{aligned}$$

**0, 1, 5, 20, 84, 409, 2365, 16064, 125664, 1112073, 10976173, 119481284, 1421542628, 18348340113, 255323504917, 3809950976992, 60683990530208, 1027542662934897**

**Réf.** JCT B24 208 78.
**HIS2 A6234** **Approximants de Padé**
**HIS1** **Fraction rationnelle**

$$\frac{2z - 1}{(1 - 3z)^2}$$

**1, 4, 15, 54, 189, 648, 2187, 7290, 24057, 78732, 255879, 826686, 2657205, 8503056, 27103491, 86093442, 272629233, 860934420, 2711943423, 8523250758, 26732013741**

---

## Complexity of doubled cycle

**Réf.** JCT B24 208 78.
**HIS2 A6235** **Approximants de Padé**
**HIS1** **Fraction rationnelle**

$$\frac{1 + 2z - 10z^2 + 2z^3 + z^4}{(z - 1)^2 (1 - 4z + z^2)^2}$$

**1, 12, 75, 384, 1805, 8100, 35287, 150528, 632025, 2620860, 10759331, 43804800, 177105253, 711809364, 2846259375, 11330543616, 44929049777, 177540878700, 699402223099**

## Triangular hex numbers

**Réf.** GA88 19. jos.

**HIS2 A6244** **Approximants de Padé**

**HIS1** **Fraction rationnelle**

$$\frac{1 - 8 z + z^2}{(1 - z) (z^2 - 98 z + 1)}$$

**1, 91, 8911, 873181, 85562821, 8384283271, 821574197731, 80505887094361, 7888755361049641, 773017519495770451, 75747828155224454551, 7422514141692500775541**

## Stacking bricks

**Réf.** GKP 360.

**HIS2 A6253** **Approximants de Padé**

**HIS1** **Fraction rationnelle**

$$\frac{1 - z}{(1 + z) (1 - 4 z + z^2)}$$

**1, 2, 9, 32, 121, 450, 1681, 6272, 23409, 87362, 326041, 1216800, 4541161, 16947842, 63250209, 236052992, 880961761, 3287794050, 12270214441, 45793063712, 170902040409**

**Réf.** MIS 4(3) 32 75.

**HIS2 A6261** **Approximants de Padé**

**HIS1** **Fraction rationnelle**

$$\frac{(1 - z + z^2)(1 - 3z + 3z^2)}{(z - 1)^6}$$

**1, 2, 4, 8, 16, 32, 63, 120, 219, 382, 638, 1024, 1586, 2380, 3473, 4944, 6885, 9402, 12616, 16664, 21700, 27896, 35443, 44552, 55455, 68406, 83682, 101584, 122438, 146596, 174437**

---

**Rooted genus-2 maps with n edges**

**Réf.** WA71. JCT 13 215 72.

**HIS2 A6298** **Hypergéométrique** **Suite P-récurrente**

**HIS1** **algébrique**

$$\frac{21 z (1 + z)}{(1 - 4z)^{11/2}}$$

**21, 483, 6468, 66066, 570570, 4390386, 31039008, 205633428, 1293938646, 7808250450, 45510945480**

## Royal paths in a lattice

**Réf.** CRO 20 12 73.

**HIS2 A6318** **Inverse fonctionnel** **Suite P-récurrente**

**HIS1** **algébrique**

**n a(n) = (6 n - 9) a(n - 1) + (- n + 3) a(n - 2)**

$$1/2 - 1/2\, z - 1/2\, (1 - 6\, z + z^2)^{1/2}$$

**1, 2, 6, 22, 90, 394, 1806, 8558, 41586, 206098, 1037718, 5293446, 27297738, 142078746, 745387038, 3937603038, 20927156706, 111818026018, 600318853926, 3236724317174**

## Royal paths in a lattice

**Réf.** CRO 20 18 73.

**HIS2 A6319** **Inverse fonctionnel** **Suite P-récurrente**

**HIS1** **algébrique**

**(n + 1) a(n) = (n - 4) a(n - 3) + (7 n - 4) a(n - 1) + (- 7 n + 17) a(n - 2)**

**S(z) est son propre inverse fonctionnel**

$$(1/2 - 1/2\, z - 1/2\, (1 - 6\, z + z^2)^{1/2})^2$$

**1, 4, 16, 68, 304, 1412, 6752, 33028, 164512, 831620, 4255728, 22004292, 114781008, 603308292, 3192216000, 16989553668, 90890869312, 488500827908, 2636405463248**

## Royal paths in a lattice

**Réf.** CRO 20 18 73.

**HIS2 A6320** **Inverse fonctionnel** **Suite P-récurrente**

**HIS1** **algébrique**

$(n + 2)\, a(n) = (9 n - 30)\, a(n - 3) + (- n + 5)\, a(n - 4) + (9 n + 3)\, a(n - 1)$

$$\left(1/2 - 1/2\, z - 1/2\, (1 - 6 z + z^2)^{1/2}\right)^3$$

**1, 6, 30, 146, 714, 3534, 17718, 89898, 461010, 2386390, 12455118, 65478978, 346448538, 1843520670, 9859734630, 52974158938, 285791932578, 1547585781414, 8408765223294**

## Royal paths in a lattice

**Réf.** CRO 20 18 73.

**HIS2 A6321** **LLL** **Suite P-récurrente**

**HIS1** **algébrique**

$(n + 3)\, a(n) = n\, a(n - 5) + (36 n - 88)\, a(n - 3) + (- 11 n + 47)\, a(n - 4) + (11 n + 14)\, a(n - 1) + (- 36 n + 20)\, a(n - 2) - 6\, a(n - 5)$

$$\left(1/2 - 1/2\, z - 1/2\, (1 - 6 z + z^2)^{1/2}\right)^4$$

**1, 8, 48, 264, 1408, 7432, 39152, 206600, 1093760, 5813000, 31019568, 166188552, 893763840, 4823997960, 26124870640, 141926904328, 773293020928, 4224773978632**

## Total preorders

**Réf.** MSH 53 20 76.

**HIS2 A6327** **Approximants de Padé** **Conjecture**

**HIS1** **Fraction rationnelle**

$$\frac{2 + z}{(1 - z)\,(1 - z - z^2)}$$

**2, 5, 10, 18, 31, 52, 86**

## From the enumeration of corners

**Réf.** CRO 6 82 65.

**HIS2 A6331** **Approximants de Padé**

**HIS1** **Fraction rationnelle**

$$\frac{2\,(1 + z)}{(z - 1)^4}$$

**2, 10, 28, 60, 110, 182, 280, 408, 570, 770, 1012, 1300, 1638, 2030, 2480, 2992, 3570, 4218, 4940, 5740, 6622, 7590, 8648, 9800, 11050, 12402, 13860, 15428, 17110, 18910, 20832**

## From the enumeration of corners

**Réf.** CRO 6 82 65.

**HIS2 A6332** **Approximants de Padé**

**HIS1** **Fraction rationnelle**

$$\frac{2\,(1+z)\,(1+6\,z+z^2)}{(1-z)^7}$$

**2, 28, 168, 660, 2002, 5096, 11424, 23256, 43890, 77924, 131560, 212940, 332514, 503440, 742016, 1068144, 1505826, 2083692, 2835560, 3801028, 5026098, 6563832, 8475040**

## From the enumeration of corners

**Réf.** CRO 6 82 65.

**HIS2 A6333** **Approximants de Padé**

**HIS1** **Fraction rationnelle**

$$\frac{z^5+20\,z^4+75\,z^3+75\,z^2+20\,z+1}{(z-1)^{10}}$$

**2, 60, 660, 4290, 20020, 74256, 232560, 639540, 1586310, 3617900, 7696260, 15438150, 29451240, 53796160, 94607040, 160908264, 265670730, 427156860, 670609940, 1030350090**

## From the enumeration of corners

**Réf.** CRO 6 82 65.

**HIS2 A6334** **hypergéométrique**

**HIS1** **Fraction rationnelle**

$$\frac{(z^7 + 42 z^6 + 364 z^5 + 1001 z^4 + 1001 z^3 + 364 z^2 + 42 z + 1) z}{(1 - z)^{13}}$$

**2, 110, 2002, 20020, 136136, 705432, 2984520, 10786908, 34370050, 98768670, 260390130, 638110200, 1468635168, 3200871520, 6650874912, 13248113736, 25415833170**

---

**Réf.** CRO 6 99 65.

**HIS2 A6335** **P-récurrences** **Suite P-récurrente**

**HIS1** **algébrique 3è degré**

$$- (2 n - 1) n a(n) =$$

$$- 6 (3 n - 4) (3 n - 5) a(n - 1)$$

**1, 2, 16, 192, 2816, 46592, 835584, 15876096, 315031552, 6466437120, 136383037440, 2941129850880, 64614360416256, 1442028424527872, 32619677465182208**

## Coloring a circuit with 4 colors

**Réf.** TAMS 60 355 46. BE74.

**HIS2 A6342** **Approximants de Padé**

**HIS1** **Fraction rationnelle**

$$\frac{2 z - 1}{(z - 1) (1 - 3 z) (1 + z)}$$

**1, 1, 4, 10, 31, 91, 274, 820, 2461**

---

## Related to series-parallel networks

**Réf.** AAP 4 127 72.

**HIS2 A6351** **Inverse fonctionnel**

**HIS1** **exponentielle** **f.g. exponentielle**

**S(z) est l'inverse fonctionnel de 2 ln(1 + z) - z**

$$- 1 - 2 W(- 1/2 \exp(- 1/2 + 1/2 z))$$

**1, 2, 8, 52, 472, 5504, 78416, 1320064, 25637824, 564275648, 13879795712, 377332365568, 11234698041088, 363581406419456, 12707452084972544, 477027941930515456**

## Distributive lattices

**Réf.** MSH 53 19 76. MSG 121 121 76.

**HIS2 A6356** **Approximants de Padé**

**HIS1** **Fraction rationnelle**

$$\frac{z^2 - z - 1}{z^3 - z^2 - 2z + 1}$$

**1, 3, 6, 14, 31, 70, 157, 353, 793, 1782, 4004**

## Distributive lattices

**Réf.** MSH 53 19 76. MSG 121 121 76.

**HIS2 A6357** **Approximants de Padé**

**HIS1** **Fraction rationnelle**

$$\frac{1 - z^2 + 2z - z^3}{(1 + z)(z^3 - 3z + 1)}$$

**1, 4, 10, 30, 85, 246, 707, 2037, 5864, 16886, 48620**

## Distributive lattices

**Réf.** MSH 53 19 76. MSG 121 121 76.

**HIS2 A6358** **Approximants de Padé**

**HIS1** **Fraction rationnelle**

$$\frac{(z - 1)(z^3 - 3z - 1)}{1 - 3z - 3z^2 + 4z^3 + z^4 - z^5}$$

**1, 5, 15, 55, 190, 671, 2353, 8272, 29056, 102091, 358671**

---

**Réf.** UPNT E17. jhc.

**HIS2 A6368** **Approximants de Padé**

**HIS1** **Fraction rationnelle**

$$\frac{1 + 3z + z^2 + 3z^3 + z^4}{(1 + z^2)(z - 1)^2 (1 + z)^2}$$

**1, 3, 2, 6, 4, 9, 5, 12, 7, 15, 8, 18, 10, 21, 11, 24, 13, 27, 14, 30, 16, 33, 17, 36, 19, 39, 20, 42, 22, 45, 23, 48, 25, 51, 26, 54, 28, 57, 29, 60, 31, 63, 32, 66, 34, 69, 35, 72, 37, 75, 38, 78, 40**

**Réf.** UPNT E17. jhc.
**HIS2 A6369** **Approximants de Padé**
**HIS1** **Fraction rationnelle**

$$\frac{(1 + z^2)(z^2 + 3 z + 1)}{(z - 1)^2 (z^2 + z + 1)^2}$$

**1, 3, 2, 5, 7, 4, 9, 11, 6, 13, 15, 8, 17, 19, 10, 21, 23, 12, 25, 27, 14, 29, 31, 16, 33, 35, 18, 37, 39, 20, 41, 43, 22, 45, 47, 24, 49, 51, 26, 53, 55, 28, 57, 59, 30, 61, 63, 32, 65, 67, 34, 69, 71**

---

**Image of n under the 3x+1 map**

**Réf.** UPNT 16.
**HIS2 A6370** **Approximants de Padé**
**HIS1** **Fraction rationnelle**

$$\frac{4 + z + 2 z^2}{(z - 1)^2 (1 + z)^2}$$

**4, 1, 10, 2, 16, 3, 22, 4, 28, 5, 34, 6, 40, 7, 46, 8, 52, 9, 58, 10, 64, 11, 70, 12, 76, 13, 82, 14, 88, 15, 94, 16, 100, 17, 106, 18, 112, 19, 118, 20, 124, 21, 130, 22, 136, 23, 142, 24, 148, 25, 154**

## Rooted nonseparable maps on the torus

**Réf.** JCT B18 241 75.

**HIS2 A6408** **Approximants de Padé**

**HIS1** **Fraction rationnelle**

$$\frac{z^2 + 11z + 4}{(z-1)^7}$$

**4, 39, 190, 651, 1792, 4242, 8988, 17490, 31812**

## Non-separable planar tree-rooted maps

**Réf.** JCT B18 243 75.

**HIS2 A6411** **Dérivée logarithmique**

**HIS1** **Fraction rationnelle**

$$\frac{2z + 3}{(1-z)^6}$$

**3, 20, 75, 210, 490, 1008, 1890, 3300, 5445, 8580, 13013**

## Non-separable toroidal tree-rooted maps

**Réf.** JCT B18 243 75.

**HIS2 A6414** **Dérivée logarithmique**

**HIS1** **Fraction rationnelle**

$$\frac{z^2 + 3z + 1}{(z - 1)^6}$$

**1, 9, 40, 125, 315, 686, 1344, 2430, 4125, 6655, 10296**

## Rooted planar maps

**Réf.** JCT B18 248 75.

**HIS2 A6416** **Approximants de Padé**

**HIS1** **Fraction rationnelle**

$$\frac{1 + 4z - 6z^2 + 2z^3}{(z - 1)^4}$$

**1, 8, 20, 38, 63, 96, 138, 190, 253, 328, 416, 518, 635**

## Rooted planar maps

**Réf.** JCT B18 248 75.

**HIS2 A6417** **Dérivée logarithmique**

**HIS1** **exponentielle**

$$\frac{\exp(z)\,(360 + 6840 z + 16560 z^2 + 8100 z^3 + 1395 z^4 + 93 z^5 + 2 z^6)}{360}$$

**1, 20, 131, 469, 1262, 2862, 5780, 10725, 18647, 30784, 48713, 74405**

## Rooted planar maps

**Réf.** JCT B18 249 75.

**HIS2 A6419** **P-récurrences** **Suite P-récurrente**

**HIS1**

$$(n + 2)\, a(n) =$$
$$(9 n + 10)\, a(n - 1)$$
$$- (24 n + 2)\, a(n - 2)$$
$$+ (16 n - 24)\, a(n - 3)$$

**1, 7, 37, 176, 794, 3473, 14893, 63004, 263950, 1097790, 4540386, 18696432, 76717268**

## Tree-rooted planar maps

**Réf.** JCT B18 256 75.

**HIS2 A6428** **Approximants de Padé**

**HIS1** **exponentielle**

$$\exp(z)\,(3 + 33 z + 33 z^{2} + 10 z^{3} + 9/8 z^{4} + 1/24 z^{5})$$

**0, 3, 36, 135, 360, 798, 1568, 2826, 4770, 7645, 11748, 17433**

## Tree-rooted planar maps

**Réf.** JCT B18 257 75.

**HIS2 A6431** **Hypergéométrique** **Suite P-récurrente**

**HIS1** **algébrique**

$$\frac{6 z^{2} - 6 z + 1 - (1 - 4 z)^{3/2}}{- 2 (1 - 4 z)^{3/2} z^{2}}$$

**0, 2, 15, 84, 420, 1980, 9009, 40040, 175032, 755820, 3233230, 13728792, 57946200**

## √n divides n

**Réf.** AMM 82 854 75. jos.

**HIS2 A6446** **Approximants de Padé**

**HIS1** **Fraction rationnelle**

$$\frac{1 + z + z^2 - z^3}{(z^2 + z + 1)^2 (z - 1)^3}$$

1, 2, 3, 4, 6, 8, 9, 12, 15, 16, 20, 24, 25, 30, 35, 36, 42, 48, 49, 56, 63, 64, 72, 80, 81, 90, 99, 100, 110, 120, 121, 132, 143, 144, 156, 168, 169, 182, 195, 196, 210, 224, 225, 240, 255, 256

## Solution to a diophantine equation

**Réf.** TR July 1973 p. 74. jos.

**HIS2 A6451** **Approximants de Padé**

**HIS1** **Fraction rationnelle**

$$\frac{z (2 + 3 z - 2 z^2 - z^3)}{(z - 1) (1 - 2 z - z^2) (z^2 - 2 z - 1)}$$

0, 2, 5, 15, 32, 90, 189, 527, 1104, 3074, 6437, 17919, 37520, 104442, 218685, 608735, 1274592, 3547970, 7428869, 20679087, 43298624, 120526554, 252362877, 702480239

## Solution to a diophantine equation

**Réf.** TR July 1973 p. 74. jos.

**HIS2 A6452** **Approximants de Padé**

**HIS1** **Fraction rationnelle**

$$\frac{(1 - z)(z^2 + 3z + 1)}{(1 - 2z - z^2)(z^2 - 2z - 1)}$$

**1, 2, 4, 11, 23, 64, 134, 373, 781, 2174, 4552, 12671, 26531, 73852, 154634, 430441, 901273, 2508794, 5253004, 14622323, 30616751, 85225144, 178447502, 496728541, 1040068261**

---

## Solution to a diophantine equation

**Réf.** TR July 1973 p. 74. jos.

**HIS2 A6454** **Approximants de Padé**

**HIS1** **Fraction rationnelle**

$$\frac{z(1 + 4z + z^2)}{3(1 - z)(z^2 - 6z + 1)(1 + 6z + z^2)}$$

**0, 3, 15, 120, 528, 4095, 17955, 139128, 609960, 4726275, 20720703, 160554240, 703893960, 5454117903, 23911673955, 185279454480, 812293020528, 6294047334435**

## Number of elements in Z[i] whose "smallest algorithm" is <= n

**Réf.** JALG 19 290 71. hwl.

**HIS2 A6457** **Approximants de Padé**

**HIS1** **Fraction rationnelle**

$$\frac{1 + z + 2 z^3}{(2 z - 1) (1 - 2 z^2) (1 - z)^2}$$

**1, 5, 17, 49, 125, 297, 669, 1457, 3093, 6457, 13309, 27201, 55237, 111689, 225101, 452689, 908885, 1822809, 3652701, 7315553, 14645349, 29311081, 58650733, 117342321, 234741877**

## Number of elements in Z[Ω] whose "smallest algorithm" is <= n

**Réf.** JALG 19 290 71. hwl.

**HIS2 A6458** **Approximants de Padé**

**HIS1** **Fraction rationnelle**

$$\frac{1 + 2 z + z^2 + 2 z^4 + 6 z^5}{(- 1 + 3 z) (2 z^3 + 2 z^2 - 1) (z - 1)^2}$$

**1, 7, 31, 115, 391, 1267, 3979, 12271, 37423, 113371, 342091, 1029799, 3095671, 9298147, 27914179, 83777503, 251394415, 754292827, 2263072411, 6789560412**

## Rooted planar maps

**Réf.** JCT B18 249 75.

**HIS2 A6468** **Approximants de Padé**

**HIS1** **Fraction rationnelle**

$$\frac{z^3 - 4 z^2 + 2 z + 5}{(z - 1)^7}$$

**5, 37, 150, 449, 1113, 2422, 4788, 8790, 15213, 25091, 39754, 60879**

## Rooted planar maps

**Réf.** JCT B18 251 75.

**HIS2 A6469** **Approximants de Padé**

**HIS1** **Fraction rationnelle**

$$\frac{3 z^2 - 9 z - 10}{(z - 1)^7}$$

**10, 79, 340, 1071, 2772, 6258, 12768, 24090, 42702, 71929**

## Rooted planar maps

**Réf.** JCT B18 257 75.

**HIS2 A6471** **Approximants de Padé**

**HIS1** **Fraction rationnelle**

$$\frac{4 z^3 + 35 z^2 + 34 z + 5}{(z - 1)^{10}}$$

**5, 84, 650, 3324, 13020, 42240, 118998, 300300, 693693, 1490060, 3011580**

---

**Réf.** JSCS 12 122 81.

**HIS2 A6472** **hypergéométrique** **f.g. exponentielle double**

**HIS1** **Fraction rationnelle**

**2 a(n) = (n - 1) n a(n - 1)**

$$\frac{4}{(z - 2)^2}$$

**1, 1, 3, 18, 180, 2700, 56700, 1587600, 57153600, 2571912000, 141455160000, 9336040560000, 728211163680000, 66267215894880000, 6958057668962400000**

**Réf.** BIT 13 93 73.
**HIS2 A6478** **Approximants de Padé**
**HIS1** **Fraction rationnelle**

$$\frac{1}{(z-1)\,(1-z-z^2)^2}$$

**1, 3, 8, 18, 38, 76, 147, 277, 512, 932, 1676, 2984, 5269, 9239, 16104, 27926, 48210, 82900, 142055, 242665, 413376, 702408, 1190808, 2014608, 3401833, 5734251, 9650312**

---

**From variance of Fibonacci search**

**Réf.** BIT 13 93 73.
**HIS2 A6479** **Approximants de Padé**
**HIS1** **Fraction rationnelle**

$$\frac{z^3\,(z^2+z+1)}{(1-z)\,(1-z-z^2)^3}$$

**0, 0, 0, 1, 5, 18, 52, 134, 318, 713, 1531, 3180, 6432, 12732, 24756, 47417, 89665, 167694, 310628, 570562, 1040226, 1883953, 3391799, 6073848, 10824096, 19204536, 33936456**

**Réf.** dsk.
**HIS2 A6483** **Approximants de Padé**
**HIS1** **Fraction rationnelle**

$$\frac{1 - 6 z^2}{(z - 1) (4 z^2 + 2 z - 1)}$$

**1, 3, 5, 17, 49, 161, 513, 1665, 5377, 17409, 56321, 182273, 589825, 1908737, 6176769, 19988481, 64684033, 209321985, 677380097, 2192048129, 7093616641, 22955425793**

**Réf.** dsk.
**HIS2 A6484** **Approximants de Padé**
**HIS1** **Fraction rationnelle**

$$\frac{1 - 2 z + 5 z^2}{(1 - z)^5}$$

**1, 3, 10, 30, 75, 161, 308, 540, 885, 1375, 2046, 2938, 4095, 5565, 7400, 9656, 12393, 15675, 19570, 24150, 29491, 35673, 42780, 50900, 60125, 70551, 82278, 95410, 110055, 126325**

## Generalized Lucas numbers

**Réf.** FQ 15 252 77.

**HIS2 A6490** **Approximants de Padé**

**HIS1** **Fraction rationnelle**

$$\frac{1 - 2z + 2z^2}{(1 - z - z^2)^3}$$

**1, 0, 3, 4, 10, 18, 35, 64, 117, 210, 374, 660, 1157, 2016, 3495, 6032, 10370, 17766, 30343, 51680, 87801, 148830, 251758, 425064, 716425, 1205568, 2025675, 3399004, 5696122**

## Generalized Lucas numbers

**Réf.** FQ 15 252 77.

**HIS2 A6491** **Approximants de Padé**

**HIS1** **Fraction rationnelle**

$$\frac{(1 - 2z + 2z^2)(z - 1)}{(1 - z - z^2)^3}$$

**1, 0, 4, 5, 15, 28, 60, 117, 230, 440, 834, 1560, 2891, 5310, 9680, 17527, 31545, 56468, 100590, 178395, 315106, 554530, 972564, 1700400, 2964325, 5153868, 8938300, 15465497**

## Generalized Lucas numbers

**Réf.** FQ 15 252 77.

**HIS2 A6492** **Approximants de Padé**

**HIS1** **Fraction rationnelle**

$$\frac{(1 - 2z + 2z^2)(z - 1)^2}{(1 - z - z^2)^4}$$

**1, 0, 5, 6, 21, 40, 93, 190, 396, 796, 1586, 3108, 6025, 11552, 21947, 41346, 77311, 143580, 265013, 486398, 888122, 1613944, 2920100, 5261880, 9445905, 16897328, 30127665**

## Generalized Lucas numbers

**Réf.** FQ 15 252 77.

**HIS2 A6493** **Approximants de Padé**

**HIS1** **Fraction rationnelle**

$$\frac{(1 - 2z + 2z^2)(z - 1)^3}{(1 - z - z^2)^5}$$

**1, 0, 6, 7, 28, 54, 135, 286, 627, 1313, 2730, 5565, 11212, 22304, 43911, 85614, 165490, 317373, 604296, 1143054, 2149074, 4017950, 7473180, 13832910, 25490115, 46774448**

**Réf.** FQ 15 292 77.
**HIS2 A6497** **Approximants de Padé**
**HIS1** **Fraction rationnelle**

$$\frac{2 - 3z}{1 - 3z - z^2}$$

**2, 3, 11, 36, 119, 393, 1298, 4287, 14159, 46764, 154451, 510117, 1684802, 5564523, 18378371, 60699636, 200477279, 662131473, 2186871698, 7222746567, 23855111399**

---

## Restricted combinations

**Réf.** FQ 16 113 78.
**HIS2 A6498** **Approximants de Padé**
**HIS1** **Fraction rationnelle**

$$\frac{1 + z + 2z^2 + z^3}{(1 - z - z^2)(1 + z^2)}$$

**1, 2, 4, 6, 9, 15, 25, 40, 64, 104, 169, 273, 441, 714, 1156, 1870, 3025, 4895, 7921, 12816, 20736, 33552, 54289, 87841, 142129, 229970, 372100, 602070, 974169, 1576239, 2550409**

## Restricted circular combinations

**Réf.** FQ 16 115 78.

**HIS2 A6499** **Approximants de Padé**

**HIS1** **Fraction rationnelle**

$$\frac{1 + 2z + 6z^2 + 2z^3}{(1 - z - z^2)(1 + z^2)}$$

**1, 3, 9, 12, 16, 28, 49, 77, 121, 198, 324, 522, 841, 1363, 2209, 3572, 5776, 9348, 15129, 24477, 39601, 64078, 103684, 167762, 271441, 439203, 710649, 1149852, 1860496**

## Restricted combinations

**Réf.** FQ 16 116 78.

**HIS2 A6500** **Approximants de Padé**

**HIS1** **Fraction rationnelle**

$$\frac{z^7 + 2z^6 + z^5 - z^4 - 3z^3 - z^2 - z - 1}{(z^6 - z^3 - 1)(1 - z - z^2)}$$

**1, 2, 4, 8, 12, 18, 27, 45, 75, 125, 200, 320, 512, 832, 1352, 2197, 3549, 5733, 9261, 14994, 24276, 39304, 63580, 102850, 166375, 269225, 435655, 704969, 1140624, 1845504, 2985984**

**Réf.** FQ 16 116 78.

**HIS2 A6501** **Approximants de Padé**

**HIS1** **Fraction rationnelle**

$$\frac{1 + z^2}{(z^2 + z + 1)^2 (z - 1)^4}$$

**1, 2, 4, 8, 12, 18, 27, 36, 48, 64, 80, 100, 125, 150, 180, 216, 252, 294, 343, 392, 448, 512, 576, 648, 729, 810, 900, 1000, 1100, 1210, 1331, 1452, 1584, 1728, 1872, 2028, 2197, 2366**

---

**Réf.** FQ 14 43 76.

**HIS2 A6503** **Approximants de Padé**

**HIS1** **Fraction rationnelle**

$$\frac{2 z - 3}{(1 - z)^4}$$

**3, 10, 22, 40, 65, 98, 140, 192, 255, 330, 418, 520, 637, 770, 920, 1088, 1275, 1482, 1710, 1960, 2233, 2530, 2852, 3200, 3575**

**Réf.** FQ 14 43 76.
**HIS2 A6504** **Approximants de Padé**
**HIS1** **Fraction rationnelle**

$$\frac{5 - 5z + z^2}{(1 - z)^5}$$

**5, 20, 51, 105, 190, 315, 490, 726, 1035, 1430, 1925, 2535, 3276, 4165, 5220, 6460, 7905, 9576, 11495, 13685, 16170, 18975, 22126, 25650, 29575**

---

**Réf.** FQ 14 69 76.
**HIS2 A6505** **équations différentielles** **Formule de B. Salvy**
**HIS1** **exponentielle**

$$\exp(\exp(z) - z - 1/2\, z^2 - 1)$$

**1, 0, 0, 1, 1, 1, 11, 36, 92, 491, 2557, 11353, 60105, 362506, 2169246, 13580815, 91927435, 650078097, 4762023647, 36508923530, 292117087090, 2424048335917, 20847410586719**

**Réf.** HO73 113.

**HIS2 A6516** **Approximants de Padé**

**HIS1** **Fraction rationnelle**

$$\frac{1}{(1 - 2 z) (1 - 4 z)}$$

**1, 6, 28, 120, 496, 2016, 8128, 32640, 130816, 523776, 2096128, 8386560, 33550336, 134209536, 536854528, 2147450880, 8589869056, 34359607296, 137438691328, 549755289600**

---

**Réf.** HO73 102.

**HIS2 A6522** **Approximants de Padé**

**HIS1** **Fraction rationnelle**

$$\frac{1 - z + z^2}{(z - 1)^5}$$

**1, 4, 11, 25, 50, 91, 154, 246, 375, 550, 781, 1079, 1456, 1925, 2500, 3196, 4029, 5016, 6175, 7525, 9086, 10879, 12926, 15250, 17875, 20826, 24129, 27811, 31900, 36425, 41416**

**Réf.** GA66 246.
**HIS2 A6527** **Approximants de Padé**
**HIS1** **Fraction rationnelle**

$$\frac{z(1+z^2)}{(z-1)^4}$$

**0, 1, 4, 11, 24, 45, 76, 119, 176, 249, 340, 451, 584, 741, 924, 1135, 1376, 1649, 1956, 2299, 2680, 3101, 3564, 4071, 4624, 5225, 5876, 6579, 7336, 8149, 9020, 9951, 10944, 12001**

---

**Réf.** GA66 246.
**HIS2 A6528** **Approximants de Padé**
**HIS1** **Fraction rationnelle**

$$\frac{z(1+z+4z^2)}{(1-z)^5}$$

**0, 1, 6, 24, 70, 165, 336, 616, 1044, 1665, 2530, 3696, 5226, 7189, 9660, 12720, 16456, 20961, 26334, 32680, 40110, 48741, 58696, 70104, 83100, 97825, 114426, 133056, 153874**

## Cubes with sides of n colors

**Réf.** GA66 246.

**HIS2 A6529** **Approximants de Padé**

**HIS1** **Fraction rationnelle**

$$\frac{z\,(1 + 5\,z + 17\,z^2 + 77\,z^3)}{(1 - z)^5}$$

**0, 1, 10, 57, 272, 885, 2226, 4725, 8912, 15417, 24970, 38401, 56640, 80717, 111762, 151005, 199776, 259505, 331722, 418057, 520240, 640101, 779570, 940677, 1125552, 1336425**

---

## C(n , 3) C(n - 1, 3) / 4

**Réf.**

**HIS2 A6542** **Approximants de Padé**

**HIS1** **Fraction rationnelle**

$$\frac{1 + 3\,z + z^2}{(1 - z)^7}$$

**1, 10, 50, 175, 490, 1176, 2520, 4950, 9075, 15730, 26026, 41405, 63700, 95200, 138720, 197676, 276165, 379050, 512050, 681835, 896126, 1163800, 1495000, 1901250, 2395575**

## n-coloring a cube

**Réf.** C1 254.

**HIS2 A6550** **Approximants de Padé**

**HIS1** **Fraction rationnelle**

$$\frac{1 + 3 z + 8 z^2 + 10 z^3 + 14 z^4 - 6 z^5}{(1 - z)^7}$$

1, 10, 57, 234, 770, 2136, 5180, 11292, 22599, 42190, 74371, 124950, 201552, 313964, 474510, 698456, 1004445, 1414962, 1956829, 2661730, 3566766, 4715040, 6156272, 7947444

## Icosahedral numbers

**Réf.**

**HIS2 A6564** **Approximants de Padé**

**HIS1** **Fraction rationnelle**

$$\frac{1 + 8 z + 6 z^2}{(z - 1)^4}$$

1, 12, 48, 124, 255, 456, 742, 1128, 1629, 2260, 3036, 3972, 5083, 6384, 7890, 9616, 11577, 13788, 16264, 19020, 22071, 25432, 29118, 33144, 37525, 42276, 47412, 52948, 58899

## Colored hexagons

Réf.

HIS2 A6565 Approximants de Padé

HIS1 Fraction rationnelle

$$\frac{1 + 7z + 53z^2 + 49z^3 + 10z^4}{(1 - z)^7}$$

1, 14, 130, 700, 2635, 7826, 19684, 43800, 88725, 166870, 295526, 498004, 804895, 1255450, 1899080, 2796976, 4023849, 5669790, 7842250, 10668140, 14296051, 18898594

## Dodecahedral numbers

Réf.

HIS2 A6566 Approximants de Padé

HIS1 Fraction rationnelle

$$\frac{1 + 16z + 10z^2}{(1 - z)^4}$$

1, 20, 84, 220, 455, 816, 1330, 2024, 2925, 4060, 5456, 7140, 9139, 11480, 14190, 17296, 20825, 24804, 29260, 34220, 39711, 45760, 52394, 59640, 67525, 76076, 85320, 95284

**Réf.** mlb.

**HIS2 A6578** **Approximants de Padé**

**HIS1** **Fraction rationnelle**

**a(n) = ∑ max(n,n-k), k=1...n-1**

$$\frac{1 + 2z}{(1 + z)(1 - z)^3}$$

**1, 4, 8, 14, 21, 30, 40, 52, 65, 80, 96, 114, 133, 154, 176, 200, 225, 252, 280, 310, 341, 374, 408, 444, 481, 520, 560, 602, 645, 690, 736, 784, 833, 884, 936, 990, 1045, 1102, 1160**

---

## Generalized Fibonacci numbers

**Réf.** LNM 622 186 77.

**HIS2 A6603** **LLL** **Suite P-récurrente**

**HIS1** **algébrique**

$$n\, a(n) = -n\, a(n-5) + (7n-9)\, a(n-1) + (-8n+12)\, a(n-2) + (6n-12)\, a(n-3) + (5n-6)\, a(n-4) + 3\, a(n-5)$$

$$1/2\ \frac{1 - z - 2z^2 - (1 - 6z + z^2)^{1/2}}{2z - z^2 + z^3 + z^4}$$

**1, 2, 7, 26, 107, 468, 2141, 10124, 49101, 242934, 1221427, 6222838, 32056215, 166690696, 873798681, 4612654808, 24499322137, 130830894666, 702037771647, 3783431872018**

## Generalized Fibonacci numbers

**Réf.** LNM 622 186 77.

**HIS2 A6604** **LLL** **Suite P-récurrente**

**HIS1** **algébrique**

**n a(n) = (- 1/2 n + 3/2) a(n - 5) + (7/2 n - 6) a(n - 4) + (13/2 n - 9) a(n - 1) + (- 7/2 n + 15/2) a(n - 2) + (- 3 n + 3) a(n - 3)**

$$\frac{1}{2}\,\frac{1 + z - 2 z^2 - (1 - 6 z + z^2)^{1/2}}{2 z - z^2 - z^3 + z^4}$$

**1, 1, 4, 13, 53, 228, 1037, 4885, 23640, 116793, 586633, 2986616, 15377097, 79927913, 418852716, 2210503285, 11738292397, 62673984492, 336260313765**

## Modes of connections of 2n points

**Réf.** LNM 686 326 78.

**HIS2 A6605** **LLL** **Suite P-récurrente**

**HIS1** **algébrique** **P-récurrence du 3è degré**

S(z) satisfait à

$$\frac{1 - S(z) + S(z)^2 z + S(z)^4 z^2}{z^2}$$

**1, 1, 3, 11, 46, 207, 979, 4797, 24138, 123998, 647615, 3428493, 18356714, 99229015, 540807165, 2968468275, 16395456762, 91053897066, 508151297602, 2848290555562**

## From generalized Catalan numbers

**Réf.** LNM 952 279 82.

**HIS2 A6629** **LLL** **La F.G. est algébrique du 3è degré et**
**HIS1** **algébrique** **prend trop de place.**

$$3F2([2, 5/3, 4/3],[3, 5/2],27 z/4)$$

**1, 4, 18, 88, 455, 2448, 13566, 76912, 444015, 2601300, 15426840, 92431584, 558685348, 3402497504, 20858916870, 128618832864, 797168807855, 4963511449260, 31032552351570**

## From generalized Catalan numbers

**Réf.** LNM 952 279 82.

**HIS2 A6630** **Hypergéométrique** **La F.G. est algébrique du 3è degré et**
**HIS1** **algébrique** **prend trop de place.**

$$3F2([2, 8/3, 7/3],[4, 7/2],27 z/4)$$

**1, 6, 33, 182, 1020, 5814, 33649, 197340, 1170585, 7012200, 42364476, 257854776, 1579730984, 9734161206, 60290077905, 375138262520, 2343880406595, 14699630061270**

## From generalized Catalan numbers

**Réf.** LNM 952 279 82.

**HIS2 A6631** **LLL** **Suite P-récurrente**

**HIS1** **algébrique**

**La F.G. est algébrique du 3è degré et prend trop de place.**

$$ {}_3F_2([3, 8/3, 10/3],[5, 9/2],27 z/4) $$

**1, 8, 52, 320, 1938, 11704, 70840, 430560, 2629575, 16138848, 99522896, 616480384, 3834669566, 23944995480, 150055305008, 943448717120, 5949850262895, 37628321318280**

---

## From generalized Catalan numbers

**Réf.** LNM 952 280 82.

**HIS2 A6632** **Hypergéométrique** **Suite P-récurrente**

**HIS1** **algébrique** **Inverse de A2293**

$$ \frac{1}{1 + z\, {}_4F_3([1, 7/4, 5/4, 3/2], [2, 5/3, 7/3],256 z/ 27)} $$

**1, 3, 15, 91, 612, 4389, 32890, 254475, 2017356, 16301164, 133767543, 1111731933, 9338434700, 79155435870, 676196049060, 5815796869995, 50318860986108**

## From generalized Catalan numbers

**Réf.** LNM 952 280 82.

**HIS2 A6633** **Hypergéométrique** **Suite P-récurrente**

**HIS1** **algébrique**

$${}_4F_3\left([2, 9/4, 3/2, 7/4], [3, 8/3, 7/3], 256 z / 27\right)$$

**1, 6, 39, 272, 1995, 15180, 118755, 949344, 7721604, 63698830, 531697881, 4482448656, 38111876530, 326439471960, 2814095259675, 24397023508416, 212579132600076**

## From generalized Catalan numbers

**Réf.** LNM 952 280 82.

**HIS2 A6634** **Hypergéométrique** **Suite P-récurrente**

**HIS1** **algébrique**

$${}_4F_3\left([3, 9/4, 5/2, 11/4], [4, 10/3, 11/3], 256 z / 27\right)$$

**1, 9, 72, 570, 4554, 36855, 302064, 2504304, 20974005, 177232627, 1509395976, 12943656180, 111676661460, 968786892675, 8445123522144, 73940567860896,**

## From generalized Catalan numbers

**Réf.** LNM 952 280 82.

**HIS2 A6635** **Hypergéométrique** **Suite P-récurrente**

**HIS1** **algébrique**

$$ {}_4F_3\left([3, 7/2, 15/4, 13/4], [5, 14/3, 13/3], 256 z / 27\right) $$

**1, 12, 114, 1012, 8775, 75516, 649264, 5593068, 48336171, 419276660, 3650774820, 31907617560, 279871768995, 2463161027292, 21747225841440, 19257567355**

## Closed meanders

**Réf.** SFCA 292.

**HIS2 A6659** **Hypergéométrique** **Suite P-récurrente**

**HIS1** **algébrique**

$$ \frac{32}{(1 - 4 z)^{1/2} \left(1 + (1 - 4 z)^{1/2}\right)^4} $$

**2, 12, 56, 240, 990, 4004**

**Planted binary phylogenetic trees with n labels**

**Réf.** LNM 884 196 81.

**HIS2 A6677** **Inverse fonctionnel** **erreurs dans la suite**

**HIS1** **exponentielle (algébrique)**

$$1 - (3 - 2 \exp(z))^{1/2}$$

**1, 2, 7, 41, 346, 3797, 51157, 816356, 15050581, 34459425**

---

**Planted binary phylogenetic trees with n labels**

**Réf.** LNM 884 196 81.

**HIS2 A6678** **Inverse fonctionnel**

**HIS1** **algébrique**

$$\frac{1 - (1 - 2 z - 2 z^2)^{1/2}}{1 + z}$$

**1, 1, 6, 39, 390, 4815, 73080, 1304415, 26847450, 625528575**

## Planted binary phylogenetic trees with n labels

**Réf.** LNM 884 196 81.

**HIS2 A6679** **Inverse fonctionnel**

**HIS1** **exponentielle (algébrique)**

$$\frac{1}{\exp(z)} + \frac{(1 + 2\exp(z) - 2\exp(z)^2)^{1/2}}{\exp(z)}$$

**1, 2, 10, 83, 946, 13772, 244315, 5113208, 123342166, 3369568817**

---

**Réf.** R1 38. sls.

**HIS2 A6790** **Recoupements**

**HIS1** **exponentielle**

$$\frac{\exp(z)}{2 - \exp(z)}$$

**1, 2, 6, 26, 150, 1082, 9366, 94586, 1091670, 14174522, 204495126, 3245265146, 56183135190, 1053716696762, 21282685940886, 460566381955706, 10631309363962710**

## Extreme points of set of n x n symmetric doubly-stochastic matrices

**Réf.** JCT 8 422 70. EJC 1 180 80.

**HIS2 A6847** **Dérivée logarithmique** **Suite P-récurrente**

**HIS1** **exponentielle (algébrique)**

$$a(n) = n^3\ a(n-1) + (4n^3 - 4n^2 + n)\ a(n-2) + (-3n^3 + 5/2\ n^2 - 1/2\ n)\ a(n-3) + (24n^3 - 26n^2 + 9n - 1)\ a(n-4)$$

$$\frac{(z+1)^{1/4}\ \exp(1/2\ z\ (z+1))}{(z-1)^{1/4}}$$

**1, 1, 2, 5, 14, 58, 238, 1516, 9020, 79892, 635984, 7127764, 70757968, 949723600, 11260506056, 175400319992, 2416123951952, 42776273847184, 671238787733920**

## Extreme points of set of n x n symmetric doubly-substochastic matrices

**Réf.** EJC 1 180 80.

**HIS2 A6848** **Dérivée logarithmique**

**HIS1** **exponentielle (algébrique)**

$$\frac{(z+1)^{1/4}\ \exp\left(\frac{z\ (z^3 - z + 2z^2 - 3)}{2\ (z-1)\ (z+1)}\right)}{(z-1)^{1/4}}$$

**1, 2, 5, 18, 75, 414, 2643, 20550, 180057, 1803330, 19925541, 242749602, 3218286195, 46082917278, 710817377715, 11689297807734, 205359276208113, 3812653265319810**